\documentclass{article}
\usepackage{amssymb}
\usepackage{graphicx}
\usepackage{amsmath}
\numberwithin{equation}{section}
\newtheorem{theorem}{Theorem}[subsection]

\newtheorem{corollary}[theorem]{Corollary}

\newtheorem{definition}[theorem]{Definition}

\newtheorem{lemma}[theorem]{Lemma}

\newtheorem{proposition}[theorem]{Proposition}
\newtheorem{remark}[theorem]{Remark}

\begin{document}

\title{Quantitative estimates of unique continuation for parabolic
equations, determination of unknown time-varying boundaries and optimal
stability estimates}
\author{Sergio Vessella\thanks{%
This work was partially supported by MIUR, PRIN, number 2006014115.}}
\date{}
\maketitle

\begin{abstract}
In this paper we will review the main results concerning the issue of
stability for the determination unknown boundary portion of a thermic
conducting body from Cauchy data for parabolic equations. We give detailed
and selfcontained proofs. We prove that such problems are severely ill-posed
in the sense that under a priori regularity assumptions on the unknown
boundaries, up to any finite order of differentiability, the continuous
dependence of unknown boundary from the measured data is, at best, of
logarithmic type.

\ \ We review the main results concerning quantitative estimates of unique
continuation for solutions to second order parabolic equations. We give a
detailed proof of \ a Carleman estimate crucial for the derivation of the
stability estimates.
\end{abstract}

\section{Introduction\label{Sec1}}

In this paper we will review the main results concerning the quantitative
estimates of unique continuation for solutions to second order parabolic
equations (with real coefficients) and the issue of the stability for
various type of inverse parabolic problems with unknown boundaries. We want
to stress that we are going to be quite sketchy on the applications and the
methods of numerical reconstruction of the solution of such problems. Such a
choice is due to the wideness of the results on this subject.

Let $\left\{ \Omega \left( t\right) \right\} _{t\in \left[
0,T\right] }$ be a family of bounded domains in $\mathbb{R}^{n}$,
where $T$ is a given positive number. We shall suppose that the
boundary of $\Omega \left( \left( 0,T\right) \right)
:=\bigcup\limits_{t\in \left( 0,T\right) }\Omega \left(
t\right) \times \left\{ t\right\} $ is sufficiently smooth, but for every $%
t\in \left( 0,T\right) $ a part $I\left( t\right) $\ of $\partial \Omega
\left( t\right) $ is not known. In many applications $\Omega \left( t\right)
$ represents (at a fixed time $t$) a thermic conducting body or specimen.
The portion $I\left( t\right) $ can be, for any $t\in \left( 0,T\right) $,
some interior component of $\partial \Omega \left( t\right) $ or some
inaccessible portion of the exterior component of $\partial \Omega \left(
t\right) $. The inverse problems consist in determining $I\left( t\right) $,
for every $t\in \left( 0,T\right) $, by means of thermal measurements on the
accessible part $A\left( t\right) :=\partial \Omega \left( t\right)
\smallsetminus I\left( t\right) $. Throughout the paper, for the sake of
simplicity, we shall assume that $A\left( t\right) $ is not time varying. In
order to give a sufficiently general mathematical formulation of the
problems, let us consider a boundary linear operator $\mathcal{B}$ on $%
\bigcup\limits_{t\in \left[ 0,T\right] }\partial \Omega \left( t\right)
\times \left\{ t\right\} $ which we shall specify later and let $f$ be a
nontrivial function, we consider the following initial-boundary value
problem (the direct problem), where we set $I\left( \left( 0,T\right]
\right) :=\bigcup\limits_{t\in \left( 0,T\right] }I\left( t\right) \times
\left\{ t\right\} $.
\begin{equation}
\left\{
\begin{array}{c}
div\left( \kappa \nabla u\right) -\partial _{t}u=0\text{, \ \ in }\Omega
\left( \left( 0,T\right) \right) \text{,} \\
\mathcal{B}u=f\text{, \ \ \ \ \ \ \ \ \ \ \ on }A\times \left( 0,T\right]
\text{,\ } \\
\mathcal{B}u=0\text{, \ \ \ \ \ \ \ \ \ \ \ \ \ on }I\left( \left( 0,T\right]
\right) \text{, \ \ } \\
u\left( .,0\right) =u_{0}\text{, \ \ \ \ \ \ \ \ \ \ in }\Omega \left(
0\right) .%
\end{array}%
\right.  \label{In1}
\end{equation}%
Here $\kappa =\left\{ \kappa ^{ij}\left( x,t,u\right) \right\} _{i,j=1}^{n}$
denotes the known symmetric thermal conductivity tensor which satisfies a
hypothesis of uniform ellipticity. The boundary operator $\mathcal{B}$ may
be one of the following operators

\noindent (a) $\mathcal{B}u=u$ (boundary Dirichlet operator),

\noindent (b) $\mathcal{B}u=\kappa \nabla u\cdot \nu $ (boundary Neumann
operator),

\noindent (c) $\mathcal{B}u=\kappa \nabla u\cdot \nu +b_{0}u$, $%
b_{0}\not\equiv 0$ (boundary Robin operator),

\noindent where $\nu $ denotes the unit exterior normal to $\Omega \left(
t\right) $.

\noindent Given an open portion $\Sigma $ of $\partial \Omega \left(
t\right) $ such that $\Sigma \subset A$, we consider the inverse problems of
determinig $I\left( t\right) $, $t\in \left( 0,T\right] $, from the
knowledge, in case (a) of $\kappa \nabla u\cdot \nu $ on $\Sigma \times %
\left[ 0,T\right] $ and, in case (b) and (c), from the knowledge of $u$ on $%
\Sigma \times \left[ 0,T\right] $. We shall refer to the inverse problems
introduced above as the \textit{inverse Dirichlet problem, }the \textit{%
inverse Neumann problem }and\textit{\ }the\textit{\ inverse Robin problem }%
respectively (for short: \textit{Dirichlet, Neumann} or \textit{Robin }case
respectively).

In many applications $I\left( t\right) $ represents:

\noindent (i) the boundary of a cavity or a corroded part of $\Omega \left(
t\right) $, \cite{BiCFIn}- \cite{BrC05}, \cite{CrWi}, \cite{CroFKT}, \cite%
{GriBiM}, \cite{HohRS}, \cite{Mal}, \cite{VaGrBi}, \cite{ZPM}

or

\noindent (ii) a privileged isothermal surface, such as a solidification
front of $\Omega \left( t\right) $, that is not accessible to direct
inspection, \cite{ApAtk}, \cite{Bin}, \cite{BinEnVe}, \cite{CEn}, \cite%
{EnLaMa}, \cite{Gre}, \cite{MaVe}, \cite{RiTo}, \cite{SaGaFrSp}.

The boundary operator $\mathcal{B}$ is, in case (i), either the boundary
Neumann operator \cite{AnAbJ}, \cite{BrC96}, \cite{KaSaVo} or the boundary
Robin operator \cite{Ing1}, \cite{BiCFIn}, \cite{BiFIn04}. In case (ii), $%
\mathcal{B}$ is the boundary Dirichlet operator. Most of the problems of (i)
arise in thermal imaging. In such problems, and in all the mentioned papers
of (i), the unknown boundary is assumed not time varying. Therefore our
formulation with time varying boundary seems a mere, although quite natural,
mathematical generalization of the original problems. On the contrary, the
assumption of time varying boundary, as well as the assumption of dependence
of $\kappa $ on $u$, seems particularly appropriate for the problems of case
(ii).

We want to stress that in the present paper we study in details the
stability issue for the Dirichlet case only. The stability issue for
problems in the type (i) (Dirichlet and Neumann case) is studied in \cite%
{Ve1}, \cite{CRoVe1}, \cite{CRoVe2}, \cite{DcRVe}. The stability issue for
problems of type (ii) is studied in \cite{BinEnVe}, \cite{MaVe}, \cite{Ve0}.
We shall illustrate below extensively the main problems that arise in the
stability issue. Actually, concerning the uniqueness and the stability, the
Neumann case, when the unknown boundary is not time varying, is treated
similarly to the Dirichlet case. In the case of time varying unknown
boundaries, uniqueness and stability for the inverse Neumann problem is an
open question. It is quite remarkable that uniqueness and stability for the
Robin case is an open question even if the unknown boundaries are not time
varying. Concerning such a problem we recall the papers \cite{BiFIn04}, \cite%
{Is5}. In \cite{Is5} the uniqueness for such a case has been proved under
the (essential) hypotheses that $\kappa $ depends on $x$ only and the
unknown boundary is not time varying. In \cite{BiFIn04} the uniqueness for
the Robin case is proved in dimension two under the assumption that the
unknown boundary is a small perturbation of a graph. Concerning the methods
of reconstruction of unknown boundaries, we refer to the papers \cite{Da},
\cite{Ike1}, \cite{Ike2}, \cite{Nak}, \cite{KaMoTs}. Here the so-called
probe method and enclosure method seem to be effective and promising, at
least for parabolic equations with constant coefficients, for a constructive
determination of unknown boundaries from boundary measurements. The issue of
numerical reconstruction of the unknown boundaries is studied in \cite%
{BiCFIn}, \cite{BiFIn04}, \cite{Br2}, \cite{BrC96}, \cite{BrC98}, \cite%
{BrC05}, \cite{CKY98}, \cite{CKY99}.

The inverse problems studied in the present paper are severely ill-posed in
the sense that under a priori smoothness assumptions on the unknown
boundaries $I\left( t\right) $, $t\in \left[ 0,T\right] $, up to any finite
order of differentiability, the continuous dependence of $I\left( t\right) $
from the measured data is, at best, of logarithmic type. Such a severely
ill-posed character of the above inverse problems has been proved in \cite%
{DcRVe}, see also \cite{Dc}, when $\kappa $ is constant, $u_{0}=0$ and the
unknown boundary \textit{is not} time varying , see Section \ref{INSTNT} for
a sketch of the proof. In the present paper, (see Section \ref{INSTT}), we
prove the (severely) ill-posedness when, $\kappa $ is a constant, $u_{0}=0$
and the unknown boundary \textit{is} time-varying. It is interesting to
observe that the severely ill-posed character has been proved, in \cite%
{DcRVe} and Section \ref{INSTT}, when instead of a single measurement the
whole Dirichlet-to-Neumann map is known.

The first results concerning this type of inverse problems have been
obtained for elliptic equations, see \cite{AnAbJ}, \cite{FIng1}, \cite{FIng2}%
, \cite{KaSa}, \cite{KaSaVo}, \cite{IngSa}, \cite{Ing1}, \cite{IngMa}, \cite%
{Sin}, \cite{VaGrBi}, \cite{ZPM}. We will describe first the main steps of
the proofs of stability and unique continuation in this case because many of
these can be extended to the parabolic case. The mathematical formulation of
the inverse elliptic problem is the following one. Let $\Omega $ be a domain
in $\mathbb{R}^{n}$ with a sufficiently smooth boundary $\partial \Omega $ a
part of which $I$ is unknown. Let $\mathcal{B}$ be a boundary linear
operator on $\partial \Omega $ and let $f$ be a nontrivial function. Let $u$
be the solution to the following boundary value problem.
\begin{equation}
\left\{
\begin{array}{c}
div\left( \sigma \nabla u\right) =0\text{, \ \ in }\Omega \text{,} \\
\mathcal{B}u=f\text{, \ \ \ \ \ \ \ \ \ \ \ on }A\text{,\ } \\
\mathcal{B}u=0\text{, \ \ \ \ \ \ \ \ \ \ \ \ \ on }I\text{, \ \ }%
\end{array}%
\right.  \label{Ellpr}
\end{equation}%
(with a possible normalization condition, such as $\int\nolimits_{\Omega
}u=0 $, if $\mathcal{B}$ is the boundary Neumann operator). Here $\sigma $
denotes the known electrical conductivity (symmetric) tensor satisfying a
hypothesis of uniform ellipticity. Given an open portion $\Sigma $ of $%
\partial \Omega $ such that $\Sigma \subset A$ we are interested in
determining $I$ from the knowledge of $\kappa \nabla u\cdot \nu $ on $\Sigma
$, in the Dirichlet case, and from the knowledge of $u$ on $\Sigma $, in the
Neumann and in the Robin case. The first stability results for the just
mentioned problems (in the Dirichlet and Neumann case) have been proved,
when $\sigma $ is the identity matrix, in \cite{BeVe} for $n=2$ and in \cite%
{CheHY}, \cite{BukCY1}, \cite{BukCY2} for $n=3$. In the last three papers
the unknown boundary is a graph. In all the cited papers the stability
results, proved under a priori regularity assumption on $I$ up to a finite
order of differentiability, are of logarithmic type. Such a type of
stability is optimal, \cite{Al1}, \cite{DcR}. Logarithmic stability
estimates for $n=2$, when $\sigma $ is inhomogeneous, discontinuous and
under relaxed a priori assumption on $I$, have been proved in \cite{AlR},
\cite{R} and, for $n\geq 3$, in \cite{AlBRVe1} (see also \cite{AlBRVe2}). In
\cite{BaVe} a Lipschitz stability estimate has been proved (in Dirichlet
case) under the a priori assumption that $I$ is a polygonal line, $\sigma $
is homogeneous, $n=2$.

In the above mentioned papers \cite{BeVe}, \cite{CheHY}, \cite{BukCY1}, \cite%
{BukCY2}, \cite{AlR}, \cite{R} the proof of stability results is based on
arguments related to complex analytic methods which do not carry over the
case $n\geq 3$. In \cite{AlBRVe1}, where the case $n\geq 3$ is studied, the
proof of the results is based on some quantitative versions of strong unique
continuation properties of solutions to elliptic equations. The methods in
\cite{AlBRVe1} has been used to prove stability for many inverse problems
(not only for second order ellipic equations \cite{MoRo1}, \cite{MoRo2})
with unknown boundaries. In particular, for the parabolic inverse problems
with unknown boundaries, such arguments have been used in Section \ref{Sec4}
of the present paper for the case of time-varying boundaries and in papers
\cite{CRoVe1}, \cite{CRoVe2}, \cite{DcRVe} for the case of not time varying
boundaries. In \cite{CRoVe1}, \cite{CRoVe2} $\kappa $ depends on $x$ only,
in \cite{DcRVe} $\kappa $ depends on $x$ ad $t$.

The strong unique continuation property in the \textit{interior} for
elliptic equations asserts that if $u$ is a solution to $div\left( \sigma
\nabla u\right) =0$, in $\Omega $, and $u$ doesn't vanish in $\Omega $ then,
for every $x_{0}\in \Omega $, there exist $C>0$ and $K\geq 0$ such that for
every $r$ sufficiently small we have%
\begin{equation}
\int\nolimits_{B_{r}\left( x_{0}\right) }u^{2}dx\geq Cr^{K}\text{.}
\label{Ellpr1}
\end{equation}%
Likewise, the strong unique continuation property at the \textit{boundary}
asserts that given a solution $u$ to the equation $div\left( \sigma \nabla
u\right) =0$, in $\Omega $, and, for instance, $u=0$ on $I\subset \partial
\Omega $, such that $u$ doesn't vanish identically in $\Omega $ then, for
very $x_{0}\in I$, there exists $C>0$ and $K\geq 0$ such that, for any
sufficiently small $r$ we have
\begin{equation}
\int\nolimits_{B_{r}\left( x_{0}\right) \cap \Omega }u^{2}dx\geq Cr^{K}%
\text{.}  \label{Ellpr2}
\end{equation}%
We emphasize that in (\ref{Ellpr1}) and (\ref{Ellpr2}) the quantities $C$
and $K$ may depend on $u$, but they do not depend on $r$. In other words, (%
\ref{Ellpr1}) and (\ref{Ellpr2}) tell us that if $u$ doesn't vanish in $%
\Omega $ then it cannot have a zero of infinite order in a point of $\Omega
\cup I$. \

The interior strong unique continuation property for second order elliptic
equations is known since 1956 by Aronzajn and Cordes papers, \cite{Ar}, \cite%
{Cor}, and since 1963, with optimal assumption regarding the regularity of $%
\sigma $, by N. Aronszajn, A. Krzywicki and J. Szarski paper \cite{AKS}. Two
suitable quantitative versions of the strong unique continuation property in
\textit{the interior} are the three sphere inequality, \cite{Ku1}, \cite%
{Lan2}, \cite{Lan3} and the doubling inequality \cite{GaLi}. The strong
unique continuation properties at \textit{the boundary} for second order
elliptic equations has been studied in the 90's by Adolfsson, Escauriaza,
Kenig, Kukavica and Wang, \cite{AE}, \cite{AEK}, \cite{KeWa}, \cite{Ku2}.

Now we sketch a crucial and significant step of the proof of \cite{AlBRVe1}
when $\mathcal{B}$ in (\ref{Ellpr}) is the boundary Dirichlet operator. Let $%
\Omega _{1}$, $\Omega _{2}$ be two domains with sufficiently smooth
boundaries. For the sake of simplicity let us assume that $\Omega _{1}$, $%
\Omega _{2}$ are two convex domains and $\left\vert \Omega _{i}\right\vert
\leq 1$, $i=1,2$, where $\left\vert \Omega _{i}\right\vert $ denotes the
Lebesgue measure of $\Omega _{i}$. Assume $\partial \Omega _{i}=A\cup I_{i}$%
, Int$_{\partial \Omega }\left( I_{i}\right) \cap $Int$_{\partial \Omega
}\left( A\right) =\varnothing $ , $I_{i}\cap A=\partial A=\partial I_{i}$ $%
i=1,2$. Let $u_{i}$ be the solution to (\ref{Ellpr}) when $\Omega =\Omega
_{i}$, $i=1,2$. Assume that $\left\Vert \sigma \nabla u_{1}\cdot \nu -\sigma
\nabla u_{2}\cdot \nu \right\Vert _{L^{2}\left( \Sigma \times \left(
0,T\right) \right) }\leq \varepsilon $.

\noindent By standard estimates for the Cauchy problem and propagation
smallness estimates we have
\begin{equation}
\left\Vert u_{1}-u_{2}\right\Vert _{L^{\infty }\left( G\right) }\leq \eta
\left( \varepsilon \right) \text{,}  \label{Ellpr3}
\end{equation}%
where $G=\Omega _{1}\cap \Omega _{2}$ such that $\Sigma \subset \overline{G}$
and $\eta \left( \varepsilon \right) $ is an infinitesimal function as $%
\varepsilon $ tend to $0$, $\eta \left( \varepsilon \right) $ depends on $%
\left\Vert f\right\Vert _{H^{1/2\left( A\right) }}$ and on the a priori
information on $\Omega _{i}$, $i=1,2$. By the maximum principle and (\ref%
{Ellpr3}) we have%
\begin{equation}
\left\Vert u_{i}\right\Vert _{L^{2}\left( \Omega _{i}\backslash G\right)
}\leq \eta \left( \varepsilon \right) \text{, }i=1,2\text{.}  \label{Ellpr4}
\end{equation}%
Let $d$ be the Hausdorff distance between $\Omega _{1}$ and $\Omega _{2}$.
Without any restriction we may assume that there exists $x_{0}\in I_{1}$
such that $d=dist\left( x_{0},\Omega _{2}\right) $. Now, by (\ref{Ellpr4})
and by using a quantitative version of (\ref{Ellpr2}) we get
\begin{equation}
d\leq \left( \frac{\eta \left( \varepsilon \right) }{C}\right) ^{1/K}\text{.}
\label{Ellpr5}
\end{equation}%
Finally, the quantities $C$ and $K$ which appear in inequality (\ref{Ellpr5}%
) can be estimated respectively from below and from above in terms of the
quantity $\frac{\left\Vert f\right\Vert _{H^{1/2\left( A\right) }}}{%
\left\Vert f\right\Vert _{L^{2}\left( A\right) }}$.

In the sequel of the introduction we shall outline the main steps to derive
some quantitative estimates of unique continuation for parabolic equations
useful to prove stability results for inverse parabolic problems with
unknown boundaries.

Such estimates can be divided in two large classes:

\noindent (i) Stability estimates for noncharacteristic Cauchy problems and
quantitative versions of weak unique continuation properties for solutions
to parabolic equations;

\noindent (ii) Quantitative versions of strong unique continuation
properties at the interior and at the boundary, for solutions to parabolic
equations.

To make clear the exposition we introduce some notation. First we assume
that the leading coefficients of equations does not depend on $u$, to
emphasize such a condition we replace, in the parabolic operator of (\ref%
{In1}), the matrix $\kappa $ by $a$ and we assume that $a$ depends on $x$
and $t$ only, $a$ is a real symmetric matrix and satisfies a uniform
ellipticity condition. We denote by $b$ and $c$ two measurable functions
from $\mathbb{R}^{n+1}$ with value in $\mathbb{R}^{n}$ and $\mathbb{R}$
respectively. We assume that $b$ and $c$ are bounded. We denote by $L$ the
parabolic operator
\begin{equation*}
Lu=div\left( a\left( x,t\right) \nabla u\right) -\partial _{t}u+b\left(
x,t\right) \cdot \nabla u+c\left( x,t\right) u.
\end{equation*}%
Let $D$ be a bounded domain in $\mathbb{R}^{n}$ and let $\Gamma $ be a
sufficiently smooth portion of the boundary $\partial D$. Let $g_{1}$, $%
g_{2} $ be given functions defined on $\Gamma \times \left( 0,T\right) $. As
above, $T$ is a positive number. The \textit{noncharacteristic Cauchy problem%
} for $Lu=0$ can be formulated as follows. Determine $u$ such that
\begin{equation}
\left\{
\begin{array}{c}
Lu=0\text{, in }D\times \left( 0,T\right) \text{,} \\
u=g_{1}\text{, on }\Gamma \times \left( 0,T\right) \text{,} \\
a\nabla u\cdot \nu =g_{2}\text{, on }\Gamma \times \left( 0,T\right) \text{,}%
\end{array}%
\right.  \label{In5}
\end{equation}%
where $\nu $ is the exterior unit normal to $D$.

\noindent We say that the Cauchy problem (\ref{In5}) enjoies the uniqueness
property if $u$ vanishes whenever $g_{1}$, $g_{2}$ do. It is well known,
\cite{Ha}, \cite{LRS}, \cite{Pu1}, that the noncharacteristic Cauchy problem
is a severely ill posed problem, as a small error on the data $g_{1}$, $%
g_{2} $ may have uncontrollable effects on the solution $u$ of (\ref{In5}).
Therefore it is very important for the applications to have a \textit{%
stability estimate} for the solutions of (\ref{In5}) whenever such solutions
belong to a certain class of functions. There exists a very large
literature, for parabolic and other types of equations, on the issue of
stability, we refer to the books \cite{Is2} and \cite{LRS} for a first
introduction on the subject. The uniqueness for the Cauchy problem (\ref{In5}%
) is equivalent to, \cite{Nir}, the so called \textit{weak unique
continuation property} for operator $L$. Such a property asserts that,
denoting by $\omega $ an open subset of $D$, if $u$ is a solution to $Lu=0$
in $D\times \left( 0,T\right) $ such that $u=0$ in $\omega \times \left(
0,T\right) $ then $u=0$ in $D\times \left( 0,T\right) $. If $\Gamma $ is
smooth enough, say $C^{1,1}$, then stability estimates for the Cauchy
problem (\ref{In5}) can be derived by a quantitative version of the above
mentioned weak unique continuation property, see Section \ref{Subs3.3} of
the present paper for details. In turns, such a weak unique continuation
property, is a consequence of (but not equivalent to) the following \textit{%
spacelike strong unique continuation property} \cite{AlVe} for the operator $%
L$. Such a property asserts that, if $u$ is a solution to $Lu=0$ in $D\times
\left( 0,T\right) $, $t_{0}\in \left( 0,T\right) $ and $u\left(
.,t_{0}\right) $ doesn't vanish identically in $D$ then, for every $x_{0}\in
D$, there exist $C>0$ and $K\geq 0$ such that for every $r$, $0<r<dist\left(
x_{0},\partial D\right) $, we have
\begin{equation}
\int\nolimits_{B_{r}\left( x_{0}\right) }u^{2}\left( x,t_{0}\right) dx\geq
Cr^{K}\text{.}  \label{In10}
\end{equation}

\noindent In other words, if the the operator $L$ enjoies such a spacelike
strong unique continuation property then $u\left( .,t_{0}\right) $ either
vanishes in $D$ or it cannot have a zero of infinite order in a point of $D$%
. We emphasize that in (\ref{In10}), $C$ and $K$ may depend on $u$. The
natural quantitative versions of spacelike strong unique continuation
property are the doubling inequality for $u\left( .,t_{0}\right) $, \cite%
{EsFeVe}, and the two-sphere one-cylinder inequality. In a rough form the
latter (of which we make an extensive use in this paper) has the following
form
\begin{equation}
\left\Vert u\left( .,t_{0}\right) \right\Vert _{L^{2}\left( B_{\rho }\left(
x_{0}\right) \right) }\leq C\left\Vert u\left( .,t_{0}\right) \right\Vert
_{L^{2}\left( B_{r}\left( x_{0}\right) \right) }^{\theta }\left\Vert
u\right\Vert _{L^{2}\left( B_{R}\left( x_{0}\right) \times \left(
t_{0}-R^{2},t_{0}\right) \right) }^{1-\theta }\text{,}  \label{In11}
\end{equation}%
where $\theta =\left( C\log \frac{R}{Cr}\right) ^{-1}$, $0<r<\rho <R$, the
cylinder $B_{R}\left( x_{0}\right) \times \left( t_{0}-R^{2},t_{0}\right) $
is contained in $D\times \left( 0,T\right) $ and $C$ depends neither on $u$
nor on $r$. See Theorem \ref{Th3.4} for a precise statement.

Before proceeding, we believe of interest to dwell a bit on the regularity
of the leading coefficients of $L$ that guarantees the uniqueness and the
stability for solutions to problem (\ref{In5}). Indeed, by some examples of
Miller \cite{Mi} and Pli\v{s} \cite{Pli} in the elliptic case, it is clear
that the minimal regularity of $a$ with respect to the space variables has
to be Lipschitz continuity. On the other hand, if $a$ depends on $x$ and $t$
optimal regularity assumptions (to the author knowledge) are not known. Even
the cases $n=1$ or $n=2$ present open issues of optimal regularity for $a$.

Now we attempt to give a view of the literature and of the main results
concerning quantitative estimates of unique continuation for parabolic
equations. To this aim we collect the contributions on the subject in the
following way.

(A) \textit{Pioneering works},

(B) \textit{The case} $n=1$,

(C) \textit{Weak unique continuation property and Cauchy problem for }$n>1$,

(D) \textit{Spacelike strong unique property continuation.}

\noindent Concerning (A) and (B) we make no claim of bibliographical
completeness. A good source of references on such points is the book \cite%
{Ca}.

(A) \textit{Pioneering works. }Such works concern mainly the case where $L$
is the heat operator $\partial _{t}-\Delta $ and those paper that are
related to the investigation on the regularity properties of solutions to $%
\partial _{t}u-\Delta u=0$, \cite{Ge}, \cite{Hol}, \cite{Pu2}. In this group
of contributions we have to enclose the classical papers and books in which
the ill posed character of the Cauchy problem for parabolic equation was
investigated for the first time \cite{Ha}, \cite{Jo}, \cite{Pu1}.

(B) \textit{The case} $n=1$. This a special, but neverthless interesting
case. Of course much more than in the general case ($n\geq 1$), operational
transformation, particular tricks and even resolution formulae are available
in the case $n=1$, \cite{Ca1}, \cite{Ca2}, \cite{Col}, \cite{Hi}. In such a
case the regularity assumption, especially concerning the leading
coefficient, can be releaxed \cite{Hao}, \cite{Lan1}, \cite{KnVe}.

(C) \textit{Weak unique continuation property and Cauchy problem for }$n>1$.
The most parts of the papers of this group share the Carleman estimates
technique for proving the uniqueness and the stability estimates. The basic
idea of such a technique has been introduced in the paper \cite{Car} to
prove the uniqueness of solution of a Cauchy problem for elliptic systems in
two variables with nonanalytic coefficients. Nowadays the general theory of
Carleman estimates is presented in several books and papers \cite{FI}, \cite%
{Ho63}, \cite{Ho85}, \cite{Is1}, \cite{Is2}, \cite{Is4}, \cite{KlTi}, \cite%
{LRS}, \cite{ShAAt}, \cite{Tr1}, \cite{Zu}. In particular in \cite{Is1} a
general theory of Carleman estimates for anisotropic operators (of which $L$
is an example) has been developed. Nirenberg \cite{Nir} has proved, in a
very general context, the uniqueness for the Cauchy problem (and the weak
unique continuation property) when the entries of matrix $a$ are constants.
It is remarkable that in such a paper Nirenberg has posed the question
whether for a solution $u$ of the equation $Lu=0$ in $D\times \left(
0,T\right) $, $D$ as above, it is true that $u\left( .,t_{0}\right) =0$ in $%
\omega $, $\omega $ open subset of $D$, it implies $u\left( .,t_{0}\right)
=0 $ in $D$. In other words, the question posed by Nirenberg is whether a
weak form of the above mentioned spacelike strong unique continuation holds
true. John \cite{Jo} has proved some stability estimate for the Cauchy
problem for the operator $L$ under the same hypotheses of \cite{Nir}. Lees
and Protter, \cite{LP}, \cite{Pr}, have proved uniqueness for the Cauchy
problem, when the entries of the matrix $a$ are assumed twice continuously
differentiable. To the author's knowledge a stability estimate of H\"{o}lder
type for the Cauchy problem for the operator $L,$ under the same hypotheses
of \cite{LP}, \cite{Pr}, was proved for the first time in \cite{AmSh}, see
also \cite{Am}. Many authors have contributed to reduce the regularity
assumption on $a$ in order to obtain the weak unique continuation for
operator $L$, among those we recall \cite{FI}, \cite{SS}, the above
mentioned \cite{Is1} and the references therein. In the context of a
quantitative version of weak unique continuation properties for operator $L$
we mention the so called \textit{three cylinder inequalities}. Roughly
speaking these are estimates of type
\begin{equation}
\left\Vert u\right\Vert _{L^{2}\left( B_{\rho }\left( x_{0}\right) \times
\left( t_{0},T-t_{0}\right) \right) }\leq C\left\Vert u\right\Vert
_{L^{2}\left( B_{r}\left( x_{0}\right) \times \left( 0,T\right) \right)
}^{\gamma }\left\Vert u\right\Vert _{L^{2}\left( B_{R}\left( x_{0}\right)
\times \left( 0,T\right) \right) }^{1-\gamma }  \label{In15}
\end{equation}%
where $\gamma \in \left( 0,1\right) $, $0<r<\rho <R$, $t_{0}\in \left(
0,T/2\right) $, the cylinder $B_{R}\left( x_{0}\right) \times \left(
0,T\right) $ is contained in $D\times \left( 0,T\right) $ and $C$ depend
neither on $u$ nor on $r$, but it may depend on $t_{0}$. Such an estimate
has been proved in \cite{Gla}, \cite{VA} when $a\in C^{3}$. More recently,
\cite{EsVe}, \cite{Ve3}, the regularity assumption on $a$ (up to Lipschitz
continuity with respect to $x$ and $t$) has been reduced and and (\ref{In15}%
) has been proved with an optimal exponent $\gamma $, that is $\gamma
=\left( C\log \frac{R}{Cr}\right) ^{-1}$. Inequality (\ref{In15}) and its
version at the boundary has been used in \cite{DcRVe} to prove optimal
stability estimates for the inverse parabolic problem (Dirichlet case) with
unknown boundary in the case of not time-varying boundary.

(D) \textit{Spacelike strong unique property continuation. }Such a property
and the inequality (\ref{In11}) have been proved for the first time in 1974
by Landis and Oleinik \cite{LanO} in the case where all the coefficients of
operator $L$ do not depend on $t$. In \cite{LanO} the authors found out a
method, named by the authors "elliptic continuation tecnique", to derive
unique continuation properties and their quantitative versions by properties
and inequalities which hold true in the elliptic context. They have employed
this method also for parabolic equation of order higher than two and for
systems. The elliptic continuation technique depends strongly on the time
independence of the coefficients of the equation. Some fundamental ideas of
such technique can be traced back to the pioneering work of Ito and Yamabe
\cite{ItYa}. Roughly speaking the above method consists of the following
idea. Let $u$ be a solution to the parabolic equation
\begin{equation}
div\left( a\left( x\right) \nabla u\right) -\partial _{t}u+b\left( x\right)
\cdot \nabla u+c\left( x\right) u=0  \label{In20}
\end{equation}%
and let $t_{0}$ be fixed, then $u\left( .,t_{0}\right) $ can be continued to
a solution $U\left( x,y\right) $, $y\in \left( -\delta ,\delta \right) $, $%
\delta >0$, of an elliptic equation in the variables $x$ and $y$. In this
way some unique continuation properties (and inequalities) of solutions to
elliptic equations can be transfered to solutions to equation (\ref{In20}).
In \cite{LanO} the regularity assumption on the matrix $a$ are very strong,
but in 1990 Lin \cite{Lin} has employed the same technique to prove the
spacelike strong unique continuation properties for solution to (\ref{In20})
assuming $a$ merely Lipschitz continuous. In \cite{CRoVe1}, \cite{CRoVe2}
the elliptic continuation technique has been used to prove some two-sphere
one-cylinder inequality at the boundary (with Dirichlet or Neumann
condition) for solution to (\ref{In20}).

A very important step towards the proof of the spacelike strong unique
continuation property for parabolic equation with time dependent
coefficients consists in proving the following strong unique continuation
properties: let $u$ be a solution of $Lu=0$ in $D\times \left( 0,T\right) $,
and $\left( x_{0},t_{0}\right) \in D\times \left( 0,T\right) $ then $u\left(
.,t_{0}\right) $ vanishes whenever
\begin{equation*}
u\left( x,t\right) =O\left( \left( \left\vert x-x_{0}\right\vert
^{2}+\left\vert t-t_{0}\right\vert \right) ^{N/2}\right) \text{, for every }%
N\in \mathbb{N}\text{.}
\end{equation*}%
Such a properties has been proved for the first time by Poon \cite{Poon}
when $Lu=\partial _{t}u-\Delta u+b\left( x,t\right) \cdot \nabla u+c\left(
x,t\right) u$, $D=\mathbb{R}^{n}$ and $u$ satisfies some growth condition at
infinity. Hence in \cite{Poon} the strong unique continuation property has
been not yet proved as a local property. Such a local strong unique
continuation property has been proved, by using Carleman estimates
techniques, in \cite{EsFe} for the operator $Lu=\partial _{t}u-div\left(
a\left( x,t\right) \nabla u\right) +b\left( x,t\right) \cdot \nabla
u+c\left( x,t\right) u$ when $a$ is Lipschitz continuous with respect to the
parabolic distance and $b$, $c$ are bounded. In \cite{AlVe} the spacelike
strong unique continuation has been derived as a consequence of the
mentioned theorem of \cite{EsFe} and of the local behaviour of solutions to
parabolic equation proved in \cite{AlVe0}. In \cite{Fe} the same spacelike
strong unique continuation property has been proved as a consequence of a
refined version of the Carleman estimate proved in \cite{EsFe}. The natural
quantitative versions, that is the two-sphere one cylinder inequality and
doubling inequality on characteristic hyperplane, have been proved in \cite%
{EsFeVe}. In the present paper, Subsections \ref{Subsrq2.0} and \ref%
{SubsRq2.1}\textit{\ }are devoted to a complete proof of the two-sphere
one-cylinder inequality (at the interior and at the boundary) and a
simplified version of the proof of the Carleman estimate proved in \cite%
{EsFe} and \cite{Fe}. Quite recently the spacelike strong unique
continuation property has been proved in \cite{KoTa} with some improvement
on the regularity assumption on coefficients of operator $L$.

The two-sphere one-cylinder inequalities are the essential tools to prove
the stability estimate of logarithmic type for the inverse parabolic
problems presented at the beginning of his paper, see Section \ref{Sec4}
and, in particular Theorem \ref{Th4.1}, for details.

Finally let us give the plan of the paper.

\textit{Section \ref{Sec2}}:\textit{\ }Main Notations and Definitions.

\textit{Section \ref{Sec3}}:\textit{\ }Quantitative estimate of unique
continuation.

\noindent Such a Section is subdivided in the following five Subsections.

\noindent \textit{Subsection \ref{Subsrq2.0}}: Parabolic equations with time
independent coefficients,

\noindent \textit{Subsection \ref{SubsRq2.1}}:\textit{\ }Carleman estimate,

\noindent \textit{Subsection \ref{Subs3.1}}: Two-sphere one-cylinder
inequalities,

\noindent \textit{Subsection \ref{Subs3.2}}:\textit{\ }Smallness propagation
estimate,

\noindent \textit{Subsection \ref{Subs3.3}}:Stability estimates from Cauchy
data.

\textit{Section \ref{Sec4}}:\textit{\ }Stability estimates for Dirichlet
inverse problem with unknown time-varying boundaries.

\noindent Such a Section is subdivided in the following four Subsections.

\noindent \textit{Subsection \ref{Subs4.1}}:\textit{\ }Proof of Theorem \ref%
{Th4.1} (it s the main Theorem of \ref{Sec4}),

\noindent \textit{Subsection \ref{Subs4.2}}:\textit{\ }Proofs of
Propositions \ref{Pr4.2}, \ref{Pr4.3}, \ref{Pr4.4}, \ref{Pr4.5}, \ref{Pr4.6}
(the statements of such propositions are given in Subsection \ref{Subs4.1}),

\noindent \textit{Subsection \ref{Subs4.3}}:\textit{\ }Some extensions of
Theorem \ref{Th4.1},

\textit{Section \ref{INST}}:\textit{\ }Exponential instability,

\noindent such a Section \ref{INST} is subdivided in the following two
Subsections,

\noindent \textit{Subsubsection \ref{INSTT}}: Exponential instability of the
Dirichlet inverse problem with time-varying unknown boundary,

\noindent \textit{Subsubsection \ref{INSTT}}: Stability properties of the
Dirichlet inverse problem with unknown boundary independent on time.

\section{Main Notations and Definitions\label{Sec2}}

For every $x\in \mathbb{R}^{n}$, $n\geq 2$, $x=\left( x_{1},...,x_{n}\right)
$, we shall set $x=\left( x^{\prime },x_{n}\right) $, where $x^{\prime
}=\left( x_{1},...,x_{n-1}\right) \in \mathbb{R}^{n-1}$. We shall use $%
X=\left( x,t\right) $ to denote a point in $\mathbb{R}^{n+1}$, where $x\in
\mathbb{R}^{n}$ and $t\in \mathbb{R}$. For every $x\in \mathbb{R}^{n}$, $%
X=\left( x,t\right) \in \mathbb{R}^{n+1}$ we shall set $\left\vert
x\right\vert =\left( \sum\limits_{i=1}^{n}x_{i}^{2}\right) ^{1/2}$ and $%
\left\vert X\right\vert =\left( \left\vert x\right\vert ^{2}+\left\vert
t\right\vert \right) ^{1/2}$. Let $r$ be a positive number. For every $%
x_{0}\in \mathbb{R}^{n}$ we shall denote by $B_{r}\left( x_{0}\right)
=\left\{ x\in \mathbb{R}^{n}:\text{ }\left\vert x-x_{0}\right\vert
<r\right\} $ (the $n$-dimensional open ball of radius $r$ centered at $x_{0}$%
) and $B_{r}^{\prime }\left( x_{0}^{\prime }\right) =\left\{ x\in \mathbb{R}%
^{n-1}:\text{ }\left\vert x^{\prime }-x_{0}^{\prime }\right\vert <r\right\} $
(the $\left( n-1\right) $-dimensional open ball of radius $r$ centered at $%
x_{0}^{\prime }$). We set generally $B_{r}=B_{r}\left( 0\right) $ and $%
B_{r}^{\prime }=B_{r}^{\prime }\left( 0\right) $. We shall denote by $%
B_{r}^{+}=\left\{ x\in B_{r}:x_{n}>0\right\} $. For every $X_{0}\in \mathbb{R%
}^{n+1}$ we shall set $Q_{r}\left( X_{0}\right) =B_{r}\left( x_{0}\right)
\times \left( t_{0}-r^{2},t_{0}\right) $.

Let $I$ be an interval of $\mathbb{R}$ and let $\left\{ D\left( t\right)
\right\} _{t\in I}$ be a family of subsets of $\mathbb{R}^{n}$, we shall
denote $D\left( I\right) =\bigcup\limits_{t\in I}D\left( t\right) \times
\left\{ t\right\} $.

Given a sufficiently smooth function $u$ of $x$ and $t$, we shall denote $%
\partial _{i}u=\dfrac{\partial u}{\partial x_{i}}$, $\partial _{ij}^{2}u=%
\dfrac{\partial ^{2}u}{\partial x_{i}\partial x_{j}}$, $i,j=1,...,n$ and $%
\partial _{t}u=\dfrac{\partial u}{\partial t}$. For a multi-index $\beta
=\left( \beta _{1},...,\beta _{n}\right) $, $\beta _{i}\in \mathbb{N\cup }%
\left\{ 0\right\} $, $i=1,...,n$ and $k\in \mathbb{N\cup }\left\{ 0\right\} $%
, we shall denote, as usual, $\partial _{x}^{\beta }\partial _{t}^{k}u=%
\dfrac{\partial ^{\left\vert \beta \right\vert +k}u}{\partial x_{1}^{\beta
_{1}}...\partial x_{n}^{\beta _{n}}\partial t^{k}}$, where $\left\vert \beta
\right\vert =\sum\limits_{i=1}^{n}\beta _{i}$. Also, we shall write $\nabla
\left( .\right) =\nabla _{x}\left( .\right) $, $div\left( .\right)
=div_{x}\left( .\right) $.

Let $D$ be a subset of $\mathbb{R}^{n+1}$, $f$ a function defined on $D$ and
$\alpha \in \left( 0,1\right] $, we shall set
\begin{equation*}
\left[ f\right] _{\alpha ;D}=\sup \left\{ \frac{\left\vert f\left(
x,t\right) -f\left( y,\tau \right) \right\vert }{\left( \left\vert
x-y\right\vert ^{2}+\left\vert t-\tau \right\vert \right) ^{\alpha /2}}:%
\text{ }\left( x,t\right) ,\left( y,\tau \right) \in D\text{, }\left(
x,t\right) \neq \left( y,\tau \right) \right\} \text{.}
\end{equation*}%
If $\alpha \in \left( 0,2\right] $ we shall set
\begin{equation*}
\left\langle f\right\rangle _{\alpha ;D}=\sup \left\{ \frac{\left\vert
f\left( x,t\right) -f\left( x,\tau \right) \right\vert }{\left\vert t-\tau
\right\vert ^{\alpha /2}}:\text{ }\left( x,t\right) ,\left( x,\tau \right)
\in D\text{, }t\neq \tau \right\} \text{.}
\end{equation*}%
Let $k$ be a positive integer number, $D$ an open subset of $\mathbb{R}%
^{n+1} $, $f$ a sufficiently smooth function and $\alpha \in \left( 0,1%
\right] $. We shall denote by
\begin{equation*}
\left[ f\right] _{k+\alpha ;D}=\sum\limits_{\left\vert \beta \right\vert
+2j=k}\left[ \partial _{x}^{\beta }\partial _{t}^{j}f\right] _{\alpha ;D}%
\text{ ,}
\end{equation*}%
\begin{equation*}
\left\langle f\right\rangle _{k+\alpha ;D}=\sum\limits_{\left\vert \beta
\right\vert +2j=k-1}\left\langle \partial _{x}^{\beta }\partial
_{t}^{j}f\right\rangle _{1+\alpha ;D}\text{.}
\end{equation*}%
If $\alpha \in \left( 0,1\right] $ and $\left[ f\right] _{\alpha ;D}$ is
finite, we shall say that $f$ is \textit{H\"{o}lder continuous} (\textit{%
Lipschitz continuous} whenever $\alpha =1$) in $D$ or that $f$ belongs to $%
C^{0,\alpha }$. Let $k$ be a positive integer number, $\alpha \in \left( 0,1%
\right] $ and let $D$ be an open subset of $\mathbb{R}^{n+1}$, we shall say
that $f$ belongs to the class $C^{k,\alpha }$ whenever for every multi-index
$\beta $ and every nonnegative integer number $j$ such that $\left\vert
\beta \right\vert +2j\leq k$ there exist the derivatives $\partial
_{x}^{\beta }\partial _{t}^{j}f$ and the quantities $\sup\limits_{D}\left%
\vert \partial _{x}^{\beta }\partial _{t}^{j}f\right\vert $, $\left[ f\right]
_{k+\alpha ;D}$ and $\left\langle f\right\rangle _{k+\alpha ;D}$ are finite.
If $f$ is a function which doesn't depend explicitly on $t$ we shall
continue to use the definition above, more precisely for a function $%
f:\Omega \rightarrow \mathbb{R}$, where $\Omega \subset \mathbb{R}^{n}$, we
shall say that $f\in C^{k,\alpha }\left( \Omega \right) $ whenever,
considering the function $\widetilde{f}:\Omega \times \mathbb{R}\rightarrow
\mathbb{R}$, $\widetilde{f}\left( x,t\right) =f\left( x\right) $ for every $%
\left( x,t\right) \in \Omega \times \mathbb{R}$, we have $\widetilde{f}\in
C^{k,\alpha }\left( \Omega \times \mathbb{R}\right) $.

Let $X_{0}=\left( x_{0},t_{0}\right) $, $Y_{0}=\left( y_{0},\tau _{0}\right)
$ be two points of $\mathbb{R}^{n+1}$, we shall say that $S:\mathbb{R}%
^{n+1}\rightarrow \mathbb{R}^{n+1}$ is a \textit{a rigid transformation of
space coordinates under which we have} $X_{0}\equiv Y_{0}$ if $S\left(
X\right) =\left( \sigma \left( x\right) ,t-t_{0}+\tau _{0}\right) $ where $%
\sigma $ is an isometry of $\mathbb{R}^{n}$ such that $\sigma \left(
x_{0}\right) =y_{0}$.

\begin{definition}
\label{Def2.1}Let $\Omega $ be a domain in $\mathbb{R}^{n+1}$. Given a
positive integer number $k$ and $\alpha \in \left( 0,1\right] $, we shall
say that a portion $\Gamma $ of $\partial \Omega $ is of class $C^{k,\alpha
} $ with constants $\rho _{0}$, $E>0$ if for any $X_{0}\in \Gamma $, there
exists a rigid tranformation of space coordinates under which we have $%
X_{0}\equiv 0$ and
\begin{equation*}
\Omega \cap \left( B_{\rho _{0}}^{\prime }\left( 0\right) \times \left(
-\rho _{0}^{2},\rho _{0}^{2}\right) \right) =\left\{ X\in B_{\rho
_{0}}^{\prime }\left( 0\right) \times \left( -\rho _{0}^{2},\rho
_{0}^{2}\right) :\text{ }x_{n}>\varphi \left( x^{\prime },t\right) \right\}
\text{,}
\end{equation*}%
where $\varphi \in C^{k,\alpha }\left( B_{\rho _{0}}^{\prime }\left(
0\right) \times \left( -\rho _{0}^{2},\rho _{0}^{2}\right) \right) $
satisfying
\begin{equation*}
\varphi \left( 0,0\right) =\left\vert \nabla _{x^{\prime }}\varphi \left(
0,0\right) \right\vert =0
\end{equation*}%
and
\begin{equation*}
\left\Vert \varphi \right\Vert _{C^{k,\alpha }\left( B_{\rho _{0}}^{\prime
}\left( 0\right) \times \left( -\rho _{0}^{2},\rho _{0}^{2}\right) \right)
}\leq E\rho _{0}.
\end{equation*}
\end{definition}

\begin{remark}
\label{Rem2.1}We have chosen to normalize all norms in such a way that their
terms are dimensional homogeneous and coincide with the standard definition
when $\rho _{0}=1$. For instance, for any $\varphi \in C^{k,\alpha }\left(
B_{\rho _{0}}^{\prime }\left( 0\right) \times \left( -\rho _{0}^{2},\rho
_{0}^{2}\right) \right) $ we set
\begin{eqnarray*}
\left\Vert \varphi \right\Vert _{C^{k,\alpha }\left( B_{\rho _{0}}^{\prime
}\left( 0\right) \times \left( -\rho _{0}^{2},\rho _{0}^{2}\right) \right) }
&=&\sum\limits_{l=0}^{k}\rho _{0}^{k}\sum\limits_{\left\vert \beta
\right\vert +2j=l}\left\Vert \partial _{x}^{\beta }\partial _{t}^{j}\varphi
\right\Vert _{L^{\infty }\left( B_{\rho _{0}}^{\prime }\left( 0\right)
\times \left( -\rho _{0}^{2},\rho _{0}^{2}\right) \right) }+ \\
&&\rho _{0}^{k+\alpha }\left( \left\langle \varphi \right\rangle _{k+\alpha
;B_{\rho _{0}}^{\prime }\left( 0\right) \times \left( -\rho _{0}^{2},\rho
_{0}^{2}\right) }+\left[ \varphi \right] _{k+\alpha ;B_{\rho _{0}}^{\prime
}\left( 0\right) \times \left( -\rho _{0}^{2},\rho _{0}^{2}\right) }\right)
\text{ .}
\end{eqnarray*}
Similarly we shall set
\begin{equation*}
\left\Vert u\right\Vert _{L^{2}\left( \Omega \right) }=\rho _{0}^{-\left(
n+1\right) /2}\left( \int\nolimits_{\Omega }u^{2}dX\right) ^{1/2}\text{,}
\end{equation*}
where $dX=dxdt$.
\end{remark}

For any matrix $M$ in $G\left( \left( 0,T\right) \right) $, we shall denote
its transposed by $M^{\ast }$ and its trace by tr$\left( M\right) $.

We shall fix the space dimension $n\geq 1$ throughout the paper. Therefore
we shall omit the dependence of various quantities on $n$

We shall use the the letters $C$, $C_{0}$, $C_{1}$ ... to denote constants.
The value of the constants may change from line to line, but we shall
specified their dependence everywhere they appear.

\section{Quantitative estimate of unique continuation\label{Sec3}}

In Section \ref{Subsrq2.0} we present the so-called elliptic continuation
technique due to Landis and Oleinik, \cite{LanO}. In Section \ref{SubsRq2.1}
we state and prove the Carleman estimate (Theorem \ref{ThRq9} below) for
parabolic operators proved in \cite{EsFe}, \cite{Fe}, see also \cite{EsFeVe}%
. In order to prove such an inequality we adapt to the case of variable
coefficients the approach used in \cite{EsSeSv} for the heat operator. In
Section \ref{Subs3.1}\ we apply the Carleman estimate to prove the
two-sphere one-cylinder inequality at the interior and at the (time varying)
boundary. In Section \ref{Subs3.2} we shall prove a smallness propagation
estimates for solutions to parabolic equations on a characteristic
hyperplane up to a time varying portion of a Lipschitz boundary. Finally, in
Section \ref{Subs3.3} we prove some sharp stability estimates for the Cauchy
problem for parabolic equations.

\subsection{Parabolic equations with time independent coefficients\label%
{Subsrq2.0}}

In this section we present some quantitative estimates of unique
continuation for solutions to parabolic equations whose coefficients do not
depend on $t$. The method for deriving such estimates has been introduced in
\cite{LanO} and, concerning the smoothness assumption on the coefficients,
has been improved in \cite{Lin}. In what follows we give an outline of the
above method and we limit ourselves to present detailed proofs only of the
main points of the method, we omit the most technical proofs for which we
refer to the quoted literature and especially to \cite{CRoVe1}. Throughout
this section, for any positive number $r$, we shall denote by $\widetilde{B}%
_{r}$ the $\left( n+1\right) $-dimensional open ball of radius $r$ centered
at $0$ and we shall denote by $\widetilde{B}_{r}^{+}=\left\{ x\in \widetilde{%
B}_{r}:x_{n+1}>0\right\} $.

In order to simplify the exposition we consider the equation
\begin{equation}
Lu:=div\left( a\left( x\right) \nabla u\right) -\partial _{t}u=0\text{ , in }%
B_{1}\times \left( 0,1\right] \text{,}  \label{rq.1.10}
\end{equation}%
where $a\left( x\right) =\left\{ a^{ij}\left( x\right) \right\} _{i,j=1}^{n}$
is a symmetric $n\times n$ matrix whose entries are real valued functions.
When $\xi \in \mathbb{R}^{n}$ and $x$, $y\in \mathbb{R}^{n}$ we assume that
\begin{equation}
\lambda \left\vert \xi \right\vert ^{2}\leq
\sum\limits_{i,j=1}^{n}a^{ij}\left( x\right) \xi _{i}\xi _{j}\leq \lambda
^{-1}\left\vert \xi \right\vert ^{2}\text{ }  \label{rq.1.15}
\end{equation}%
and
\begin{equation}
\left( \sum\limits_{i,j=1}^{n}\left( a^{ij}\left( x\right) -a^{ij}\left(
y\right) \right) ^{2}\right) ^{1/2}\leq \Lambda \left\vert x-y\right\vert
\text{ ,}  \label{rq.1.20}
\end{equation}%
where $\lambda $ and $\Lambda $ are positive numbers with $\lambda \in
\left( 0,1\right] $.

Let $u\in H^{2,1}\left( B_{1}\times \left( 0,1\right) \right) $ be a
solution to (\ref{rq.1.10}).Denoting by $\overline{u}$ the following
extension of $u$%
\begin{equation}
\overline{u}=\left\{
\begin{array}{c}
u\left( x,t\right) \text{, if }\left( x,t\right) \in B_{1}\times \left(
0,1\right) \text{,} \\
-3u\left( x,2-t\right) +4u\left( x,\frac{3}{2}-\frac{1}{2}t\right) \text{,
if }\left( x,t\right) \in B_{1}\times \left[ 1,2\right) \text{,}%
\end{array}%
\right.  \label{rq.1.25}
\end{equation}%
we have $\overline{u}\in H^{2,1}\left( B_{1}\times \left( 0,2\right) \right)
$ and, \cite{Ev},
\begin{equation}
\left\Vert \overline{u}\right\Vert _{H^{2,1}\left( B_{1}\times \left(
0,2\right) \right) }\leq C\left\Vert u\right\Vert _{H^{2,1}\left(
B_{1}\times \left( 0,1\right) \right) }\text{,}  \label{rq.1.30}
\end{equation}%
where $C$ is an absolute constant.

\noindent Let $\eta \in C^{1}\left( \left[ 0,+\infty \right) \right) $ be a
function satisfying $0\leq \eta \leq 1$, $\eta =1$ in $\left[ 0,1\right] $, $%
\eta =0$ in $\left[ 2,+\infty \right) $ and $\left\vert \eta ^{\prime
}\right\vert \leq c$, in $\left[ 1,2\right] $. Moreover, denoting by $%
\widetilde{u}$ the trivial extension of $\eta \left( t\right) \overline{u}%
\left( x,t\right) $ to $\left( 0,+\infty \right) $ (i.e. $\widetilde{u}%
\left( x,t\right) =0$ if $t\geq 2$), let us denote by $u_{1}$, $u_{2}$ the
weak solutions to the following initial-boundary value problems respectively
\begin{equation}
\left\{
\begin{array}{c}
Lu_{1}=0\text{ , in }B_{1}\times \left( 0,+\infty \right) \text{,} \\
u_{1}=0\text{ , on }\partial B_{1}\times \left( 0,+\infty \right) \text{,}
\\
u_{1}\left( .,0\right) =u\left( .,0\right) \text{ , in }B_{1}%
\end{array}
\right.  \label{rq.1.35-45}
\end{equation}
and
\begin{equation}
\left\{
\begin{array}{c}
Lu_{2}=0\text{ , in }B_{1}\times \left( 0,+\infty \right) \text{,} \\
u_{2}=\widetilde{u}\text{ , on }\partial B_{1}\times \left( 0,+\infty
\right) \text{,} \\
u_{2}\left( .,0\right) =0\text{ , in }B_{1}\text{ .}%
\end{array}
\right.  \label{rq.1.50-60}
\end{equation}
Since $u_{1}+u_{2}=u$ on $\partial B_{1}\times \left( 0,1\right] $ and $%
\left( u_{1}+u_{2}\right) \left( .,0\right) =u\left( .,0\right) $ in $B_{1}$%
, by uniqueness theorem for initial-boundary value problem for parabolic
equations we have
\begin{equation}
u_{1}\left( x,1\right) +u_{2}\left( x,1\right) =u\left( x,1\right) \text{ ,
for every }x\in B_{1}\text{.}  \label{rq.1.65}
\end{equation}

In what follows we shall denote by $C_{\mathcal{P}}$ the constant appearing
in the following Poincar\'{e} inequality
\begin{equation*}
\int\nolimits_{B_{1}}f^{2}\left( x\right) dx\leq C_{\mathcal{P}%
}\int\nolimits_{B_{1}}\left\vert \nabla f\left( x\right) \right\vert ^{2}dx%
\text{ , for every }f\in H_{0}^{1}\left( B_{1}\right) \text{ ,}
\end{equation*}
where we recall that $C_{\mathcal{P}}\leq 1$, \cite{GTr}, and $C_{\mathcal{P}%
}=\frac{1}{k_{0}^{2}}$, where $k_{0}$ is the smallest positive root of the
Bessel function of first kind $J_{\frac{n-2}{2}}$ , \cite{CHi}. Let us
denote
\begin{equation}
b_{1}=\frac{\lambda }{C_{\mathcal{P}}}\text{ ,}  \label{rq.1.70}
\end{equation}
and
\begin{equation}
H=\sup\limits_{t\in \left[ 0,1\right] }\left\Vert u\left( .,t\right)
\right\Vert _{H^{1}\left( B_{1}\right) }\text{ .}  \label{rq.1.75}
\end{equation}

\begin{proposition}
\label{Prrq.1}Let $u$ and $u_{2}$ be as above. We have
\begin{equation}
\left\Vert u_{2}\left( .,t\right) \right\Vert _{H^{1}\left( B_{1}\right)
}\leq C_{1}He^{-b_{1}\left( t-2\right) _{+}}\text{ , for every }t\in \left(
0,+\infty \right) \text{,}  \label{rq.1.80}
\end{equation}
($\left( t-2\right) _{+}$ is equal to $t-2$ if $t\geq 2$ and it is equal to $%
0$ if $t<2$) where $C_{1}$, $C_{1}>1$, depends on $\lambda $ and $\Lambda $
only.
\end{proposition}

\textbf{Proof.} Let
\begin{equation*}
v=u_{2}-\widetilde{u}\text{ ,}
\end{equation*}%
where $\widetilde{u}$ is the function defined above. Denoting $F=-L%
\widetilde{u}$, we have
\begin{equation}
\left\{
\begin{array}{c}
Lv=F\text{ , in }B_{1}\times \left( 0,+\infty \right) \text{,} \\
v=0\text{ , on }\partial B_{1}\times \left( 0,+\infty \right) \text{,} \\
v\left( .,0\right) =-u\left( .,0\right) \text{ , in }B_{1}\text{.}%
\end{array}%
\right.  \label{rq.1.85-95}
\end{equation}%
We claim
\begin{equation}
\left\Vert u_{2}\left( .,t\right) \right\Vert _{L^{2}\left( B_{1}\right)
}\leq CH\text{ , for every }t\in \left[ 0,2\right] \text{,}  \label{rq.1.100}
\end{equation}%
\begin{equation}
\left\Vert \nabla u_{2}\left( .,t\right) \right\Vert _{L^{2}\left(
B_{1}\right) }\leq CH\text{ , for every }t\in \left[ 0,2\right] \text{,}
\label{rq.1.105}
\end{equation}%
where $C$ depends on $\lambda $ and $\Lambda $ only.

In order to prove (\ref{rq.1.100}) let us multiply the first equation in (%
\ref{rq.1.85-95}) by $v$ and integrate over $B_{1}\times \left( 0,t\right) $%
. We get, for every $t\in \left[ 0,2\right] $,
\begin{equation*}
\int\nolimits_{B_{1}}v^{2}\left( x,t\right) dx\leq 2\left\Vert F\right\Vert
_{L^{2}\left( B_{1}\times \left( 0,2\right) \right)
}^{2}+\int\nolimits_{B_{1}}u^{2}\left( x,0\right)
dx+2\int\nolimits_{0}^{t}d\tau \int\nolimits_{B_{1}}v^{2}\left( x,\tau
\right) dx\text{.}
\end{equation*}%
By applying Gronwall inequality we have
\begin{equation}
\int\nolimits_{B_{1}}v^{2}\left( x,t\right) dx\leq e^{4}\left( 2\left\Vert
F\right\Vert _{L^{2}\left( B_{1}\times \left( 0,2\right) \right)
}^{2}+H^{2}\right) \text{ , for every }t\in \left[ 0,2\right] \text{.}
\label{rq.1.110}
\end{equation}%
Now by (\ref{rq.1.30}) and regularity estimate for parabolic equations we
have
\begin{equation}
\left\Vert F\right\Vert _{L^{2}\left( B_{1}\times \left( 0,2\right) \right)
}\leq CH,  \label{rq.1.115}
\end{equation}%
where $C$ depends on $\lambda $ and $\Lambda $ only. From this inequality
and (\ref{rq.1.110}) we obtain (\ref{rq.1.100}).

\noindent In order to prove (\ref{rq.1.105}) let us multiply the first
equation in (\ref{rq.1.85-95}) by $\partial _{t}v$ and integrate over $%
B_{1}\times \left( 0,t\right) $. We obtain
\begin{equation*}
\lambda \int\nolimits_{B_{1}}\left\vert \nabla v\left( x,t\right)
\right\vert ^{2}dx\leq \lambda ^{-1}\int\nolimits_{B_{1}}\left\vert \nabla
u\left( x,0\right) \right\vert ^{2}dx+2\left\Vert F\right\Vert _{L^{2}\left(
B_{1}\times \left( 0,2\right) \right) }^{2}\text{,}
\end{equation*}%
so, by (\ref{rq.1.75}) and (\ref{rq.1.115}) we have (\ref{rq.1.105}). Claims
are proved.

Now let us denote by $\mu _{k}$, $k\in \mathbb{N}$, the decreasing sequence
of eigenvalues associated to the problem
\begin{equation}
\left\{
\begin{array}{c}
div\left( a\nabla \varphi \right) =\mu \varphi \text{, in }B_{1}\text{,} \\
\varphi =0\text{, on }\partial B_{1}%
\end{array}
\right.  \label{rq.1.116}
\end{equation}
and by $\varphi _{k}$, $k\in \mathbb{N}$, the corresponding eigenfunctions
normalized by
\begin{equation*}
\int\nolimits_{B_{1}}\varphi _{k}^{2}\left( x\right) dx=1\text{, }k\in
\mathbb{N}\text{.}
\end{equation*}
We have
\begin{equation}
0>-b_{1}\geq \mu _{1}\geq \mu _{2}\geq ...\geq \mu _{k}\geq ...\text{ .}
\label{rq.1.120}
\end{equation}
Since $u_{2}=0$ , on $\partial B_{1}\times \left[ 2,+\infty \right) $, we
have
\begin{equation}
u_{2}\left( x,t\right) =\sum\limits_{k=1}^{\infty }\beta _{k}\varphi
_{k}\left( x\right) e^{\mu _{k}\left( t-2\right) }\text{ , for every }t\geq 2%
\text{,}  \label{rq.1.125}
\end{equation}
where
\begin{equation*}
\beta _{k}=\int\nolimits_{B_{1}}u_{2}\left( x,2\right) \varphi _{k}\left(
x\right) dx\text{, }k\in \mathbb{N}\text{.}
\end{equation*}
Since $u_{2}\left( .,2\right) =v\left( .,2\right) $ in $B_{1}$, by (\ref%
{rq.1.100}), (\ref{rq.1.120}) and (\ref{rq.1.125}) we have
\begin{equation}
\int\nolimits_{B_{1}}u_{2}^{2}\left( x,t\right)
dx=\sum\limits_{k=1}^{\infty }\beta _{k}^{2}e^{2\mu _{k}\left( t-2\right) }
\label{rq.1.126}
\end{equation}
\begin{equation*}
\leq e^{-2b_{1}\left( t-2\right) }\int\nolimits_{B_{1}}u_{2}^{2}\left(
x,2\right) dx\leq CH^{2}e^{-2b_{1}\left( t-2\right) }\text{, for every }%
t\geq 2\text{,}
\end{equation*}
where $C$ depends on $\lambda $ and $\Lambda $ only.

\noindent Moreover, for every $t\geq 2$, we have
\begin{equation}
\int\nolimits_{B_{1}}a\left( x\right) \nabla u_{2}\left( x,t\right) \cdot
\nabla u_{2}\left( x,t\right) dx  \label{rq.1.130}
\end{equation}
\begin{equation*}
=-\int\nolimits_{B_{1}}\partial _{t}u_{2}\left( x,t\right) u_{2}\left(
x,t\right) dx=\sum\limits_{k=1}^{\infty }\left\vert \mu _{k}\right\vert
\beta _{k}^{2}e^{2\mu _{k}\left( t-2\right) }\text{.}
\end{equation*}
By choosing $t=2$ in (\ref{rq.1.130}) and using (\ref{rq.1.105}) we get
\begin{equation}
\sum\limits_{k=1}^{\infty }\left\vert \mu _{k}\right\vert \beta
_{k}^{2}\leq CH^{2}\text{, }  \label{rq.1.135}
\end{equation}
where $C$ depends on $\lambda $ and $\Lambda $ only.

\noindent By (\ref{rq.1.120}), (\ref{rq.1.130}) and (\ref{rq.1.135}) we
derive
\begin{equation}
\int\nolimits_{B_{1}}\left\vert \nabla u_{2}\left( x,t\right) \right\vert
^{2}dx\leq CH^{2}e^{-2b_{1}\left( t-2\right) }\text{, for every }t\geq 2%
\text{,}  \label{rq.1.140}
\end{equation}
where $C$ depends on $\lambda $ and $\Lambda $ only.

Finally, from (\ref{rq.1.100}), (\ref{rq.1.105}), (\ref{rq.1.126}) and (\ref%
{rq.1.140}) we get (\ref{rq.1.80}).$\blacksquare $

\bigskip

Let us still denote by $u_{2}$ the extension by $0$ of $u_{2}$ to $%
B_{1}\times \mathbb{R}$ and let us consider the Fourier tranform of $u_{2}$
with respect to the $t$ variable
\begin{equation}
\widehat{u}_{2}\left( x,\mu \right) =\int\nolimits_{-\infty }^{+\infty
}e^{-i\mu t}u_{2}\left( x,t\right) dt=\int\nolimits_{0}^{+\infty }e^{-i\mu
t}u_{2}\left( x,t\right) dt\text{ , }\mu \in \mathbb{R}\text{ .}
\label{rq.1.145}
\end{equation}
We have that $\widehat{u}_{2}$ satisfies
\begin{equation}
div\left( a\left( x\right) \nabla \widehat{u}_{2}\left( x,\mu \right)
\right) -i\mu \widehat{u}_{2}\left( x,\mu \right) =0\text{ , if }\left(
x,\mu \right) \in B_{1}\times \mathbb{R}\text{.}  \label{rq.1.150}
\end{equation}

\begin{proposition}
\label{Prrq.2}Let $\widehat{u}_{2}$ be as above. We have, for every $\mu \in
\mathbb{R}$,
\begin{equation}
\left\Vert \widehat{u}_{2}\left( .,\mu \right) \right\Vert _{L^{2}\left(
B_{1/2}\right) }\leq cC_{1}\lambda ^{-1}H\left( 2+\frac{1}{b_{1}}\right)
e^{-\left\vert \mu \right\vert ^{1/2}\delta }\text{ ,}  \label{rq.1.155}
\end{equation}
where $c$ is an absolute constant, $C_{1}$ is the constant that appears in (%
\ref{rq.1.80}) and $\delta $ is given by
\begin{equation}
\delta =\frac{\lambda }{8\pi e}\text{ .}  \label{rq.1.160}
\end{equation}
\end{proposition}

\textbf{Proof.} Let us denote, for every $\mu ,\xi \in \mathbb{R}$, $x\in
B_{1}$%
\begin{equation*}
v\left( x,\xi ;\mu \right) =e^{i\left\vert \mu \right\vert ^{1/2}\xi }%
\widehat{u}_{2}\left( x,\mu \right) \text{ .}
\end{equation*}%
For every $\mu \in \mathbb{R\smallsetminus }\left\{ 0\right\} $, the
function $v\left( .,.;\mu \right) $ solves the uniformly elliptic equation
\begin{equation}
div\left( a\left( x\right) \nabla v\left( x,\xi ;\mu \right) \right) +i\text{%
sgn}\left( \mu \right) \partial _{\xi }^{2}v\left( x,\xi ;\mu \right) =0%
\text{, in }B_{1}\times \mathbb{R}\text{.}  \label{rq.1.165}
\end{equation}%
Let $k\in \mathbb{N}$ and denote by $r_{m}=1-\frac{m}{2k}$, for every $m\in
\left\{ 0,1,...,k\right\} $. Moreover set
\begin{equation*}
h_{m}\left( s\right) =\left\{
\begin{array}{c}
0\text{, if }\left\vert s\right\vert >r_{m}\text{,} \\
\frac{1}{2}\left( 1+\cos 2\pi k\left( r_{m+1}-s\right) \right) \text{ , if }%
r_{m+1}\leq \left\vert s\right\vert \leq r_{m}\text{,} \\
1\text{, }\left\vert s\right\vert <r_{m+1}\text{,}%
\end{array}%
\right.
\end{equation*}%
\begin{equation*}
\eta _{m}\left( x,\xi \right) =h_{m}\left( \left\vert x\right\vert \right)
h_{m}\left( \left\vert \xi \right\vert \right)
\end{equation*}%
and
\begin{equation*}
v_{m}\left( x,\xi ;\mu \right) =\partial _{\xi }^{m}v\left( x,\xi ;\mu
\right) \text{ , }m\in \left\{ 0,1,...,k\right\} \text{.}
\end{equation*}%
We have that $v_{m}\left( .,.;\mu \right) $ satisfies
\begin{equation}
div\left( a\left( x\right) \nabla v_{m}\left( x,\xi ;\mu \right) \right) +i%
\text{sgn}\left( \mu \right) \partial _{\xi }^{2}v_{m}\left( x,\xi ;\mu
\right) =0\text{, in }B_{1}\times \mathbb{R}\text{.}  \label{rq.1.170}
\end{equation}%
Multiplying equation (\ref{rq.1.170}) by $\overline{v}_{m}\left( x,\xi ;\mu
\right) \eta _{m}^{2}\left( x,\xi \right) $ (here $\overline{v}_{m}$ denotes
the complex conjugate of $v_{m}$) and integrating over $D_{m}=B_{r_{m}}%
\times \left( -r_{m},r_{m}\right) $, we obtain
\begin{equation*}
\int\nolimits_{D_{m}}\left( a\nabla v_{m}\cdot \nabla \overline{v}%
_{m}\right) \eta _{m}^{2}dxd\xi +\int\nolimits_{D_{m}}\left\vert \partial
_{\xi }v_{m}\right\vert ^{2}\eta _{m}^{2}dxd\xi \leq
C_{2}m^{2}\int\nolimits_{D_{m}}\left\vert v_{m}\right\vert ^{2}\eta
_{m}^{2}dxd\xi \text{.}
\end{equation*}%
where $C_{2}=\frac{8\sqrt{2}\pi ^{2}}{\lambda }$.

\noindent Therefore, for every $m\in \left\{ 0,1,...,k\right\} $, we have
\begin{equation}
\int\nolimits_{D_{m+1}}\left\vert \partial _{\xi }^{m+1}v\right\vert
^{2}dxd\xi \leq \frac{C_{2}m^{2}}{\lambda }\int\nolimits_{D_{m}}\left\vert
\partial _{\xi }^{m}v\right\vert ^{2}dxd\xi \text{ .}  \label{rq.1.175}
\end{equation}
By iteration of (\ref{rq.1.175}) for $m=0,1,...,k-1$, we get
\begin{equation}
\int\nolimits_{B_{1/2}\times \left( -1/2,1/2\right) }\left\vert \partial
_{\xi }^{k}v\right\vert ^{2}dxd\xi \leq 2\left( \frac{C_{2}k^{2}}{\lambda }%
\right) ^{k}\int\nolimits_{B_{1}}\left\vert \widehat{u}_{2}\left( x,\mu
\right) \right\vert ^{2}dx\text{ .}  \label{rq.1.180}
\end{equation}

\noindent Now, let us estimate the integral at the right-hand side of (\ref%
{rq.1.180}). By (\ref{rq.1.80}) and Schwarz inequality we have
\begin{equation*}
\int\nolimits_{B_{1}}\left\vert \widehat{u}_{2}\left( x,\mu \right)
\right\vert ^{2}dx=\int\nolimits_{B_{1}}\left\vert
\int\nolimits_{0}^{+\infty }e^{-i\mu t}u_{2}\left( x,t\right) dt\right\vert
^{2}dx
\end{equation*}
\begin{equation*}
\leq \int\nolimits_{B_{1}}dx\left( \int\nolimits_{0}^{+\infty
}e^{-b_{1}\left( t-2\right) _{+}}dt\right) \left(
\int\nolimits_{0}^{+\infty }e^{-b_{1}\left( t-2\right) _{+}}u_{2}^{2}\left(
x,t\right) dt\right) \leq cC_{1}^{2}\left( 2+\frac{1}{b_{1}}\right) ^{2}%
\text{,}
\end{equation*}
where $c$ is an an absolute constant and $C_{1}$ is the constant that
appears in (\ref{rq.1.80}).

\noindent By the just obtained inequality and by (\ref{rq.1.180}) we get,
for every $k\in \mathbb{N}$,
\begin{equation}
\int\nolimits_{B_{1/2}\times \left( -1/2,1/2\right) }\left\vert \partial
_{\xi }^{k}v\right\vert ^{2}dxd\xi \leq cC_{1}^{2}H^{2}\left( 2+\frac{1}{%
b_{1}}\right) ^{2}\left( \frac{C_{2}k^{2}}{\lambda }\right) ^{k}\text{ ,}
\label{rq.1.185}
\end{equation}
where $c$ is an absolute constant.

\noindent For fixed $\mu \in \mathbb{R\smallsetminus }\left\{ 0\right\} $
and $\psi \in L^{2}\left( B_{1/2},\mathbb{C}\right) $, let us denote by $%
\Psi $ the function
\begin{equation*}
\Psi \left( \xi \right) =\int\nolimits_{B_{1/2}}v\left( x,\xi ;\mu \right)
\overline{\psi \left( x\right) }dx\text{, }\xi \in \left( -\frac{1}{2},\frac{%
1}{2}\right) \text{ .}
\end{equation*}
Recall now the inequality
\begin{equation}
\left\Vert f\right\Vert _{L^{\infty }\left( J\right) }\leq c\left(
\left\vert J\right\vert \left\Vert f\right\Vert _{L^{2}\left( J\right)
}^{2}+\left\vert J\right\vert ^{-1}\left\Vert f^{\prime }\right\Vert
_{L^{2}\left( J\right) }^{2}\right) ^{1/2}\text{,}  \label{A.1}
\end{equation}
where $J$ is a bounded interval of $\mathbb{R}$, $\left\vert J\right\vert $
is the lenght of $J$ and $c$ is an absolute constant.

\noindent By (\ref{rq.1.185}) and (\ref{A.1}), we have that for every $k\in
\mathbb{N}\cup \left\{ 0\right\} $ and $\xi \in \left( -\frac{1}{2},\frac{1}{%
2}\right) $,
\begin{equation}
\left\vert \Psi ^{\left( k\right) }\left( \xi \right) \right\vert \leq
cC_{1}H\left\Vert \psi \right\Vert _{L^{2}\left( B_{1/2}\right) }\left( 2+%
\frac{1}{b_{1}}\right) \left( \left( C_{2}\lambda ^{-1}\right)
^{1/2}k\right) ^{k}\text{ .}  \label{rq.1.190}
\end{equation}%
By using inequality (\ref{rq.1.190}) and the power series of $\Psi $ at any
point $\xi _{0}$ such that \textit{Re}$\left( \xi _{0}\right) \in \left( -%
\frac{1}{2},\frac{1}{2}\right) $, \textit{Re}$\left( \xi _{0}\right) =0$, we
have that the function $\Psi $ can be analitically extended in the rectangle
\begin{equation*}
Q=\left\{ \xi \in \mathbb{C}:\text{\textit{Re}}\left( \xi \right) \in \left(
-\frac{1}{2},\frac{1}{2}\right) \text{, \textit{Im}}\left( \xi \right) \in
\left( -2\delta ,2\delta \right) \right\} \text{ ,}
\end{equation*}%
where $\delta =\dfrac{\lambda }{8\pi e}$. Denoting by $\widetilde{\Psi }$
such an analitic extension of $\Psi $ to $Q$ we have
\begin{equation}
\left\vert \widetilde{\Psi }\left( -i\delta \right) \right\vert \leq
c\lambda ^{-1}C_{1}H\left\Vert \psi \right\Vert _{L^{2}\left( B_{1/2}\right)
}\left( 2+\frac{1}{b_{1}}\right) \text{ ,}  \label{rq.1.195}
\end{equation}%
where $c$ is an absolute constant.

\noindent On the other side, by the definition of $v$ we have
\begin{equation*}
\widetilde{\Psi }\left( -i\delta \right)
=\int\nolimits_{B_{1/2}}e^{\left\vert \mu \right\vert ^{1/2}\delta }%
\widehat{u}_{2}\left( x,\mu \right) \overline{\psi \left( x\right) }\text{ ,}
\end{equation*}
so that, by (\ref{rq.1.195}), we obtain (\ref{rq.1.155}).$\blacksquare $

\bigskip

Estimate (\ref{rq.1.155}) allows us to define, for every $x\in B_{1/2}$ and $%
y\in \left( -\sqrt{2}\delta ,\sqrt{2}\delta \right) $, the function
\begin{equation}
w_{2}\left( x,y\right) =\frac{1}{2\pi }\int\nolimits_{-\infty }^{+\infty
}e^{i\mu }\widehat{u}_{2}\left( x,\mu \right) \cosh \left( \sqrt{-i\mu }%
y\right) d\mu \text{ ,}  \label{rq.1.198}
\end{equation}%
where
\begin{equation*}
\sqrt{-i\mu }=\left\vert \mu \right\vert ^{1/2}e^{-\frac{\pi }{4}i\text{sgn}%
\left( \mu \right) }\text{.}
\end{equation*}%
It turns out that $w_{2}$ satisfies the following elliptic equation
\begin{equation}
div\left( a\left( x\right) \nabla w_{2}\left( x,y\right) \right) +\partial
_{y}^{2}w_{2}\left( x,y\right) =0\text{, in }B_{1/2}\times \left( -\sqrt{2}%
\delta ,\sqrt{2}\delta \right)  \label{rq.1.199}
\end{equation}%
and the following conditions at $y=0$%
\begin{equation}
w_{2}\left( x,0\right) =u_{2}\left( x,1\right) \text{ , for every }x\in
B_{1/2}\text{ ,}  \label{rq.1.200}
\end{equation}%
\begin{equation}
\partial _{y}w_{2}\left( x,0\right) =0\text{ , for every }x\in B_{1/2}\text{.%
}  \label{rq.1.205}
\end{equation}%
Moreover, notice that $w_{2}$ is an even function with respect to the
variable $y$.

\bigskip

Now let us consider the function $u_{1}$. Since such a function is the
solution to initial-boundary value problem (\ref{rq.1.35-45}) we have
\begin{equation}
u_{1}\left( x,t\right) =\sum\limits_{k=1}^{\infty }\alpha _{k}e^{\mu
_{k}t}\varphi _{k}\left( x\right) \text{ ,}  \label{rq.1.206}
\end{equation}
where $\mu _{k}$ and $\varphi _{k}\left( x\right) $, for $k\in \mathbb{N}$,
are respectively the decreasing sequence of eigenvalues and the
corresponding eigenfunctions associated to the problem
\begin{equation*}
\left\{
\begin{array}{c}
div\left( a\nabla \varphi \right) =\mu \varphi \text{, in }B_{1}\text{,} \\
\varphi =0\text{, on }\partial B_{1}%
\end{array}
\right.
\end{equation*}
and
\begin{equation*}
\alpha _{k}=\int\nolimits_{B_{1}}u_{1}\left( x,0\right) \varphi _{k}\left(
x\right) dx\text{, }k\in \mathbb{N}\text{.}
\end{equation*}
Let us define
\begin{equation}
w_{1}\left( x,y\right) =\sum\limits_{k=1}^{\infty }\alpha _{k}e^{\mu
_{k}}\varphi _{k}\left( x\right) \cosh \left( \sqrt{\left\vert \mu
_{k}\right\vert }y\right) \text{ ,}  \label{rq.1.210}
\end{equation}
It is easy to check that $w_{1}$ is a solution to equation (\ref{rq.1.199})
and satisfies the following conditions
\begin{equation}
w_{1}\left( x,0\right) =u_{1}\left( x,1\right) \text{ , for every }x\in B_{1}%
\text{ ,}  \label{rq.1.215}
\end{equation}
\begin{equation}
\partial _{y}w_{1}\left( x,0\right) =0\text{ , for every }x\in B_{1}\text{.}
\label{rq.1.220}
\end{equation}
Moreover, notice that $w_{1}$ is an even function with respect to the
variable $y$.

\noindent Now, let us define
\begin{equation}
w=w_{1}+w_{2}\text{, in }B_{1/2}\times \left( -\sqrt{2}\delta ,\sqrt{2}%
\delta \right) \text{ ,}  \label{rq.1.225}
\end{equation}%
we have that $w$ is again a solution to equation (\ref{rq.1.199}), it is an
even function with respect to the variable $y$ and, as a consequence of (\ref%
{rq.1.65}), (\ref{rq.1.200}), (\ref{rq.1.215}), we have that $w\left(
x,0\right) =u\left( x,1\right) $ for every $x\in B_{1/2}$. In addition by (%
\ref{rq.1.205}) and (\ref{rq.1.220}) we get $\partial _{y}w\left( x,0\right)
=0$ , for every $x\in B_{1/2}$. Also, observe that by (\ref{rq.1.155}), (\ref%
{rq.1.198}) and (\ref{rq.1.210}) we derive
\begin{equation}
\left\Vert w\left( .,y\right) \right\Vert _{L^{2}\left( B_{1/2}\right) }\leq
c\frac{C_{1}H}{\lambda \delta }\left( 1+\frac{1}{b_{1}}\right) \text{ , for
every }y\in \left( -\frac{\delta }{2},\frac{\delta }{2}\right) \text{ .}
\label{rq.1.225b}
\end{equation}%
Indeed, if $y\in \left( -\frac{\delta }{2},\frac{\delta }{2}\right) $ then,
by Schwarz inequality, we have
\begin{equation*}
\left\Vert w_{2}\left( .,y\right) \right\Vert _{L^{2}\left( B_{1/2}\right)
}^{2}=\int\nolimits_{B_{1/2}}dx\left\vert \frac{1}{2\pi }%
\int\nolimits_{-\infty }^{+\infty }e^{i\mu }\widehat{u}_{2}\left( x,\mu
\right) \cosh \left( \sqrt{-i\mu }y\right) d\mu \right\vert ^{2}
\end{equation*}%
\begin{equation*}
\leq \frac{1}{\left( 2\pi \right) ^{2}}\int\nolimits_{B_{1/2}}\left(
\int\nolimits_{-\infty }^{+\infty }\left\vert \widehat{u}_{2}\left( x,\mu
\right) \right\vert ^{2}e^{2\delta \left\vert \mu \right\vert ^{1/2}}d\mu
\right) \left( \int\nolimits_{-\infty }^{+\infty }e^{-\delta \left( 2-\sqrt{%
2}\right) \left\vert \mu \right\vert ^{1/2}}d\mu \right)
\end{equation*}%
\begin{equation*}
\leq c\left( \frac{C_{1}H}{\lambda \delta }\left( 1+\frac{1}{b_{1}}\right)
\right) ^{2}\text{ ,}
\end{equation*}%
where $c$ is an absolute constant and $C_{1}$ is the constant that appears
in (\ref{rq.1.80}).

\noindent Moreover, for every $y\in \mathbb{R}$ we have
\begin{equation*}
\left\Vert w_{1}\left( .,y\right) \right\Vert _{L^{2}\left( B_{1}\right)
}^{2}=\sum\limits_{k=1}^{\infty }\alpha _{k}^{2}e^{2\mu _{k}}\cosh
^{2}\left( \sqrt{\left\vert \mu _{k}\right\vert }y\right)
\end{equation*}
\begin{equation*}
\leq \sum\limits_{k=1}^{\infty }\alpha _{k}^{2}=\left\Vert u\left(
.,0\right) \right\Vert _{L^{2}\left( B_{1}\right) }^{2}\leq H^{2}\text{.}
\end{equation*}
By the inequalities obtained above we get (\ref{rq.1.225b}).

In the following proposition we summarize the results obtained above

\begin{proposition}
\label{Prrq.3}Let $w$ be the function defined in (\ref{rq.1.225}). Then $w$
is a solution to the equation
\begin{equation}
div\left( a\left( x\right) \nabla w\left( x,y\right) \right) +\partial
_{y}^{2}w\left( x,y\right) =0\text{, in }B_{1/2}\times \left( -\sqrt{2}%
\delta ,\sqrt{2}\delta \right) \text{ }  \label{rq.1.230}
\end{equation}
and $w$ satisfies the following conditions
\begin{equation}
w\left( x,0\right) =u\left( x,1\right) \text{ , for every }x\in B_{1}\text{ ,%
}  \label{rq.1.235}
\end{equation}
\begin{equation}
\partial _{y}w\left( x,0\right) =0\text{ , for every }x\in B_{1}\text{.}
\label{rq.1.240}
\end{equation}
Moreover $w$ is an even function with respect to the variable $y$ and $w$
satisfies inequality (\ref{rq.1.225b}).
\end{proposition}

The following lemma is nothing else than a stability estimate for the
elliptic Cauchy problem (\ref{rq.1.230})-(\ref{rq.1.240}). A detailed proof
of such a lemma can be found in Lemma 3.1.5 of \cite{CRoVe1}.

\begin{lemma}
\label{Lerq.4}Let $w$ be the function defined in (\ref{rq.1.225}). Then
there exist constants $C$, $C>1$, $\gamma _{0}$, $\beta $, $\gamma
_{0},\beta \in \left( 0,1\right) $, depending on $\lambda $ only such that
for every $r\leq \frac{1}{2}\gamma _{0}$ the following estimate holds true
\begin{equation}
\int\nolimits_{\widetilde{B}_{r}}w^{2}dxdy\leq C\left(
\int\nolimits_{B_{\rho /2}}w^{2}\left( x,0\right) dx\right) ^{\beta }\left(
\int\nolimits_{\widetilde{B}_{\widetilde{\rho }}}w^{2}dxdy\right) ^{1-\beta
}\text{ ,}  \label{rq.1.245}
\end{equation}%
where $\rho =\frac{8\sqrt{2}e\pi }{3\lambda }r$ and $\widetilde{\rho }=2%
\sqrt{2}\rho $.
\end{lemma}

Now we recall the three sphere inequality for elliptic equations \cite{Ku1}.
In what follows we denote by $\widetilde{a}\left( x,y\right) =\left\{
\widetilde{a}^{ij}\left( x,y\right) \right\} _{i,j=1}^{n+1}$ a symmetric $%
\left( n+1\right) \times \left( n+1\right) $ matrix whose entries are real
valued functions. When $\xi \in \mathbb{R}^{n+1}$ and $\left( x,y\right) \in
\mathbb{R}^{n+1}$ we assume that
\begin{equation}
\widetilde{\lambda }\left\vert \xi \right\vert ^{2}\leq
\sum\limits_{i,j=1}^{n+1}\widetilde{a}^{ij}\left( x,y\right) \xi _{i}\xi
_{j}\leq \widetilde{\lambda }^{-1}\left\vert \xi \right\vert ^{2}\text{ }
\label{rq.1.250}
\end{equation}%
and, for every $\left( x_{1},y_{1}\right) ,\left( x_{2},y_{2}\right) \in
\mathbb{R}^{n+1}$,
\begin{equation}
\left( \sum\limits_{l,j=1}^{n+1}\left( \widetilde{a}^{ij}\left(
x_{1},y_{1}\right) -\widetilde{a}^{ij}\left( x_{2},y_{2}\right) \right)
^{2}\right) ^{1/2}\leq \widetilde{\Lambda }\left( \left\vert
x_{1}-x_{2}\right\vert +\left\vert y_{1}-y_{2}\right\vert \right) \text{ ,}
\label{rq.1.255}
\end{equation}%
where $\widetilde{\lambda }$ and $\widetilde{\Lambda }$ are positive numbers
with $\widetilde{\lambda }\in \left( 0,1\right] $.

\begin{lemma}
\label{Lerq.5}\textbf{(three sphere inequality for elliptic equations)}. Let
$\widetilde{a}\left( x,y\right) =\left\{ \widetilde{a}^{lj}\left( x,y\right)
\right\} _{l,j=1}^{n+1}$ satisfy (\ref{rq.1.250}), (\ref{rq.1.255}). Let $%
\widetilde{w}\in H^{1}\left( \widetilde{B}_{1}\right) $ be a weak solution
to
\begin{equation}
div\left( \widetilde{a}\nabla \widetilde{w}\right) =0\text{ , in }\widetilde{%
B}_{1}\text{.}  \label{rq.1.260}
\end{equation}%
Then there exist constants $\gamma _{1}$, $\gamma _{1}\in \left( 0,1\right) $%
, $C$, $C>1$, depending on $\widetilde{\lambda }$ and $\widetilde{\Lambda }$
only such that for every $r_{1}$, $r_{2}$, $r_{3}$ that satisfy $0<r_{1}\leq
r_{2}\leq \frac{\lambda }{2}r_{3}\leq \gamma _{1}$ we have
\begin{equation}
\int\nolimits_{\widetilde{B}_{r_{2}}}\widetilde{w}^{2}dxdy\leq C\left(
\frac{r_{3}}{r_{2}}\right) ^{C}\left( \int\nolimits_{\widetilde{B}_{r_{1}}}%
\widetilde{w}^{2}dxdy\right) ^{\vartheta _{0}}\left( \int\nolimits_{%
\widetilde{B}_{r_{3}}}\widetilde{w}^{2}dxdy\right) ^{1-\vartheta _{0}}\text{%
, }  \label{rq.1.265}
\end{equation}%
where
\begin{equation}
\vartheta _{0}:=\vartheta _{0}\left( r_{1},r_{2},r_{3}\right) =\frac{\log
\left( \frac{1}{2}+\frac{\widetilde{\lambda }r_{3}}{2r_{2}}\right) }{\log
\left( \frac{1}{2}+\frac{\widetilde{\lambda }r_{3}}{2r_{2}}\right) +C\log
\frac{2r_{2}}{\widetilde{\lambda }r_{1}}}\text{ .}  \label{rq.1.270}
\end{equation}
\end{lemma}

\begin{theorem}
\label{Thrq.6}\textbf{(two-sphere one cylinder inequality for equation (\ref%
{rq.1.10})). }Let $u\in H^{2,1}\left( B_{1}\times \left( 0,1\right) \right) $
be a solution to equation (\ref{rq.1.10}) and let (\ref{rq.1.15}), (\ref%
{rq.1.20}) be satisfied. Then there exist constants $\gamma _{2}$, $\gamma
_{2}\in \left( 0,1\right) $, $C$, $C>1$, depending on $\lambda $ and $%
\Lambda $ only such that for every $r_{1}$, $r_{2}$, $r_{3}$ that satisfy $%
0<r_{1}\leq r_{2}\leq \gamma _{2}r_{3}\leq \gamma _{2}^{2}$ we have
\begin{equation}
\int\nolimits_{B_{r_{2}}}u^{2}\left( x,1\right) dx\leq K\left(
r_{2},r_{3}\right) H^{2\left( 1-\beta \vartheta _{1}\right) }\left(
\int\nolimits_{B_{r_{1}}}u^{2}\left( x,1\right) dx\right) ^{\beta \vartheta
_{1}}\text{ ,}  \label{rq.1.275}
\end{equation}%
where
\begin{equation*}
K\left( r_{2},r_{3}\right) =C\frac{r_{3}}{r_{3}-r_{2}}\left( \frac{r_{3}}{%
r_{2}}\right) ^{C}\text{ ,}
\end{equation*}%
\begin{equation}
\vartheta _{1}=\frac{\log \left( \frac{1}{2}+\frac{\lambda r_{3}}{2%
\widetilde{r}_{2}}\right) }{\log \left( \frac{1}{2}+\frac{\lambda r_{3}}{2%
\widetilde{r}_{2}}\right) +C\log \frac{2\widetilde{r}_{2}}{\lambda br_{1}}}%
\text{ ,}  \label{rq.1.276}
\end{equation}%
\begin{equation}
\widetilde{r}_{2}=\left( 1-s_{\lambda }\right) r_{2}+s_{\lambda }r_{3}\text{
, }s_{\lambda }=\frac{\lambda }{4-\lambda }\text{ , }b=\frac{3\lambda }{4%
\sqrt{2}e\pi }\text{ }  \label{notatrq.6}
\end{equation}%
and $\beta $ is the same exponent that appears in inequality (\ref{rq.1.245})
\end{theorem}

\textbf{Proof.} Let us consider the function $w$ introduced in (\ref%
{rq.1.225}). Recall that $w$ satisfies (\ref{rq.1.230}), (\ref{rq.1.235})
and (\ref{rq.1.240}). Let $r_{1}$, $r_{2}$, $r_{3}$ satisfy
\begin{equation*}
0<r_{1}\leq r_{2}\leq \left( 4\max \left\{ \sqrt{2},\lambda ^{-1}\right\}
\right) ^{-1}r_{3}\leq \gamma _{\ast }\text{ ,}
\end{equation*}%
where $\gamma _{\ast }=\min \left\{ \gamma _{0},\frac{\gamma _{1}\lambda
\delta }{\sqrt{2}}\right\} $, $\gamma _{0}$ and $\gamma _{1}$ have been
defined in Lemma \ref{Lerq.4} and Lemma \ref{Lerq.5} respectively, $\delta $
is given by (\ref{rq.1.160}). Let $\widetilde{r}_{2}$ and $b$ be defined by (%
\ref{notatrq.6}), we have $0<br_{1}<r_{1}\leq r_{2}\leq \widetilde{r}%
_{2}\leq \frac{\lambda }{2}r_{3}\leq \gamma _{1}$. By applying Lemma \ref%
{Lerq.5} to the triplet of radii $br_{1}$, $\widetilde{r}_{2}$, $r_{3}$ we
get
\begin{equation}
\int\nolimits_{\widetilde{B}_{\widetilde{r}_{2}}}w^{2}dxdy\leq C\left(
\frac{r_{3}}{r_{2}}\right) ^{C}\left( \int\nolimits_{\widetilde{B}%
_{br_{1}}}w^{2}dxdy\right) ^{\vartheta _{1}}\left( \int\nolimits_{%
\widetilde{B}_{r_{3}}}w^{2}dxdy\right) ^{1-\vartheta _{1}}\text{ ,}
\label{rq.1.280}
\end{equation}%
where $\vartheta _{1}=\vartheta _{0}\left( br_{1},r_{2}^{\prime
},r_{3}\right) $, with $\vartheta _{0}$ defined by (\ref{rq.1.276}) and $C$
depending on $\lambda $ and $\Lambda $ only.

Now, let us recall the following trace inequality. Given $r$, $\rho $, $%
0<\rho <r$, $f\in H^{1}\left( B_{r}^{+}\right) $, we have
\begin{equation}
\int\nolimits_{B_{\rho }}f^{2}\left( x,0\right) dx\leq c\left( \frac{r}{%
r^{2}-\rho ^{2}}\int\nolimits_{\widetilde{B}_{r}^{+}}f^{2}\left( x,y\right)
dxdy+r\int\nolimits_{\widetilde{B}_{r}^{+}}\left( \partial _{y}f\left(
x,y\right) \right) ^{2}dxdy\right) \text{ ,}  \label{rq.1.285}
\end{equation}%
where $c$ is an absolute constant.

\noindent Set $r_{2}^{\ast }=\frac{1}{2}\left( r_{2}+\widetilde{r}%
_{2}\right) $. By Caccioppoli inequality and (\ref{rq.1.285}) we have
\begin{equation*}
\int\nolimits_{\widetilde{B}_{\widetilde{r}_{2}}}w^{2}dxdy=\frac{1}{2}%
\int\nolimits_{\widetilde{B}_{\widetilde{r}_{2}}}w^{2}dxdy+\frac{1}{2}%
\int\nolimits_{\widetilde{B}_{\widetilde{r}_{2}}}w^{2}dxdy
\end{equation*}%
\begin{equation*}
\geq \frac{1}{2}\int\nolimits_{\widetilde{B}_{r_{2}^{\ast
}}}w^{2}dxdy+C\left( r_{2}-\widetilde{r}_{2}\right) ^{2}\int\nolimits_{%
\widetilde{B}_{r_{2}^{\ast }}}\left( \left\vert \nabla _{x}w\right\vert
^{2}+\left( \partial _{y}w\right) ^{2}\right) dxdy
\end{equation*}%
\begin{equation*}
\geq C^{\prime }\left( r_{3}-r_{2}\right) \left( \frac{r_{2}^{\ast }}{%
r_{2}^{\ast 2}-r_{2}^{2}}\int\nolimits_{\widetilde{B}_{r_{2}^{\ast
}}}w^{2}dxdy+r_{2}^{\ast }\int\nolimits_{\widetilde{B}_{r_{2}^{\ast
}}}\left( \partial _{y}w\right) ^{2}dxdy\right)
\end{equation*}%
\begin{equation*}
\geq C^{\prime \prime }\left( r_{3}-r_{2}\right)
\int\nolimits_{B_{r_{2}}}w^{2}\left( x,0\right) dx=C^{\prime \prime }\left(
r_{3}-r_{2}\right) \int\nolimits_{B_{r_{2}}}u^{2}\left( x,1\right) dx\text{
,}
\end{equation*}%
where $C$, $C^{\prime }$, $C^{\prime \prime }$ depend on $\lambda $ only.

\noindent By the just obtained inequality and recalling (\ref{rq.1.225}) and
(\ref{rq.1.280}) we obtain (\ref{rq.1.275}).$\blacksquare $

\begin{corollary}
\label{Cor}\textbf{(Spacelike strong unique continuation for equation (\ref%
{rq.1.10})).} Let $u\in H^{2,1}\left( B_{1}\times \left( 0,1\right) \right) $
satisfy equation (\ref{rq.1.10}).

\noindent If for every $k\in \mathbb{N}$ we have
\begin{equation}
\int\nolimits_{B_{r}}u^{2}\left( x,1\right) dx=O\left( r^{2k}\right) \text{
, as }r\rightarrow 0\text{,}  \label{bigO}
\end{equation}
then
\begin{equation*}
u\left( .,1\right) =0\text{ , in }B_{1}\text{ .}
\end{equation*}
\end{corollary}

\textbf{Proof.} Let us fix $\rho \in \left( 0,\gamma _{2}^{2}\right] $,
where $\gamma _{2}$ has been introduced in Theorem \ref{Thrq.6}. By applying
Lemma \ref{Lerq.5} to the triplet of radii $r_{1}=r$, $r_{2}=\rho $ and $%
r_{3}=\gamma _{2}$, by (\ref{bigO}) and passing to the limit as $r$ tends to
$0$, we get
\begin{equation}
\int\nolimits_{B_{\rho }}u^{2}\left( x,t_{0}\right) dx\leq Ce^{-C_{1}k}%
\text{ , for every }k\in \mathbb{N}\text{, }  \label{n1}
\end{equation}
where $C$ depends on $\lambda $, $\Lambda $ and $\left\Vert u\right\Vert
_{L^{2}\left( B_{1}\times \left( 0,1\right) \right) }$ only and $C_{1}$
depends on $\lambda $ and $\Lambda $, only. Passing to the limit as $%
k\rightarrow \infty $, (\ref{n1}) yields $u\left( .,1\right) =0$ in $B_{\rho
}$. By iteration the thesis follows.$\blacksquare $

\subsection{Carleman estimate\label{SubsRq2.1}}

In the present section it is convenient to carry out the calculations for
the backward parabolic operator. Since in other parts of the paper\ we deal
with the forward parabolic operator, for the reader convenience, when we
consider the backward operator, we rename the time variable by $s$. Thus we
denote by
\begin{equation}
Pu=\partial _{i}\left( g^{ij}\left( x,s\right) \partial _{j}u\right)
+\partial _{s}u\text{,}  \label{Rq2.100}
\end{equation}%
where $\left\{ g^{ij}\left( x,s\right) \right\} _{i,j=1}^{n}$ is a symmetric
$n\times n$ matrix whose entries are real functions. When $\xi \in \mathbb{R}%
^{n}$ and $\left( x,s\right) ,\left( y,\tau \right) \in \mathbb{R}^{n+1}$ we
assume that
\begin{equation}
\lambda \left\vert \xi \right\vert ^{2}\leq
\sum\limits_{i,j=1}^{n}g^{ij}\left( x,s\right) \xi _{i}\xi _{j}\leq \lambda
^{-1}\left\vert \xi \right\vert ^{2}\text{ }  \label{Rq2.105}
\end{equation}%
and
\begin{equation}
\left( \sum\limits_{i,j=1}^{n}\left( g^{ij}\left( x,s\right) -g^{ij}\left(
y,\tau \right) \right) ^{2}\right) ^{1/2}\leq \Lambda \left( \left\vert
x-y\right\vert ^{2}+\left\vert s-\tau \right\vert \right) ^{1/2}\text{ ,}
\label{Rq2.110}
\end{equation}%
where $\lambda $ and $\Lambda $ are positive numbers with $\lambda \in
\left( 0,1\right] $.

In order to simplify the exposition and the calculations in the proof of
Theorem \ref{ThRq9} below we use some of the standard notations of
Riemannian geometry, but always dropping the corresponding volume element in
the definition of the Laplace-Beltrami operator associated to a Riemannian
metric. More precisely, letting $g\left( x,s\right) =\left\{ g_{ij}\left(
x,s\right) \right\} _{i,j=1}^{n}$ to denote the inverse matrix of the matrix
of coefficients of $P$, we have $g^{-1}\left( x,s\right) =\left\{
g^{ij}\left( x,s\right) \right\} _{i,j=1}^{n}$ and we use the following
notation when considering either a function $v$ or two vector fields $\xi $
and $\eta $:

1. $\left( \xi \mid \eta \right) =g_{ij}\left( x,s\right) \xi _{i}\eta _{j}$
, $\left\vert \xi \right\vert _{g}^{2}=g_{ij}\left( x,s\right) \xi _{i}\xi
_{j}$ ,

2. $\nabla _{x}v=\left( \partial _{1}v,...\partial _{n}v\right) $ , $\nabla
_{g}v\left( x,s\right) =g^{-1}\left( x,s\right) \nabla _{x}v\left(
x,s\right) $ , $div\left( \xi \right) =\sum\limits_{i=1}^{n}\partial
_{i}\xi _{i}$, $\Delta _{g}=div\left( \nabla _{g}\right) $ .

With this notations the following formulae hold true when $u$, $v$ and $w$
are smooth functions
\begin{equation}
Pu=\Delta _{g}u+\partial _{s}u\text{ , }\Delta _{g}\left( v^{2}\right)
=2v\Delta _{g}v+2\left\vert \nabla _{g}v\right\vert _{g}^{2}  \label{Rq2.115}
\end{equation}%
and
\begin{equation}
\int v\Delta _{g}wdx=\int w\Delta _{g}vdx=-\int \left( \nabla _{g}v\mid
\nabla _{g}w\right) dx\text{ .}  \label{Rq2.120}
\end{equation}%
Also, we shall use the following Rellich-Ne\v{c}as-Pohozaev identity
\begin{eqnarray}
2\left( \beta \mid \nabla _{g}v\right) \Delta _{g}v &=&2div\left( \left(
\left( \beta \mid \nabla _{g}v\right) \right) \nabla _{g}v\right) -div\left(
\beta \left\vert \nabla _{g}v\right\vert _{g}^{2}\right) +  \label{Rq2.125}
\\
&&div\left( \beta \right) \left\vert \nabla _{g}v\right\vert
_{g}^{2}-2\partial _{i}\beta ^{k}g^{ij}\partial _{j}v\partial _{k}v+\beta
^{k}\partial _{k}g^{ij}\partial _{i}v\partial _{j}v\text{ ,}  \notag
\end{eqnarray}%
where $\beta =\left( \beta ^{1},...,\beta ^{n}\right) $ is a smooth vector
field.

\begin{lemma}
\label{LeRq7}Assume that $\theta :\left( 0,1\right) \rightarrow \left(
0,+\infty \right) $ satisfies
\begin{equation}
0\leq \theta \leq C_{0}\text{ , }\left\vert s\theta ^{\prime }\left(
s\right) \right\vert \leq C_{0}\theta \left( s\right) \text{ , }%
\int\nolimits_{0}^{1}\left( 1+\log \frac{1}{s}\right) \frac{\theta \left(
s\right) }{s}ds\leq C_{0}\text{ ,}  \label{Rq2.130}
\end{equation}
for some constant $C_{0}$. Let $\gamma \geq 1$ and
\begin{equation}
\sigma \left( s\right) =s\exp \left( -\int\nolimits_{0}^{\gamma s}\left(
1-\exp \left( -\int\nolimits_{0}^{t}\frac{\theta \left( \eta \right) }{\eta
}d\eta \right) \right) \frac{dt}{t}\right) \text{ .}  \label{Rq2.135}
\end{equation}
Then $\sigma $ is the solution to the Cauchy problem
\begin{equation}
\frac{d}{ds}\log \left( \frac{\sigma }{s\sigma ^{\prime }}\right) =\frac{%
\theta \left( \gamma s\right) }{s}\text{ , }\sigma \left( 0\right) =0\text{,
}\sigma ^{\prime }\left( 0\right) =1  \label{Rq2.140}
\end{equation}
and satisfies the following properties when $0\leq \gamma s\leq 1$%
\begin{equation}
se^{-C_{0}}\leq \sigma \left( s\right) \leq s\text{ ,}  \label{Rq2.145}
\end{equation}
\begin{equation}
e^{-C_{0}}\leq \sigma ^{\prime }\left( s\right) \leq 1\text{ .}
\label{Rq2.150}
\end{equation}
\end{lemma}

\textbf{Proof.} The proof is straightforward.$\blacksquare $

\begin{lemma}
\label{LeRq8}For every positive number $\mu $ there exists a constant $C$
such that for all $y\geq 0$ and $\varepsilon \in \left( 0,1\right) $,
\begin{equation}
y^{\mu }e^{-y}\leq C\left( \varepsilon +\left( \log \frac{1}{\varepsilon }%
\right) ^{\mu }e^{-y}\right) \text{.}  \label{Rq2.155}
\end{equation}
\end{lemma}

\textbf{Proof.} Consider the function $\varphi \left( y\right)
=y-\varepsilon e^{y}$ on $\left[ 0,+\infty \right) $. Since the maximum of $%
\varphi $ is $\log \dfrac{1}{\varepsilon }-1$, $y-\varepsilon e^{y}<\log
\dfrac{1}{\varepsilon }$ and we get (\ref{Rq2.155}) when $\mu =1$. If $\mu
>1 $, we use the just proved inequality and the convexity of $y^{\mu }$ to
get
\begin{equation*}
\left( \frac{y}{\mu }\right) ^{\mu }\leq \left( \varepsilon ^{1/\mu
}e^{y/\mu }+\log \frac{1}{\varepsilon ^{1/\mu }}\right) ^{\mu }\leq 2^{\mu
-1}\left( \varepsilon e^{y}+\left( \frac{1}{\mu }\log \frac{1}{\varepsilon }%
\right) ^{\mu }\right) \text{ ,}
\end{equation*}%
that gives (\ref{Rq2.155}).

\noindent If $0<\mu <1$ we obtain (\ref{Rq2.155}) by the inequality for $\mu
=1$ and the inequality $\left( a+b\right) ^{\mu }\leq a^{\mu }+b^{\mu }$.$%
\blacksquare $

\bigskip

In the sequel of this section we shall denote by%
\begin{equation}
\theta \left( s\right) =s^{1/2}\left( \log \frac{1}{s}\right) ^{3/2}\text{ ,
}s\in \left( 0,1\right] \text{ .}  \label{Rq2.160}
\end{equation}%
It is easy to check that $\theta $ satisfies (\ref{Rq2.130}) of Lemma \ref%
{LeRq7}. We denote by $\sigma $ the function defined in (\ref{Rq2.135})
where $\theta $ is given by (\ref{Rq2.160}) and $\gamma =\dfrac{\alpha }{%
\delta ^{2}}$ with $\alpha \geq 1$, $\delta \in \left( 0,1\right) $. In
addition, for a given number $a>0$, we shall denote by $\sigma _{a}\left(
s\right) =\sigma \left( s+a\right) $.

\begin{theorem}
\label{ThRq9}Let $P$ be the operator defined in (\ref{Rq2.100}) and assume
that (\ref{Rq2.105}) and (\ref{Rq2.110}) are satisfied. Assume that $%
g^{ij}\left( 0,0\right) =\delta ^{ij}$, $i,j=1,...,n$. Then there exist
constants $C$, $C\geq 1$, $\eta _{0}\in \left( 0,1\right) $ and $\delta
_{1}\in \left( 0,1\right) $ depending on $\lambda $ only $\Lambda $ only,
such that for every $\alpha $, $\alpha \geq 2$, $a$, $0<a\leq \dfrac{\delta
^{2}}{4\alpha }$, $\delta \in \left( 0,\delta _{1}\right] $ and $u\in
C_{0}^{\infty }\left( \mathbb{R}^{n}\times \left[ 0,+\infty \right) \right) $%
, with supp$\,u\subset B_{\eta _{0}}\times \left[ 0,\dfrac{\delta ^{2}}{%
2\alpha }\right) $, the following inequality holds true
\begin{equation}
\alpha \int\nolimits_{\mathbb{R}_{+}^{n+1}}\sigma _{a}^{-2\alpha -2}\left(
s\right) \theta \left( \gamma \left( s+a\right) \right) u^{2}e^{-\frac{%
\left\vert x\right\vert ^{2}}{4\left( s+a\right) }}dX  \label{Rq2.165}
\end{equation}%
\begin{equation*}
+\int\nolimits_{\mathbb{R}_{+}^{n+1}}\left( s+a\right) \sigma
_{a}^{-2\alpha -2}\left( s\right) \theta \left( \gamma \left( s+a\right)
\right) \left\vert \nabla _{g}u\right\vert _{g}^{2}e^{-\frac{\left\vert
x\right\vert ^{2}}{4\left( s+a\right) }}dX
\end{equation*}%
\begin{equation*}
\leq C\int\nolimits_{\mathbb{R}_{+}^{n+1}}\left( s+a\right) ^{2}\sigma
_{a}^{-2\alpha -2}\left( s\right) \left\vert Pu\right\vert ^{2}e^{-\frac{%
\left\vert x\right\vert ^{2}}{4\left( s+a\right) }}dX+C^{\alpha }\alpha
^{\alpha }\int\nolimits_{\mathbb{R}_{+}^{n+1}}\left( u^{2}+\left(
s+a\right) \left\vert \nabla _{g}u\right\vert _{g}^{2}\right) dX
\end{equation*}%
\begin{equation*}
+C\alpha \sigma ^{-2\alpha -1}\left( a\right) \int\nolimits_{\mathbb{R}%
^{n}}u^{2}\left( x,0\right) e^{-\frac{\left\vert x\right\vert ^{2}}{4a}}dx-%
\frac{\sigma ^{-2\alpha }\left( a\right) }{C}\int\nolimits_{\mathbb{R}%
^{n}}\left\vert \left( \nabla _{g}u\right) \left( x,0\right) \right\vert
_{g\left( .,0\right) }^{2}e^{-\frac{\left\vert x\right\vert ^{2}}{4a}}dx%
\text{.}
\end{equation*}
\end{theorem}

\textbf{Proof.} We divide the proof into two steps. In the first step we
prove that if the matrix $g$ does not depend on $s$ then there exists $%
\delta _{0}\in \left( 0,1\right) $ such that, for every $\delta \in \left(
0,\delta _{0}\right] $, inequality (\ref{Rq2.165}) holds true. In the second
step we prove inequality (\ref{Rq2.165}), for a suitable $\delta _{1}$, in
the general case. For the sake of brevity in what follows we prove the
inequality when the matrix $g\left( x,s\right) $ and the function $u\left(
x,s\right) $ are replaced by $\widetilde{g}\left( x,s\right) :=g\left(
x,s-a\right) $ and $\widetilde{u}\left( x,s\right) :=u\left( x,s-a\right) $
respectively, so that $\widetilde{g}^{ij}\left( 0,a\right) =\delta ^{ij}$
and $\widetilde{u}$ is a function belonging to $C_{0}^{\infty }\left(
\mathbb{R}^{n}\times \left[ a,+\infty \right) \right) $. Also, the sign " $%
\widetilde{}$ " over $g$ and $u$ will be dropped. Finally, we shall denote
by $\int \left( .\right) dX$ the integral $\int\nolimits_{\mathbb{R}%
^{n}\times \left[ a,+\infty \right) }\left( .\right) dX$.

\underline{STEP 1}. Set $g\left( x\right) =g\left( x,a\right) $ and denote
\begin{equation*}
P_{0}u=\partial _{i}\left( g^{ij}\left( x\right) \partial _{j}u\right)
+\partial _{s}u\text{ ,}
\end{equation*}
\begin{equation*}
\phi \left( x,s\right) =-\frac{\left\vert x\right\vert ^{2}}{8s}-\left(
\alpha +1\right) \log \sigma \left( s\right) \text{ ,}
\end{equation*}
\begin{equation*}
v=e^{\phi }u\text{ ,}
\end{equation*}
\begin{equation*}
Lv=e^{\phi }P_{0}\left( e^{-\phi }v\right) \text{ .}
\end{equation*}
Denoting by $S$ and $A$ the symmetric and skew-symmetric part of the
operator $sL$ respectively, we have
\begin{equation*}
Sv=s\left( \Delta _{g}v+\left( \left\vert \nabla _{g}\phi \right\vert
_{g}^{2}-\partial _{s}\phi \right) v\right) -\frac{1}{2}v\text{ ,}
\end{equation*}
\begin{equation*}
Av=\frac{1}{2}\left( \partial _{s}\left( sv\right) +s\partial _{s}v\right)
-s\left( v\Delta _{g}\phi +2\nabla _{g}v\cdot \nabla _{g}\phi \right)
\end{equation*}
and
\begin{equation*}
sLv=Sv+Av\text{.}
\end{equation*}
Furthermore, noticing that
\begin{equation*}
\Delta _{g}\phi =-\frac{n}{4s}+\psi _{0}\left( x,s\right) \text{ ,}
\end{equation*}
where
\begin{equation*}
\psi _{0}\left( x,s\right) =\frac{1}{4s}\left( n-tr\left( g^{-1}\right)
-x_{j}\partial _{i}g^{ij}\right) \text{ ,}
\end{equation*}
we write
\begin{equation}
sLv=Sv+A_{0}v-s\psi _{0}v\text{ ,}  \label{Rq2.180}
\end{equation}
where
\begin{equation*}
A_{0}v=\frac{1}{2}\left( \partial _{s}\left( sv\right) +s\partial
_{s}v\right) -s\left( -\frac{n}{4s}v+2\left( \nabla _{g}v\mid \nabla
_{g}\phi \right) \right) \text{ .}
\end{equation*}
Noticing that $\left\vert \psi _{0}\right\vert \leq \dfrac{C\Lambda
\left\vert x\right\vert }{s}$, where $C$ is an absolute constant, we have by
(\ref{Rq2.180})
\begin{equation}
2\int s^{2}\left\vert Lv\right\vert ^{2}dX\geq \int \left\vert
Sv+A_{0}v\right\vert ^{2}dX-4\int s^{2}\psi _{0}^{2}v^{2}dX
\label{Rq2.185b}
\end{equation}

\begin{equation*}
\geq \int \left\vert Sv\right\vert ^{2}dX+\int \left\vert
A_{0}v\right\vert ^{2}dX+2\int SvA_{0}vdX-4C^{2}\Lambda ^{2}\int
\left\vert x\right\vert ^{2}v^{2}dX\text{.}
\end{equation*}
Let us consider the term $I:=2\int SvA_{0}vdX$ at the right-hand side of (%
\ref{Rq2.185b}). We have
\begin{equation}
I=-2\int s^{2}\left( \Delta _{g}v+\left( \left\vert \nabla _{g}\phi
\right\vert _{g}^{2}-\partial _{s}\phi \right) v\right) \left( -\frac{n}{4s}%
v+2\left( \nabla _{g}v\mid \nabla _{g}\phi \right) \right) dX
\label{Rq2.186}
\end{equation}

\begin{equation*}
+\int s\left( \Delta _{g}v+\left( \left\vert \nabla _{g}\phi \right\vert
_{g}^{2}-\partial _{s}\phi \right) v\right) \left( \partial _{s}\left(
sv\right) +s\partial _{s}v\right) dX
\end{equation*}
\begin{equation*}
+\int \left\{ -\frac{1}{2}v\left( \partial _{s}\left( sv\right) +s\partial
_{s}v\right) +sv\left( -\frac{n}{4s}v+2\left( \nabla _{g}v\mid \nabla
_{g}\phi \right) \right) \right\} dX
\end{equation*}
\begin{equation*}
:=I_{1}+I_{2}+I_{3}\text{ .}
\end{equation*}

\bigskip

\textit{Evaluation of} $I_{1}$.

\noindent We have
\begin{equation}
I_{1}=-4\int s^{2}\Delta _{g}v\left( \nabla _{g}v\mid \nabla _{g}\phi
\right) dX-4\int s^{2}\left( \left\vert \nabla _{g}\phi \right\vert
_{g}^{2}-\partial _{s}\phi \right) v\left( \nabla _{g}v\mid \nabla _{g}\phi
\right) dX  \label{Rq2.187}
\end{equation}%
\begin{equation*}
+\frac{n}{2}\int sv\Delta _{g}vdX+\frac{n}{2}\int s\left( \left\vert
\nabla _{g}\phi \right\vert _{g}^{2}-\partial _{s}\phi \right)
v^{2}dX:=I_{11}+I_{12}+I_{13}+I_{14}\text{ .}
\end{equation*}%
By using identity (\ref{Rq2.125}) with the vector field $\beta =\nabla
_{g}\phi $ we have
\begin{eqnarray}
I_{11} &=&\int \left\{ 4s^{2}\partial _{i}\left( g^{kj}\partial _{j}\phi
\right) g^{ih}\partial _{h}v\partial _{k}v-2s^{2}\Delta _{g}\phi \left\vert
\nabla _{g}v\right\vert _{g}^{2}\right.  \label{Rq2.189} \\
&&\left. -2s^{2}\left( g^{kj}\partial _{j}\phi \right) \left( \partial
_{k}g^{hl}\partial _{h}v\partial _{l}v\right) \right\} dX\text{ .}  \notag
\end{eqnarray}%
Concerning the term $I_{12}$, by using the identity $v\nabla _{g}v=\frac{1}{2%
}\nabla _{g}\left( v^{2}\right) $ and by integration parts, we get
\begin{equation*}
I_{12}=2\int s^{2}v^{2}\left\{ \Delta _{g}\phi \left( \left\vert \nabla
_{g}\phi \right\vert _{g}^{2}-\partial _{s}\phi \right) +g^{ij}\partial
_{i}\phi \partial _{j}\left( \left\vert \nabla _{g}\phi \right\vert
_{g}^{2}-\partial _{s}\phi \right) \right\} dX\text{ .}
\end{equation*}%
Now we have
\begin{eqnarray*}
&&g^{ij}\partial _{i}\phi \partial _{j}\left( \left\vert \nabla _{g}\phi
\right\vert _{g}^{2}-\partial _{s}\phi \right) \\
&=&\left( 2\partial _{jk}^{2}\phi g^{kh}\partial _{h}\phi +\partial
_{j}g^{hk}\partial _{k}\phi \partial _{h}\phi \right) g^{ij}\partial
_{i}\phi -\frac{1}{2}\partial _{s}\left( \left\vert \nabla _{g}\phi
\right\vert _{g}^{2}\right) \text{ ,}
\end{eqnarray*}%
hence
\begin{equation}
I_{12}=2\int s^{2}v^{2}\Delta _{g}\phi \left( \left\vert \nabla _{g}\phi
\right\vert _{g}^{2}-\partial _{s}\phi \right) dX  \label{Rq2.191}
\end{equation}%
\begin{equation*}
+4\int s^{2}v^{2}\left( \partial _{jk}^{2}\phi g^{kh}\partial _{h}\phi +%
\frac{1}{2}\partial _{j}g^{hk}\partial _{k}\phi \partial _{h}\phi \right)
g^{ij}\partial _{i}\phi dX-\int s^{2}v^{2}\partial _{s}\left( \left\vert
\nabla _{g}\phi \right\vert _{g}^{2}\right) dX\text{ .}
\end{equation*}%
Concerning $I_{13}$, by the identity $v\Delta _{g}v=div\left( v\nabla
_{g}v\right) -\left\vert \nabla _{g}v\right\vert _{g}^{2}$ and by
integration by parts we obtain
\begin{equation}
I_{13}=-\frac{n}{2}\int s\left\vert \nabla _{g}v\right\vert _{g}^{2}dX\text{
.}  \label{Rq2.193}
\end{equation}%
Now, let us introduce the notation
\begin{equation*}
\psi _{1}=\left( \frac{n}{2}s+2s^{2}\Delta _{g}\phi \right) \left(
\left\vert \nabla _{g}\phi \right\vert _{g}^{2}-\partial _{s}\phi \right)
+2s^{2}\partial _{j}g^{hk}\partial _{k}\phi \partial _{h}\phi g^{ij}\partial
_{i}\phi \text{ ,}
\end{equation*}%
\begin{equation*}
Q\left( \nabla _{g}v\right) =-2s^{2}\psi _{0}\left\vert \nabla
_{g}v\right\vert _{g}^{2}-2s^{2}\left( g^{kj}\partial _{j}\phi \right)
\left( \partial _{k}g^{hl}\partial _{h}v\partial _{l}v\right) \text{ .}
\end{equation*}%
By (\ref{Rq2.187}), (\ref{Rq2.189}), (\ref{Rq2.191}), (\ref{Rq2.193}) we
have
\begin{eqnarray}
I_{1} &=&4\int s^{2}\left\{ \partial _{i}\left( g^{kj}\partial _{j}\phi
\right) g^{ih}\partial _{h}v\partial _{k}v+\partial _{jk}^{2}\phi
g^{kh}\partial _{h}\phi g^{ij}\partial _{i}\phi v^{2}\right\} dX
\label{Rq2.200} \\
&&-\int s^{2}v^{2}\partial _{s}\left( \left\vert \nabla _{g}\phi
\right\vert _{g}^{2}\right) dX+\int \psi _{1}v^{2}dX+\int Q\left( \nabla
_{g}v\right) dX\text{ .}  \notag
\end{eqnarray}

\bigskip

\textit{Evaluation of} $I_{2}$.

\noindent We have
\begin{equation}
I_{2}=\int sv\Delta _{g}vdX+\int 2s^{2}\Delta _{g}v\partial _{s}vdX
\label{Rq2.205}
\end{equation}
\begin{equation*}
+\int 2s^{2}\left( \left\vert \nabla _{g}\phi \right\vert _{g}^{2}-\partial
_{s}\phi \right) v\partial _{s}vdX+\int s\left( \left\vert \nabla _{g}\phi
\right\vert _{g}^{2}-\partial _{s}\phi \right) v^{2}dX\text{ .}
\end{equation*}
Now, we use the identities $v\Delta _{g}v=div\left( v\nabla _{g}v\right)
-\left\vert \nabla _{g}v\right\vert _{g}^{2}$, $2\partial _{s}v\Delta
_{g}v=2div\left( \partial _{s}v\nabla _{g}v\right) -\partial _{s}\left(
\left\vert \nabla _{g}v\right\vert _{g}^{2}\right) $ and $2v\partial
_{s}v=\partial _{s}\left( v^{2}\right) $ in the first, the second and the
third integral at the right-hand side of (\ref{Rq2.205}) respectively, then
we integrate by parts and we obtain
\begin{eqnarray}
I_{2} &=&\int s\left\vert \nabla _{g}v\right\vert _{g}^{2}dX-\int s\left(
\left\vert \nabla _{g}\phi \right\vert _{g}^{2}-\partial _{s}\phi \right)
v^{2}-\int s^{2}\partial _{s}\left( \left\vert \nabla _{g}\phi \right\vert
_{g}^{2}-\partial _{s}\phi \right) v^{2}  \label{Rq2.210} \\
&&+a^{2}\int\nolimits_{\mathbb{R}^{n}}\left\vert \nabla _{g}v\left(
x,a\right) \right\vert _{g}^{2}dx-a^{2}\int\nolimits_{\mathbb{R}^{n}}\left(
\left\vert \nabla _{g}\phi \left( x,a\right) \right\vert _{g}^{2}-\partial
_{s}\phi \left( x,a\right) \right) v^{2}\left( x,a\right) dx\text{ .}  \notag
\end{eqnarray}

\bigskip

\textit{Evaluation of} $I_{3}$

\noindent We have
\begin{equation}
I_{3}=-\int \left( \frac{1}{2}+\frac{n}{4}\right) v^{2}dX-\int sv\partial
_{s}vdX+\int 2sv\left( \nabla _{g}v\mid \nabla _{g}\phi \right) dX\text{ .}
\label{Rq2.215}
\end{equation}
We use the identities $v\partial _{s}v=\dfrac{1}{2}\partial _{s}\left(
v^{2}\right) $ and $2v\left( \nabla _{g}v\mid \nabla _{g}\phi \right)
=div\left( v^{2}\nabla _{g}\phi \right) -v^{2}\Delta _{g}\phi $ respectively
in the second and the third integral at the right-hand side of (\ref{Rq2.215}%
) then we integrate by parts and we obtain
\begin{equation}
I_{3}=-\int s\psi _{0}v^{2}dX+\frac{a}{2}\int\limits_{\mathbb{R}%
^{n}}v^{2}\left( x,a\right) dx\text{ .}  \label{Rq2.220}
\end{equation}

Now, denote by
\begin{eqnarray}
J_{1} &=&\int 4s^{2}\partial _{jk}^{2}\phi g^{kh}\partial _{h}\phi
g^{ij}\partial _{i}\phi v^{2}dX  \label{Rq2.221} \\
&&-\int \left\{ s^{2}\partial _{s}\left( 2\left\vert \nabla _{g}\phi
\right\vert _{g}^{2}-\partial _{s}\phi \right) v^{2}+s\left( \left\vert
\nabla _{g}\phi \right\vert _{g}^{2}-\partial _{s}\phi \right) v^{2}\right\}
dX\text{ ,}  \notag
\end{eqnarray}
\begin{equation}
J_{2}=\int \left\{ 4s^{2}\partial _{i}\left( g^{kj}\partial _{j}\phi
\right) g^{ih}\partial _{h}v\partial _{k}v+s\left\vert \nabla
_{g}v\right\vert _{g}^{2}+Q\left( \nabla _{g}v\right) \right\} dX\text{ ,}
\label{Rq2.222}
\end{equation}
\begin{equation}
J_{3}=\int \left( \psi _{1}-s\psi _{0}\right) v^{2}dX\text{ ,}
\label{Rq2.223}
\end{equation}
\begin{equation}
\mathcal{B}_{1}=a^{2}\int\nolimits_{\mathbb{R}^{n}}\left\vert \nabla
_{g}v\left( x,a\right) \right\vert _{g}^{2}dx+\frac{a}{2}\int\nolimits_{%
\mathbb{R}^{n}}v^{2}\left( x,a\right) dx  \label{Rq2.224}
\end{equation}
\begin{equation*}
-a^{2}\int\nolimits_{\mathbb{R}^{n}}\left( \left\vert \nabla _{g}\phi
\left( x,a\right) \right\vert _{g}^{2}-\partial _{s}\phi \left( x,a\right)
\right) v^{2}\left( x,a\right) dx\text{.}
\end{equation*}
By (\ref{Rq2.186}), (\ref{Rq2.200}), (\ref{Rq2.210}), (\ref{Rq2.215}) we
have
\begin{equation}
I=J_{1}+J_{2}+J_{3}+\mathcal{B}_{1}\text{ .}  \label{Rq2.225}
\end{equation}

Let us now estimate from below the terms $J_{k}$, $k=1,2,3,$ at the
right-hand side of (\ref{Rq2.225}). We easily have
\begin{equation}
\partial _{i}\phi =-\frac{x_{i}}{4s}\text{ , }\partial _{ij}^{2}\phi =-\frac{%
\delta _{ij}}{4s}\text{ , }\left\vert \nabla _{g}\phi \right\vert
_{g}^{2}=g^{ij}\left( x\right) \frac{x_{i}x_{j}}{16s^{2}}\text{ ,}
\label{Rq2.230}
\end{equation}%
\begin{equation}
\Delta _{g}\phi =-\frac{tr\left( g^{-1}\left( x\right) \right)
+x_{j}\partial _{i}g^{ij}\left( x\right) }{4s}\text{ , }\partial _{s}\phi =%
\frac{\left\vert x\right\vert ^{2}}{8s^{2}}-\left( \alpha +1\right) \frac{%
\sigma ^{\prime }\left( s\right) }{\sigma \left( s\right) }\text{ ,}
\label{Rq2.231}
\end{equation}%
\begin{equation}
\partial _{s}^{2}\phi =-\frac{\left\vert x\right\vert ^{2}}{4s^{3}}-\left(
\alpha +1\right) \left( \frac{\sigma ^{\prime }\left( s\right) }{\sigma
\left( s\right) }\right) ^{\prime }\text{ , }\partial _{s}\left( \left\vert
\nabla _{g}\phi \right\vert _{g}^{2}\right) =-g^{ij}\left( x\right) \frac{%
x_{i}x_{j}}{8s^{3}}\text{ .}  \label{Rq2.232}
\end{equation}%
Observing that if $g^{ij}\left( x\right) =\delta ^{ij}$, $i,j=1,2,...n$, $%
x\in \mathbb{R}^{n}$, then the term under the integral sign at the
right-hand side of (\ref{Rq2.221}) is equal to $\left( \alpha +1\right) v^{2}%
\dfrac{\sigma ^{\prime }}{\sigma }\dfrac{d}{ds}\left( \log \dfrac{\sigma }{%
s\sigma ^{\prime }}\right) $, it is easy to check that
\begin{equation}
J_{1}\geq \left( \alpha +1\right) \int s^{2}v^{2}\dfrac{\sigma ^{\prime }}{%
\sigma }\dfrac{d}{ds}\left( \log \dfrac{\sigma }{s\sigma ^{\prime }}\right)
dX-C\Lambda \int \frac{\left\vert x\right\vert ^{3}}{s}v^{2}dX\text{ ,}
\label{Rq2.235}
\end{equation}%
where $C$ is an absolute constant.

\noindent Likewise, observing that if $g^{ij}\left( x\right) =\delta ^{ij}$,
$i,j=1,2,...n$, $\in \mathbb{R}^{n}$, then the term under the integral sign
at the right-hand side of (\ref{Rq2.222}) vanishes, so that the following
inequality holds true
\begin{equation}
J_{2}\geq -C\Lambda \int \left\vert x\right\vert s\left\vert \nabla
_{g}v\right\vert _{g}^{2}dX\text{ ,}  \label{Rq2.240}
\end{equation}
where $C$ is an absolute constant.

\noindent Now, let us consider $J_{3}$. Let $\delta \in \left( 0,1\right) $
be a number that we shall choose later, let $\sigma $ be the function (\ref%
{Rq2.135}) where $\gamma =\dfrac{\alpha }{\delta ^{2}}$, $\alpha \geq 1$ and
$\theta $ is given by (\ref{Rq2.160}). We have
\begin{equation*}
\left\vert \psi _{1}-s\psi _{0}\right\vert \leq \left\vert \psi
_{1}\right\vert +s\left\vert \psi _{0}\right\vert \leq C\left( \left\vert
x\right\vert +\frac{\left\vert x\right\vert ^{3}}{s}+\left( \alpha +1\right)
s\left\vert x\right\vert \frac{\sigma ^{\prime }}{\sigma }\right) \text{ ,}
\end{equation*}%
where $C$ depends on $\lambda $, $\Lambda $ only. By (\ref{Rq2.145}) and (%
\ref{Rq2.150}) we have
\begin{equation}
J_{3}\geq -C\int \left( \frac{\left\vert x\right\vert ^{3}}{s}+\left(
\alpha +1\right) \left\vert x\right\vert \right) v^{2}dX\text{ ,}
\label{Rq2.245}
\end{equation}%
for every $v\in C_{0}^{\infty }\left( \mathbb{R}^{n}\times \left[ a,+\infty
\right) \right) $ such that \textit{supp\thinspace }$v\subset B_{1}\times %
\left[ a,\frac{1}{\gamma }\right) $, where $C$ depends on $\lambda $ and $%
\Lambda $ only.

Now by Lemma \ref{LeRq7} and by (\ref{Rq2.235}), (\ref{Rq2.240}), (\ref%
{Rq2.245}) we get
\begin{eqnarray}
I &\geq &\frac{\left( \alpha +1\right) }{e^{C_{0}}}\int \theta \left(
\gamma s\right) v^{2}dX-C\int \left( \left( \alpha +1\right) \left\vert
x\right\vert +\frac{\left\vert x\right\vert ^{3}}{s}\right) v^{2}dX
\label{Rq2.250} \\
&&-C\int \left\vert x\right\vert s\left\vert \nabla _{g}v\right\vert
_{g}^{2}dX+\mathcal{B}_{1}\text{ ,}  \notag
\end{eqnarray}%
for every $v\in C_{0}^{\infty }\left( \mathbb{R}^{n}\times \left[ a,+\infty
\right) \right) $ such that \textit{supp\thinspace }$v\subset B_{1}\times %
\left[ a,\frac{1}{\gamma }\right) $, $\alpha \geq 1$, where $C_{0}$ is
defined in Lemma \ref{LeRq7} and $C$ depends on $\lambda $ and $\Lambda $
only.

\bigskip Let $\varepsilon _{0}\in \left( 0,1\right] $ be a number that we
shall choose later. We get
\begin{equation}
\int \left\vert Sv\right\vert ^{2}dX\geq \varepsilon _{0}\int \left\vert
Sv\right\vert ^{2}dX=\varepsilon _{0}\int \left\vert \left( Sv+\theta
\left( \gamma s\right) v\right) -\theta \left( \gamma s\right) v\right\vert
^{2}dX  \label{Rq2.255}
\end{equation}
\begin{equation*}
\geq -2\varepsilon _{0}\int s\theta \left( \gamma s\right) v\Delta
_{g}vdX-2\varepsilon _{0}\int \left( s\left( \left\vert \nabla _{g}\phi
\right\vert _{g}^{2}-\partial _{s}\phi \right) -\frac{1}{2}+\theta \left(
\gamma s\right) \right) \theta \left( \gamma s\right) v^{2}dX
\end{equation*}
\begin{equation*}
=2\varepsilon _{0}\int s\theta \left( \gamma s\right) \left( \left\vert
\nabla _{g}v\right\vert _{g}^{2}+\left\vert \nabla _{g}\phi \right\vert
_{g}^{2}v^{2}\right) dX
\end{equation*}
\begin{equation*}
-2\varepsilon _{0}\int \left( s\left( 2\left\vert \nabla _{g}\phi
\right\vert _{g}^{2}-\partial _{s}\phi \right) -\frac{1}{2}+\theta \left(
\gamma s\right) \right) \theta \left( \gamma s\right) v^{2}dX\text{ .}
\end{equation*}
Now by (\ref{Rq2.145}), (\ref{Rq2.150}) and (\ref{Rq2.230})-(\ref{Rq2.232})
we obtain
\begin{eqnarray}
&&\left\vert \left( s\left( 2\left\vert \nabla _{g}\phi \right\vert
_{g}^{2}-\partial _{s}\phi \right) -\frac{1}{2}+\theta \left( \gamma
s\right) \right) \theta \left( \gamma s\right) \right\vert  \label{Rq2.260}
\\
&\leq &C\left( \frac{\left\vert x\right\vert ^{3}}{s}+\left( \alpha
+1\right) \theta \left( \gamma s\right) \right) \text{ , for every }0<\gamma
s\leq \frac{1}{2}\text{, }\alpha \geq 1\text{ ,}  \notag
\end{eqnarray}
where $C$ depends on $\lambda $ and $\Lambda $ only.

\noindent By (\ref{Rq2.255}) and (\ref{Rq2.260}) we have
\begin{eqnarray}
\int \left\vert Sv\right\vert ^{2}dX &\geq &2\varepsilon _{0}\int s\theta
\left( \gamma s\right) \left( \left\vert \nabla _{g}v\right\vert
_{g}^{2}+\left\vert \nabla _{g}\phi \right\vert _{g}^{2}v^{2}\right) dX
\label{Rq2.264} \\
&&-C_{1}\varepsilon _{0}\int \left( \frac{\left\vert x\right\vert ^{3}}{s}%
+\left( \alpha +1\right) \theta \left( \gamma s\right) \right) v^{2}dX\text{
,}  \notag
\end{eqnarray}%
for every $v\in C_{0}^{\infty }\left( \mathbb{R}^{n}\times \left[ a,+\infty
\right) \right) $ such that \textit{supp\thinspace }$v\subset B_{1}\times %
\left[ a,\frac{1}{2\gamma }\right) $, $\alpha \geq 1$, where $C_{1}$ depends
on $\lambda $ and $\Lambda $ only.

\noindent By the just obtained inequality and (\ref{Rq2.250}) we get
\begin{equation}
I+\int \left\vert Sv\right\vert ^{2}dX\geq \left( \alpha +1\right) \left(
\frac{1}{e^{C_{0}}}-C_{1}\varepsilon _{0}\right) \int \theta \left( \gamma
s\right) v^{2}dX  \label{Rq2.265}
\end{equation}%
\begin{equation*}
+2\varepsilon _{0}\int s\theta \left( \gamma s\right) \left( \left\vert
\nabla _{g}v\right\vert _{g}^{2}+\left\vert \nabla _{g}\phi \right\vert
_{g}^{2}v^{2}\right) dX-C\int \left( \left( \alpha +1\right) \left\vert
x\right\vert +\frac{\left\vert x\right\vert ^{3}}{s}\right) v^{2}dX
\end{equation*}%
\begin{equation*}
-C\int \left\vert x\right\vert s\left( \left\vert \nabla _{g}v\right\vert
_{g}^{2}+\left\vert \nabla _{g}\phi \right\vert _{g}^{2}v^{2}\right) +%
\mathcal{B}_{1}\text{ ,}
\end{equation*}%
for every $v\in C_{0}^{\infty }\left( \mathbb{R}^{n}\times \left[ a,+\infty
\right) \right) $ such that \textit{supp}$\,v\subset B_{1}\times \left[ a,%
\frac{1}{2\gamma }\right) $, $\alpha \geq 1$, where $C$ depends on $\lambda $
and $\Lambda $ only and $C_{1}$ is the same constant that appears in (\ref%
{Rq2.264}).

\noindent Now we choose $\varepsilon _{0}=\dfrac{1}{2C_{1}e^{C_{0}}}$ at the
right-hand side of (\ref{Rq2.265}), then, recalling (\ref{Rq2.185b}) and
coming back to the function $u$, we have
\begin{equation}
\int s^{2}\left\vert P_{0}u\right\vert ^{2}e^{2\phi }dX\geq \frac{\left(
\alpha +1\right) }{C}\int \theta \left( \gamma s\right) u^{2}e^{2\phi }dX
\label{Rq2.270}
\end{equation}%
\begin{equation*}
+\frac{1}{C}\int s\theta \left( \gamma s\right) \left\vert \nabla
_{g}u\right\vert _{g}^{2}e^{2\phi }dX+\mathcal{B}_{1}
\end{equation*}%
\begin{equation*}
-C\int \left( \left( \alpha +1\right) \left\vert x\right\vert +\frac{%
\left\vert x\right\vert ^{3}}{s}\right) u^{2}e^{2\phi }dX-C\int \left\vert
x\right\vert s\left\vert \nabla _{g}u\right\vert _{g}^{2}e^{2\phi }dX\text{ ,%
}
\end{equation*}%
for every $u\in C_{0}^{\infty }\left( \mathbb{R}^{n}\times \left[ a,+\infty
\right) \right) $ such that \textit{supp}$\,u\subset B_{1}\times \left[ a,%
\frac{1}{2\gamma }\right) $, $\alpha \geq 1$, where $C$, $C>1$, depends on $%
\lambda $ and $\Lambda $ only.

\noindent Now, we use Lemma \ref{LeRq8} to prove the following estimates
\begin{equation}
\int \left( \left( 1+\alpha \right) \left\vert x\right\vert +\frac{%
\left\vert x\right\vert ^{3}}{s}\right) u^{2}e^{2\phi }dX\leq C\left(
e^{C_{0}}\gamma \right) ^{2\alpha +\frac{5}{2}}\int u^{2}dX  \label{Rq2.275}
\end{equation}%
\begin{equation*}
+C\alpha \delta \int \theta \left( \gamma s\right) u^{2}e^{2\phi }dX\text{ ,%
}
\end{equation*}%
\begin{equation}
\int \left\vert x\right\vert s\left\vert \nabla _{g}u\right\vert
_{g}^{2}e^{2\phi }dX\leq C\left( e^{C_{0}}\gamma \right) ^{2\alpha +\frac{3}{%
2}}\int s\left\vert \nabla _{g}u\right\vert _{g}^{2}dX  \label{Rq2.280}
\end{equation}%
\begin{equation*}
+C\delta \int \theta \left( \gamma s\right) s\left\vert \nabla
_{g}u\right\vert _{g}^{2}e^{2\phi }dX\text{ ,}
\end{equation*}%
where $C$ is an absolute constant.

\noindent By applying Lemma \ref{LeRq8} with $\mu =\dfrac{1}{2}$, $y=\frac{%
\left\vert x\right\vert ^{2}}{4s}$, $\varepsilon =\left( \gamma s\right) ^{%
\frac{3}{2}+2\alpha }$ and by using (\ref{Rq2.145}) we have
\begin{equation}
\left\vert x\right\vert e^{2\phi }=\left( 4s\right) ^{1/2}\left( \frac{%
\left\vert x\right\vert ^{2}}{4s}\right) ^{1/2}e^{-\frac{\left\vert
x\right\vert ^{2}}{4s}}\sigma ^{-2\left( \alpha +1\right) }  \label{Rq2.285}
\end{equation}

\begin{equation*}
\leq C\left( \left( \gamma s\right) ^{\alpha +\frac{3}{2}}s^{1/2}\sigma
^{-2\left( \alpha +1\right) }+s^{1/2}\left( \left( \frac{3}{2}+2\alpha
\right) \log \frac{1}{\gamma s}\right) ^{1/2}\sigma ^{-2\left( \alpha
+1\right) }e^{-\frac{\left\vert x\right\vert ^{2}}{4s}}\right)
\end{equation*}
\begin{equation*}
\leq C^{\prime }\left( \left( e^{C_{0}}\gamma \right) ^{2\alpha +\frac{5}{2}%
}+\delta \theta \left( \gamma s\right) e^{2\phi }\right) \text{ ,}
\end{equation*}
hence (\ref{Rq2.280}) follows.

\noindent Likewise, by applying Lemma \ref{LeRq8} with $\mu =\dfrac{3}{2}$, $%
y=\frac{\left\vert x\right\vert ^{2}}{4s}$, $\varepsilon =\left( \gamma
s\right) ^{\frac{3}{2}+2\alpha }$ we have
\begin{equation}
\frac{\left\vert x\right\vert ^{3}}{s}e^{2\phi }=4\left( 4s\right)
^{1/2}\left( \frac{\left\vert x\right\vert ^{2}}{4s}\right) ^{3/2}e^{-\frac{%
\left\vert x\right\vert ^{2}}{4s}}\sigma ^{-2\left( \alpha +1\right) }
\label{Rq2.290b}
\end{equation}

\begin{equation*}
\leq C\left( \left( e^{C_{0}}\gamma \right) ^{2\alpha +2}+\alpha \delta
\theta \left( \gamma s\right) e^{2\phi }\right) \text{ .}
\end{equation*}
By (\ref{Rq2.285}) and (\ref{Rq2.290b}) we get (\ref{Rq2.275}).

Finally, by (\ref{Rq2.270}), (\ref{Rq2.275}) and (\ref{Rq2.280}) we have
that there exists $\delta _{0}\in \left( 0,1\right] $ such that for every $%
\delta \in \left( 0,\delta _{0}\right] $, $\alpha \geq 1$, $u\in
C_{0}^{\infty }\left( \mathbb{R}^{n}\times \left[ a,+\infty \right) \right) $
such that \textit{supp}$\,u\subset B_{1}\times \left[ a,\frac{1}{2\gamma }%
\right) $, the following estimate holds true
\begin{eqnarray}
&&\int s^{2}\left\vert P_{0}u\right\vert ^{2}e^{2\phi }dX+C\left(
e^{C_{0}}\gamma \right) ^{2\alpha +\frac{5}{2}}\int \left(
u^{2}+s\left\vert \nabla _{g}u\right\vert _{g}^{2}\right) e^{2\phi }dX
\label{Rq2.295} \\
&\geq &\frac{\alpha +1}{C}\int \theta \left( \gamma s\right) u^{2}e^{2\phi
}dX+\frac{1}{C}\int \theta \left( \gamma s\right) s\left\vert \nabla
_{g}u\right\vert _{g}^{2}e^{2\phi }dX+\mathcal{B}_{1}\text{,}  \notag
\end{eqnarray}%
where $C\geq 1$ depends on $\lambda $ and $\Lambda $ only.

\underline{STEP 2.} It is simple to check the following identity, for $k\in
\left\{ 1,...,n\right\} $%
\begin{eqnarray*}
P_{0}\left( \left( \partial _{k}u\right) ^{2}\right) &=&2\partial _{k}\left(
\partial _{k}uP_{0}u\right) -2\partial _{i}\left( \partial _{k}u\left(
\partial _{k}g^{ij}\left( x,a\right) \right) \partial _{j}u\right) \\
&&-2\partial _{kk}^{2}uP_{0}u+2\left( \partial _{k}g^{ij}\left( x,a\right)
\right) \partial _{j}u\partial _{ik}^{2}u+2g^{ij}\left( x,a\right) \partial
_{ik}^{2}u\partial _{jk}^{2}u\text{ ,}
\end{eqnarray*}%
(recall that $\left( P_{0}u\right) \left( x,s\right) =\partial _{j}\left(
g^{ij}\left( x,a\right) \partial _{i}u\left( x,s\right) \right) +\partial
_{s}u\left( x,s\right) $). Now we multiply both sides of the above
inequality by $\sigma ^{1-2\alpha }e^{-\frac{\left\vert x\right\vert ^{2}}{4s%
}}$ and we integrate by parts. We obtain, for any $k\in \left\{
1,...,n\right\} $,
\begin{equation}
2\int \left( g^{ij}\left( x,a\right) \partial _{ik}^{2}u\partial
_{jk}^{2}u\right) \sigma ^{1-2\alpha }e^{-\frac{\left\vert x\right\vert ^{2}%
}{4s}}dX  \label{Rq2.300}
\end{equation}%
\begin{equation*}
\leq \int \left( \partial _{k}u\right) ^{2}\left\vert P_{0}^{\ast }\left(
\sigma ^{1-2\alpha }e^{-\frac{\left\vert x\right\vert ^{2}}{4s}}\right)
\right\vert dX+2\int \left\vert \partial _{k}uP_{0}u\right\vert \left\vert
\partial _{k}\left( \sigma ^{1-2\alpha }e^{-\frac{\left\vert x\right\vert
^{2}}{4s}}\right) \right\vert dX
\end{equation*}%
\begin{equation*}
+2\int \left\vert \partial _{k}g^{ij}\left( x,a\right) \partial
_{j}u\partial _{k}u\right\vert \left\vert \partial _{i}\left( \sigma
^{1-2\alpha }e^{-\frac{\left\vert x\right\vert ^{2}}{4s}}\right) \right\vert
dX+2\int \left\vert \partial _{kk}^{2}uP_{0}u\right\vert \sigma ^{1-2\alpha
}e^{-\frac{\left\vert x\right\vert ^{2}}{4s}}dX
\end{equation*}%
\begin{equation*}
+2\int \left\vert \partial _{k}g^{ij}\left( x,a\right) \partial
_{j}u\partial _{ik}^{2}u\right\vert \sigma ^{1-2\alpha }e^{-\frac{\left\vert
x\right\vert ^{2}}{4s}}dX-\int\nolimits_{\mathbb{R}^{n}}\left( \partial
_{k}u\left( x,a\right) \right) ^{2}\sigma ^{1-2\alpha }\left( a\right) e^{-%
\frac{\left\vert x\right\vert ^{2}}{4a}}dx\text{ ,}
\end{equation*}%
where $P_{0}^{\ast }=\partial _{j}\left( g^{ij}\left( x,a\right) \partial
_{i}u\left( x,s\right) \right) -\partial _{s}u\left( x,s\right) $.

\noindent We have
\begin{eqnarray*}
\left\vert P_{0}^{\ast }\left( \sigma ^{1-2\alpha }e^{-\frac{\left\vert
x\right\vert ^{2}}{4s}}\right) \right\vert &=&\left\vert P_{0}^{\ast }\left(
\left( \sigma ^{1-2\alpha }s^{n/2}\right) \left( s^{-n/2}e^{-\frac{%
\left\vert x\right\vert ^{2}}{4s}}\right) \right) \right\vert \\
&\leq &C\left( \left( \frac{\left\vert x\right\vert }{s}+\frac{\left\vert
x\right\vert ^{3}}{s^{2}}\right) \sigma ^{1-2\alpha }+\alpha \sigma
^{-2\alpha }\right) e^{-\frac{\left\vert x\right\vert ^{2}}{4s}}
\end{eqnarray*}
and
\begin{equation*}
\left\vert \nabla _{x}\left( \sigma ^{1-2\alpha }e^{-\frac{\left\vert
x\right\vert ^{2}}{4s}}\right) \right\vert \leq C^{\prime }\frac{\left\vert
x\right\vert }{s}\sigma ^{1-2\alpha }e^{-\frac{\left\vert x\right\vert ^{2}}{%
4s}}\text{ ,}
\end{equation*}
where $C$ depends on $\lambda $ and $\Lambda $ only and $C^{\prime }$ is an
absolute constant.

\noindent By the above inequalities and (\ref{Rq2.300}) we get
\begin{equation}
\int \left\vert D^{2}u\right\vert ^{2}\sigma ^{1-2\alpha }e^{-\frac{%
\left\vert x\right\vert ^{2}}{4s}}dX\leq C\alpha \int \sigma ^{-2\alpha
}\left\vert \nabla _{g}u\right\vert _{g}^{2}e^{-\frac{\left\vert
x\right\vert ^{2}}{4s}}dX  \label{Rq2.305}
\end{equation}
\begin{equation*}
+C\int \left( \frac{\left\vert x\right\vert }{s}+\frac{\left\vert
x\right\vert ^{3}}{s^{2}}\right) \sigma ^{1-2\alpha }\left\vert \nabla
_{g}u\right\vert _{g}^{2}e^{-\frac{\left\vert x\right\vert ^{2}}{4s}%
}dX+C\int \frac{\left\vert x\right\vert }{s}\sigma ^{1-2\alpha }\left\vert
\nabla _{g}u\right\vert _{g}\left\vert P_{0}u\right\vert e^{-\frac{%
\left\vert x\right\vert ^{2}}{4s}}dX
\end{equation*}
\begin{equation*}
+C\int \left\vert D^{2}u\right\vert \left\vert P_{0}u\right\vert \sigma
^{1-2\alpha }e^{-\frac{\left\vert x\right\vert ^{2}}{4s}}dX+C\int
\left\vert D^{2}u\right\vert \left\vert \nabla _{g}u\right\vert _{g}\sigma
^{1-2\alpha }e^{-\frac{\left\vert x\right\vert ^{2}}{4s}}dX
\end{equation*}
\begin{equation*}
-\sum\limits_{k=1}^{n}\int\nolimits_{\mathbb{R}^{n}}\sigma ^{1-2\alpha
}\left( a\right) \left( \partial _{k}u\left( x,a\right) \right) ^{2}e^{-%
\frac{\left\vert x\right\vert ^{2}}{4a}}dx=I^{\left( 1\right) }+I^{\left(
2\right) }+I^{\left( 3\right) }+I^{\left( 4\right) }+I^{\left( 5\right) }-%
\mathcal{B}_{2}\text{ ,}
\end{equation*}
(here $\left\vert D^{2}u\right\vert ^{2}=\sum\limits_{i,j=1}^{n}\left(
\partial _{ij}^{2}u\right) ^{2}$) where $C$ depends on $\lambda $ and $%
\Lambda $ only.

\noindent By the trivial inequality $\dfrac{\theta \left( s\right) }{s}\geq
\dfrac{1}{C}$, for every $s\in \left( 0,\dfrac{1}{2}\right] $, $C\geq 1$, we
have $C\delta ^{2}\dfrac{\theta \left( \gamma s\right) }{s}\geq \alpha $,
whenever $0<\gamma s\leq \dfrac{1}{2}$. By the last inequality and (\ref%
{Rq2.145}) we have
\begin{equation}
I^{\left( 1\right) }\leq C\delta ^{2}\int s\theta \left( \gamma s\right)
\sigma ^{-2-2\alpha }\left\vert \nabla _{g}u\right\vert _{g}^{2}e^{-\frac{%
\left\vert x\right\vert ^{2}}{4s}}dX\text{ ,}  \label{Rq2.310}
\end{equation}%
for every $\alpha \geq 1$, $u\in C_{0}^{\infty }\left( \mathbb{R}^{n}\times %
\left[ a,+\infty \right) \right) $ such that \textit{supp}$\,u\subset
B_{1}\times \left[ a,\frac{1}{2\gamma }\right) $, where $C$ depends on $%
\lambda $ and $\Lambda $ only.

\noindent To estimate $I^{\left( 2\right) }$ we observe that (\ref{Rq2.135})
gives $\sigma \left( s\right) \leq \dfrac{C}{\gamma }$, whenever $0<\gamma
s\leq 1$. By this inequality, (\ref{Rq2.145}) and Lemma \ref{LeRq8} we have
\begin{equation}
I^{\left( 2\right) }\leq \dfrac{C}{\gamma }\int \left( \left\vert
x\right\vert +\frac{\left\vert x\right\vert ^{3}}{s}\right) \sigma
^{-1-2\alpha }\left\vert \nabla _{g}u\right\vert _{g}^{2}e^{-\frac{%
\left\vert x\right\vert ^{2}}{4s}}dX  \label{Rq2.320}
\end{equation}%
\begin{equation*}
\leq \dfrac{C}{\gamma }\left( e^{C_{0}}\gamma \right) ^{2\alpha +\frac{3}{2}%
}\int s\left\vert \nabla _{g}u\right\vert _{g}^{2}dX+C\delta ^{3}\int
s\theta \left( \gamma s\right) \sigma ^{-2-2\alpha }\left\vert \nabla
_{g}u\right\vert _{g}^{2}e^{-\frac{\left\vert x\right\vert ^{2}}{4s}}dX\text{
}
\end{equation*}%
for every $\alpha \geq 1$, $u\in C_{0}^{\infty }\left( \mathbb{R}^{n}\times %
\left[ a,+\infty \right) \right) $ such that \textit{supp}$\,u\subset
B_{1}\times \left[ a,\frac{1}{2\gamma }\right) $, where $C$ depends on $%
\lambda $ and $\Lambda $ only.

In order to estimate $I^{\left( 3\right) }$, $I^{\left( 4\right) }$ and $%
I^{\left( 5\right) }$ we use the inequalities $2ab\leq a^{2}+b^{2}$, $\sigma
\left( s\right) \leq \dfrac{C}{\gamma }$, $s\leq \theta \left( \gamma
s\right) $, when $0<\gamma s\leq \frac{1}{2}$, and Lemma \ref{LeRq8} and we
obtain, for every $\alpha \geq 1$, $u\in C_{0}^{\infty }\left( \mathbb{R}%
^{n}\times \left[ a,+\infty \right) \right) $ such that \textit{supp}$%
\,u\subset B_{1}\times \left[ a,\frac{1}{2\gamma }\right) $

\begin{eqnarray}
I^{\left( 3\right) } &\leq &C\left( e^{C_{0}}\gamma \right) ^{2\alpha
+1}\int s\left\vert \nabla _{g}u\right\vert _{g}^{2}dX  \label{Rq2.325} \\
&&+C\delta ^{2}\int s\theta \left( \gamma s\right) \sigma ^{-2-2\alpha
}\left\vert \nabla _{g}u\right\vert _{g}^{2}e^{-\frac{\left\vert
x\right\vert ^{2}}{4s}}dX  \notag \\
&&+C\delta ^{2}\int s^{2}\sigma ^{-2-2\alpha }\left\vert P_{0}u\right\vert
^{2}e^{-\frac{\left\vert x\right\vert ^{2}}{4s}}dX,  \notag
\end{eqnarray}
\begin{eqnarray}
I^{\left( 4\right) } &\leq &\frac{1}{4}\int \left\vert D^{2}u\right\vert
^{2}\sigma ^{1-2\alpha }e^{-\frac{\left\vert x\right\vert ^{2}}{4s}}dX
\label{Rq2.330} \\
&&+C\delta ^{2}\int s^{2}\sigma ^{-2-2\alpha }\left\vert P_{0}u\right\vert
^{2}e^{-\frac{\left\vert x\right\vert ^{2}}{4s}}dX\text{ ,}  \notag
\end{eqnarray}
\begin{eqnarray}
I^{\left( 5\right) } &\leq &\frac{1}{4}\int \left\vert D^{2}u\right\vert
^{2}\sigma ^{1-2\alpha }e^{-\frac{\left\vert x\right\vert ^{2}}{4s}}dX
\label{Rq2.335} \\
&&+C\delta ^{2}\int s\theta \left( \gamma s\right) \sigma ^{-2-2\alpha
}\left\vert \nabla _{g}u\right\vert _{g}^{2}e^{-\frac{\left\vert
x\right\vert ^{2}}{4s}}dX\text{ ,}  \notag
\end{eqnarray}
where $C$ depends on $\lambda $ and $\Lambda $ only.

\noindent By (\ref{Rq2.305})-(\ref{Rq2.335}) we have
\begin{equation}
\frac{1}{2\delta ^{2}}\int \left\vert D^{2}u\right\vert ^{2}\sigma
^{1-2\alpha }e^{-\frac{\left\vert x\right\vert ^{2}}{4s}}dX\leq \dfrac{C}{%
\delta ^{2}}\left( e^{C_{0}}\gamma \right) ^{2\alpha +\frac{3}{2}}\int
s\left\vert \nabla _{g}u\right\vert _{g}^{2}dX  \label{Rq2.340}
\end{equation}%
\begin{equation*}
+C\int s\theta \left( \gamma s\right) \sigma ^{-2-2\alpha }\left\vert
\nabla _{g}u\right\vert _{g}^{2}e^{-\frac{\left\vert x\right\vert ^{2}}{4s}%
}dX+C\int s^{2}\sigma ^{-2-2\alpha }\left\vert P_{0}u\right\vert ^{2}e^{-%
\frac{\left\vert x\right\vert ^{2}}{4s}}dX\text{ ,}
\end{equation*}%
for every $\alpha \geq 1$, $u\in C_{0}^{\infty }\left( \mathbb{R}^{n}\times %
\left[ a,+\infty \right) \right) $ such that \textit{supp\thinspace }$%
u\subset B_{1}\times \left[ a,\frac{1}{2\gamma }\right) $, where $C$ depends
on $\lambda $ and $\Lambda $ only.

\noindent Now, by the inequality
\begin{equation*}
\left\vert P_{0}u\right\vert \leq \left\vert Pu\right\vert +C\left\vert
\nabla _{g}u\right\vert _{g}+C\sqrt{s}\left\vert D^{2}u\right\vert \text{ ,}
\end{equation*}%
where $C$ depends on $\lambda $ and $\Lambda $ only and by (\ref{Rq2.295})
and (\ref{Rq2.340}) we obtain
\begin{equation}
\mathcal{B}_{1}+\frac{1}{C\delta ^{2}}\int \left\vert D^{2}u\right\vert
^{2}\sigma ^{1-2\alpha }e^{-\frac{\left\vert x\right\vert ^{2}}{4s}}dX
\label{Rq2.345}
\end{equation}%
\begin{equation*}
+\frac{1}{C}\int s\theta \left( \gamma s\right) \sigma ^{-2-2\alpha
}\left\vert \nabla _{g}u\right\vert _{g}^{2}e^{-\frac{\left\vert
x\right\vert ^{2}}{4s}}dX+\frac{\left( \alpha +1\right) }{C}\int \theta
\left( \gamma s\right) \sigma ^{-2-2\alpha }u^{2}e^{-\frac{\left\vert
x\right\vert ^{2}}{4s}}dX
\end{equation*}%
\begin{equation*}
\leq \dfrac{C}{\delta ^{2}}\left( e^{C_{0}}\gamma \right) ^{2\alpha +\frac{5%
}{2}}\int \left( u^{2}+s\left\vert \nabla _{g}u\right\vert _{g}^{2}\right)
dX+C\int s^{2}\sigma ^{-2-2\alpha }\left\vert Pu\right\vert ^{2}e^{-\frac{%
\left\vert x\right\vert ^{2}}{4s}}dX
\end{equation*}%
\begin{equation*}
C\int \left\vert D^{2}u\right\vert ^{2}\sigma ^{1-2\alpha }e^{-\frac{%
\left\vert x\right\vert ^{2}}{4s}}dX+C\int s^{2}\sigma ^{-2-2\alpha
}\left\vert \nabla _{g}u\right\vert _{g}^{2}e^{-\frac{\left\vert
x\right\vert ^{2}}{4s}}dX\text{ ,}
\end{equation*}%
for every $\alpha \geq 1$, $u\in C_{0}^{\infty }\left( \mathbb{R}^{n}\times %
\left[ a,+\infty \right) \right) $ such that \textit{supp\thinspace }$%
u\subset B_{1}\times \left[ a,\frac{1}{2\gamma }\right) $, where $C$, $C\geq
1$, depends on $\lambda $ and $\Lambda $ only.

\noindent Observing that $\theta \left( \gamma s\right) \geq \dfrac{\gamma s%
}{C}$, whenever $0<\gamma s\leq \dfrac{1}{2}$, we obtain that, if $\delta $
is small enough and $\alpha \geq 1$, the third and the fourth term at the
right-hand side of (\ref{Rq2.345}) are absorbed by the second and the third
term at the left-hand side of (\ref{Rq2.345}). Therefore there exists $%
\delta _{1}\in \left( 0,\delta _{0}\right] $ such that for every $\delta \in
\left( 0,\delta _{1}\right] $, $\alpha \geq 1$ and $u\in C_{0}^{\infty
}\left( \mathbb{R}^{n}\times \left[ a,+\infty \right) \right) $ such that
\textit{supp}$u\subset B_{1}\times \left[ a,\frac{1}{2\gamma }\right) $, the
following inequality holds true
\begin{equation}
\frac{1}{C}\int s\theta \left( \gamma s\right) \sigma ^{-2-2\alpha
}\left\vert \nabla _{g}u\right\vert _{g}^{2}e^{-\frac{\left\vert
x\right\vert ^{2}}{4s}}dX  \label{Rq2.350}
\end{equation}
\begin{equation*}
+\frac{\left( \alpha +1\right) }{C}\int \theta \left( \gamma s\right)
\sigma ^{-2-2\alpha }u^{2}e^{-\frac{\left\vert x\right\vert ^{2}}{4s}}dX
\end{equation*}
\begin{equation*}
\leq \dfrac{C}{\delta ^{2}}\left( e^{C_{0}}\gamma \right) ^{2\alpha +2}\int
\left( u^{2}+s\left\vert \nabla _{g}u\right\vert _{g}^{2}\right) dX+C\int
s^{2}\sigma ^{-2-2\alpha }\left\vert Pu\right\vert ^{2}e^{-\frac{\left\vert
x\right\vert ^{2}}{4s}}dX-\mathcal{B}_{1}\text{ ,}
\end{equation*}
where $C\geq 1$ depends on $\lambda $ and $\Lambda $ only.

Now for a fixed $\delta \in \left( 0,\delta _{1}\right] $ and $a\in \left( 0,%
\dfrac{1}{4\gamma }\right] $ we estimate from above the term $-\mathcal{B}%
_{1}$ at the right-hand side of (\ref{Rq2.350}).

\noindent We have
\begin{eqnarray*}
-\mathcal{B}_{1} &=&-a^{2}\int\nolimits_{\mathbb{R}^{n}}\left\vert \left(
\nabla _{g}u\right) \left( x,a\right) +\left( \nabla _{g}\phi \right) \left(
x,a\right) u\left( x,a\right) \right\vert _{g\left( .,a\right) }^{2}e^{2\phi
\left( x,a\right) }dx \\
&&+a^{2}\int\nolimits_{\mathbb{R}^{n}}\left( \left\vert \left( \nabla
_{g}\phi \right) \left( x,a\right) \right\vert _{g\left( .,a\right)
}^{2}-\partial _{s}\phi \left( x,a\right) -\frac{1}{2a}\right) u^{2}\left(
x,a\right) e^{2\phi \left( x,a\right) }dx\text{ .}
\end{eqnarray*}
Let $\varepsilon \in \left( 0,1\right) $ be a number that we shall choose
later. We get
\begin{eqnarray*}
-\mathcal{B}_{1} &\leq &-\varepsilon a^{2}\int\nolimits_{\mathbb{R}%
^{n}}\left\vert \left( \nabla _{g}u\right) \left( x,a\right) \right\vert
_{g\left( .,a\right) }^{2}e^{2\phi \left( x,a\right) }dx \\
&&+a^{2}\int\nolimits_{\mathbb{R}^{n}}\left( \frac{1}{1-\varepsilon }%
\left\vert \left( \nabla _{g}\phi \right) \left( x,a\right) \right\vert
_{g\left( .,a\right) }^{2}-\partial _{s}\phi \left( x,a\right) -\frac{1}{2a}%
\right) u^{2}\left( x,a\right) e^{2\phi \left( x,a\right) }dx\text{.}
\end{eqnarray*}
Denote by $\eta _{0}=\min \left\{ 1,\frac{1}{2\Lambda }\right\} $. Let us
choose $\varepsilon =\frac{1}{4}$, we have, for every $x\in B_{\eta _{0}}$,
\begin{equation*}
\frac{1}{1-\varepsilon }\left\vert \left( \nabla _{g}\phi \right) \left(
x,a\right) \right\vert _{g\left( .,a\right) }^{2}-\partial _{s}\phi \left(
x,a\right) -\frac{1}{2a}
\end{equation*}
\begin{equation*}
\leq \frac{1}{3}\left( 2\Lambda \left\vert x\right\vert -1\right) \frac{%
\left\vert x\right\vert ^{2}}{8a^{2}}+\frac{\left( \alpha +1\right) \sigma
^{\prime }\left( a\right) }{\sigma \left( a\right) }\leq \frac{\left( \alpha
+1\right) }{e^{C_{0}}\sigma \left( a\right) }\text{.}
\end{equation*}

\noindent Therefore
\begin{equation*}
-\mathcal{B}_{1}\leq -\frac{a^{2}}{4}\int\nolimits_{\mathbb{R}%
^{n}}\left\vert \left( \nabla _{g}\right) \left( x,a\right) \right\vert
_{g\left( .,a\right) }^{2}e^{2\phi \left( x,a\right) }dx+\frac{\left( \alpha
+1\right) a^{2}}{e^{C_{0}}\sigma \left( a\right) }\int\nolimits_{\mathbb{R}%
^{n}}u^{2}\left( x,a\right) e^{2\phi \left( x,a\right) }dx\text{ .}
\end{equation*}
for every $\delta \in \left( 0,\delta _{1}\right] $, $a\in \left( 0,\dfrac{1%
}{4\gamma }\right] $, $\alpha \geq 1$ and $u\in C_{0}^{\infty }\left(
\mathbb{R}^{n}\times \left[ a,+\infty \right) \right) $ such that \textit{%
supp}$u\subset B_{\eta _{0}}\times \left[ a,\frac{1}{2\gamma }\right) $.
This inequality and (\ref{Rq2.350}) give (\ref{Rq2.165}).$\blacksquare $

\subsection{Two-sphere one-cylinder inequalities\label{Subs3.1}}

In the present section and in the next two sections we shall prove some
quantitative estimates of unique continuation for the forward parabolic
operator
\begin{equation}
Lu=\partial _{i}\left( g^{ij}\left( x,t\right) \partial _{j}u\right)
-\partial _{t}u\text{ .}  \label{L1}
\end{equation}%
Since such estimates are consequence of Theorem \ref{ThRq9}, for the sake of
simplicity it is better to present their proofs in terms of the "backward"
parabolic operators (\ref{Rq2.100}). On the contrary, since in the
applications to the stability of the inverse problem of Section \ref{Sec4}
we shall use the quantitative estimates of unique continuation for the
"forward" operator (\ref{L1}), we state such estimates with the forward
notation.

In order to prove next lemma we need some properties of the fundamental
solution $\Gamma \left( x,s;y,\tau \right) $ of the adjoint operator $%
P^{\ast }=\partial _{i}\left( g^{ij}\left( x,s\right) \partial _{j}u\right)
-\partial _{s}u$ of operator $P$ appearing in (\ref{Rq2.100}). We refer to
\cite{Aro} for the definition and the proofs of the properties of function $%
\Gamma \left( x,s;y,\tau \right) $. In what follows we recall some
properties of $\Gamma \left( x,s;y,\tau \right) $ that we shall use later on.

\noindent i) For every $\left( y,\tau \right) \in \mathbb{R}^{n+1}$, the
function $\Gamma \left( .,.;y,\tau \right) $ is a solution to the equation
\begin{equation}
P^{\ast }\left( \Gamma \left( .,.;y,\tau \right) \right) =0\text{ , in }%
\mathbb{R}^{n}\times \left( \tau ,+\infty \right) \text{ .}  \label{3.60}
\end{equation}
ii) For every $\left( x,s\right) ,\left( y,\tau \right) \in \mathbb{R}^{n+1}$
such that $\left( x,s\right) \neq \left( y,\tau \right) $ we have
\begin{equation}
\Gamma \left( x,s;y,\tau \right) \geq 0  \label{3.70}
\end{equation}
and
\begin{equation}
\Gamma \left( x,s;y,\tau \right) =0\text{ , for every }s<\tau \text{ .}
\label{3.80}
\end{equation}
iii) There exists a constant $C>1$ depending on $\lambda $ (and $n$) only
such that, for every $\left( x,s\right) ,\left( y,\tau \right) \in \mathbb{R}%
^{n+1}$, $s>\tau $ we have
\begin{equation}
\frac{C^{-1}}{\left( s-\tau \right) ^{n/2}}e^{-\frac{C\left\vert
x-y\right\vert ^{2}}{\left( s-\tau \right) }}\leq \Gamma \left( x,s;y,\tau
\right) \leq \frac{C}{\left( s-\tau \right) ^{n/2}}e^{-\frac{\left\vert
x-y\right\vert ^{2}}{C\left( s-\tau \right) }}\text{ .}  \label{3.90}
\end{equation}
iv) If $u_{0}$ is a function on $\mathbb{R}^{n}$ continuous at a point $y\in
\mathbb{R}^{n}$ which satisfies $e^{-\mu \left\vert x\right\vert
^{2}}u_{0}\in L^{2}\left( \mathbb{R}^{2}\right) $, for a positive number $%
\mu $, then
\begin{equation}
\lim\limits_{s\rightarrow 0^{+}}\int\nolimits_{\mathbb{R}^{n}}\Gamma \left(
x,s;y,0\right) u_{0}\left( x\right) dx=u_{0}\left( y\right) \text{ .}
\label{3.100}
\end{equation}
Moreover by standard regularity results \cite{LSU}, \cite{Li} we have that,
by (\ref{Rq2.105}) and (\ref{Rq2.110}), for every $\left( y,\tau \right) \in
\mathbb{R}^{n+1}$ the function $\Gamma \left( .,.;y,\tau \right) $ belongs
to $H_{loc}^{2,1}\left( \mathbb{R}^{n+1}\smallsetminus \left\{ \left( y,\tau
\right) \right\} \right) $.

\begin{lemma}
\label{Le3.2}Let $P$ be operator (\ref{Rq2.100}) whose coefficients satisfy (%
\ref{Rq2.105}) and (\ref{Rq2.110}). Assume that $u\in H^{2,1}\left(
B_{1}\times \left( 0,1\right) \right) $ satisfies the inequality
\begin{equation}
\left\vert Pu\right\vert \leq \Lambda \left( \left\vert \nabla u\right\vert
+\left\vert u\right\vert \right) \text{ , in }B_{1}\times \left[ 0,1\right)
\text{ .}  \label{3.110}
\end{equation}
Then there exists a constant $C$, $C>1$, depending on $\lambda $ and $%
\Lambda $ only such that for every $\rho _{0},\rho _{1},\rho _{2},T\in
\left( 0,\dfrac{1}{2}\right] $ satisfying $\rho _{0}<\rho _{1}<\rho _{2}$
the following inequality holds true
\begin{equation}
\int\nolimits_{B_{\rho _{2}}}u^{2}\left( x,s\right) dx\geq \frac{1}{C}%
\int\nolimits_{B_{\rho _{0}}}u^{2}\left( x,0\right) dx\text{ , for every }%
s\in \left[ 0,s_{0}\right] \text{ ,}  \label{3.111}
\end{equation}
where
\begin{equation}
s_{0}=\min \left\{ T,\frac{\left( \rho _{1}-\rho _{0}\right) ^{2}}{C}\left[
\left( \log \left( \widetilde{C}N\left( u\right) \right) \right) _{+}\right]
^{-1}\right\} \text{ ,}  \label{3.112}
\end{equation}
(here $x_{+}=\max \left\{ x,0\right\} $),
\begin{equation}
N\left( u\right) =\frac{\int\nolimits_{B_{\rho _{2}}\times \left(
0,2T\right) }u^{2}dX}{\int\nolimits_{B_{\rho _{0}}}u^{2}\left( x,0\right) dx%
}\text{ ,}  \label{3.112b}
\end{equation}
and
\begin{equation}
\widetilde{C}=\frac{C\left( \rho _{2}^{n}-\rho _{1}^{n}\right) \rho _{0}^{n}%
}{T\left( \rho _{1}-\rho _{0}\right) ^{n-1}\left( \rho _{2}-\rho _{1}\right)
^{2}}\text{ .}  \label{3.113}
\end{equation}
\end{lemma}

\textbf{Proof.} Let $\rho _{1},\rho _{2}$ satisfy $0<\rho _{1}<\rho _{2}\leq
1$ and let $\varphi $ be a function belonging to $C_{0}^{2}\left( \mathbb{R}%
^{n}\right) $ and satisfying $0\leq \varphi \leq 1$ in $\mathbb{R}^{n}$, $%
\varphi =1$ in $B_{\rho _{1}}$, $\varphi =0$ in $\mathbb{R}%
^{n}\smallsetminus B_{\rho _{2}}$ and
\begin{equation}
\left( \rho _{2}-\rho _{1}\right) \left\vert \nabla \varphi \right\vert
+\left( \rho _{2}-\rho _{1}\right) ^{2}\left\vert D^{2}\varphi \right\vert
\leq C\text{ , in }B_{\rho _{2}}\smallsetminus B_{\rho _{1}}\text{ ,}
\label{3.120}
\end{equation}
where $C$ is an absolute constant.

Denote
\begin{equation*}
v\left( x,t\right) =u\left( x,t\right) \varphi \left( x\right) \text{.}
\end{equation*}
By (\ref{Rq2.105}), (\ref{Rq2.110}), (\ref{3.110}) and (\ref{3.120}) we have
\begin{equation}
\left\vert Pv\right\vert \leq C\left( \left\vert \nabla v\right\vert
+\left\vert v\right\vert +E\chi _{B_{\rho _{2}}\smallsetminus B_{\rho
_{1}}}\right) \text{ ,}  \label{3.130}
\end{equation}
where
\begin{equation}
E=\left( \rho _{2}-\rho _{1}\right) ^{-2}\left\Vert u\right\Vert _{L^{\infty
}\left( B_{\rho _{2}}\times \left( 0,T\right) \right) }+\left( \rho
_{2}-\rho _{1}\right) ^{-1}\left\Vert \nabla u\right\Vert _{L^{\infty
}\left( B_{\rho _{2}}\times \left( 0,T\right) \right) }\text{ ,}
\label{3.140}
\end{equation}
where $\chi _{B_{\rho _{2}}\setminus B_{\rho _{1}}}$ is the characteristic
function of $B_{\rho _{2}}\smallsetminus B_{\rho _{1}}$ and $C$ depends on $%
\lambda $ and $\Lambda $ only.

Now, let $\rho _{0}\in \left( 0,\rho _{1}\right) $ and let $y$ be a fixed
point of $B_{\rho _{0}}$ and denote by
\begin{equation}
H\left( s\right) =\int\nolimits_{\mathbb{R}^{n}}v^{2}\left( x,s\right)
\Gamma \left( x,s;y,0\right) dx\text{ , }s>0\text{ .}  \label{3.150}
\end{equation}
By differentiating $H$ we get
\begin{equation}
\frac{dH\left( s\right) }{ds}=2\int\nolimits_{\mathbb{R}^{n}}\partial
_{s}v\left( x,s\right) v\left( x,s\right) \Gamma \left( x,s;y,0\right)
dx+\int\nolimits_{\mathbb{R}^{n}}v^{2}\left( x,s\right) \partial _{s}\Gamma
\left( x,s;y,0\right) dx  \label{3.160}
\end{equation}
\begin{equation*}
=2\int\nolimits_{\mathbb{R}^{n}}\left( Pv\left( x,s\right) \right) v\left(
x,s\right) \Gamma \left( x,s;y,0\right) dx+\int\nolimits_{\mathbb{R}%
^{n}}v^{2}\left( x,s\right) v\left( x,s\right) \partial _{s}\Gamma \left(
x,s;y,0\right) dx
\end{equation*}
\begin{equation*}
-2\int\nolimits_{\mathbb{R}^{n}}\left( \Delta _{g}v\right) \left(
x,s\right) v\left( x,s\right) \Gamma \left( x,s;y,0\right) dx\text{.}
\end{equation*}
By the identity $2\left( \Delta _{g}v\right) v=\Delta _{g}\left(
v^{2}\right) -2\left\vert \nabla _{g}v\right\vert _{g}^{2}$ and integrating
by parts we obtain
\begin{equation*}
2\int\nolimits_{\mathbb{R}^{n}}\left( \Delta _{g}v\right) \left( x,s\right)
v\left( x,s\right) \Gamma \left( x,s;y,0\right) dx
\end{equation*}
\begin{equation*}
=\int\nolimits_{\mathbb{R}^{n}}v^{2}\left( x,s\right) \Delta _{g}\Gamma
\left( x,s;y,0\right) dx-2\int\nolimits_{\mathbb{R}^{n}}\left\vert \nabla
_{g}v\left( x,s\right) \right\vert _{g}^{2}\Gamma \left( x,s;y,0\right) dx%
\text{.}
\end{equation*}
By plugging the just obtained identity at the right-hand side of (\ref{3.160}%
) and by using (\ref{3.60}) we have
\begin{equation}
\frac{dH\left( s\right) }{ds}=2\int\nolimits_{\mathbb{R}^{n}}\left(
Pv\left( x,s\right) \right) v\left( x,s\right) \Gamma \left( x,s;y,0\right)
dx+2\int\nolimits_{\mathbb{R}^{n}}\left\vert \nabla _{g}v\left( x,s\right)
\right\vert _{g}^{2}\Gamma \left( x,s;y,0\right) dx\text{ .}  \label{3.170}
\end{equation}
By using (\ref{3.130}) in the identity (\ref{3.170}) we get
\begin{eqnarray}
\frac{dH\left( s\right) }{ds} &\geq &+2\int\nolimits_{\mathbb{R}%
^{n}}\left\vert \nabla _{g}v\left( x,s\right) \right\vert _{g}^{2}\Gamma
\left( x,s;y,0\right) dx  \label{3.180} \\
&&-C\int\nolimits_{\mathbb{R}^{n}}\left\vert v\left( x,s\right) \right\vert
\left( \left\vert \nabla _{g}v\left( x,s\right) \right\vert _{g}+\left\vert
v\left( x,s\right) \right\vert +E\chi _{B_{\rho _{2}}\setminus B_{\rho
_{1}}}\right) \Gamma \left( x,s;y,0\right) dx\text{,}  \notag
\end{eqnarray}
where $C$ depends on $\lambda $ and $\Lambda $ only.

\noindent By using the inequality $2ab\leq \dfrac{a^{2}}{\varepsilon }%
+\varepsilon b^{2}$ to estimate from below the right-hand side of (\ref%
{3.180}), we have
\begin{equation}
\frac{dH\left( s\right) }{ds}\geq -CH\left( s\right) -E^{2}\int\nolimits_{%
\mathbb{R}^{n}}\chi _{B_{\rho _{2}}\setminus B_{\rho _{1}}}\Gamma \left(
x,s;y,0\right) dx\text{ ,}  \label{3.190}
\end{equation}
where $C$ depends on $\lambda $ and $\Lambda $ only.

\noindent Since $y\in B_{\rho _{0}}$, by (\ref{3.90}) and (\ref{3.190}) we
have
\begin{equation}
\frac{dH\left( s\right) }{ds}\geq -CH\left( s\right) -C_{1}E^{2}\left( \rho
_{2}^{n}-\rho _{1}^{n}\right) \frac{e^{-\frac{\left( \rho _{1}-\rho
_{0}\right) ^{2}}{4C_{1}s}}}{s^{n/2}}\text{ ,}  \label{3.200}
\end{equation}
where $C$ depends on $\lambda $ and $\Lambda $ only and $C_{1}$ depends on $%
\lambda $ only.

\noindent Now let us multiply both sides of (\ref{3.200}) by $e^{Cs}$ and
integrate the obtained inequality. We get, for every $\varepsilon \in \left(
0,T\right) $ and $s\in \left[ \varepsilon ,T\right] $,
\begin{equation}
H\left( s\right) e^{Cs}\geq H\left( \varepsilon \right) e^{C\varepsilon
}-C_{1}E^{2}\left( \rho _{2}^{n}-\rho _{1}^{n}\right)
\int\nolimits_{\varepsilon }^{s}\frac{e^{C\tau -\frac{\left( \rho _{1}-\rho
_{0}\right) ^{2}}{4C_{1}\tau }}}{\tau ^{n/2}}d\tau \text{ .}  \label{3.210}
\end{equation}%
By using (\ref{3.100}) and passing to the limit as $\varepsilon $ tends to $%
0 $ in the inequality (\ref{3.210}), we get
\begin{eqnarray}
H\left( s\right) e^{Cs} &\geq &u^{2}\left( y,0\right) -C_{1}E^{2}\left( \rho
_{2}^{n}-\rho _{1}^{n}\right) \int\nolimits_{\varepsilon }^{s}\frac{%
e^{C\tau -\frac{\left( \rho _{1}-\rho _{0}\right) ^{2}}{4C_{1}\tau }}}{\tau
^{n/2}}d\tau  \label{3.210b} \\
&\geq &u^{2}\left( y,0\right) -C_{1}E^{2}\frac{\left( \rho _{2}^{n}-\rho
_{1}^{n}\right) }{\left( \rho _{1}-\rho _{0}\right) ^{n-1}}e^{-\frac{\left(
\rho _{1}-\rho _{0}\right) ^{2}}{8C_{1}s}}\text{ , for every }s\in \left(
0,T_{0}\right] \text{ ,}  \notag
\end{eqnarray}%
where $T_{0}=\min \left\{ \frac{\left( \rho _{1}-\rho _{0}\right) ^{2}}{%
16C_{1}},T\right\} $ and $C_{1}$ depend on $\lambda $ only.

\noindent Now, integrating over $B_{\rho _{0}}$ the inequality (\ref{3.210b}%
) we obtain
\begin{eqnarray}
&&\int\nolimits_{B_{\rho _{0}}}dy\int\nolimits_{\mathbb{R}^{n}}v^{2}\left(
x,s\right) \Gamma \left( x,s;y,0\right) dx  \label{3.220} \\
&\geq &\frac{1}{C}\left( \int\nolimits_{B_{\rho _{0}}}u^{2}\left(
y,0\right) dy-CE^{2}\frac{\left( \rho _{2}^{n}-\rho _{1}^{n}\right) \rho
_{0}^{n}}{\left( \rho _{1}-\rho _{0}\right) ^{n-1}}e^{-\frac{\left( \rho
_{1}-\rho _{0}\right) ^{2}}{Cs}}\right) \text{ ,}  \notag
\end{eqnarray}
where $C$, $C>1$, depends on $\lambda $ and $\Lambda $ only.

\noindent On the other hand, by (\ref{3.90}) we have
\begin{eqnarray}
&&\int\nolimits_{B_{\rho _{0}}}dy\int\nolimits_{\mathbb{R}^{n}}v^{2}\left(
x,s\right) \Gamma \left( x,s;y,0\right) dx  \label{3.230} \\
&\leq &\int\nolimits_{\mathbb{R}^{n}}\left( v^{2}\left( x,s\right)
\int\nolimits_{\mathbb{R}^{n}}\frac{Ce^{-\frac{\left\vert x-y\right\vert
^{2}}{4Cs}}}{s^{n/2}}dy\right) dx\leq C^{\prime }\int\nolimits_{B_{\rho
_{2}}}u^{2}\left( x,s\right) dx  \notag
\end{eqnarray}%
and by standard regularity estimates we have
\begin{equation}
E^{2}\leq \frac{C}{\left( \rho _{2}-\rho _{1}\right) ^{2}}\frac{1}{\rho _{2}T%
}\int\nolimits_{B_{\rho _{2}}\times \left( 0,2T\right) }u^{2}dX\text{ ,}
\label{3.240}
\end{equation}%
where $C$ depends on $\lambda $ and $\Lambda $ only.

\noindent By (\ref{3.220}), (\ref{3.230}) and (\ref{3.240}) we obtain, for
every $s\in \left( 0,T_{0}\right] $,
\begin{equation}
\int\nolimits_{B_{\rho _{2}}}u^{2}\left( x,s\right) dx\geq \frac{1}{C}%
\int\nolimits_{B_{\rho _{0}}}u^{2}\left( x,0\right) dx-C_{3}e^{-\frac{%
\left( \rho _{1}-\rho _{0}\right) ^{2}}{Cs}}\int\nolimits_{B_{\rho
_{2}}\times \left( 0,2T\right) }u^{2}dX\text{,}  \label{3.240c}
\end{equation}
where
\begin{equation*}
C_{3}=\frac{C\left( \rho _{2}^{n}-\rho _{1}^{n}\right) \rho _{0}^{n}}{%
T\left( \rho _{1}-\rho _{0}\right) ^{n-1}\left( \rho _{2}-\rho _{1}\right)
^{2}}
\end{equation*}
and $C$, $C>1$, depends on $\lambda $ and $\Lambda $ only. By (\ref{3.240c})
we get (\ref{3.111}).$\blacksquare $

\bigskip

Now let us introduce a notation which we shall use in Lemma \ref{Pr3.3} and
in Theorem \ref{Th3.4} below.

Let $\alpha \geq 1$, $a\geq 0$, $k\geq 0$, $\rho >0$ be given numbers.
Denote
\begin{equation}
D_{\rho }^{\left( a\right) }=\left\{ \left( x,s\right) \in \mathbb{R}%
^{n}\times \left( 0,+\infty \right) :\frac{\left\vert x\right\vert ^{2}}{%
4\alpha \left( s+a\right) }+\log \left( s+a\right) \leq \log \frac{\rho ^{2}%
}{4\alpha }\right\}  \label{3.250}
\end{equation}
and, for $\rho >\left( 4\alpha a\right) ^{1/2}$,
\begin{equation}
L_{\rho }^{\left( a\right) }=\sup \left\{ \left( s+a\right) ^{-\left( \alpha
+k\right) }e^{-\frac{\left\vert x\right\vert ^{2}}{4\left( s+a\right) }%
}:\left( x,s\right) \in D_{2\rho }^{\left( a\right) }\smallsetminus D_{\rho
}^{\left( a\right) }\right\} \text{ .}  \label{3.260}
\end{equation}

\begin{lemma}
\label{Pr3.3}For any positive numbers $\alpha $, $a$, $k$ and $\rho >0$ such
that $\alpha \geq 1$ and $\rho >\left( 4\alpha a\right) ^{1/2}$ we have
\begin{equation}
L_{\rho }^{\left( a\right) }\leq c_{k}4^{\alpha }\left( \frac{\alpha }{\rho
^{2}}\right) ^{\alpha +k}\text{ ,}  \label{3.270}
\end{equation}
where $c_{k}$ depends on $k$ only.
\end{lemma}

\textbf{Proof.} First notice that $L_{\rho }^{\left( a\right) }\leq L_{\rho
}^{\left( 0\right) }$ for every $a\geq 0$, hence it is enough to prove (\ref%
{3.270}) for $a=0$.

Now, for any fixed $x\in \mathbb{R}^{n}$, the function $s\rightarrow
s^{-\alpha -k}e^{-\frac{\left\vert x\right\vert ^{2}}{4s}}$ attains the
maximum when $s=\frac{\left\vert x\right\vert ^{2}}{4\left( \alpha +k\right)
}$. Therefore in order to estimate $L_{\rho }^{\left( 0\right) }$ from above
we need to consider the intersction $\Gamma $ of the paraboloid $\left\{
\left( x,s\right) \in \mathbb{R}^{n+1}:s=\frac{\left\vert x\right\vert ^{2}}{%
4\left( \alpha +k\right) }\right\} $ with the set $D_{2\rho }^{\left(
a\right) }\smallsetminus D_{\rho }^{\left( a\right) }$. It is easy to check
that the projection of $\Gamma $ on the $s$-axis is equal to the interval $%
I=\left( \frac{\rho ^{2}}{4\alpha e^{1+\left( k/\alpha \right) }},\frac{%
\left( 2\rho \right) ^{2}}{4\alpha e^{1+\left( k/\alpha \right) }}\right] $.
Hence we have
\begin{eqnarray*}
L_{\rho }^{\left( 0\right) } &=&\sup \left\{ s^{-\alpha -k}e^{-\alpha
-k}:s\in I\right\} \\
&\leq &e^{k\left( k+1\right) }\left( 4\alpha \rho ^{-2}\right) ^{\alpha +k}%
\text{ .}
\end{eqnarray*}
$\blacksquare $

\bigskip

\begin{theorem}
\label{Th3.4}(\textbf{two-sphere one-cylinder inequality in the interior}).
Let $\lambda $, $\Lambda $ and $R$ be positive numbers, with $\lambda \in
\left( 0,1\right] $ and $t_{0}\in \mathbb{R}$. Let $L$ be the parabolic
operator
\begin{equation*}
L=\partial _{i}\left( g^{ij}\left( x,t\right) \partial _{j}\right) -\partial
_{t}\text{ ,}
\end{equation*}%
where $\left\{ g^{ij}\left( x,t\right) \right\} _{i,j=1}^{n}$ is a real
symmetric $n\times n$ matrix. When $\xi \in \mathbb{R}^{n}$, $\left(
x,t\right) ,\left( y,\tau \right) \in \mathbb{R}^{n+1}$ assume that
\begin{equation}
\lambda \left\vert \xi \right\vert ^{2}\leq
\sum\limits_{i,j=1}^{n}g^{ij}\left( x,t\right) \xi _{i}\xi _{j}\leq \lambda
^{-1}\left\vert \xi \right\vert ^{2}  \label{Rq2.105b}
\end{equation}%
and
\begin{equation}
\left( \sum\limits_{i,j=1}^{n}\left( g^{ij}\left( x,t\right) -g^{ij}\left(
y,\tau \right) \right) ^{2}\right) ^{1/2}\leq \frac{\Lambda }{R}\left(
\left\vert x-y\right\vert ^{2}+\left\vert t-\tau \right\vert \right) ^{1/2}%
\text{.}  \label{Rq2.110b}
\end{equation}%
Let $u$ be a function belonging to $H^{2,1}\left( B_{R}\times \left(
t_{0}-R^{2},t_{0}\right) \right) $ which satisfies the inequality
\begin{equation}
\left\vert Lu\right\vert \leq \Lambda \left( \frac{\left\vert \nabla
u\right\vert }{R}+\frac{\left\vert u\right\vert }{R^{2}}\right) \text{ , in }%
B_{R}\times \left( t_{0}-R^{2},t_{0}\right] \text{ .}  \label{3.280}
\end{equation}%
Then there exist constants $\eta _{1}\in \left( 0,1\right) $ and $C$, $C\geq
1$, depending on $\lambda $ and $\Lambda $ only such that for every $r$ and $%
\rho $ such that $0<r\leq \rho \leq \eta _{1}R$ we have
\begin{equation}
\int\nolimits_{B_{\rho }}u^{2}\left( x,t_{0}\right) dx\leq \frac{CR}{\rho }%
\left( R^{-2}\int\nolimits_{B_{R}\times \left( t_{0}-R^{2},t_{0}\right)
}u^{2}dX\right) ^{1-\theta _{1}}\left( \int\nolimits_{B_{r}}u^{2}\left(
x,t_{0}\right) dx\text{ }\right) ^{\theta _{1}}\text{,}  \label{3.281}
\end{equation}%
where
\begin{equation}
\theta _{1}=\frac{1}{C\log \frac{R}{r}}\text{ .}  \label{3.282}
\end{equation}
\end{theorem}

\textbf{Proof.} Denoting $\widetilde{u}\left( x,s\right) =u\left(
x,t_{0}-s\right) $, $\widetilde{g}^{-1}\left( x,s\right) =g^{-1}\left(
x,t_{0}-s\right) $, $P=\partial _{i}\left( \widetilde{g}^{ij}\left(
x,s\right) \partial _{j}\right) +\partial _{s}$, by (\ref{3.280}) it turns
out that $\widetilde{u}$ belongs to $H^{2,1}\left( B_{R}\times \left(
0,R^{2}\right) \right) $ and satisfies the differential inequality
\begin{equation*}
\left\vert P\widetilde{u}\right\vert \leq \Lambda \left( \frac{\left\vert
\nabla \widetilde{u}\right\vert }{R}+\frac{\left\vert \widetilde{u}%
\right\vert }{R^{2}}\right) \text{ , in }B_{R}\times \left[ 0,R^{2}\right)
\text{.}
\end{equation*}

By scaling we have that inequality (\ref{3.281}) is equivalent to
\begin{equation*}
\int\limits_{B_{\rho }}\widetilde{u}^{2}\left( x,0\right) dx\leq \frac{C}{%
\rho }\left( \int\limits_{B_{1}\times \left( 0,1\right) }\widetilde{u}%
^{2}dX\right) ^{1-\theta _{1}}\left( \int\limits_{B_{r}}\widetilde{u}%
^{2}\left( x,0\right) dx\text{ }\right) ^{\theta _{1}}\text{ ,}
\end{equation*}%
where $\theta _{1}=\frac{1}{C\log \frac{1}{r}}$. In what follows we drop the
sign " $\widetilde{}$ " over $u$ and the matrix $g^{-1}$. First we consider
the case $g^{ij}\left( 0,0\right) =\delta ^{ij}$ and we write the Carleman
estimate (\ref{Rq2.165}) in a slight different form which is more suitable
for our purposes. It is easy to check that, by (\ref{Rq2.145}) and by $%
C\delta ^{2}\dfrac{\theta \left( \gamma s\right) }{s}\geq \alpha $, whenever
$0<\gamma s\leq \dfrac{1}{2}$ and by (\ref{Rq2.165}), we have

\begin{equation}
\alpha ^{2}\int\nolimits_{\mathbb{R}_{+}^{n+1}}\sigma _{a}^{-\alpha
}u^{2}e^{-\frac{\left\vert x\right\vert ^{2}}{4\left( s+a\right) }}dX+\alpha
\int\nolimits_{\mathbb{R}_{+}^{n+1}}\sigma _{a}^{-\alpha }\left\vert \nabla
_{g}u\right\vert _{g}^{2}e^{-\frac{\left\vert x\right\vert ^{2}}{4\left(
s+a\right) }}dX  \label{Rq2.165b}
\end{equation}%
\begin{equation*}
\leq C\int\nolimits_{\mathbb{R}_{+}^{n+1}}\left( s+a\right) ^{2}\sigma
_{a}^{1-\alpha }\left\vert Pu\right\vert ^{2}e^{-\frac{\left\vert
x\right\vert ^{2}}{4\left( s+a\right) }}dX+C^{\alpha }\alpha ^{\alpha
}\int\nolimits_{\mathbb{R}_{+}^{n+1}}\left( u^{2}+\left( s+a\right)
\left\vert \nabla _{g}u\right\vert _{g}^{2}\right) dX
\end{equation*}%
\begin{equation*}
+C\alpha \sigma ^{-\alpha }\left( a\right) \int\nolimits_{\mathbb{R}%
^{n}}u^{2}\left( x,0\right) e^{-\frac{\left\vert x\right\vert ^{2}}{4a}}dx%
\text{ .}
\end{equation*}%
for every $\alpha \geq 2$, $0<a\leq \dfrac{T_{1}}{\alpha }$ and $u\in
C_{0}^{\infty }\left( \mathbb{R}^{n}\times \left[ 0,+\infty \right) \right) $
such that \textit{supp\thinspace }$u\subset B_{\eta _{0}}\times \left[ 0,%
\dfrac{3T_{1}}{\alpha }\right) $, where $T_{1}=\dfrac{\delta _{1}^{2}}{3}$,
and $\eta _{0}$, $\eta _{0}\in \left( 0,1\right) $, $C$, $C>1$, depend on $%
\lambda $ and $\Lambda $ only (recall that $\eta _{0}$ is defined in Theorem %
\ref{ThRq9}).

By using Friedrichs density theorem we can apply Carleman estimate (\ref%
{Rq2.165b}) to the function $v=u\varphi $, where $\varphi $ is a function
belonging to $C_{0}^{\infty }\left( B_{\eta _{0}}\times \left[ 0,\dfrac{%
3T_{1}}{\alpha }\right) \right) $ that we are going to define.

\noindent Let $\overline{R}_{0}=\min \left\{ \sqrt{T_{1}},\sqrt{e}\eta
_{0}\right\} $, $R_{1}\in \left( 0,\overline{R}_{0}\right] $. For any $%
\alpha \geq 2$ denote
\begin{equation*}
d_{1}=\log \frac{\left( R_{1}/2\right) ^{2}}{4\alpha }\text{ , }d_{2}=\log
\frac{R_{1}^{2}}{4\alpha }\text{ .}
\end{equation*}
Let $\psi _{1}$ be the function which is equal to $0$ in $\mathbb{%
R\smallsetminus }\left[ d_{1},d_{2}\right] $ and such that $\psi _{1}\left(
\tau \right) =\exp \frac{\left( d_{2}-d_{1}\right) ^{2}}{\left( d_{1}-\tau
\right) \left( \tau -d_{2}\right) }$, if $\tau \in \left( d_{1},d_{2}\right)
$. Denote by $\psi _{2}$ the function such that
\begin{equation*}
\psi _{2}\left( \tau \right) =\frac{\int\nolimits_{\tau }^{d_{2}}\psi
_{1}\left( \xi \right) d\xi }{\int\nolimits_{d_{1}}^{d_{2}}\psi _{1}\left(
\xi \right) d\xi }\text{ , if }\tau \in \mathbb{R}\text{ .}
\end{equation*}
It is easy to check that $\psi _{2}\left( \tau \right) =1$, for every $\tau
\in \left( -\infty ,d_{1}\right) $, $\psi _{2}\left( \tau \right) =0$, for
every $\tau \in \left( d_{2},+\infty \right) $ and
\begin{equation}
\left\vert \psi _{2}^{\prime }\left( \tau \right) \right\vert \leq \frac{C}{%
d_{2}-d_{1}}\text{ , }\left\vert \psi _{2}^{\prime \prime }\left( \tau
\right) \right\vert \leq \frac{C}{\left( d_{2}-d_{1}\right) ^{2}}\text{, if }%
\tau \in \left[ d_{1},d_{2}\right] \text{ ,}  \label{3.290}
\end{equation}
where $C$ is an absolute constant. Notice that the right-hand side of
inequalities does not depend on $\alpha $ and $R_{1}$.

\noindent Now, let us define
\begin{equation}
\varphi \left( x,s\right) =\psi _{2}\left( \frac{\left\vert x\right\vert ^{2}%
}{4\alpha \left( s+a\right) }+\log \left( s+a\right) \right) \text{ ,}
\label{3.310}
\end{equation}
for every $a\in \left( 0,\dfrac{T_{1}}{\alpha }\right] $ and $\alpha \geq 2$.

\noindent It is easy to check that $\varphi \in C_{0}^{\infty }\left(
B_{\eta _{0}}\times \left[ 0,\dfrac{3T_{1}}{\alpha }\right) \right) $, $%
\varphi =1$ for every $\left( x,s\right) \in D_{R_{1}/2}^{\left( a\right) }$
and $\varphi =0$ for every $\left( x,s\right) \in \left( B_{\eta _{0}}\times %
\left[ 0,\dfrac{3T_{1}}{\alpha }\right) \right) \smallsetminus
D_{R_{1}}^{\left( a\right) }$. Moreover by (\ref{3.280}) we have that the
function $v=u\varphi $ satisfies the inequality
\begin{equation}
\left\vert Pv\right\vert \leq \Lambda \left( \left\vert \nabla u\right\vert
+\left\vert u\right\vert \right) \chi _{D_{R_{1}}^{\left( a\right)
}}+C\left( \frac{\left\vert x\right\vert \left\vert \nabla _{g}u\right\vert
_{g}}{\alpha \left( s+a\right) }+\frac{\left\vert u\right\vert }{\left(
s+a\right) ^{2}}\right) \chi _{D_{R_{1}}^{\left( a\right) }\smallsetminus
D_{R_{1}/2}^{\left( a\right) }}\text{ ,}  \label{3.320}
\end{equation}
where $C$ depends on $\lambda $ and $\Lambda $ only.

Now we apply inequality (\ref{Rq2.165b}) to the function $v$. By (\ref{3.320}%
) we obtain, for every $\alpha \geq 2$,
\begin{equation}
\alpha ^{2}\int\nolimits_{D_{R_{1}/2}^{\left( a\right) }}\sigma
_{a}^{-\alpha }u^{2}e^{-\frac{\left\vert x\right\vert ^{2}}{4\left(
s+a\right) }}dX+\alpha \int\nolimits_{D_{R_{1}/2}^{\left( a\right) }}\sigma
_{a}^{-\alpha }\left\vert \nabla _{g}u\right\vert _{g}^{2}e^{-\frac{%
\left\vert x\right\vert ^{2}}{4\left( s+a\right) }}dX  \label{3.330}
\end{equation}
\begin{equation*}
\leq C\int\nolimits_{D_{R_{1}/2}^{\left( a\right) }}\sigma _{a}^{1-\alpha
}u^{2}e^{-\frac{\left\vert x\right\vert ^{2}}{4\left( s+a\right) }%
}dX+C\int\nolimits_{D_{R_{1}/2}^{\left( a\right) }}\sigma _{a}^{1-\alpha
}\left\vert \nabla _{g}u\right\vert _{g}^{2}e^{-\frac{\left\vert
x\right\vert ^{2}}{4\left( s+a\right) }}dX
\end{equation*}
\begin{equation*}
+C\int\nolimits_{D_{R_{1}}^{\left( a\right) }\smallsetminus
D_{R_{1}/2}^{\left( a\right) }}\sigma _{a}^{1-\alpha }\left( \frac{%
\left\vert x\right\vert ^{2}\left\vert \nabla _{g}u\right\vert _{g}^{2}}{%
\alpha ^{2}\left( s+a\right) ^{2}}+\frac{u^{2}}{\left( s+a\right) ^{4}}%
\right) e^{-\frac{\left\vert x\right\vert ^{2}}{4\left( s+a\right) }}dX+I%
\text{ ,}
\end{equation*}
where
\begin{equation}
I=C^{\alpha }\alpha ^{\alpha }\int\nolimits_{D_{R_{1}/2}^{\left( a\right)
}}\left\{ \left( \frac{\left\vert x\right\vert ^{2}}{\alpha ^{2}\left(
s+a\right) ^{2}}+1\right) u^{2}+\left\vert \nabla _{g}u\right\vert
_{g}^{2}\right\} dX  \label{3.331}
\end{equation}
\begin{equation*}
+C\alpha \left( \sigma \left( a\right) \right) ^{-\alpha
}\int\nolimits_{B_{1}}v^{2}\left( x,0\right) e^{-\frac{\left\vert
x\right\vert ^{2}}{4a}}dx
\end{equation*}
and $C$, $C>1$, depends on $\lambda $ and $\Lambda $ only.

\noindent Now, for $\alpha $ large enough, the first and the second integral
at the right-hand side of (\ref{3.330}) obtained above can be absorbed by
the left-hand side of (\ref{3.330}) and we have
\begin{eqnarray}
&&\alpha ^{2}\int\nolimits_{D_{R_{1}/2}^{\left( a\right) }}\sigma
_{a}^{-\alpha }u^{2}e^{-\frac{\left\vert x\right\vert ^{2}}{4\left(
s+a\right) }}dX  \label{3.340} \\
&\leq &C\int\nolimits_{D_{R_{1}}^{\left( a\right) }\smallsetminus
D_{R_{1}/2}^{\left( a\right) }}\sigma _{a}^{1-\alpha }\left( \frac{%
\left\vert x\right\vert ^{2}\left\vert \nabla _{g}u\right\vert _{g}^{2}}{%
\alpha ^{2}\left( s+a\right) ^{2}}+\frac{u^{2}}{\left( s+a\right) ^{4}}%
\right) e^{-\frac{\left\vert x\right\vert ^{2}}{4\left( s+a\right) }}dX+I%
\text{ ,}  \notag
\end{eqnarray}
for every $\alpha \geq C$, where $C$ depends on $\lambda $ and $\Lambda $
only.

\noindent In order to estimate from above the first integral at the
right-hand side of (\ref{3.340}) we use (\ref{Rq2.145}) and Lemma \ref{Pr3.3}%
, so we get
\begin{equation}
\int\nolimits_{D_{R_{1}}^{\left( a\right) }\smallsetminus
D_{R_{1}/2}^{\left( a\right) }}\sigma _{a}^{1-\alpha }\left( \frac{%
\left\vert x\right\vert ^{2}\left\vert \nabla _{g}u\right\vert _{g}^{2}}{%
\alpha ^{2}\left( s+a\right) ^{2}}+\frac{u^{2}}{\left( s+a\right) ^{4}}%
\right) e^{-\frac{\left\vert x\right\vert ^{2}}{4\left( s+a\right) }}dX
\label{3.350}
\end{equation}
\begin{eqnarray*}
&\leq &C^{\alpha }\int\nolimits_{D_{R_{1}}^{\left( a\right) }\smallsetminus
D_{R_{1}/2}^{\left( a\right) }}\left( s+a\right) ^{-\alpha -3}\left(
\left\vert \nabla _{g}u\right\vert _{g}^{2}+u^{2}\right) e^{-\frac{%
\left\vert x\right\vert ^{2}}{4\left( s+a\right) }}dX \\
&\leq &\frac{C^{\prime \alpha }\alpha ^{\alpha }}{R_{1}^{2\left( \alpha
+3\right) }}\int\nolimits_{D_{R_{1}}^{\left( a\right) }\smallsetminus
D_{R_{1}/2}^{\left( a\right) }}\left( \left\vert \nabla _{g}u\right\vert
_{g}^{2}+u^{2}\right) dX\text{ ,}
\end{eqnarray*}
for every $\alpha \geq C$ and every $a\in \left( 0,\frac{R_{1}^{2}}{16\alpha
}\right) $, where $C$, $C>1$, and $C^{\prime }$, $C^{\prime }>1$, depend on $%
\lambda $ and $\Lambda $ only.

\noindent By the inequality obtained in (\ref{3.350}) and by (\ref{3.340})
we have
\begin{eqnarray}
&&\alpha ^{2}\int\nolimits_{D_{R_{1}/2}^{\left( a\right) }}\sigma
_{a}^{-\alpha }u^{2}e^{-\frac{\left\vert x\right\vert ^{2}}{4\left(
s+a\right) }}dX  \label{3.360} \\
&\leq &\frac{C^{\alpha }\alpha ^{\alpha }}{R_{1}^{2\left( \alpha +3\right) }}%
\int\nolimits_{D_{R_{1}}^{\left( a\right) }\smallsetminus
D_{R_{1}/2}^{\left( a\right) }}\left( \left\vert \nabla _{g}u\right\vert
_{g}^{2}+u^{2}\right) dX+I\text{ ,}  \notag
\end{eqnarray}
for every $\alpha \geq C$ and every $a\in \left( 0,\dfrac{R_{1}^{2}}{%
16\alpha }\right) $, where $C>1$ depends on $\lambda $ and $\Lambda $ only.

\noindent Now we estimate from above the term $I$ (defined in (\ref{3.331}))
at the right-hand side of (\ref{3.360}). Concerning the first integral at
the right-hand side of formula (\ref{3.331}) it is simple to check that
\begin{eqnarray}
&&\int\nolimits_{D_{R_{1}}^{\left( a\right) }}\left\{ \left( \frac{%
\left\vert x\right\vert ^{2}}{\alpha ^{2}\left( s+a\right) ^{2}}+1\right)
u^{2}+\left\vert \nabla _{g}u\right\vert _{g}^{2}\right\} dX  \label{3.370}
\\
&\leq &C\left\Vert u\right\Vert _{L^{\infty }\left( D_{R_{1}}^{\left(
a\right) }\right) }^{2}+\int\nolimits_{D_{R_{1}}^{\left( a\right)
}}\left\vert \nabla _{g}u\right\vert _{g}^{2}dX\text{ ,}  \notag
\end{eqnarray}
for every $\alpha \geq 1$ , where $C$ is an absolute constant. Concerning
the second integral at the right-hand side of (\ref{3.331}) we have
\begin{equation}
\int\nolimits_{B_{1}}v^{2}\left( x,0\right) e^{-\frac{\left\vert
x\right\vert ^{2}}{4a}}dx\leq \int\nolimits_{B_{r\left( a\right)
}}u^{2}\left( x,0\right) dx\text{ ,}  \label{3.371}
\end{equation}
where
\begin{equation*}
r\left( a\right) =\left( 4\alpha a\log \frac{R_{1}^{2}}{4\alpha a}\right)
^{1/2}\text{ .}
\end{equation*}
By (\ref{3.360}), (\ref{3.370}) and (\ref{3.371}) we have
\begin{equation}
\alpha ^{2}\int\nolimits_{D_{R_{1}/2}^{\left( a\right) }}\sigma
_{a}^{-\alpha }u^{2}e^{-\frac{\left\vert x\right\vert ^{2}}{4\left(
s+a\right) }}dX  \label{3.380}
\end{equation}
\begin{equation*}
\leq \frac{C^{\alpha }\alpha ^{\alpha }}{R_{1}^{2\left( \alpha +3\right) }}%
\left( \left\Vert u\right\Vert _{L^{\infty }\left( D_{R_{1}}^{\left(
a\right) }\right) }^{2}+\int\nolimits_{D_{R_{1}}^{\left( a\right)
}}\left\vert \nabla _{g}u\right\vert _{g}^{2}dX\right) +\frac{C^{\alpha }}{%
a^{\alpha }}\int\nolimits_{B_{r\left( a\right) }}u^{2}\left( x,0\right) dx%
\text{ ,}
\end{equation*}
for every $\alpha \geq C$ and every $a\in \left( 0,\dfrac{R_{1}^{2}}{%
16\alpha }\right) $, where $C$, $C>1$, depends on $\lambda $ and $\Lambda $
only.

\noindent Let $r$ be a number in the interval $\left( 0,e^{-1/2}R_{1}\right)
$ and, for every $\alpha \geq C$, let $\overline{a}$ belong to $\left(
0,e^{-1}\dfrac{R_{1}^{2}}{4\alpha }\right) $ and such that
\begin{equation}
r\left( \overline{a}\right) =\left( 4\alpha \overline{a}\log \frac{R_{1}^{2}%
}{4\alpha \overline{a}}\right) ^{1/2}=r\text{ .}  \label{3.390}
\end{equation}%
By asymptotic estimates of $\overline{a}$ we have
\begin{equation}
\frac{1}{4e\alpha }r^{2}\left( \log \frac{R_{1}^{2}}{r^{2}}\right) ^{-1}\leq
\overline{a}\leq \frac{1}{4\alpha }r^{2}\left( \log \frac{R_{1}^{2}}{r^{2}}%
\right) ^{-1}\text{ .}  \label{3.391}
\end{equation}%
Furthermore, using standard regularity estimates \cite{Li} we estimate from
above the first and the second integral at the right-hand side of (\ref%
{3.380}) and we get
\begin{equation}
\alpha ^{2}\int\nolimits_{D_{R_{1}/2}^{\left( \overline{a}\right) }}\sigma
_{\overline{a}}^{-\alpha }u^{2}e^{-\frac{\left\vert x\right\vert ^{2}}{%
4\left( s+\overline{a}\right) }}dX  \label{3.400}
\end{equation}%
\begin{equation*}
\leq \frac{C^{\alpha }\alpha ^{\alpha }}{R_{1}^{2\left( \alpha +5\right) }}%
\int\nolimits_{B_{R_{1}}\times \left( 0,R_{1}^{2}\right) }u^{2}dX+\frac{%
C^{\alpha }}{\overline{a}^{\alpha }}\int\nolimits_{B_{r}}u^{2}\left(
x,0\right) dx\text{ ,}
\end{equation*}%
for every $\alpha \geq C$, where $C$, $C>1$, depends on $\lambda $ and $%
\Lambda $ only.

\noindent Now we estimate from below the left-hand side of (\ref{3.400}).
Let $\rho \in \left( 0,\dfrac{R_{1}}{2}\right) $, by (\ref{Rq2.145}) we have
\begin{equation}
\int\nolimits_{D_{R_{1}/2}^{\left( \overline{a}\right) }}\sigma _{\overline{%
a}}^{-\alpha }u^{2}e^{-\frac{\left\vert x\right\vert ^{2}}{4\left( s+%
\overline{a}\right) }}dX  \label{3.410}
\end{equation}
\begin{equation*}
\geq \int\nolimits_{D_{\rho }^{\left( \overline{a}\right) }}\sigma _{%
\overline{a}}^{-\alpha }u^{2}e^{-\frac{\left\vert x\right\vert ^{2}}{4\left(
s+\overline{a}\right) }}dX\geq \frac{1}{C^{\alpha }}\left( \frac{\rho ^{2}}{%
4\alpha }\right) ^{-\alpha }\int\nolimits_{D_{\rho }^{\left( \overline{a}%
\right) }}u^{2}dX\text{ ,}
\end{equation*}
where $C$, $C>1$, depends on $\lambda $ and $\Lambda $ only. By the
inequalities obtained in (\ref{3.391}), (\ref{3.400}) and (\ref{3.410}) we
have
\begin{equation}
\alpha ^{2}\int\nolimits_{D_{\rho }^{\left( \overline{a}\right)
}}u^{2}dX\leq \frac{\left( C\rho ^{2}\right) ^{\alpha }}{R_{1}^{2\left(
\alpha +5\right) }}\int\nolimits_{B_{R_{1}}\times \left( 0,R_{1}^{2}\right)
}u^{2}dX  \label{3.420}
\end{equation}
\begin{equation*}
+C^{\alpha }\rho ^{2\alpha }\left( r^{2}\left( \log \frac{R_{1}^{2}}{r^{2}}%
\right) ^{-1}\right) ^{-\alpha }\int\nolimits_{B_{r}}u^{2}\left( x,0\right)
dx\text{ ,}
\end{equation*}
for every $\alpha \geq C$, where $C$, $C>1$, depends on $\lambda $ and $%
\Lambda $ only.

\noindent Denote
\begin{equation*}
t_{1}=\frac{\rho ^{2}e^{-1}}{8\alpha }-\overline{a}\text{ , }t_{2}=\frac{%
\rho ^{2}e^{-1}}{8\alpha }\left( 1+e^{-1}\right) -\overline{a}\text{ ,}
\end{equation*}%
\begin{equation*}
Q=B_{\frac{\rho }{2e^{1/2}}}\times \left( t_{1},t_{2}\right) \text{ ,}
\end{equation*}%
\begin{equation*}
E=\left( \int\nolimits_{B_{R_{1}}\times \left( 0,R_{1}^{2}\right)
}u^{2}dX\right) ^{1/2}\text{ , }\varepsilon =\left(
\int\nolimits_{B_{r}}u^{2}\left( x,0\right) dx\right) ^{1/2}\text{ .}
\end{equation*}%
Observing that for every $\rho \in \left( \dfrac{r}{e^{1/2}},\dfrac{R_{1}}{2}%
\right) $ we have $t_{2}>t_{1}>0$ and $Q\subset D_{\rho }^{\left( \overline{a%
}\right) }$, inequality (\ref{3.420}) gives
\begin{equation}
\alpha ^{2}\int\nolimits_{Q}u^{2}dX\leq \frac{\left( C\rho \right)
^{2\alpha }}{R_{1}^{2\left( \alpha +5\right) }}E^{2}+\left( \frac{C\rho }{%
r\left( \log R_{1}^{2}/r^{2}\right) ^{-1/2}}\right) ^{2\alpha }\varepsilon
^{2}\text{ ,}  \label{3.430}
\end{equation}%
for every $\alpha \geq C$, where $C$, $C>1$, depends on $\lambda $ and $%
\Lambda $ only.

\noindent In order to estimate from below the left-hand side of inequality (%
\ref{3.430}) we apply Lemma \ref{Le3.2}. In doing so we choose $\rho _{0}=%
\dfrac{e^{-1/2}\rho }{4}$, $\rho _{1}=\dfrac{3e^{-1}\rho }{8}$, $\rho _{2}=%
\dfrac{e^{-1/2}\rho }{2}$ and $T=\rho _{2}^{2}$. Let us denote
\begin{equation*}
K=\left( \int\nolimits_{B_{\rho _{0}}}u^{2}\left( x,0\right) dx\right)
^{1/2}\text{.}
\end{equation*}%
We have that there exists a constant $C_{0}$, $C_{0}>1$, depending on $%
\lambda $ and $\Lambda $ only such that if
\begin{equation}
\alpha \geq \alpha _{0}:=C_{0}\max \left\{ 1,\log \left( \frac{\rho
^{n-4}E^{2}}{K^{2}}\right) \right\}  \label{3.440}
\end{equation}%
then
\begin{equation}
\int\nolimits_{Q}u^{2}dX=\int\nolimits_{t_{1}}^{t_{2}}dt\int%
\nolimits_{B_{\rho _{2}}}u^{2}\left( x,t\right) dx\geq \frac{t_{2}-t_{1}}{C}%
\int\nolimits_{B_{\rho _{0}}}u^{2}\left( x,0\right) dx=\frac{\rho ^{2}}{%
8\alpha C}K^{2}\text{ .}  \label{3.450}
\end{equation}%
By inequality (\ref{3.430}) and (\ref{3.450}) we get
\begin{equation}
\rho ^{2}K^{2}\leq \frac{\left( C\rho \right) ^{2\alpha }}{R_{1}^{2\left(
\alpha +5\right) }}E^{2}+\left( \frac{C\rho }{r\left( \log
R_{1}^{2}/r^{2}\right) ^{-1/2}}\right) ^{2\alpha }\varepsilon ^{2}\text{ ,}
\label{3.460}
\end{equation}%
for every $\alpha \geq \alpha _{0}$ and every $r,R_{1},\rho $ satisfying the
relations $R_{1}\in \left( 0,\overline{R}_{0}\right] $, $r\in \left(
0,e^{-1/2}R_{1}\right] $, $\rho \in \left( e^{1/2}r,\dfrac{R_{1}}{2}\right] $%
, where $C>1$ depends on $\lambda $ and $\Lambda $ only.

\noindent Denote by
\begin{equation*}
\alpha _{1}=\frac{\log \left( E^{2}\varepsilon ^{-2}\right) }{\log \left(
R_{1}^{2}r^{-2}\log \left( R_{1}^{2}r^{-2}\right) \right) }\text{ .}
\end{equation*}
the following cases occur: i) $\alpha _{1}\geq \alpha _{0}$, ii) $\alpha
_{1}<\alpha _{0}$. If case i) occurs then we choose in (\ref{3.460}) $\alpha
=\alpha _{1}$ and we have
\begin{equation}
\rho ^{2}K^{2}\leq \frac{2}{R_{1}^{10}}E^{2\left( 1-\theta _{0}\right)
}\varepsilon ^{2\theta _{0}}\text{ ,}  \label{3.470}
\end{equation}
where
\begin{equation}
\theta _{0}=\frac{\log \left( E^{2}\rho ^{-2}/C\right) }{\log \left(
R_{1}^{2}r^{-2}\log \left( R^{2}r^{-2}\right) \right) }\text{ ,}
\label{3.480}
\end{equation}
where $C$, $C>1$, depends on $\lambda $ and $\Lambda $ only.

\noindent Now consider case ii). We have either
\begin{equation}
\alpha _{0}=C_{0}\log \left( \frac{\rho ^{n-4}E^{2}}{K^{2}}\right) \geq C_{0}%
\text{ ,}  \label{3.490}
\end{equation}
or
\begin{equation}
\alpha _{0}=C_{0}\leq C_{0}\log \left( \frac{\rho ^{n-4}E^{2}}{K^{2}}\right)
\text{ , }  \label{3.495}
\end{equation}
where $C_{0}$ is the same constant which appears in (\ref{3.440}). Let us
introduce the notation
\begin{equation*}
\widetilde{\theta }_{0}=\frac{1}{C_{0}\log \left( R_{1}^{2}r^{-2}\log \left(
R_{1}^{2}r^{-2}\right) \right) }\text{ , }\theta _{1}=\frac{1}{C_{0}\log
\frac{1}{r}}
\end{equation*}
If (\ref{3.490}) occurs then, recalling that $\alpha _{1}<\alpha $, we have
trivially
\begin{equation*}
\theta _{1}\log \left( \frac{E^{2}}{\varepsilon ^{2}}\right) \leq \log
\left( \frac{\rho ^{n-4}E^{2}}{K^{2}}\right) \text{ ,}
\end{equation*}
hence
\begin{equation}
K^{2}\leq E^{2\left( 1-\widetilde{\theta }_{0}\right) }\varepsilon ^{2%
\widetilde{\theta }_{0}}\text{ .}  \label{3.500}
\end{equation}
If (\ref{3.495}) occurs then we have trivially
\begin{equation*}
\widetilde{\theta }_{0}\log \left( \frac{E^{2}}{\varepsilon ^{2}}\right)
\leq 1\text{ ,}
\end{equation*}
hence
\begin{equation}
E^{2\theta _{1}}\leq e\varepsilon ^{2\theta _{1}}\text{ .}  \label{3.510}
\end{equation}
On the other side, by standard regularity estimates for parabolic equations
we have
\begin{equation}
K^{2}\leq \frac{C_{1}}{R_{1}^{2}}E^{2}\text{ ,}  \label{3.520}
\end{equation}
where $C_{1}$ depends on $\lambda $ only. Therefore (\ref{3.510}) and (\ref%
{3.520}) yield
\begin{equation}
K^{2}\leq \frac{C_{1}}{R_{1}^{2}}E^{2}=\frac{C_{1}}{R_{1}^{2}}E^{2\left( 1-%
\widetilde{\theta }_{0}\right) }E^{2\widetilde{\theta }_{0}}\leq \frac{C_{1}%
}{R_{1}^{2}}E^{2\left( 1-\widetilde{\theta }_{0}\right) }\varepsilon ^{2%
\widetilde{\theta }_{0}}\text{ .}  \label{3.530}
\end{equation}
Now there exists a constant $C$, $C>1$, depending on $\lambda $ and $\Lambda
$ only such that if $\rho ^{2}\leq \frac{R_{1}^{2}}{e^{1/C}C}$ then we have $%
\widetilde{\theta }_{0}\leq \theta _{0}$. Therefore by (\ref{3.470}), (\ref%
{3.500}), (\ref{3.520}) and (\ref{3.530}) we get
\begin{equation}
K^{2}\leq \frac{C_{2}}{R_{1}^{5}\rho }E^{2\left( 1-\widetilde{\theta }%
_{0}\right) }\varepsilon ^{2\widetilde{\theta }_{0}}\text{ ,}  \label{3.540}
\end{equation}
for every $R_{1}\in \left( 0,\overline{R}_{0}\right] $, $r\in \left(
0,e^{-1/2}R_{1}\right] $, $\rho \in \left( e^{1/2}r,\dfrac{R_{1}}{%
Ce^{1/\left( 2C\right) }}\right] $, where $C_{2}$, $C_{2}>1$, depends on $%
\lambda $ only.

In the case $g^{ij}\left( 0,0\right) \neq \delta ^{ij}$ we can consider a
linear transformation $S:\mathbb{R}^{n}\rightarrow \mathbb{R}^{n}$, $%
Sx=\left\{ S_{j}^{i}x_{j}\right\} _{i=1}^{n}$, such that, denoting $%
\widetilde{g}^{ij}\left( y,t\right) =\frac{1}{\left\vert \det S\right\vert }%
S_{k}^{i}g^{kl}\left( S^{-1}y,t\right) S_{l}^{j}$ we have $\widetilde{g}%
^{ij}\left( 0,0\right) =\delta ^{ij}$ and
\begin{equation*}
B_{\varrho /\sqrt{\lambda }}\subset S\left( B_{\varrho }\right) \subset B_{%
\sqrt{\lambda }\varrho }\text{ , for every }\varrho >0\text{.}
\end{equation*}
Therefore, by a change of variables, by using the inequality (\ref{3.540})
and the trivial inequality $r^{2}\left( \log \frac{1}{r^{2}}\right)
^{-1}>r^{3}$, for every $r\in \left( 0,1\right) $ the thesis follows.$%
\blacksquare $

\begin{corollary}
\label{Cor3.5}\textbf{(Spacelike strong unique continuation).} Let $u\in
H^{2,1}\left( B_{R}\times \left( t_{0}-R^{2},t_{0}\right) \right) $ satisfy
inequality (\ref{3.280}).

\noindent If for every $k\in \mathbb{N}$ we have
\begin{equation}
\int\nolimits_{B_{r}}u^{2}\left( x,t_{0}\right) dx=O\left( r^{2k}\right)
\text{ , as }r\rightarrow 0\text{,}  \label{3.550}
\end{equation}
then
\begin{equation*}
u\left( .,t_{0}\right) =0\text{ , in }B_{R}\text{ .}
\end{equation*}
\end{corollary}

\textbf{Proof.} Is the same of the proof of Corollary \ref{Cor}.$%
\blacksquare $

\bigskip

In order to state the theorem below we need to introduce some new notations.
Let $E$ and $R$ be positive numbers, $t_{0}\in \mathbb{R}$ and let $\varphi :%
\mathbb{R}^{n-1}\times \mathbb{R\rightarrow R}$ be a function satisfying the
conditions
\begin{equation}
\varphi \left( 0,t_{0}\right) =\left\vert \nabla _{x^{\prime }}\varphi
\left( 0,t_{0}\right) \right\vert =0\text{ ,}  \label{3.560}
\end{equation}
\begin{equation}
\left\Vert \varphi \right\Vert _{L^{\infty }\left( Q_{R}^{t_{0}\prime
}\right) }+R\left\Vert \nabla _{x^{\prime }}\varphi \right\Vert _{L^{\infty
}\left( Q_{R}^{t_{0}\prime }\right) }+R^{2}\left\Vert D_{x^{\prime
}}^{2}\varphi \right\Vert _{L^{\infty }\left( Q_{R}^{t_{0}\prime }\right)
}+R^{2}\left\Vert \partial _{s}\varphi \right\Vert _{L^{\infty }\left(
Q_{R}^{t_{0}\prime }\right) }\leq ER\text{ ,}  \label{3.570}
\end{equation}
where $Q_{R}^{t_{0}\prime }=B_{R}^{\prime }\times \left( t_{0}-R^{2},t_{0}%
\right] $.

\noindent For any numbers $\rho \in \left( 0,R\right] $, $\tau \in \left(
t_{0}-R^{2},t_{0}\right] $, denote by
\begin{equation*}
Q_{\rho ,\varphi }^{t_{0}}=\left\{ \left( x,s\right) \in B_{\rho }\times
\left( t_{0}-\rho ^{2},t_{0}\right) :x_{n}>\varphi \left( x^{\prime
},s\right) \right\} \text{,}
\end{equation*}
\begin{equation*}
Q_{\rho ,\varphi }\left( \tau \right) =\left\{ x\in B_{\rho }:x_{n}>\varphi
\left( x^{\prime },\tau \right) \right\} \text{,}
\end{equation*}
\begin{equation*}
\Gamma _{\rho ,\varphi }^{t_{0}}=\left\{ \left( x,s\right) \in B_{\rho
}\times \left( t_{0}-\rho ^{2},t_{0}\right) :x_{n}=\varphi \left( x^{\prime
},s\right) \right\} \text{.}
\end{equation*}

\begin{theorem}
\label{Th3.6}\textbf{(two-sphere one-cylinder inequality at the boundary). }%
Let $L$ be a second order parabolic operator as in Theorem \ref{Th3.4}. Let $%
\varphi :\mathbb{R}^{n-1}\times \mathbb{R\rightarrow R}$ be a function
satisfying (\ref{3.560}) and (\ref{3.570}). Assume that $u\in H^{2,1}\left(
Q_{R,\varphi }\right) $ satisfies the inequality
\begin{equation}
\left\vert Lu\right\vert \leq \Lambda \left( \frac{\left\vert \nabla
u\right\vert }{R}+\frac{\left\vert u\right\vert }{R^{2}}\right) \text{ , in }%
Q_{R,\varphi }  \label{3.580}
\end{equation}
and
\begin{equation}
u\left( x,t\right) =0\text{, for every }\left( x,t\right) \in \Gamma
_{R,\varphi }^{t_{0}}\text{ .}  \label{3.590}
\end{equation}
Then there exist constants $\eta _{2}\in \left( 0,1\right) $ and $C$, $C>1$,
depending on $\lambda $, $\Lambda $ and $E$ only such that for every $r$, $%
\rho $, $0<r\leq \rho \leq \eta _{2}R$ we have
\begin{equation}
\int\nolimits_{Q_{\rho ,\varphi }\left( t_{0}\right) }u^{2}\left(
x,t_{0}\right) dx\leq \frac{CR^{2}}{\rho ^{2}}\left(
R^{-2}\int\nolimits_{Q_{R,\varphi }^{t_{0}}}u^{2}dX\right) ^{1-\theta
_{2}}\left( \int\nolimits_{Q_{r,\varphi }\left( t_{0}\right) }u^{2}\left(
x,t_{0}\right) dx\right) ^{\theta _{2}}\text{ ,}  \label{3.600}
\end{equation}
where
\begin{equation}
\theta _{2}=\frac{1}{C\log \frac{R}{r}}\text{ .}  \label{3.605}
\end{equation}
\end{theorem}

For the sake of simplicity, here we give a proof of Theorem \ref{Th3.6} with
the following additional condition on the matrix $\left\{ g^{ij}\left(
x,t\right) \right\} _{i,j=1}^{n}$:
\begin{equation}
\left( \sum\limits_{i,j=1}^{n}\left( g^{ij}\left( x,t\right) -g^{ij}\left(
x,\tau \right) \right) ^{2}\right) ^{1/2}\leq \frac{\Lambda }{R^{2}}%
\left\vert t-\tau \right\vert \text{ ,}  \label{3.610}
\end{equation}%
for every $\left( x,t\right) ,\left( x,\tau \right) \in \mathbb{R}^{n+1}$.

\noindent Actually the additional condition (\ref{3.610}) is really used to
prove Theorem \ref{Th4.1}\ (the stability result for the inverse problem).
An outline of the proof of Theorem \ref{Th3.6} (without condition (\ref%
{3.610})) is given in \cite{EsFeVe}.

We need the following proposition proved \cite{Ve3}.

\begin{proposition}
\label{Pr3.7}Let $\lambda _{0}$ and $\Lambda _{0}$ be positive numbers, with
$\lambda _{0}\in \left( 0,1\right] $. Let $A\left( t\right) $ be a $n\times
n $ symmetric real matrix whose entries are Lipschitz continuous in $\mathbb{%
R} $. Assume that
\begin{equation}
\lambda _{0}\left\vert \xi \right\vert ^{2}\leq A\left( t\right) \xi \cdot
\xi \leq \lambda _{0}^{-1}\left\vert \xi \right\vert ^{2}\text{ , for every }%
\xi \in \mathbb{R}^{n}\text{ and }s\in \mathbb{R}\text{ ,}  \label{3.620}
\end{equation}
\begin{equation}
\left\vert \frac{d}{ds}A\left( t\right) \right\vert \leq \Lambda _{0}\text{,
a.e. in }\mathbb{R}\text{.}  \label{3.630}
\end{equation}
Denote by $\sqrt{A\left( t\right) }$ the positive square root of $A\left(
t\right) $. We have
\begin{equation}
\left\vert \frac{d}{ds}\sqrt{A\left( t\right) }\right\vert \leq
2^{3/2}\lambda _{0}^{-7/2}\Lambda _{0}\text{, a.e. in }\mathbb{R}\text{.}
\label{3.640}
\end{equation}
\end{proposition}

\bigskip

\textbf{Proof of Theorem \ref{Th3.6} with condition }(\ref{3.610})\textbf{.}
Without loss of generality we assume that $t_{0}=0$. Let us introduce the
following notations. Let $\psi _{1}:\mathbb{R}^{n}\mathbb{\rightarrow R}^{n}$
and $\Phi _{1}:\mathbb{R}^{n+1}\mathbb{\rightarrow R}^{n+1}$ be respectively
the maps $\psi _{1}\left( y,\tau \right) =\left( y^{\prime },y_{n}+\varphi
\left( y^{\prime },\tau \right) \right) $ and $\Phi _{1}\left( y,\tau
\right) =\left( \psi _{1}\left( y,\tau \right) ,\tau \right) $. We have $%
\Phi _{1}\left( y^{\prime },0,\tau \right) =\left( y^{\prime },\varphi
\left( y^{\prime },\tau \right) ,\tau \right) $, for every $y^{\prime }\in
\mathbb{R}^{n-1},\tau \in \mathbb{R}$. Furthermore, denoting $R_{1}=\dfrac{R%
}{1+2E}$, by (\ref{3.570}) we have
\begin{equation}
\Phi _{1}\left( B_{1}^{+}\times \left( -R_{1}^{2},0\right] \right) \subset
Q_{R,\varphi }^{0}\text{ .}  \label{3.650}
\end{equation}%
Denote
\begin{equation*}
g_{1}^{-1}\left( y,\tau \right) =\left( \frac{\partial \psi _{1}\left(
y,\tau \right) }{\partial y}\right) ^{-1}g^{-1}\left( \Phi _{1}\left( y,\tau
\right) \right) \left( \left( \frac{\partial \psi _{1}\left( y,\tau \right)
}{\partial y}\right) ^{-1}\right) ^{\ast }\text{ ,}
\end{equation*}%
(here $\frac{\partial \psi _{1}\left( y,\tau \right) }{\partial y}$ is the
Jacobian matrix of $\psi _{1}$),
\begin{equation*}
u_{1}\left( y,\tau \right) =u\left( \Phi _{1}\left( y,\tau \right) \right)
\end{equation*}%
and
\begin{equation*}
\left( L_{1}u_{1}\right) \left( y,\tau \right) =div\left( g_{1}^{-1}\left(
y,\tau \right) \nabla _{y}u_{1}\right) -\partial _{\tau }u_{1}\text{.}
\end{equation*}%
We have
\begin{equation*}
\left( Lu\right) \left( \Phi _{1}\left( y,\tau \right) \right) =\left(
L_{1}u_{1}\right) \left( y,\tau \right) +\left( \frac{\partial \psi
_{1}\left( y,\tau \right) }{\partial \tau }\right) ^{\ast }\left( \left(
\frac{\partial \psi _{1}\left( y,\tau \right) }{\partial y}\right)
^{-1}\right) ^{\ast }\nabla _{y}u_{1}\text{.}
\end{equation*}%
By (\ref{Rq2.105b}), (\ref{Rq2.110b}), (\ref{3.570}), (\ref{3.580}), (\ref%
{3.610}) and (\ref{3.650}) we have that there exist $\lambda _{1}\in \left(
0,1\right] $, $\Lambda _{1}>0$, $\lambda _{1}$ depending on $\lambda $ and $%
E $ only and $\Lambda _{1}$ depending on $\lambda $, $\Lambda $ and $E$
only, such that, when $\left( y,\tau \right) ,\left( z,s\right) \in \mathbb{R%
}^{n+1}$, $\xi \in \mathbb{R}^{n}$, we have
\begin{equation}
\lambda _{1}\left\vert \xi \right\vert ^{2}\leq
\sum\limits_{i,j=1}^{n}g_{1}^{ij}\left( y,\tau \right) \xi _{i}\xi _{j}\leq
\lambda _{1}^{-1}\left\vert \xi \right\vert ^{2}\text{ ,}  \label{3.660}
\end{equation}%
\begin{equation}
\left( \sum\limits_{i,j=1}^{n}\left( g_{1}^{ij}\left( y,\tau \right)
-g_{1}^{ij}\left( z,s\right) \right) ^{2}\right) ^{1/2}\leq \frac{\Lambda
_{1}}{R_{1}}\left( \left\vert y-z\right\vert ^{2}+\left\vert \tau
-s\right\vert \right) ^{1/2}\text{ ,}  \label{3.670}
\end{equation}%
\begin{equation}
\left( \sum\limits_{i,j=1}^{n}\left( g_{1}^{ij}\left( y,\tau \right)
-g_{1}^{ij}\left( y,s\right) \right) ^{2}\right) ^{1/2}\leq \frac{\Lambda
_{1}}{R_{1}^{2}}\left\vert \tau -s\right\vert  \label{3.680}
\end{equation}%
and
\begin{equation}
\left\vert L_{1}u_{1}\right\vert \leq \Lambda _{1}\left( \frac{\left\vert
\nabla u_{1}\right\vert }{R_{1}}+\frac{\left\vert u_{1}\right\vert }{%
R_{1}^{2}}\right) \text{ , in }B_{1}^{+}\times \left( -R_{1}^{2},0\right]
\text{ .}  \label{3.690}
\end{equation}%
Furthermore
\begin{equation}
u_{1}\left( y^{\prime },0,\tau \right) =0\text{, for every }y^{\prime }\in
B_{R_{1}}^{\prime }\text{, }\tau \in \left( -R_{1}^{2},0\right] \text{ .}
\label{3.700}
\end{equation}%
Now let $H\left( \tau \right) $ and $S\left( \tau \right) $ be respectively
the positive square root of $g_{1}\left( 0,\tau \right) $ and the Lipschitz
continuous rotation of $\mathbb{R}^{n}$ such that $S\left( \tau \right)
e_{n}=vers\left( H\left( \tau \right) e_{n}\right) $, for every $\tau \in
\mathbb{R}$. Denote
\begin{equation*}
K\left( \tau \right) =\left( H\left( \tau \right) \right) ^{-1}S\left( \tau
\right) \text{ , }g_{2}^{-1}\left( \widetilde{y},\tau \right) =\left(
K\left( \tau \right) \right) ^{-1}g_{1}^{-1}\left( K\left( \tau \right)
\widetilde{y},\tau \right) \left( \left( K\left( \tau \right) \right)
^{-1}\right) ^{\ast }\text{ ,}
\end{equation*}%
\begin{equation*}
u_{2}\left( \widetilde{y},\tau \right) =u_{1}\left( K\left( \tau \right)
\widetilde{y},\tau \right) \text{ , }\left( L_{2}u_{2}\right) \left(
\widetilde{y},\tau \right) =div\left( g_{2}^{-1}\left( \widetilde{y},\tau
\right) \nabla _{\widetilde{y}}u_{2}\right) -\partial _{\tau }u_{2}\text{ .}
\end{equation*}%
We have
\begin{equation*}
g_{2}^{-1}\left( 0,\tau \right) =I_{n}\text{ , for every }\tau \in \mathbb{R}%
\text{ ,}
\end{equation*}%
\begin{equation*}
\left( L_{1}u_{1}\right) \left( K\left( \tau \right) \widetilde{y},\tau
\right) =\left( L_{2}u_{2}\right) \left( \widetilde{y},\tau \right) +\left(
K^{\prime }\left( \tau \right) \widetilde{y}\right) K^{-1}\left( \tau
\right) \nabla _{\widetilde{y}}u_{2}\text{ .}
\end{equation*}%
By (\ref{3.660})-(\ref{3.700}) and by Proposition \ref{Pr3.7} we have that
there exist $\lambda _{2}\in \left( 0,1\right] $, $\Lambda _{2}>0$ and $C>1$%
, $\lambda _{2}$ depending on $\lambda $ and $E$ only and $\Lambda _{2}$, $C$
depending on $\lambda $, $\Lambda $ and $E$ only such that we have, setting $%
R_{2}=\frac{R_{1}}{C}$, for every $\left( \widetilde{y},\tau \right) ,\left(
\widetilde{y}_{\ast },s\right) \in \mathbb{R}^{n+1}$, $\xi \in \mathbb{R}%
^{n} $,
\begin{equation}
\lambda _{2}\left\vert \xi \right\vert ^{2}\leq
\sum\limits_{i,j=1}^{n}g_{2}^{ij}\left( \widetilde{y},\tau \right) \xi
_{i}\xi _{j}\leq \lambda _{2}^{-1}\left\vert \xi \right\vert ^{2}\text{ ,}
\label{3.710}
\end{equation}%
\begin{equation}
\left( \sum\limits_{i,j=1}^{n}\left( g_{2}^{ij}\left( \widetilde{y},\tau
\right) -g_{2}^{ij}\left( \widetilde{y}_{\ast },s\right) \right) ^{2}\right)
^{1/2}\leq \frac{\Lambda _{2}}{R_{2}}\left( \left\vert \widetilde{y}-%
\widetilde{y}_{\ast }\right\vert ^{2}+\left\vert \tau -s\right\vert \right)
^{1/2}\text{ ,}  \label{3.720}
\end{equation}%
\begin{equation}
\left( \sum\limits_{i,j=1}^{n}\left( g_{2}^{ij}\left( \widetilde{y},\tau
\right) -g_{2}^{ij}\left( \widetilde{y}_{\ast },s\right) \right) ^{2}\right)
^{1/2}\leq \frac{\Lambda _{2}}{R_{2}^{2}}\left\vert \tau -s\right\vert \text{
,}  \label{3.730}
\end{equation}%
\begin{equation}
\left\vert L_{2}u_{2}\right\vert \leq \Lambda _{2}\left( \frac{\left\vert
\nabla u_{2}\right\vert }{R_{2}}+\frac{\left\vert u_{2}\right\vert }{%
R_{2}^{2}}\right) \text{ , in }B_{2}^{+}\times \left( -R_{2}^{2},0\right]
\text{ }  \label{3.740}
\end{equation}%
and
\begin{equation}
u_{2}\left( \widetilde{y}^{\prime },0,\tau \right) =0\text{, for every }%
\widetilde{y}^{\prime }\in B_{R_{2}}^{\prime }\text{, }\tau \in \left(
-R_{2}^{2},0\right] \text{ .}  \label{3.750}
\end{equation}%
For the sake of brevity we shall omit the sign " $\widetilde{}$ " over $y$.

Now we adapt ideas in \cite{AE} to the case of time varying coefficients.
Let $\zeta $ be a function belonging to $C_{0}^{\infty }\left( \mathbb{R}%
^{n-1}\right) $ and such that \textit{supp\thinspace }$\zeta \subset
B_{1}^{\prime }$, $\zeta \geq 0$ and $\int\limits_{\mathbb{R}^{n-1}}\zeta
\left( y^{\prime }\right) dy^{\prime }=1$. Denote by $\zeta _{\left(
n\right) }$ the function $y_{n}^{1-n}\zeta \left( \frac{y^{\prime }}{y_{n}}%
\right) $. Denote
\begin{equation*}
\beta \left( z^{\prime },\tau \right) =-g_{2}^{-1}\left( z^{\prime },0,\tau
\right) e_{n}\text{ ,}
\end{equation*}%
\begin{equation*}
w_{j}\left( z^{\prime },z_{n},\tau \right) =z_{j}-z_{n}\left( \zeta _{\left(
n\right) }\ast \beta _{j}\left( .,\tau \right) \right) \left( z^{\prime
}\right) \text{ , }j=1,...,n-1\text{ ,}
\end{equation*}%
\begin{equation*}
w_{n}\left( z^{\prime },z_{n},\tau \right) =-z_{n}\left( \zeta _{\left(
n\right) }\ast \beta _{n}\left( .,\tau \right) \right) \left( z^{\prime
}\right) \text{ .}
\end{equation*}%
With $\psi _{2}\left( .,\tau \right) :\mathbb{R}^{n}\mathbb{\rightarrow R}%
^{n}$ we denote the map whose components are defined as follows, for each $%
j=1,...,n$%
\begin{equation*}
\left( \psi _{2}\right) _{j}\left( z^{\prime },z_{n},\tau \right) =\left\{
\begin{array}{c}
w_{j}\left( z^{\prime },z_{n},\tau \right) \text{ , for }z_{n}\geq 0\text{,}
\\
4w_{j}\left( z^{\prime },z_{n},\tau \right) -3w_{j}\left( z^{\prime
},z_{n},\tau \right) \text{, for }z_{n}<0\text{.}%
\end{array}%
\right.
\end{equation*}%
It turns out that $\psi _{2}$ is a $C^{1,1}$ function with respect to $z$
and it is a Lipschitz continuous function with respect to $\tau \in \left[
0,R_{2}^{2}\right) $. Moreover there exist constants $C_{1},C_{2},C_{3}\in %
\left[ 1,+\infty \right) $ depending on $\lambda $, $\Lambda $ and $E$ only
such that, setting $\rho _{1}=\frac{R_{2}}{C_{1}}$, $\rho _{2}=\frac{R_{2}}{%
C_{2}}$, we have, for every $\tau \in \left( -R_{2}^{2},0\right] $,

i) $\psi _{2}\left( z,\tau \right) \in B_{\rho _{1}}$, for every $z\in
B_{\rho _{2}}$,

ii) $\psi _{2}\left( z^{\prime },0,\tau \right) =\left( z^{\prime },0\right)
$, for every $z^{\prime }\in B_{\rho _{2}}^{\prime }$,

iii) $\psi _{2}\left( z,\tau \right) \in B_{\rho _{1}}^{+}$, for every $z\in
B_{\rho _{2}}^{+}$,

iv) $C_{3}^{-1}\left\vert z-z_{\ast }\right\vert \leq \left\vert \psi
_{2}\left( z,\tau \right) -\psi _{2}\left( z_{\ast },\tau \right)
\right\vert \leq C_{3}\left\vert z-z_{\ast }\right\vert $, for every $%
z,z_{\ast }\in B_{\rho _{2}}$,

v) $\left\vert \frac{\partial ^{2}\psi _{2}\left( z,\tau \right) }{\partial
z_{i}\partial z_{j}}\right\vert \leq C_{3}$, for every $z\in B_{\rho _{2}}$,

vi) $C_{3}^{-1}\leq \left\vert \det \frac{\partial \psi _{2}\left( z,\tau
\right) }{\partial z}\right\vert \leq C_{3}$, for every $z\in B_{\rho _{2}}$.

\noindent Let us denote
\begin{equation*}
u_{3}\left( z,\tau \right) =u_{2}\left( \psi _{2}\left( z,\tau \right) ,\tau
\right) \text{ , }B\left( u_{3}\right) =\left( \frac{\partial \psi
_{2}\left( z,\tau \right) }{\partial \tau }\right) ^{\ast }\left( \left(
\frac{\partial \psi _{2}\left( z,\tau \right) }{\partial z}\right)
^{-1}\right) ^{\ast }\nabla _{z}u_{3}\text{ ,}
\end{equation*}
\begin{equation*}
g_{3}^{-1}\left( z,\tau \right) =\left( \frac{\partial \psi _{2}\left(
z,\tau \right) }{\partial y}\right) ^{-1}g_{2}^{-1}\left( \psi _{2}\left(
z,\tau \right) ,\tau \right) \left( \left( \frac{\partial \psi _{2}\left(
z,\tau \right) }{\partial z}\right) ^{-1}\right) ^{\ast }\text{ , }J\left(
z,\tau \right) =\left\vert \det \frac{\partial \psi _{2}\left( z,\tau
\right) }{\partial z}\right\vert \text{.}
\end{equation*}
We have, for every $\tau \in \left[ 0,R_{2}^{2}\right) $ and $z^{\prime }\in
B_{\rho _{2}}^{\prime }$,
\begin{equation*}
J\left( 0,\tau \right) g_{3}^{-1}\left( 0,\tau \right) =I_{n}\text{ ,}
\end{equation*}
\begin{equation*}
g_{3}^{nj}\left( z^{\prime },0,\tau \right) =g_{3}^{jn}\left( z^{\prime
},0,\tau \right) =0\text{ , for }j=1,...,n-1\text{.}
\end{equation*}
Moreover by (\ref{3.710})-(\ref{3.730}) we have that there exist $\lambda
_{3}\in \left( 0,1\right] $, $\Lambda _{3}>0$ depending on $\lambda $, $%
\Lambda $ and $E$ only such that for every $\left( z,\tau \right) ,\left(
z_{\ast },s\right) \in B_{\rho _{2}}^{+}\times \left[ 0,\rho _{2}^{2}\right)
$ and every $\xi \in \mathbb{R}^{n}$ we have
\begin{equation*}
\lambda _{3}\left\vert \xi \right\vert ^{2}\leq
\sum\limits_{i,j=1}^{n}g_{3}^{ij}\left( z,\tau \right) \xi _{i}\xi _{j}\leq
\lambda _{3}^{-1}\left\vert \xi \right\vert ^{2}\text{ ,}
\end{equation*}
\begin{equation*}
\left( \sum\limits_{i,j=1}^{n}\left( g_{3}^{ij}\left( z,\tau \right)
-g_{3}^{ij}\left( z_{\ast },s\right) \right) ^{2}\right) ^{1/2}\leq \frac{%
\Lambda _{3}}{R_{3}}\left( \left\vert z-z_{\ast }\right\vert ^{2}+\left\vert
\tau -s\right\vert \right) ^{1/2}\text{ ,}
\end{equation*}
\begin{equation}
\left( L_{2}u_{2}\right) \left( \psi _{2}\left( z,\tau \right) ,\tau \right)
=\left( L_{3}u_{3}\right) \left( z,\tau \right) +B\left( u_{3}\right) \text{
,}  \label{3.741}
\end{equation}
where
\begin{equation*}
\left( L_{3}u_{3}\right) \left( z,\tau \right) =\frac{1}{J\left( z,\tau
\right) }div\left( J\left( z,\tau \right) g_{3}^{-1}\left( z,\tau \right)
\nabla _{z}u_{3}\right) -\partial _{\tau }u_{3}\text{ .}
\end{equation*}
By (\ref{3.740}) and (\ref{3.741}) we have
\begin{equation*}
\left\vert div\left( g_{3}^{-1}\nabla _{z}u_{3}\right) -\partial _{\tau
}u_{3}\right\vert \leq C\left( \frac{\left\vert \nabla u_{3}\right\vert }{%
\rho _{2}}+\frac{\left\vert u_{3}\right\vert }{\rho _{2}^{2}}\right) \text{
, in }B_{\rho _{2}}^{+}\times \left[ 0,\rho _{2}^{2}\right) \text{ ,}
\end{equation*}
where $C$ depends on $\lambda $, $\Lambda $ and $E$ only.

\noindent For every $\left( z,\tau \right) \in B_{\rho _{2}}\times \left(
-\rho _{2}^{2},0\right] $ , let us denote by $\overline{g}_{3}^{-1}\left(
z,\tau \right) =\left\{ \overline{g}_{3}^{ij}\left( z,\tau \right) \right\}
_{i,j=1}^{n}$ the symmetric matrix whose entries are given by
\begin{equation*}
\overline{g}_{3}^{ij}\left( z^{\prime },z_{n},\tau \right) =g_{3}^{ij}\left(
z^{\prime },\left\vert z_{n}\right\vert ,\tau \right) \text{ , if either }%
1\leq i,j\leq n-1\text{ or }i=j=n\text{,}
\end{equation*}
\begin{equation*}
\overline{g}_{3}^{nj}\left( z^{\prime },z_{n},\tau \right) =\overline{g}%
_{3}^{jn}\left( z^{\prime },z_{n},\tau \right) =sgn\left( z_{n}\right)
g_{3}^{jn}\left( z^{\prime },\left\vert z_{n}\right\vert ,\tau \right) \text{
, if }1\leq j\leq n-1\text{.}
\end{equation*}
It turns out that $\overline{g}_{3}^{-1}$ satisfies the same ellipticity
condition and Lipschitz condition as $g_{3}^{-1}$.

\noindent Denoting
\begin{equation*}
v\left( z,\tau \right) =sgn\left( z_{n}\right) u_{3}\left( z^{\prime
},\left\vert z_{n}\right\vert ,\tau \right) \text{ , for every }\left(
z,\tau \right) \in B_{\rho _{2}}\times \left[ 0,\rho _{2}^{2}\right) \text{,}
\end{equation*}
we have that $v$ belongs to $H^{2,1}\left( B_{\rho _{2}}\times \left( -\rho
_{2}^{2},0\right) \right) $ and satisfies the inequality
\begin{equation*}
\left\vert div\left( \overline{g}_{3}^{-1}\nabla _{z}u_{3}\right) -\partial
_{\tau }u_{3}\right\vert \leq C\left( \frac{\left\vert \nabla
u_{3}\right\vert }{\rho _{2}}+\frac{\left\vert u_{3}\right\vert }{\rho
_{2}^{2}}\right) \text{ , in }B_{\rho _{2}}\times \left( -\rho _{2}^{2},0%
\right] \text{ ,}
\end{equation*}
where $C$ depends on $\lambda $, $\Lambda $ and $E$ only.

\noindent By Theorem \ref{Th3.4} we have there exist constants $\eta _{2}\in
\left( 0,1\right) $ and $C>1$ depending on $\lambda $, $\Lambda $ and $E$
only such that for every $0<r\leq \rho \leq \eta _{2}\rho _{2}$ we have
\begin{equation}
\int\nolimits_{B_{\rho }}v^{2}\left( x,0\right) dx\leq \frac{C\rho _{2}^{2}%
}{\rho ^{2}}\left( \rho _{2}^{-2}\int\nolimits_{B_{\rho _{2}}\times \left(
-\rho _{2}^{2},0\right) }v^{2}dX\right) ^{1-\overline{\theta }_{1}}\left(
\int\nolimits_{B_{r}}v^{2}\left( x,0\right) dx\right) ^{\overline{\theta }%
_{1}}\text{ ,}  \label{3.800}
\end{equation}
where
\begin{equation*}
\overline{\theta }_{1}=\frac{1}{C\log \frac{\rho _{2}}{r}}\text{.}
\end{equation*}
Now, denoting by $\psi _{3}\left( z,\tau \right) =\psi _{1}\left( \Gamma
\left( \tau \right) \psi _{2}\left( z,\tau \right) ,\tau \right) $, we have
for every $\rho \in \left( 0,\rho _{2}\right) $ and $\tau \in \left( -\rho
_{2}^{2},0\right] $
\begin{equation}
Q_{C^{-1}\rho ,\varphi }\left( \tau \right) \subset \left( \psi _{3}\left(
.,\tau \right) \right) \left( B_{\rho }^{+}\right) \subset Q_{C\rho ,\varphi
}\left( \tau \right) \text{ ,}  \label{3.810}
\end{equation}
where $C$, $C>1$, depends on $\lambda $, $\Lambda $ and $E$ only. By (\ref%
{3.800}) and (\ref{3.810}), by using the change of variables $x=\psi
_{3}\left( z,\tau \right) $, $s=\tau $ we obtain (\ref{3.600}).$\blacksquare
$

\subsection{Smallness propagation estimates\label{Subs3.2}}

Let us introduce some notations that we shall use in Proposition \ref{Pr3.8}
below. Let $\alpha $, $\delta $ and $R$ be positive numbers such that $%
\alpha <\frac{\pi }{2}$ and $\delta \leq 1$. Let $\eta _{1}\in \left(
0,1\right) $ be defined in Theorem \ref{Th3.4}. Given $X_{0}=\left(
x_{0},t_{0}\right) \in \mathbb{R}^{n+1}$ and $\zeta \in \mathbb{R}^{n}$, $%
\left\vert \zeta \right\vert =1$, we denote
\begin{equation}
\mathcal{C}\left( x_{0},\zeta ,\alpha ,R\right) =\left\{ x\in B_{R}\left(
x_{0}\right) :\frac{\left( x-x_{0}\right) \cdot \zeta }{\left\vert
x-x_{0}\right\vert }>\cos \alpha \right\} \text{ ,}  \label{subs3.2.1}
\end{equation}
\begin{equation}
\mathcal{S}\left( X_{0},\zeta ,\alpha ,\delta ,R\right) =\left\{ \left(
x,t\right) \in \mathbb{R}^{n+1}:x\in \mathcal{C}\left( x_{0},\zeta ,\alpha
,R\right) \text{, }-\left( \delta \left( x-x_{0}\right) \cdot \zeta \right)
^{2}<t-t_{0}\leq 0\right\} ,  \label{subs3.2.3}
\end{equation}

\begin{equation*}
\alpha _{1}=\arcsin \left( \min \left\{ \sin \alpha ,\delta \left( 1-\sin
\alpha \right) \right\} \right) \text{ ,}
\end{equation*}
\begin{equation*}
\mu _{1}=\frac{R}{1+\sin \alpha _{1}}\text{ , }w_{1}=x_{0}+\mu _{1}\zeta
\text{ , }\rho _{1}=\frac{1}{4}\mu _{1}\eta _{1}\sin \alpha _{1}\text{ ,}
\end{equation*}
\begin{equation}
a=\frac{1-\frac{1}{4}\eta _{1}\sin \alpha _{1}}{1+\frac{1}{4}\eta _{1}\sin
\alpha _{1}}\text{.}  \label{3.829}
\end{equation}

\begin{proposition}
\label{Pr3.8}\textbf{(smallness propagation estimate)}. Let us take positive
numbers $\alpha ,\gamma ,\delta ,H$ and $R$ such that $\alpha <\frac{\pi }{2}
$, $\gamma \leq 1$ and $\delta \leq 1$. Let $L$ be the parabolic operator of
Theorem \ref{Th3.4}. Assume that $u$, $u\in H_{loc}^{2,1}\left( \mathcal{S}%
\left( X_{0},\zeta ,\alpha ,\delta ,R\right) \right) $, satisfies the
differential inequality
\begin{equation}
\left\vert Lu\right\vert \leq \Lambda \left( \frac{\left\vert \nabla
u\right\vert }{R}+\frac{\left\vert u\right\vert }{R^{2}}\right) \text{, in }%
\mathcal{S}\left( X_{0},\zeta ,\alpha ,\delta ,R\right)  \label{3.820}
\end{equation}%
and
\begin{equation}
\left\Vert u\right\Vert _{L^{\infty }\left( \mathcal{S}\left( X_{0},\zeta
,\alpha ,\delta ,R\right) \right) }+R^{\gamma }\left[ u\left( .,t_{0}\right) %
\right] _{\gamma ,\mathcal{C}\left( x_{0},\zeta ,\alpha ,R\right) }\leq H%
\text{ .}  \label{3.830}
\end{equation}%
Then, setting $\sigma _{1}=\left\Vert u\left( .,t_{0}\right) \right\Vert
_{L^{\infty }\left( B_{\rho _{1}\left( w_{1}\right) }\right) }$ we have
\begin{equation}
\left\vert u\left( x_{0},t_{0}\right) \right\vert \leq CH\left\vert \log
\frac{\sigma _{1}}{eH}\right\vert ^{-B}\text{ ,}  \label{3.840}
\end{equation}%
where $C$ depend on $\lambda ,\Lambda ,\alpha $ and $\gamma $ only and
\begin{equation}
B=\frac{\left\vert \log a\right\vert }{C_{1}}\text{ , }  \label{3.841}
\end{equation}%
where $a$ is given by (\ref{3.829}) and $C_{1}$, $C_{1}>0$, depends on $%
\lambda $ and $\Lambda $ only.
\end{proposition}

\textbf{Proof.} With no loss of generality we can assume $x_{0}=0$, $t_{0}=0$%
, $\zeta =e_{n}$.

\noindent Denote
\begin{equation*}
\mu _{k}=a^{k-1}\mu _{1}\text{ , }w_{k}=\mu _{k}e_{n}\text{ , }\rho
_{k}=a^{k-1}\rho _{1}\text{ , }d_{k}=\mu _{k}-\rho _{k}=a^{k-1}\mu
_{1}\left( 1-\eta _{1}\sin \beta \right) \text{ ,}
\end{equation*}
where $a$ is defined by (\ref{3.829}).

\noindent It is simple to check that, for every $k\geq 1$, the following
inclusions hold true
\begin{equation}
B_{\rho _{k+1}}\left( w_{k+1}\right) \subset B_{3\rho _{k}}\left(
w_{k}\right) \subset B_{4\eta _{1}^{-1}\rho _{k}}\left( w_{k}\right) \subset
\mathcal{C}\left( 0,e_{n},\beta ,R\right) \text{ ,}  \label{3.850}
\end{equation}
\begin{equation}
B_{4\eta _{1}^{-1}\rho _{k}}\left( w_{k}\right) \times \left( -4\eta
_{1}^{-1}\rho _{k},0\right] \subset \mathcal{S}\left( 0,e_{n},\alpha ,\delta
,R\right) \text{ .}  \label{3.860}
\end{equation}
For a number $r$, $r\in \left( 0,d_{1}\right] $, to be chosen later, let $%
\overline{k}$ be the smallest positive integer such that $d_{k}\leq r$. We
have
\begin{equation}
\frac{\left\vert \log \left( r/d_{1}\right) \right\vert }{\left\vert \log
a\right\vert }\leq \overline{k}-1\leq \frac{\left\vert \log \left(
r/d_{1}\right) \right\vert }{\left\vert \log a\right\vert }+1\text{.}
\label{3.870}
\end{equation}
Denote
\begin{equation*}
\sigma _{j}=\left( \int\nolimits_{B_{\rho _{j}}\left( w_{j}\right)
}u^{2}\left( x,0\right) dx\right) ^{1/2}\text{ , }j=1,2...\overline{k}\text{
.}
\end{equation*}
By Theorem \ref{Th3.4}, (\ref{3.850}), (\ref{3.860}) and since
\begin{equation*}
\sigma _{j+1}\leq \left( \int\nolimits_{B_{3\rho _{j}}\left( w_{j}\right)
}u^{2}\left( x,0\right) dx\right) ^{1/2}\text{ , }j=1,2...\overline{k}-1%
\text{ ,}
\end{equation*}
we obtain
\begin{equation}
\sigma _{j+1}^{2}\leq CH^{2\left( 1-\theta _{\ast }\right) }\sigma
_{j}^{2\theta _{\ast }}\text{, }j=1,2...\overline{k}-1\text{ ,}
\label{3.880}
\end{equation}
where $\mathcal{S}=\mathcal{S}\left( 0,e_{n},\alpha ,\delta ,R\right) $,
\begin{equation*}
\theta _{\ast }=\frac{1}{C_{0}\log \frac{4}{\eta _{1}}}
\end{equation*}
and $C$, $C_{0}$, depend on $\lambda $ and $\Lambda $ only.

\noindent By iterating (\ref{3.880}) we get
\begin{equation}
\sigma _{\overline{k}}^{2}\leq C^{\frac{1}{1-\theta _{\ast }}}H^{2\left(
1-\theta _{\ast }^{\overline{k}}\right) }\sigma _{1}^{2\theta _{\ast }^{%
\overline{k}}}\text{ .}  \label{3.890}
\end{equation}
Let us recall the following interpolation inequality
\begin{equation}
\left\Vert v\right\Vert _{L^{\infty }\left( B_{\rho }\right) }\leq C\left(
\left( \int\nolimits_{B_{\rho }}v^{2}dx\right) ^{\frac{\gamma }{2\gamma +n}}%
\left[ v\right] _{\gamma ,B_{\rho }}^{\frac{n}{2\gamma +n}}+\left( \rho
^{-n}\int\nolimits_{B_{\rho }}v^{2}dx\right) ^{1/2}\right) \text{ ,}
\label{3.900}
\end{equation}
where $C$ is an absolute constant.

By (\ref{3.830}) and (\ref{3.900}) we have
\begin{equation}
\left\Vert u\left( .,0\right) \right\Vert _{L^{\infty }\left( B_{\rho _{%
\overline{k}}}\left( w_{\overline{k}}\right) \right) }\leq CH^{\frac{n}{%
2\gamma +n}}\sigma _{\overline{k}}^{\frac{2\gamma }{2\gamma +n}}\text{ ,}
\label{3.910}
\end{equation}
where $C$ depend on $\lambda $ and $\Lambda $ only.

\noindent Let us consider the point $x_{r}=re_{n}$. Noticing that $x\in
B_{\rho _{\overline{k}}}\left( w_{\overline{k}}\right) $, by (\ref{3.830})
and (\ref{3.910}) we have
\begin{equation*}
\left\vert u\left( 0,0\right) \right\vert \leq \left\vert u\left(
x_{r},0\right) -u\left( 0,0\right) \right\vert +\left\vert u\left(
x_{r},0\right) \right\vert \leq CH\left( \left( \frac{r}{d_{1}}\right)
^{\gamma }+\left( \frac{\sigma _{\overline{k}}}{H}\right) ^{\frac{2\gamma }{%
2\gamma +n}}\right) \text{ ,}
\end{equation*}
where $C$ depend on $\lambda $ and $\Lambda $ only.

\noindent From the above obtained inequality and by (\ref{3.890}) we get
\begin{equation}
\left\vert u\left( 0,0\right) \right\vert \leq CH\left( \left( \frac{r}{d_{1}%
}\right) ^{\gamma }+\left( \frac{\sigma _{1}}{H}\right) ^{\frac{2\gamma
\theta _{\ast }^{\overline{k}}}{2\gamma +n}}\right) \text{ ,}  \label{3.920}
\end{equation}
where $C$ depend on $\lambda $, $\Lambda $ and $\alpha $ only.

\noindent Let us choose
\begin{equation*}
r=d_{1}\left\vert \log \left( e^{-1}\left( \frac{\sigma _{1}}{H}\right) ^{%
\frac{2\gamma }{2\gamma +n}}\right) \right\vert ^{-\frac{\left\vert \log
a\right\vert }{2\left\vert \log \theta _{\ast }\right\vert }}\text{ .}
\end{equation*}
By (\ref{3.870}) and (\ref{3.920}) we have
\begin{equation*}
\left\vert u\left( 0,0\right) \right\vert \leq CH\left\vert \log \left(
\frac{\sigma _{1}}{eH}\right) \right\vert ^{-\frac{\left\vert \log
a\right\vert }{2\gamma \left\vert \log \theta _{\ast }\right\vert }}
\end{equation*}
and we obtain(\ref{3.840}) with $B=\frac{\left\vert \log a\right\vert }{%
2\gamma \left\vert \log \theta _{\ast }\right\vert }$.$\blacksquare $

\subsection{Stability estimates from Cauchy data\label{Subs3.3}}

We are given positive numbers $R$, $T$, $E$, $\lambda $ and $\Lambda $ such
that $\lambda \in \left( 0,1\right] $. Let us consider the following
parabolic operator
\begin{equation}
Lu=\partial _{i}\left( g^{ij}\left( x,t\right) \partial _{j}u\right)
-\partial _{t}u+b_{i}\left( x,t\right) \partial _{i}u+c\left( x,t\right) u%
\text{ ,}  \label{3.930}
\end{equation}
where $\left\{ g^{ij}\left( x,t\right) \right\} _{i,j=1}^{n}$ is a real
symmetric $n\times n$ matrix. When $\xi \in \mathbb{R}^{n}$ $\left(
x,t\right) ,\left( y,\tau \right) \in \mathbb{R}^{n+1}$ assume that
\begin{equation}
\lambda \left\vert \xi \right\vert ^{2}\leq
\sum\limits_{i,j=1}^{n}g^{ij}\left( x,t\right) \xi _{i}\xi _{j}\leq \lambda
^{-1}\left\vert \xi \right\vert ^{2}  \label{3.931}
\end{equation}
and
\begin{equation}
\left( \sum\limits_{i,j=1}^{n}\left( g^{ij}\left( x,t\right) -g^{ij}\left(
y,\tau \right) \right) ^{2}\right) ^{1/2}\leq \frac{\Lambda }{R}\left(
\left\vert x-y\right\vert ^{2}+\left\vert t-\tau \right\vert \right) ^{1/2}%
\text{.}  \label{3.932}
\end{equation}
Concerning $b_{i}\left( x,t\right) $, $i=1,...n$ and $c\left( x,t\right) $
we assume that
\begin{equation}
R\sum\limits_{i=1}^{n}\left\Vert b_{i}\right\Vert _{L^{\infty }\left(
\mathbb{R}^{n+1}\right) }+R^{2}\left\Vert c\right\Vert _{L^{\infty }\left(
\mathbb{R}^{n+1}\right) }\leq \Lambda \text{ .}  \label{3.940}
\end{equation}
Let $\varphi :\mathbb{R}^{n-1}\times \mathbb{R\rightarrow R}$ be a function
satisfying the conditions
\begin{equation}
\varphi \left( 0\right) =\left\vert \nabla _{x^{\prime }}\varphi \left(
0\right) \right\vert =0\text{ ,}  \label{3.950}
\end{equation}
\begin{equation}
\left\Vert \varphi \right\Vert _{L^{\infty }\left( B_{R}^{\prime }\right)
}+R\left\Vert \nabla _{x^{\prime }}\varphi \right\Vert _{L^{\infty }\left(
B_{R}^{\prime }\right) }+R^{2}\left\Vert D^{2}\varphi \right\Vert
_{L^{\infty }\left( B_{R}^{\prime }\right) }\leq ER\text{ .}  \label{3.960}
\end{equation}

\noindent For any positive numbers $\rho $ and $t_{0}$ denote by
\begin{equation*}
D_{\rho ,\varphi }=\left\{ x\in B_{\rho }:x_{n}>\varphi \left( x^{\prime
}\right) \right\} \text{ , }D_{\rho ,\varphi }^{t_{0}}=D_{\rho ,\varphi
}\times \left( 0,t_{0}\right) \text{ ,}
\end{equation*}
\begin{equation*}
\Gamma _{\rho ,\varphi }=\left\{ x\in B_{\rho }:x_{n}=\varphi \left(
x^{\prime }\right) \right\} \text{ , }\Gamma _{\rho ,\varphi
}^{t_{0}}=\Gamma _{\rho ,\varphi }\times \left( 0,t_{0}\right) \text{.}
\end{equation*}
If $x=\left( x^{\prime },\varphi \left( x^{\prime }\right) \right) $, we
denote by $\nu \left( x\right) $, or simply by $\nu $, the unit vector of $%
\mathbb{R}^{n-1}$%
\begin{equation*}
\nu \left( x\right) =\frac{\left( \nabla _{x^{\prime }}\varphi \left(
x^{\prime }\right) ,-1\right) }{\sqrt{1+\left\vert \nabla _{x^{\prime
}}\varphi \left( x^{\prime }\right) \right\vert ^{2}}}\text{ .}
\end{equation*}

\begin{theorem}
\label{Th3.9}Let $L$ be the parabolic operator (\ref{3.930}). Let $u\in
H^{2,1}\left( D_{R,\varphi }^{T}\right) $ be a solution to the equation
\begin{equation}
Lu=0\text{, in }D_{R,\varphi }^{T}\text{ .}  \label{3.970}
\end{equation}
Let
\begin{equation}
\varepsilon :=\left\Vert u\right\Vert _{H^{\frac{3}{2},\frac{3}{4}}\left(
\Gamma _{R,\varphi }^{T}\right) }+R\left\Vert g^{ij}\partial _{j}u\nu
_{i}\right\Vert _{H^{\frac{1}{2},\frac{1}{4}}\left( \Gamma _{R,\varphi
}^{T}\right) }\text{ .}  \label{3.980}
\end{equation}
We have that there exist constants $C$, $\eta _{3}\in \left( 0,1\right) $, $%
s_{1}\in \left( 0,1\right) $ depending on $\lambda $, $\Lambda $ and $E$
only such that for every $r\in \left( 0,\min \left\{ R,\sqrt{T}\right\}
\right) $ we have
\begin{equation}
\left\Vert u\right\Vert _{L^{2}\left( D_{\eta _{3}r,\varphi }\times \left(
r^{2},T\right) \right) }\leq C\left( \varepsilon ^{s_{1}}\left\Vert
u\right\Vert _{L^{2}\left( D_{R,\varphi }^{T}\right) }^{1-s_{1}}+\varepsilon
\right) \text{ .}  \label{3.990}
\end{equation}
If, in addition, $u$ satisfies
\begin{equation}
u\left( x,0\right) =0\text{ , in }D_{R,\varphi }\text{ ,}  \label{3.1000}
\end{equation}
then we have, for every $r\in \left( 0,R\right] $,
\begin{equation}
\left\Vert u\right\Vert _{L^{2}\left( D_{\eta _{3}r,\varphi }^{T}\right)
}\leq C\left( \varepsilon ^{s_{1}}\left\Vert u\right\Vert _{L^{2}\left(
D_{R,\varphi }^{T}\right) }^{1-s_{1}}+\varepsilon \right) \text{,}
\label{3.1010}
\end{equation}
where $C$ depends $\lambda $, $\Lambda $ and $E$ only.
\end{theorem}

\textbf{Proof.} By using the extension theorem in Sobolev spaces, \cite{LiMa}%
, there exists $v\in H^{2,1}\left( D_{R,\varphi }^{T}\right) $ such that
\begin{equation}
v=u\text{ , }g^{ij}\partial _{j}v\nu _{i}=g^{ij}\partial _{j}u\nu _{i}\text{
, on }\Gamma _{R,\varphi }^{T}  \label{3.1020}
\end{equation}%
and
\begin{equation}
\left\Vert v\right\Vert _{H^{2,1}\left( D_{R,\varphi }^{T}\right) }\leq
C\varepsilon \text{ ,}  \label{3.1030}
\end{equation}%
where $C$ depends on $E$ only. Let $w=u-v$. We have that $w\in H^{2,1}\left(
D_{R,\varphi }^{T}\right) $ and $w$ satisfies
\begin{equation}
Lw=-Lv\text{, in }D_{R,\varphi }^{T}\text{ ,}  \label{3.1040}
\end{equation}%
\begin{equation}
w=0\text{ , }g^{ij}\partial _{j}w\nu _{i}=0\text{ , on }\Gamma _{R,\varphi
}^{T}\text{ .}  \label{3.1045}
\end{equation}%
Define
\begin{equation*}
\overline{f}=\left\{
\begin{array}{c}
-Lv\text{ , in }D_{R,\varphi }^{T}\text{ ,} \\
0\text{ , in }\left( B_{R}\times \left( 0,T\right) \right) \smallsetminus
D_{R,\varphi }^{T}%
\end{array}%
\right.
\end{equation*}%
and
\begin{equation*}
\overline{w}=\left\{
\begin{array}{c}
w\text{ , in }D_{R,\varphi }^{T}\text{ ,} \\
0\text{ , in }\left( B_{R}\times \left( 0,T\right) \right) \smallsetminus
D_{R,\varphi }^{T}\text{.}%
\end{array}%
\right.
\end{equation*}%
Since $w=g^{ij}\partial _{j}w\nu _{i}=0$ , on $\Gamma _{R,\varphi }^{T}$, we
have $\overline{w}\in H^{2,1}\left( B_{R}\times \left( 0,T\right) \right) $
and $\overline{w}$ satisfies
\begin{equation}
L\overline{w}=\overline{f}\text{ , in }B_{R}\times \left( 0,T\right) \text{ .%
}  \label{3.1050}
\end{equation}%
Let $z\in H^{2,1}\left( B_{R}\times \left( 0,T\right) \right) $ be the
solution to the following initial-boundary value problem
\begin{equation}
\left\{
\begin{array}{c}
Lz=\overline{f}\text{ , in }B_{R}\times \left( 0,T\right) \text{ ,} \\
z_{\mid \partial _{p}\left( B_{R}\times \left( 0,T\right) \right) }=0\text{ ,%
}%
\end{array}%
\right.  \label{3.1060}
\end{equation}%
where $\partial _{p}\left( B_{R}\times \left( 0,T\right) \right) =\left(
\partial B_{R}\times \left( 0,T\right] \right) \cup \left( B_{R}\times
\left\{ 0\right\} \right) $ is the parabolic boundary of $B_{R}\times \left(
0,T\right) $. We have by (\ref{3.1030}) and (\ref{3.1060})
\begin{equation}
\left\Vert z\right\Vert _{H^{2,1}\left( B_{R}\times \left( 0,T\right)
\right) }\leq C\left\Vert \overline{f}\right\Vert _{L^{2}\left( B_{R}\times
\left( 0,T\right) \right) }\leq C^{\prime }\varepsilon \text{ ,}
\label{3.1070}
\end{equation}%
where $C$ and $C^{\prime }$ depend on $\lambda $, $\Lambda $ and $E$ only.
Denote by $w_{1}$ the function $\overline{w}-z$, by (\ref{3.1050}) and (\ref%
{3.1060}) we have $w_{1}\in H^{2,1}\left( B_{R}\times \left( 0,T\right)
\right) $ and
\begin{equation}
Lw_{1}=0\text{ , in }B_{R}\times \left( 0,T\right) \text{ .}  \label{3.1080}
\end{equation}%
Let $r$ be any number belonging to $\left( 0,\min \left\{ R,\sqrt{T}\right\} %
\right] $. Let us denote by $r_{1}=r\min \left\{ \frac{1}{2},\frac{1}{E}%
\right\} $, $r_{2}=\eta _{1}r_{1}$, $r_{3}=\dfrac{\eta _{1}r_{1}}{4}$, where
$\eta _{1}$ is defined as in Theorem \ref{Th3.4}. Notice that, by (\ref%
{3.950}) and (\ref{3.960}) we have
\begin{equation}
B_{r_{3}}\left( -r_{3}e_{n}\right) \subset B_{R}\smallsetminus D_{R,\varphi }%
\text{ .}  \label{3.1090}
\end{equation}%
Since $Lw_{1}=0$ , in $B_{r_{1}}\left( -r_{3}e_{n}\right) \times \left(
0,T\right) $, by Theorem \ref{Th3.4} we have, for every $t\in \left[ r^{2},T%
\right] $,
\begin{equation}
\left\Vert w_{1}\left( .,t\right) \right\Vert _{L^{2}\left( B_{r_{2}}\left(
-r_{3}e_{n}\right) \right) }^{2}  \label{3.1100}
\end{equation}%
\begin{equation*}
\leq C\left\Vert w_{1}\left( .,t\right) \right\Vert _{L^{2}\left(
B_{r_{3}}\left( -r_{3}e_{n}\right) \right) }^{2s_{1}}\left\Vert w_{1}\left(
.,t\right) \right\Vert _{L^{2}\left( B_{r_{1}}\left( -r_{3}e_{n}\right)
\times \left( 0,T\right) \right) }^{2\left( 1-s_{1}\right) }\text{ ,}
\end{equation*}%
where $C$ and $s_{1}$, $s_{1}\in \left( 0,1\right) $, depend on $\lambda $
and $\Lambda $ only. By integrating both sides of last inequality with
respect to $t$ over $\left( r^{2},T\right) $ and by H\"{o}lder inequality,
we get%
\begin{equation}
\left\Vert w_{1}\right\Vert _{L^{2}\left( B_{r_{2}}\left( -r_{3}e_{n}\right)
\times \left( r^{2},T\right) \right) }^{2}  \label{3.1110}
\end{equation}%
\begin{equation*}
\leq C^{\prime }\left\Vert w\right\Vert _{L^{2}\left( B_{r_{3}}\left(
-r_{3}e_{n}\right) \times \left( r^{2},T\right) \right) }^{2s_{1}}\left\Vert
w_{1}\right\Vert _{L^{2}\left( B_{r_{1}}\left( -r_{3}e_{n}\right) \times
\left( 0,T\right) \right) }^{2\left( 1-s_{1}\right) }\text{ ,}
\end{equation*}%
where $C$ depends on $\lambda $ and $\Lambda $ only.

\noindent By (\ref{3.1030}), (\ref{3.1070}), the triangle inequality and
recalling the definition of $w$ we have
\begin{equation}
\left\Vert w_{1}\right\Vert _{L^{2}\left( B_{r_{1}}\left( -r_{3}e_{n}\right)
\times \left( 0,T\right) \right) }\leq \left\Vert u\right\Vert _{L^{2}\left(
D_{R,\varphi }^{T}\right) }+C\varepsilon \text{ ,}  \label{3.1120}
\end{equation}%
where $C$ depends on $\lambda $, $\Lambda $ and $E$ only.

By (\ref{3.1090}) we have $w_{1}=-z$ in $B_{r_{3}}\left( -r_{3}e_{n}\right)
\times \left( 0,T\right) $, so using (\ref{3.1070}) we have
\begin{equation}
\left\Vert w_{1}\right\Vert _{L^{2}\left( B_{r_{2}}\left( -r_{3}e_{n}\right)
\times \left( r^{2},T\right) \right) }\leq C\varepsilon \text{ ,}
\label{3.1130}
\end{equation}
where $C$ depends on $\lambda $, $\Lambda $ and $E$ only.

\noindent Since $D_{3r_{3},\varphi }\subset B_{r_{2}}\left(
-r_{3}e_{n}\right) \cap D_{R,\varphi }$, using (\ref{3.1030}), (\ref{3.1070}%
) and the triangle inequality we have
\begin{equation*}
\left\Vert w_{1}\right\Vert _{L^{2}\left( B_{r_{2}}\left( -r_{3}e_{n}\right)
\times \left( r^{2},T\right) \right) }\geq \left\Vert w_{1}\right\Vert
_{L^{2}\left( D_{3r_{3},\varphi }\times \left( r^{2},T\right) \right) }
\end{equation*}%
\begin{equation*}
=\left\Vert u-v-z\right\Vert _{L^{2}\left( D_{3r_{3},\varphi }\times \left(
r^{2},T\right) \right) }\geq \left\Vert u\right\Vert _{L^{2}\left(
D_{3r_{3},\varphi }\times \left( r^{2},T\right) \right) }-C\varepsilon \text{
,}
\end{equation*}%
therefore
\begin{equation}
\left\Vert u\right\Vert _{L^{2}\left( D_{3r_{3},\varphi }\times \left(
r^{2},T\right) \right) }\leq \left\Vert w_{1}\right\Vert _{L^{2}\left(
B_{r_{2}}\left( -r_{3}e_{n}\right) \times \left( r^{2},T\right) \right)
}+C\varepsilon \text{ ,}  \label{3.1140}
\end{equation}%
where $C$ depend on $\lambda $, $\Lambda $ and $E$ only.

\noindent By (\ref{3.1110}), (\ref{3.1120}), (\ref{3.1130}) and (\ref{3.1140}%
) we obtain
\begin{equation*}
\left\Vert u\right\Vert _{L^{2}\left( D_{3r_{3},\varphi }\times \left(
r^{2},T\right) \right) }\leq C\left( \varepsilon ^{s_{1}}\left\Vert
u\right\Vert _{L^{2}\left( D_{R,\varphi }^{T}\right) }^{1-s_{1}}+\varepsilon
\right) \text{ ,}
\end{equation*}
where $C$ depend on $\lambda $, $\Lambda $ and $E$ only, so (\ref{3.990}) is
proved.

To prove (\ref{3.1010}) we only need some slight variations of the previous
proof. Indeed, due to (\ref{3.1000}) we can choose the function $v$ in such
a way that, besides (\ref{3.1120}) and (\ref{3.1130}), $v$ satisfies $%
v\left( .,0\right) =0$ on $D_{R,\varphi }$. Therefore $w_{1}\in
H^{2,1}\left( B_{R}\times \left( 0,T\right) \right) $ and
\begin{equation*}
w_{1}\left( .,0\right) =0\text{ , in }B_{R}\text{ .}
\end{equation*}%
Now, considering the trivial extension $\widetilde{w}_{1}$ of $w_{1}$ to $%
B_{R}\times \left( -\infty ,T\right) $%
\begin{equation*}
\widetilde{w}_{1}=\left\{
\begin{array}{c}
w_{1}\text{ , in }B_{R}\times \left( 0,T\right) \text{ ,} \\
0\text{ , in }B_{R}\times \left( -\infty ,0\right) \text{,}%
\end{array}%
\right.
\end{equation*}%
we have that $\widetilde{w}_{1}\in H^{2,1}\left( B_{R}\times \left( -\infty
,T\right) \right) $ and
\begin{equation*}
L\widetilde{w}_{1}=0\text{, in }B_{R}\times \left( -\infty ,T\right) \text{.}
\end{equation*}%
Therefore applying Theorem \ref{Th3.4} to the function $\widetilde{w}_{1}$
we have that (\ref{3.1100}) holds true for every $t\in \left[ 0,T\right] $.
Hence we can replace interval $\left( r^{2},T\right) $ with interval $\left(
0,T\right) $ in (\ref{3.1110}). Finally by definition of $w_{1}$ we have
\begin{equation*}
\left\Vert u\right\Vert _{L^{2}\left( D_{3r_{3},\varphi }^{T}\right) }\leq
C\left( \varepsilon ^{s_{1}}\left\Vert u\right\Vert _{L^{2}\left(
D_{R,\varphi }^{T}\right) }^{1-s_{1}}+\varepsilon \right) \text{ ,}
\end{equation*}%
where $C$ depend on $\lambda $, $\Lambda $ and $E$ only, so (\ref{3.1010})
is proved.$\blacksquare $

\section{Stability estimates for Dirichlet inverse problem with unknown
time-varying boundaries\label{Sec4}}

In this section $\left\{ \Omega \left( t\right) \right\} _{t\in \mathbb{R}}$
will be a family of bounded domains in $\mathbb{R}^{n}$, where $T$ is a
given positive number. We shall suppose that the boundary of $\Omega \left(
\left( -\infty ,+\infty \right) \right) =\bigcup\limits_{t\in \mathbb{R}%
}\Omega \left( t\right) \times \left\{ t\right\} $ is sufficiently smooth,
but for every $t\in \mathbb{R}$ a part $I\left( t\right) $\ of $\partial
\Omega \left( t\right) $ is not known. The inverse problem we are interested
in determining $I\left( t\right) $, for every $t\in \left( 0,T\right) $, by
means of thermal measurements on the accessible part $A\left( t\right)
:=\partial \Omega \left( t\right) \smallsetminus I\left( t\right) $. In what
follows, for the sake of simplicity, we shall assume that $A\left( t\right) $
is not time varying, so we set $A\left( t\right) =A$. Given a nontrivial
function $f$ on $A\times \left( 0,T\right) $, let us consider the following
Cauchy-Dirichlet (direct) problem
\begin{equation}
\left\{
\begin{array}{c}
div\left( \kappa \left( x,t,u\right) \nabla u\right) -\partial _{t}u=0\text{%
, \ \ in }\Omega \left( \left( 0,T\right) \right) \text{,} \\
u=f\text{, \ \ \ \ \ \ \ \ \ \ \ on }A\times \left( 0,T\right] \text{,\ } \\
u=0\text{, \ \ \ \ \ \ \ \ \ \ \ \ \ on }I\left( \left( 0,T\right] \right)
\text{, \ \ } \\
u\left( .,0\right) =u_{0}\text{, \ \ \ \ \ \ \ \ \ \ in }\Omega \left(
0\right) \text{ ,}%
\end{array}%
\right.  \label{4.1-4.4}
\end{equation}%
where $\kappa \left( x,t,u\right) =\left\{ \kappa ^{ij}\left( x,t,u\right)
\right\} _{i,j=1}^{n}$ denotes a known symmetric matrix which satisfies a
hypothesis of uniform ellipticity and some smothness conditions that we
shall specify below.

Assume that $\Omega \left( 0\right) $ is known. Given an open portion $%
\Sigma $ of $\partial \Omega \left( t\right) $ such that $\Sigma \subset A$,
we consider the inverse problems of determinig $I\left( t\right) $, $t\in
\left( 0,T\right] $, from the knowledge of $\kappa \left( x,t,u\right)
\nabla u\cdot \nu $ on $\Sigma \times \left[ 0,T\right] $, where $\nu $
denotes the exterior unit normal to $\Omega \left( t\right) $, for every $%
t\in \left[ 0,T\right] $.

\bigskip

\textit{i) A priori information on the domains }$\Omega \left( t\right) $.

Given $R_{0}$, $M$ positive numbers, we assume
\begin{equation}
\left\vert \Omega \left( t\right) \right\vert \leq MR_{0}^{n}\text{ , for
every }t\in \left[ 0,T\right] \text{,}  \label{4.10}
\end{equation}
here and below $\left\vert \Omega \left( t\right) \right\vert $ denotes the (%
$n$-dimensional) Lebesgue measure of $\Omega \left( t\right) $. For every $%
t\in \left[ 0,T\right] $ we shall distinguish two nonempty parts $A\left(
t\right) =A$, $I\left( t\right) $ in $\partial \Omega \left( t\right) $ and
we assume
\begin{equation}
I\left( t\right) \cup A=\partial \Omega \left( t\right) \text{ , Int}\left(
I\left( t\right) \right) \cap \text{Int}\left( A\right) =\varnothing \text{
, }I\left( t\right) \cap A=\partial A=\partial I\left( t\right) \text{ .}
\label{4.20}
\end{equation}
Here, interior and boundaries are intented in the relative topology in $%
\partial \Omega \left( t\right) $.

Concerning the regularity of $\partial \Omega \left( \left[ 0,T\right]
\right) $ we assume that, for a number $\beta \in \left( 0,1\right) $ and a
positive number $E$ we have
\begin{equation}
\partial \Omega \left( \left[ 0,T\right] \right) \text{ is a portion of }%
\partial \Omega \left( \left( -\infty ,+\infty \right) \right)  \label{4.25}
\end{equation}%
\begin{equation*}
\text{of class }C^{2,\beta }\text{ with constants }R_{0}\text{, }E
\end{equation*}%
and also
\begin{equation}
A\times \left[ 0,T\right] \text{, }I\left( \left[ 0,T\right] \right) \text{
are portions of }\partial \Omega \left( \left( -\infty ,+\infty \right)
\right)  \label{4.26}
\end{equation}%
\begin{equation*}
\text{of class }C^{2,\beta }\text{ with constants }R_{0}\text{, }E\text{.}
\end{equation*}%
In addition we assume that there exist open portions $\Sigma $, $\widetilde{%
\Sigma }$ within $A$ such that $\Sigma \subset \widetilde{\Sigma }$ and
there exists a point $P_{0}\in \Sigma $ such that
\begin{equation}
\partial \Omega \left( t\right) \cap B_{R_{0}}\left( P_{0}\right) \subset
\Sigma \text{ , for every }t\in \left[ 0,T\right]  \label{4.30}
\end{equation}%
and
\begin{equation}
\widetilde{\Sigma }\cap \left( I\left( t\right) \right) ^{R_{0}}=\varnothing
\text{ , for every }t\in \left[ 0,T\right] \text{ ,}  \label{4.35}
\end{equation}%
where $\left( I\left( t\right) \right) ^{R_{0}}=\left\{ x\in \partial \Omega
\left( t\right) :dist\left( x,I\left( t\right) \right) <R_{0}\right\} $.

\begin{remark}
\label{Re4.1}Let $t$ be any number in $\left[ 0,T\right] $. Observe that (%
\ref{4.25}) automatically implies a lower bound on the diameter of every
connected component of $\partial \Omega \left( t\right) $. In addition, by
combining (\ref{4.10}) with (\ref{4.25}), an upper bound of the diameter of $%
\Omega \left( t\right) $. Note also that (\ref{4.10}) and (\ref{4.25})
implicitly comprise an a priori upper bound on the number of connected
components of $\partial \Omega \left( t\right) $.
\end{remark}

\textit{ii) Assumptions about the initial and boundary data.}

Let $F_{0}$ be a given positive number and let $F_{1}:\left[ 0,T\right]
\rightarrow \left[ 0,+\infty \right) $ be a given strictly increasing
continuous function on $\left[ 0,T\right] $ such that $F_{1}\left( 0\right)
=0$. We assume the following conditions on the initial data $u_{0}$ and the
Dirichlet data $f$ appearing in problem (\ref{4.1-4.4})
\begin{equation}
u_{0}\in C^{2,\beta }\left( \Omega \left( 0\right) \right) \text{ , }u_{0}=0%
\text{ on }I\left( 0\right) \text{ ,}  \label{4.38}
\end{equation}%
\begin{equation}
\left\Vert u_{0}\right\Vert _{C^{2,\beta }\left( \Omega \left( 0\right)
\right) }\leq F_{0}\text{,}  \label{4.39}
\end{equation}%
\begin{equation}
f\in C^{2,\beta }\left( A\times \left( 0,T\right) \right) \text{,}
\label{4.40}
\end{equation}%
\begin{equation}
f\left( .,0\right) =u_{0}\left( .\right) \text{ on }A\text{,}  \label{4.41}
\end{equation}%
\begin{equation}
\text{\textit{supp\thinspace }}f\left( .,t\right) \subset \widetilde{\Sigma }%
\text{, for every }t\in \left[ 0,T\right] \text{,}  \label{4.45}
\end{equation}%
\begin{equation}
\left\Vert f\right\Vert _{C^{2,\beta }\left( A\times \left( 0,T\right)
\right) }\leq F_{0}\text{,}  \label{4.50}
\end{equation}%
\begin{equation}
\left\Vert f\left( .,t\right) \right\Vert _{L^{\infty }\left( A\right) }\geq
F_{1}\left( t\right) \text{ , for every }t\in \left[ 0,T\right] \text{.}
\label{4.55}
\end{equation}

\bigskip

\textit{iii) Assumptions on the thermal conductivity }$\kappa $.

The thermal conductivity $\kappa $ is assumed to be a given function from $%
\mathbb{R}^{n}\times \mathbb{R\times R}$ with values $n\times n$ matrices
satisfying the following conditions. For given constants $\lambda _{0}$, $%
\Lambda _{0}$, $\lambda _{0}\in \left( 0,1\right] $, $\Lambda _{0}\geq 0$,
we have, for every $\xi \in \mathbb{R}^{n}$ and $\left( x,t,z\right) \in
\mathbb{R}^{n}\times \mathbb{R\times R}$,
\begin{equation}
\lambda _{0}\left\vert \xi \right\vert ^{2}\leq \kappa \left( x,t,z\right)
\xi \cdot \xi \leq \lambda _{0}^{-1}\left\vert \xi \right\vert ^{2}\text{, }
\label{4.60}
\end{equation}
\begin{equation}
R_{0}\left\Vert \nabla _{x}\kappa \right\Vert _{L^{\infty }\left( \mathbb{R}%
^{n+2}\right) }+R_{0}^{2}\left\Vert \partial _{t}\kappa \right\Vert
_{L^{\infty }\left( \mathbb{R}^{n+2}\right) }+\left\Vert \partial
_{u}^{2}\kappa \right\Vert _{L^{\infty }\left( \mathbb{R}^{n+2}\right) }\leq
\Lambda _{0}\text{ .}  \label{4.65}
\end{equation}
In the sequel we shall refer to the set $\left\{ \lambda _{0},\Lambda
_{0},E,M,F_{0},\beta ,R_{0}^{2}T^{-1},F_{1}\left( .\right) \right\} $ as to
the \textit{a priori data}.

\begin{theorem}
\label{Th4.1}Let $\left\{ \Omega _{1}\left( t\right) \right\} _{t\in \mathbb{%
R}}$, $\left\{ \Omega _{2}\left( t\right) \right\} _{t\in \mathbb{R}}$ be
two families of domains satisfying (\ref{4.10}), (\ref{4.25}). For any $%
i=1,2 $ and any $t\in \left[ 0,T\right] $ let $A_{i}\left( t\right) $, $%
I_{i}\left( t\right) $ satisfying (\ref{4.20}) and (\ref{4.26}), be the
accessible and inaccessible part of $\partial \Omega _{i}\left( t\right) $,
respectively. Assume that for any $t\in \left[ 0,T\right] $ $A_{1}\left(
t\right) =A_{2}\left( t\right) =A$ and that $\Omega _{1}\left( t\right) $
and $\Omega _{2}\left( t\right) $ lie on the same side of $A$. Let us take $%
\Sigma $, $\widetilde{\Sigma }$ such that $\Sigma \subset \widetilde{\Sigma }%
\subset A$ satisfying (\ref{4.35}). Assume that
\begin{equation}
\Omega _{1}\left( 0\right) =\Omega _{2}\left( 0\right) \text{.}  \label{4.66}
\end{equation}
Let $u_{0}\in C^{2,\beta }\left( \Omega _{1}\left( 0\right) \right) $ and $%
f\in C^{2,\beta }\left( A\times \left( 0,T\right) \right) $ satisfy (\ref%
{4.38})-(\ref{4.55}). Let us assume that (\ref{4.60}) and (\ref{4.65}) are
satisfied and let $u_{i}\in C^{2,\beta }\left( \Omega _{i}\left( \left[ 0,T%
\right] \right) \right) $ be the solution to (\ref{4.1-4.4}) when $\Omega
\left( \left[ 0,T\right] \right) =\Omega _{i}\left( \left[ 0,T\right]
\right) $, $i=1,2$.

If, given $\varepsilon >0$, we have
\begin{equation}
R_{0}\left\Vert \kappa _{1}\nabla u_{1}\cdot \nu -\kappa _{2}\nabla
u_{2}\cdot \nu \right\Vert _{L^{2}\left( \Sigma \times \left( 0,T\right)
\right) }\leq \varepsilon \text{ ,}  \label{4.70}
\end{equation}%
where $\kappa _{i}=\kappa \left( x,t,u_{i}\left( x,t\right) \right) $, $%
i=1,2 $, then we have
\begin{equation}
\sup\limits_{t\in \left[ \tau ,T\right] }d_{\mathcal{H}}\left( \overline{%
\Omega _{1}\left( t\right) },\overline{\Omega _{2}\left( t\right) }\right)
\leq R_{0}C_{1}\left( \tau \right) \left\vert \log \varepsilon \right\vert
^{-\frac{1}{C_{2}\left( \tau \right) }}\text{ ,}  \label{4.75}
\end{equation}
for every $\tau \in \left( 0,T\right] $ and $\varepsilon \in \left(
0,1\right) $, where $C_{1}\left( \tau \right) $ and $C_{2}\left( \tau
\right) $ are positive constants depending on $\tau ,\lambda _{0},\Lambda
_{0},E,M,F_{0},\beta ,R_{0}^{2}T^{-1}$ and $F_{1}\left( .\right) $ only.
\end{theorem}

Here $d_{\mathcal{H}}$ denotes the Hausdorff distance, namely for every $%
t\in \left[ 0,T\right] $%
\begin{equation*}
d_{\mathcal{H}}\left( \overline{\Omega _{1}\left( t\right) },\overline{%
\Omega _{2}\left( t\right) }\right) =\max \left\{ \sup\limits_{x\in
\overline{\Omega _{1}\left( t\right) }}dist\left( x,\overline{\Omega
_{2}\left( t\right) }\right) ,\sup\limits_{x\in \overline{\Omega _{2}\left(
t\right) }}dist\left( x,\overline{\Omega _{1}\left( t\right) }\right)
\right\} \text{.}
\end{equation*}

\subsection{Proof of Theorem \protect\ref{Th4.1}\label{Subs4.1}}

Here and in the sequel we shall denote, for any $t\in \mathbb{R}$, by $%
G\left( t\right) $ the connected component of $\Omega _{1}\left( t\right)
\cap \Omega _{2}\left( t\right) $ such that $\widetilde{\Sigma }\subset
\overline{G\left( t\right) }$.

From regularity estimates for solutions to the Cauchy-Dirichlet problem for
parabolic equations \cite{Li}, we have
\begin{equation}
\left\Vert u_{i}\right\Vert _{C^{2,\beta }\left( \Omega _{i}\left( \left[ 0,T%
\right] \right) \right) }\leq H\text{ ,}  \label{4.95}
\end{equation}%
where $H$ is a positive constant depending on $\lambda _{0},\Lambda
_{0},E,M,F_{0},\beta $ and $R_{0}^{2}T^{-1}$ only.

The proof of Theorem \ref{Th4.1} is obtained from a sequence of
propositions. In the present section we give the statement of such
propositions in the next section we shall prove them.

\begin{proposition}
\label{Pr4.2}\textbf{(Stability estimates of continuation from Cauchy data)}%
. Let the hypotheses of Theorem \ref{Th4.1} be satisfied. We have
\begin{equation}
\int\nolimits_{\Omega _{i}\left( t\right) \smallsetminus G\left( t\right)
}u_{i}^{2}\left( x,t\right) dx\leq R_{0}^{n}H^{2}\omega \left( \frac{%
\varepsilon }{H}\right) \text{ , }i=1,2\text{, }t\in \left[ 0,T\right] \text{
,}  \label{4.85}
\end{equation}
where $\omega $ is an increasing continuous function on $\left[ 0,+\infty
\right) $ which satisfies
\begin{equation}
\omega \left( s\right) \leq C\left( \log \left\vert \log s\right\vert
\right) ^{-\frac{1}{n}}\text{ , for every }s\in \left( 0,e^{-1}\right)
\label{4.90}
\end{equation}
and $C$ depends on $\lambda _{0},\Lambda _{0},E,M,F_{0},\beta $ and $%
R_{0}^{2}T^{-1}$ only.
\end{proposition}

\begin{proposition}
\label{Pr4.3}\textbf{(Improved stability estimate of continuation from
Cauchy data)}. Let the hypotheses of Proposition \ref{Pr4.2} satisfied and,
in addition let us assume that there exist $L>0$ and $\rho _{0}\in \left(
0,R_{0}\right] $ such that, for any $t\in \left[ 0,T\right] $, $\partial
G\left( t\right) $ is of Lipschitz class with constants $\rho _{0},L$. Then (%
\ref{4.85}) holds with $\omega $ satisfying
\begin{equation}
\omega \left( s\right) \leq C\left\vert \log s\right\vert ^{-\frac{1}{C}}%
\text{ , for every }s\in \left( 0,1\right) \text{,}  \label{4.185}
\end{equation}
where $C$, $C>1$, depends on $\lambda _{0},\Lambda _{0},E,M,F_{0},\beta $
and $R_{0}^{2}T^{-1}$ only.
\end{proposition}

\begin{proposition}
\label{Pr4.4}\textbf{(Lipschitz stability estimate of continuation from the
interior)}. Let $\left\{ \Omega \left( t\right) \right\} _{t\in \mathbb{R}}$
be a family of domains satisfying (\ref{4.10}), (\ref{4.20}), (\ref{4.25})
and (\ref{4.26}). Let $u_{0}\in C^{2,\beta }\left( \Omega \left( 0\right)
\right) $ satisfy (\ref{4.38}) and (\ref{4.39}). Let $f\in C^{2,\beta
}\left( A\times \left( 0,T\right) \right) $ satisfy (\ref{4.41})-(\ref{4.55}%
). Assume that (\ref{4.60}) and (\ref{4.65}) are satisfied and let $u\in
C^{2,\beta }\left( \Omega \left( \left[ 0,T\right] \right) \right) $ be the
solution to (\ref{4.1-4.4}).

For every $\overline{\rho }>0$, every $t\in \left[ 0,T\right] $ and every $%
z\in \left( \Omega \left( t\right) \right) _{2\overline{\rho }}$ we have
\begin{eqnarray}
&&\int\nolimits_{B_{\overline{\rho }}\left( z\right) }u^{2}\left(
x,t\right) dx  \label{4.220} \\
&\geq &H^{2}\overline{\rho }^{n}\exp \left( -e^{\frac{CR_{0}^{n}}{\min
\left\{ t^{n/2},\overline{\rho }^{n}\right\} }}\left( \frac{F_{1}\left(
t\right) }{H}\right) ^{\frac{CR_{0}}{\min \left\{ t^{1/2},\overline{\rho }%
\right\} }}\right) \text{ ,}  \notag
\end{eqnarray}
where $C$, $C>1$, depends on $\lambda _{0},\Lambda _{0},E,M,F_{0},\beta $
and $R_{0}^{2}T^{-1}$ only.
\end{proposition}

In what follows we need this definition.

\begin{definition}
\label{Defgrafrel}\textbf{(relative graphs)}. We shall say that two bounded
domains $\Omega _{1}$ and $\Omega _{2}$ in $\mathbb{R}^{n}$ of class $%
C^{1,\beta }$ with constants $R_{0},E$ are \textit{relative graphs} if for
any $P\in \partial \Omega _{1}$ there exists a rigid transformation of
coordinates under which we have $P\equiv 0$ and there exist $\varphi
_{P,1},\varphi _{P,2}\in C^{1,\beta }\left( B_{r_{0}}^{\prime }\left(
0\right) \right) $, where $\dfrac{r_{0}}{R_{0}}\leq 1$ depends on $E$ and $%
\beta $ only, satisfying the following conditions

i) $\varphi _{P,1}\left( 0\right) =0$ , $\left\vert \varphi _{P,2}\left(
0\right) \right\vert \leq \dfrac{r_{0}}{2}$,

ii) $\left\Vert \varphi _{P,i}\right\Vert _{C^{1,\beta }\left(
B_{r_{0}}^{\prime }\left( 0\right) \right) }\leq ER_{0}$,

iii) $\Omega _{i}\cap B_{r_{0}}\left( 0\right) =\left\{ x\in B_{r_{0}}\left(
0\right) :x_{n}>\varphi _{P,i}\left( x^{\prime }\right) \right\} $, $i=1,2$.

\noindent We shall denote
\begin{equation}
\gamma \left( \Omega _{1},\Omega _{2}\right) =\sup\limits_{P\in \partial
\Omega _{1}}\left\Vert \varphi _{P,1}-\varphi _{P,2}\right\Vert _{L^{\infty
}\left( B_{r_{0}}^{\prime }\left( 0\right) \right) }\text{ .}  \label{4.304}
\end{equation}
\end{definition}

\begin{proposition}
\label{Pr4.5}Let $\Omega _{1}$ and $\Omega _{2}$ be two domains in $\mathbb{R%
}^{n}$ of class $C^{1,\beta }$ with constants $R_{0},E$ and satisfying $%
\left\vert \Omega _{i}\right\vert \leq MR_{0}^{n}$. Assume that $\Omega _{1}$
and $\Omega _{2}$ are relative graphs. Then there exists a $C^{1,\beta }$
diffeomorphism $\Phi :\mathbb{R}^{n}\rightarrow \mathbb{R}^{n}$ such that $%
\Phi \left( \Omega _{2}\right) =\Omega _{1}$ and
\begin{equation}
\left\Vert \Phi -Id\right\Vert _{L^{\infty }\left( \mathbb{R}^{n}\right)
}\leq C\gamma \left( \Omega _{1},\Omega _{2}\right) \text{ ,}  \label{4.305}
\end{equation}%
where $C$ depends on $E,M$ and $\beta $ only.
\end{proposition}

\begin{proposition}
\label{Pr4.6} Let $\left\{ \Omega _{1}\left( t\right) \right\} _{t\in
\mathbb{R}}$ and $\left\{ \Omega _{2}\left( t\right) \right\} _{t\in \mathbb{%
R}}$ be two families of domains satisfying (\ref{4.10}), (\ref{4.25}). We
have, for every $t,\tau \in \left[ 0,T\right] $,
\begin{equation}
\left\vert d_{\mathcal{H}}\left( \overline{\Omega _{1}\left( t\right) },%
\overline{\Omega _{2}\left( t\right) }\right) -d_{\mathcal{H}}\left(
\overline{\Omega _{1}\left( \tau \right) },\overline{\Omega _{2}\left( \tau
\right) }\right) \right\vert \leq C\frac{\left\vert t-\tau \right\vert }{%
R_{0}}\text{ ,}  \label{4.310}
\end{equation}
where $C$ depends on $E,M$ and $\beta $ only.
\end{proposition}

\begin{definition}
\label{Def4.7}\textbf{(modified distance). }Let $\Omega _{1}$ and $\Omega
_{2}$ be bounded domains in $\mathbb{R}^{n}$. We call modified distance
between $\Omega _{1}$ and $\Omega _{2}$ the number
\begin{equation}
d_{m}\left( \Omega _{1},\Omega _{2}\right) =\max \left\{ \sup\limits_{x\in
\partial \Omega _{1}}dist\left( x,\overline{\Omega _{2}}\right)
,\sup\limits_{x\in \partial \Omega _{2}}dist\left( x,\overline{\Omega _{1}}%
\right) \right\} \text{ .}  \label{4.350}
\end{equation}
\end{definition}

Notice that
\begin{equation}
d_{m}\left( \Omega _{1},\Omega _{2}\right) \leq d_{\mathcal{H}}\left(
\overline{\Omega _{1}},\overline{\Omega _{2}}\right) \text{ , }
\label{4.355}
\end{equation}
but, in general, the reverse inequality does not hold as the following
example makes clear: $\Omega _{1}=B_{1}\left( 0\right) $, $\Omega
_{2}=B_{1}\left( 0\right) \smallsetminus B_{1/2}\left( 0\right) $. In this
case $d_{m}\left( \Omega _{1},\Omega _{2}\right) =0$, whereas $d_{\mathcal{H}%
}\left( \Omega _{1},\Omega _{2}\right) =\dfrac{1}{2}$. However, the
following proposition holds true, \cite{AlBRVe1}.

\begin{proposition}
\label{Pr4.8}Let $\Omega _{1}$ and $\Omega _{2}$ be bounded domains in $%
\mathbb{R}^{n}$ of class $C^{1,\beta }$ with constants $R_{0}$, $E$ and
satisfying $\left\vert \Omega _{i}\right\vert \leq MR_{0}^{n}$. There exist
numbers $d_{0}$, $\rho _{0}\in \left( 0,R_{0}\right] $ such that $\dfrac{%
d_{0}}{R_{0}}$ and $\dfrac{\rho _{0}}{R_{0}}$ depend on $\beta $ and $E$
only, such that if we have
\begin{equation}
d_{\mathcal{H}}\left( \overline{\Omega _{1}},\overline{\Omega _{2}}\right)
\leq d_{0}\text{,}  \label{4.360}
\end{equation}
then the following facts hold true

\noindent i) $\Omega _{1}$ and $\Omega _{2}$ are relative graphs and
\begin{equation}
\gamma \left( \Omega _{1},\Omega _{2}\right) \leq Cd_{\mathcal{H}}\left(
\overline{\Omega _{1}},\overline{\Omega _{2}}\right) \text{ ,}  \label{4.361}
\end{equation}
where $C$ depends $\beta $ and $E$ only,

\noindent ii) there exists an absolute positive constant $c$ such that
\begin{equation}
d_{\mathcal{H}}\left( \overline{\Omega _{1}},\overline{\Omega _{2}}\right)
\leq cd_{m}\left( \Omega _{1},\Omega _{2}\right) \text{ ,}  \label{4.365}
\end{equation}
iii) any connected component of $\Omega _{1}\cap \Omega _{2}$ has boundary
of Lipschitz class with constants $\rho _{0}$, $L$, where $\rho _{0}$ is as
above and $L>0$ depends on $E$ only.
\end{proposition}

\textbf{Proof of Theorem \ref{Th4.1}.} For the sake of simplicity we denote,
for any $t\in \left[ 0,T\right] $, $d\left( t\right) =d_{\mathcal{H}}\left(
\overline{\Omega _{1}\left( t\right) },\overline{\Omega _{2}\left( t\right) }%
\right) $ and $d_{m}\left( t\right) =d_{m}\left( \Omega _{1}\left( t\right)
,\Omega _{2}\left( t\right) \right) $. We will prove that denoting
\begin{equation}
\sigma =R_{0}^{-n}\max \left\{ \sup\limits_{t\in \left[ 0,T\right]
}\int\nolimits_{\Omega _{1}\left( t\right) \smallsetminus G\left( t\right)
}u_{1}^{2}\left( x,t\right) dx,\sup\limits_{t\in \left[ 0,T\right]
}\int\nolimits_{\Omega _{2}\left( t\right) \smallsetminus G\left( t\right)
}u_{2}^{2}\left( x,t\right) dx\right\} \text{ ,}  \label{4.370}
\end{equation}
we have
\begin{equation}
d_{m}\left( t\right) \leq \frac{CR_{0}^{n+1}}{b\left( t\right) \min \left\{
t^{n/2},R_{0}^{n}\right\} }\left( \frac{\sigma }{MH^{2}}\right) ^{\frac{1}{%
\log \frac{C}{b\left( t\right) }}}\text{ }  \label{4.375}
\end{equation}
and
\begin{equation}
d\left( t\right) \leq \frac{CR_{0}^{n+2}}{b\left( t\right) \min \left\{
t^{\left( n+1\right) /2},R_{0}^{n+1}\right\} }\left( \frac{\sigma }{MH^{2}}%
\right) ^{\frac{1}{\log \frac{C}{b\left( t\right) }}}\text{ ,}  \label{4.376}
\end{equation}
where
\begin{equation}
b\left( t\right) =\exp \left( -e^{\frac{CR_{0}^{n}}{\min \left\{
t^{n/2},R_{0}^{n}\right\} }}\left( \frac{F_{1}\left( t\right) }{H}\right) ^{-%
\frac{CR_{0}}{\min \left\{ t^{1/2},R_{0}\right\} }}\right)  \label{4.377}
\end{equation}
and $C$ is a constant depending on $\lambda _{0},\Lambda
_{0},E,M,F_{0},\beta $ and $R_{0}^{2}T^{-1}$ only.

First we prove (\ref{4.375}). Let $t\in \left[ 0,T\right] $ be fixed. Let us
assume, without loss of generality, that there exists $x_{0}\in I_{1}\left(
t\right) \subset \partial \Omega _{1}\left( t\right) $ such that $dist\left(
x_{0},\Omega _{2}\left( t\right) \right) =d_{m}\left( t\right) $. By (\ref%
{4.370}) we have
\begin{equation}
\int\nolimits_{\Omega _{1}\left( t\right) \cap B_{d_{m}\left( t\right)
}\left( x_{0}\right) }u_{1}^{2}\left( x,t\right) dx\leq R_{0}^{n}\sigma
\text{ .}  \label{4.380}
\end{equation}
Let $\eta _{2}\in \left( 0,1\right) $ be defined in Theorem \ref{Th3.4}, $%
\eta _{2}$ depends on $\lambda _{0},\Lambda _{0},E,M,F_{0},\beta $ and $%
R_{0}^{2}T^{-1}$ only. Let us denote $R\left( t\right) =\min \left\{
t^{1/2},R_{0}\right\} $. We distinguish two cases. If $d_{m}\left( t\right)
\geq \dfrac{\eta _{2}R\left( t\right) }{2}$, let $\overline{d}\left(
t\right) =\frac{\eta _{2}R\left( t\right) }{2\left( 1+\sqrt{1+E^{2}}\right) }
$, $\overline{x}\left( t\right) =x_{0}-\nu \overline{d}\left( t\right) \sqrt{%
1+E^{2}}$, where $\nu $ denotes the outer unit normal to $\Omega _{1}\left(
t\right) $ at $x_{0}$. We have
\begin{equation}
B_{\overline{d}\left( t\right) }\left( \overline{x}\left( t\right) \right)
\subset \Omega _{1}\left( t\right) \cap B_{\frac{\eta _{2}R\left( t\right) }{%
2}}\left( x_{0}\right) \text{ , }dist\left( B_{\overline{d}\left( t\right)
}\left( \overline{x}\left( t\right) \right) ,\partial \Omega _{1}\left(
t\right) \right) \geq \overline{d}\left( t\right) \text{.}  \label{4.381}
\end{equation}
By Proposition \ref{Pr4.4}, (\ref{4.370}) and (\ref{4.381}) we obtain
\begin{equation}
\sigma \geq R_{0}^{-n}\int\nolimits_{\Omega _{1}\left( t\right) \cap B_{%
\frac{\eta _{2}R\left( t\right) }{2}}\left( x_{0}\right) }u_{1}^{2}\left(
x,t\right) dx\geq \frac{H^{2}R^{n}\left( t\right) }{CR_{0}^{n}}b\left(
t\right) \ \text{,}  \label{4.385}
\end{equation}
where $b\left( t\right) $ is defined by (\ref{4.377}) and $C$, $C>1$,
depends on $\lambda _{0},\Lambda _{0},E,M,F_{0},\beta $ and $R_{0}^{2}T^{-1}$
only. Since it is evident that $d_{m}\left( t\right) \leq CR_{0}$, where $C$
depends on $E\ $and $M$ only, by (\ref{4.385}) we have
\begin{equation}
d_{m}\left( t\right) \leq \frac{CR_{0}^{n+1}}{R^{n}\left( t\right) b\left(
t\right) }\frac{\sigma }{H^{2}}\ ,  \label{4.395}
\end{equation}
where $C$, $C>1$, depends on $\lambda _{0},\Lambda _{0},E,M,F_{0},\beta $
and $R_{0}^{2}T^{-1}$ only.

Otherwise, if $d_{m}\left( t\right) <\dfrac{\eta _{2}R\left( t\right) }{2}$,
let us apply Theorem \ref{Th3.6} with $r=d_{m}\left( t\right) $, $\rho =%
\dfrac{\eta _{2}R\left( t\right) }{2}$, $R=R\left( t\right) $ and by (\ref%
{4.380}) we have
\begin{equation}
\int\nolimits_{\Omega _{1}\left( t\right) \cap B_{\frac{\eta _{2}R\left(
t\right) }{2}}\left( x_{0}\right) }u_{1}^{2}\left( x,t\right) dx\leq
CH^{2}R^{n}\left( t\right) \left( \frac{R_{0}^{n}\sigma }{R^{n}\left(
t\right) H^{2}}\right) ^{\frac{1}{C\log \frac{R\left( t\right) }{d_{m}\left(
t\right) }}}\text{ ,}  \label{4.400}
\end{equation}
and $C$, $C>1$, depends on $\lambda _{0},\Lambda _{0},E,M,F_{0},\beta $ and $%
R_{0}^{2}T^{-1}$ only.

\noindent By (\ref{4.400}) and Proposition \ref{Pr4.4}, taking into account
that $\sigma <MH^{2}$, we have
\begin{equation}
d_{m}\left( t\right) \leq MR\left( t\right) \left( \frac{R_{0}^{n}\sigma }{%
MH^{2}R_{0}^{n}}\right) ^{\frac{1}{C\log \frac{R\left( t\right) }{%
d_{m}\left( t\right) }}}\text{ ,}  \label{4.410}
\end{equation}
where $C$, $C>1$, depends on $\lambda _{0},\Lambda _{0},E,M,F_{0},\beta $
and $R_{0}^{2}T^{-1}$ only. By (\ref{4.395}) and (\ref{4.410}) we have that $%
d_{m}\left( t\right) $ satisfies (\ref{4.375}).

Now we prove (\ref{4.376}). Without loss of generality, we may assume that
there exists $y_{0}\in \overline{\Omega _{1}\left( t\right) }$ such that $%
dist\left( y_{0},\overline{\Omega _{2}\left( t\right) }\right) =d\left(
t\right) $. Denoting $\delta \left( t\right) =dist\left( y_{0},\partial
\Omega _{1}\left( t\right) \right) $, let us distinguish the following three
cases:

\noindent i) $\delta \left( t\right) \leq \frac{1}{2}d\left( t\right) $,

\noindent ii) $\frac{1}{2}d\left( t\right) <\delta \left( t\right) \leq
\frac{1}{2}d_{0}$,

\noindent iii) $\delta \left( t\right) >\max \left\{ \frac{1}{2}d\left(
t\right) ,\frac{1}{2}d_{0}\right\} $,

\noindent where $d_{0}$ is the number introduced in Proposition \ref{Pr4.8}.

If case i) occurs, let $z_{0}\in \partial \Omega _{1}\left( t\right) $ be
such that $\left\vert y_{0}-z_{0}\right\vert =\delta \left( t\right) $. We
have
\begin{equation*}
d_{m}\left( t\right) \geq dist\left( z_{0},\overline{\Omega _{2}\left(
t\right) }\right) \geq dist\left( y_{0},\overline{\Omega _{2}\left( t\right)
}\right) -\left\vert y_{0}-z_{0}\right\vert =d\left( t\right) -\delta \left(
t\right) \geq \delta \left( t\right) \text{ ,}
\end{equation*}
hence $\delta \left( t\right) \leq d_{m}\left( t\right) $ and $d\left(
t\right) -\delta \left( t\right) \leq d_{m}\left( t\right) $. Therefore $%
d\left( t\right) \leq 2d_{m}\left( t\right) $ and by (\ref{4.375}) we have (%
\ref{4.376}).

If case ii) occurs then $d\left( t\right) <d_{0}$ and Proposition \ref{Pr4.8}
applies, hence by (\ref{4.365}) we have that $d\left( t\right) \leq
cd_{m}\left( t\right) $, where $c$ is an absolute constant, therefore $%
d\left( t\right) $ satisfies (\ref{4.376}).

If case iii) occurs, let us denote $R_{1}\left( t\right) =\min \left\{ \frac{%
d_{0}}{\sqrt{2E}},R\left( t\right) \right\} $ and $d_{1}\left( t\right)
=\min \left\{ \frac{d\left( t\right) }{2},\frac{\eta _{1}R_{1}\left(
t\right) }{2}\right\} $, where $\eta _{1}\in \left( 0,1\right) $ has been
introduced in Theorem \ref{Th3.4} and depends on $\lambda _{0},\Lambda
_{0},E,M,F_{0},\beta $ and $R_{0}^{2}T^{-1}$ only.

\noindent We have
\begin{equation}
B_{d_{1}\left( t\right) }\left( y_{0}\right) \subset \Omega _{1}\left(
t\right) \smallsetminus \overline{\Omega _{2}\left( t\right) }  \label{4.415}
\end{equation}
and
\begin{equation}
B_{R_{1}\left( t\right) }\left( y_{0}\right) \times \left( t-R_{1}^{2}\left(
t\right) ,t\right] \subset \Omega _{1}\left( \left( 0,t\right] \right) \text{
.}  \label{4.420}
\end{equation}
Now let us apply Therem \ref{Th3.4} with $r=d_{1}\left( t\right) $, $\rho
=\eta _{1}R\left( t\right) $, $R=R_{1}\left( t\right) $. By (\ref{4.380}), (%
\ref{4.415}) and (\ref{4.420}) we have
\begin{equation}
\int\nolimits_{B_{\eta _{1}R\left( t\right) }\left( y_{0}\right)
}u_{1}^{2}\left( x,t\right) dx\leq CH^{2}R_{1}^{n}\left( t\right) \left(
\frac{R_{0}^{n}\sigma }{R_{1}^{n}\left( t\right) H^{2}}\right) ^{\frac{1}{%
C\log \frac{R_{1}\left( t\right) }{d_{1}\left( t\right) }}}\text{ ,}
\label{4.425}
\end{equation}
where $C$, $C>1$, depends on $\lambda _{0},\Lambda _{0},E,M,F_{0},\beta $
and $R_{0}^{2}T^{-1}$ only. By (\ref{4.425}) and Proposition \ref{Pr4.4} we
have (\ref{4.376}).

Let us introduce some notation. Let $T_{0}=\min \left\{ R_{0}^{2},T\right\} $%
,
\begin{equation*}
\psi \left( \mu \right) =\exp \left( \mu ^{-1/2}\log \frac{\left\Vert
F_{1}\right\Vert _{L^{\infty }\left( 0,T\right) }}{F\left( T_{0}\mu \right) }%
+\mu ^{-n/2}\log \left( \frac{\left\Vert F_{1}\right\Vert _{L^{\infty
}\left( 0,T\right) }}{H}+2\right) \right) \ ,\ \mu \in \left( 0,1\right]
\end{equation*}
and let $\Psi _{0}=\inf\limits_{\left( 0,1\right] }\psi $, notice $\Psi
_{0}\geq 1$.

\noindent By (\ref{4.376}) and Proposition \ref{Pr4.6} we have
\begin{equation}
\sup\limits_{t\in \left[ 0,T\right] }d\left( t\right) \leq R_{0}\left(
e^{\left( \psi \left( \mu \right) \right) ^{C_{1}}}\left( \frac{\sigma }{%
MH^{2}}\right) ^{\left( \psi \left( \mu \right) \right) ^{-C_{1}}}+\mu
\right) \ ,\ \text{for every }\mu \in \left( 0,1\right] \text{,}
\label{4.430}
\end{equation}
where $C_{1}$, $C_{1}>1$, depends on $\lambda _{0},\Lambda
_{0},E,M,F_{0},\beta $ and $R_{0}^{2}T^{-1}$ only.

\noindent Now, if $\left\vert \log \left( \frac{\sigma }{MH^{2}}\right)
\right\vert >\left( 2\Psi _{0}\right) ^{C_{1}}$, then we choose
\begin{equation*}
\mu =\overline{\mu }\left( \sigma \right) :=\left( \psi ^{-1}\left(
\left\vert \log \left( \frac{\sigma }{MH^{2}}\right) \right\vert \right)
\right) ^{\frac{1}{4C_{1}}}
\end{equation*}%
at the right-hand side of (\ref{4.430}) (here $\psi ^{-1}$ is the inverse of
function $\psi $). Otherwise the left-hand side of (\ref{4.430}) can be
easily estimated from above as follows
\begin{equation*}
\sup\limits_{t\in \left[ 0,T\right] }d_{\mathcal{H}}\left( \overline{\Omega
_{1}\left( t\right) },\overline{\Omega _{2}\left( t\right) }\right) \leq
\sup\limits_{t\in \left[ 0,T\right] }\left( \text{diam}\left( \Omega
_{1}\left( t\right) \right) +\text{diam}\left( \Omega _{2}\left( t\right)
\right) \right) \leq C_{2}e^{\left( 2\Psi _{0}\right) ^{C_{1}}}\frac{\sigma
}{MH^{2}}\text{ ,}
\end{equation*}%
(here diam$\left( \Omega _{i}\left( t\right) \right) $ denotes the diameter
of $\Omega _{i}\left( t\right) $, $i=1,2$) where $C_{2}$ depends on $\lambda
_{0},\Lambda _{0},E,M,F_{0},\beta $ and $R_{0}^{2}T^{-1}$ only. Therefore,
denoting
\begin{equation*}
a=e^{-\left( 2\Psi _{0}\right) ^{C_{1}}}MH^{2}\text{ , }m_{0}=1-\left( 2\Psi
_{0}\right) ^{-\frac{C_{1}}{2}}\text{ }
\end{equation*}%
and
\begin{equation}
\phi \left( s\right) =\left\{
\begin{array}{c}
\overline{\mu }\left( s\right) +\exp \left( -m_{0}\left\vert \log \left(
\frac{s}{MH^{2}}\right) \right\vert ^{3/4}\right) \text{, if }s\in \left(
0,a\right) \text{,} \\
C_{2}a^{-1}s\text{ , if }s\in \left[ a,+\infty \right) ,%
\end{array}%
\right.  \label{4.431}
\end{equation}%
(notice that $m_{0}>0$) we have
\begin{equation}
\sup\limits_{t\in \left[ 0,T\right] }d_{\mathcal{H}}\left( \overline{\Omega
_{1}\left( t\right) },\overline{\Omega _{2}\left( t\right) }\right) \leq
R_{0}\phi \left( \sigma \right) \text{ .}  \label{4.435}
\end{equation}%
By (\ref{4.435}) and Proposition \ref{Pr4.2} we have, for every $\varepsilon
\in \left( 0,He^{-1}\right) $,
\begin{equation}
\sup\limits_{t\in \left[ 0,T\right] }d_{\mathcal{H}}\left( \overline{\Omega
_{1}\left( t\right) },\overline{\Omega _{2}\left( t\right) }\right) \leq
R_{0}\phi \left( \frac{1}{M}\left( \log \left\vert \log \frac{\varepsilon }{H%
}\right\vert \right) ^{-1/n}\right) \text{.}  \label{4.440}
\end{equation}%
Now, let us define
\begin{equation*}
\widetilde{\varepsilon }_{0}=\sup \left\{ \tau \in \left( 0,e^{-1}\right)
:\phi \left( \frac{1}{M}\left( \log \left\vert \log \tau \right\vert \right)
^{-1/n}\right) \leq \frac{d_{0}}{R_{0}}\right\} \text{ .}
\end{equation*}%
By Proposition \ref{Pr4.8} we have that, if $\varepsilon \in \left( 0,%
\widetilde{\varepsilon }_{0}H\right] $, there exist $\rho _{0}\in \left(
0,R_{0}\right] $ and $L>0$, $\frac{\rho _{0}}{R_{0}}$ and $L$ depending on $%
E $ only, such that $\partial G\left( t\right) $ is of Lipschitz class of
constants $\rho _{0}$, $L$. Therefore by Proposition \ref{Pr4.3} and (\ref%
{4.376}) we have
\begin{equation*}
d\left( t\right) \leq \frac{CR_{0}^{n+2}}{b\left( t\right) \min \left\{
t^{\left( n+1\right) /2},R_{0}^{n+1}\right\} }\left( C\left\vert \log \frac{%
\varepsilon }{H}\right\vert ^{-1/C}\right) ^{\frac{1}{\log \frac{C}{b\left(
t\right) }}}\text{ ,}
\end{equation*}%
where $C$, $C>1$, depends on $\lambda _{0},\Lambda _{0},E,M,F_{0},\beta $
and $R_{0}^{2}T^{-1}$ only.

\noindent Otherwise, if $\varepsilon >\widetilde{\varepsilon }_{0}H$, we
have trivially
\begin{equation*}
d\left( t\right) \leq C^{\prime }\frac{\varepsilon }{\widetilde{\varepsilon }%
_{0}H}\text{ ,}
\end{equation*}
where $C^{\prime }$ depends on $E$ and $M$ only, and (\ref{4.75}) follows.$%
\blacksquare $

\subsection{Proofs of Propositions \protect\ref{Pr4.2}, \protect\ref{Pr4.3},
\protect\ref{Pr4.4}, \protect\ref{Pr4.5}, \protect\ref{Pr4.6}\label{Subs4.2}}

\textbf{Proof of Proposition \ref{Pr4.2}. }Let us introduce the following
notation
\begin{equation*}
u=u_{1}-u_{2}\text{ , in }G\left( \left( 0,T\right) \right) \text{,}
\end{equation*}
\begin{equation*}
b\left( x,t\right) =-\partial _{u}\kappa \left( x,t,u_{1}\left( x,t\right)
\right) \nabla u_{2}\left( x,t\right) \text{ , in }G\left( \left( 0,T\right)
\right) \text{,}
\end{equation*}
\begin{eqnarray*}
c\left( x,t\right) &=&-\int\nolimits_{0}^{1}\partial _{u}^{2}\kappa \left(
x,t,u_{1}\left( x,t\right) +\tau u\left( x,t\right) \right) \nabla
u_{2}\left( x,t\right) \cdot \nabla u_{2}\left( x,t\right) d\tau \\
&&-\int\nolimits_{0}^{1}\text{tr}\left( \partial _{u}\kappa \left(
x,t,u_{1}\left( x,t\right) +\tau u\left( x,t\right) \right) D^{2}u_{2}\left(
x,t\right) \right) d\tau \text{ , in }G\left( \left( 0,T\right) \right)
\text{.}
\end{eqnarray*}
Denoting $\kappa _{1}\left( x,t\right) =\kappa \left( x,t,u_{1}\left(
x,t\right) \right) $, we have that the function $u$ satisfies the equation
\begin{equation}
div\left( \kappa _{1}\left( x,t\right) \nabla u\right) -\partial
_{t}u+b\left( x,t\right) \cdot \nabla u+c\left( x,t\right) u=0\text{ , in }%
G\left( \left( 0,T\right) \right)  \label{4.100}
\end{equation}
and the following Cauchy conditions
\begin{equation}
u=0\text{ , on }\Sigma \times \left( 0,T\right) \text{, }  \label{4.102}
\end{equation}
\begin{equation}
\left\Vert \kappa _{1}\nabla u\cdot \nu \right\Vert _{L^{2}\left( \Sigma
\times \left( 0,T\right) \right) }\leq \frac{\varepsilon }{R_{0}}\text{ .}
\label{4.104}
\end{equation}
Moreover
\begin{equation}
u\left( .,0\right) =0\text{ , on }G\left( 0\right) \text{ .}  \label{4.105}
\end{equation}
Now, let us recall the interpolation inequality
\begin{equation}
\left\Vert g\right\Vert _{H^{\frac{1}{2},\frac{1}{4}}\left( \Sigma \times
\left( 0,T\right) \right) }\leq C\left\Vert g\right\Vert _{L^{2}\left(
\Sigma \times \left( 0,T\right) \right) }^{\mu }\left\Vert g\right\Vert
_{C^{0,1}\left( \Sigma \times \left( 0,T\right) \right) }^{1-\mu }\text{ ,}
\label{4.106}
\end{equation}
where $\mu =\frac{1}{1+n}$ and $C$ depends on $E$ and $R_{0}^{2}T^{-1}$ only.

\noindent By (\ref{4.95}), (\ref{4.104}) and (\ref{4.106}) we have
\begin{equation}
\left\Vert \kappa _{1}\nabla u\cdot \nu \right\Vert _{H^{\frac{1}{2},\frac{1%
}{4}}\left( \Sigma \times \left( 0,T\right) \right) }\leq C\varepsilon ^{\mu
}H^{1-\mu }\text{ ,}  \label{4.110}
\end{equation}
where $C$ depends on $\lambda _{0},\Lambda _{0},E,M,F_{0},\beta $ and $%
R_{0}^{2}T^{-1}$ only.

\noindent In order to apply Theorem \ref{Th3.9} notice that by (\ref{4.95})
we have
\begin{equation*}
R_{0}\left\Vert b\right\Vert _{L^{\infty }\left( G\left( \left( 0,T\right)
\right) \right) }+R_{0}^{2}\left\Vert c\right\Vert _{L^{\infty }\left(
G\left( \left( 0,T\right) \right) \right) }\leq C\text{ ,}
\end{equation*}
where $C$ depends on $\lambda _{0},\Lambda _{0},E,M,F_{0},\beta $ and $%
R_{0}^{2}T^{-1}$ only.

\noindent Let us denote $R_{1}=\frac{\eta _{3}R_{0}}{1+\left( 1+E^{2}\right)
^{1/2}}$, $P_{1}=P_{0}-\nu \left( 1+E^{2}\right) ^{1/2}R_{1}$, where $\eta
_{3}\in \left( 0,1\right) $, defined in Theorem \ref{Th3.9}, depends on $%
\lambda _{0},\Lambda _{0},E,M,F_{0},\beta $ and $R_{0}^{2}T^{-1}$ only.
Notice $\left( B_{R_{0}}\left( P_{0}\right) \times \left( 0,T\right) \right)
\cap G\left( \left( 0,T\right) \right) \subset B_{R_{1}}\left( P_{1}\right)
\times \left( 0,T\right) $. By (\ref{3.1010}), (\ref{4.95}), (\ref{4.100}), (%
\ref{4.102}) (\ref{4.110}) and the local boundedness estimate, \cite{Li}, we
have
\begin{equation}
\left\Vert u\right\Vert _{L^{\infty }\left( B_{R_{1}/2}\left( P_{1}\right)
\times \left( 0,T\right) \right) }\leq C\varepsilon ^{\mu s_{1}}H^{1-\mu
s_{1}}\text{ , }  \label{4.115}
\end{equation}%
where $s_{1}\in \left( 0,1\right) $, defined in Theorem \ref{Th3.9}, and $C$
depends on $\lambda _{0},\Lambda _{0},E,M,F_{0},\beta $ and $R_{0}^{2}T^{-1}$
only.

Let $\tau \in \left( 0,T\right] $, let $r_{0}\in \left( 0,\frac{R_{1}}{2}%
\right] $ and denote by%
\begin{equation*}
\widetilde{\Sigma }_{r_{0}}^{i}\left( \tau \right) =\left\{ x\in \Omega
_{i}\left( \tau \right) :dist\left( x,\widetilde{\Sigma }\right) =dist\left(
x,\partial \Omega _{i}\left( \tau \right) \right) =r_{0}\right\} \text{ , }%
i=1,2\text{.}
\end{equation*}%
Notice that $\widetilde{\Sigma }_{r_{0}}^{1}\left( \tau \right) =\widetilde{%
\Sigma }_{r_{0}}^{2}\left( \tau \right) $ and that $\widetilde{\Sigma }%
_{r_{0}}^{i}\left( \tau \right) $ does not really depend on $\tau $. Let $%
V_{r_{0}}\left( \tau \right) $ be the connected component of $\left( \Omega
_{1}\left( \tau \right) \right) _{r_{0}}\cap \left( \Omega _{2}\left( \tau
\right) \right) _{r_{0}}$ whose closure contains $\widetilde{\Sigma }%
_{r_{0}}^{1}\left( \tau \right) $, we have $B_{R_{1}}\left( P_{1}\right)
\subset V_{r_{0}}\left( \tau \right) $. Now it is easy to check that if $%
r_{0}\in \left( 0,\overline{r}_{0}\right] $, where $\overline{r}_{0}=\min
\left\{ \frac{R_{1}}{2},\frac{R_{0}^{3}E}{2T},\frac{2T}{ER_{0}}\right\} $,
then for every $x\in V_{r_{0}}\left( \tau \right) $ we have
\begin{equation}
B_{r_{0}/2}\left( x\right) \times \left( \tau -\left( \frac{r_{0}}{2}\right)
^{2},\tau \right] \subset \Omega _{1}\left( \left( -\infty ,\tau \right]
\right) \cap \Omega _{2}\left( \left( -\infty ,\tau \right] \right) \text{ .}
\label{4.120}
\end{equation}%
Let $r_{0}\in \left( 0,\overline{r}_{0}\right] $ and let $\overline{x}\in
\overline{V_{r_{0}}\left( \tau \right) }$ be such that $\left\vert u\left(
\overline{x},\tau \right) \right\vert =\max\limits_{\overline{%
V_{r_{0}}\left( \tau \right) }}\left\vert u\left( .,\tau \right) \right\vert
$. Denote by $\rho _{0}=\frac{\eta _{1}}{6}r_{0}$, where $\eta _{1}\in
\left( 0,1\right) $ is defined in Theorem \ref{Th3.4} and depends on $%
\lambda _{0},\Lambda _{0},E,M,F_{0},\beta $ and $R_{0}^{2}T^{-1}$ only. Let $%
\gamma $ be an arc in $\overline{V_{r_{0}}\left( \tau \right) }$ joining $%
\overline{x}$ to $P_{1}$. Let us define $x_{i}$, $i=1,...,m$ as follows: $%
x_{1}=P_{1}$, $x_{i+1}=\gamma \left( t_{i}\right) $, where $t_{i}=\max
\left\{ t:\left\vert \gamma \left( t\right) -x_{i}\right\vert =2\rho
_{0}\right\} $ if $\left\vert x_{i}-\overline{x}\right\vert >2\rho _{0}$,
otherwise let $i=m$ and stop the process. By construction the balls $B_{\rho
_{0}}\left( x_{i}\right) $ are pairwise disjoint, $\left\vert
x_{i+1}-x_{i}\right\vert =2\rho _{0}$, for $i=1,...,m-1$, $\left\vert x_{m}-%
\overline{x}\right\vert \leq 2\rho _{0}$. Hence we have $m\leq C_{0}\left(
\frac{R_{0}}{r_{0}}\right) ^{n}$, where $C_{0}$ depends on $\lambda
_{0},\Lambda _{0},E,M,F_{0},\beta $ and $R_{0}^{2}T^{-1}$ only. By an
iterated application of the two-sphere one-cylinder inequality (\ref{3.281})
to the trivial extension $\widetilde{u}$ of $u$ (i.e. $\widetilde{u}\left(
.,t\right) =0$ if $t<0$) with $R=\frac{r_{0}}{2}$, $\rho =3\rho _{0}$, $%
r=\rho _{0}$, over the chain of balls $B_{\rho _{0}}\left( x_{i}\right) $,
and by (\ref{4.95}), (\ref{4.115}) and (\ref{4.120}), we have
\begin{equation}
\left( \rho _{0}^{-n}\int\nolimits_{B_{\rho _{0}}\left( \overline{x}\right)
}u^{2}\left( x,\tau \right) dx\right) ^{1/2}\leq C\varepsilon ^{\mu
s_{1}s_{2}^{m}}H^{1-\mu s_{1}s_{2}^{m}}\text{ ,}  \label{4.121}
\end{equation}%
where $s_{1},s_{2}\in \left( 0,1\right) $ and $C$ depend on $\lambda
_{0},\Lambda _{0},E,M,F_{0},\beta $ and $R_{0}^{2}T^{-1}$ only.

\noindent Let us recall the following interpolation inequality
\begin{equation}
\left\Vert v\right\Vert _{L^{\infty }\left( B_{\rho }\right) }\leq c\left(
\left\Vert v\right\Vert _{L^{\infty }\left( B_{\rho }\right) }+\rho
\left\Vert \nabla v\right\Vert _{L^{\infty }\left( B_{\rho }\right) }\right)
^{\frac{n}{n+2}}\left( \rho ^{-n}\int\nolimits_{B_{\rho }}v^{2}\right) ^{%
\frac{2}{n+2}}\text{ ,}  \label{4.125}
\end{equation}
where $c$ is an absolute constant.

\noindent By (\ref{4.95}), (\ref{4.115}), (\ref{4.121}) and (\ref{4.125}) we
have
\begin{equation}
\left\Vert u\left( .\tau \right) \right\Vert _{L^{\infty }\left( B_{\rho
_{0}}\left( \overline{x}\right) \right) }\leq C\varepsilon ^{\widetilde{s}%
_{1}s_{2}^{m}}H^{1-\widetilde{s}_{1}s_{2}^{m}}\text{ ,}  \label{4.130}
\end{equation}
where $\widetilde{s}_{1}=\frac{2s_{1}\mu }{n+2}$ and $C$ depend on $\lambda
_{0},\Lambda _{0},E,M,F_{0},\beta $ and $R_{0}^{2}T^{-1}$ only.

\noindent By (\ref{4.130}) we obtain
\begin{equation}
\left\Vert u\left( .\tau \right) \right\Vert _{L^{\infty }\left(
V_{r_{0}}\left( \tau \right) \right) }\leq C\varepsilon ^{\widetilde{s}%
_{1}s_{2}^{m}}H^{1-\widetilde{s}_{1}s_{2}^{m}}\text{ , for every }\tau \in %
\left[ 0,T\right] \text{,}  \label{4.135}
\end{equation}
where $\widetilde{s}_{1}$, $s_{2}$ and $C$ are the same of (\ref{4.130}).

We have, for every $\tau \in \left[ 0,T\right] $,
\begin{equation*}
\Omega _{1}\left( \tau \right) \smallsetminus G\left( \tau \right) \subset
\left( \left( \Omega _{1}\left( \tau \right) \smallsetminus \left( \Omega
_{1}\left( \tau \right) \right) _{r_{0}}\right) \smallsetminus G\left( \tau
\right) \right) \cup \omega _{r_{0}}\left( \tau \right) \ \text{,}
\end{equation*}
where
\begin{equation*}
\omega _{r_{0}}\left( \tau \right) =\left( \Omega _{1}\left( \tau \right)
\right) _{r_{0}}\smallsetminus V_{r_{0}}\left( \tau \right)
\end{equation*}
and we denote by $\Gamma _{1,r_{0}}\left( \tau \right) $ and $\Gamma
_{2,r_{0}}\left( \tau \right) $ subsets of $\partial \left( \Omega
_{1}\left( \tau \right) \right) _{r_{0}}$ and $\partial \left( \Omega
_{2}\left( \tau \right) \right) _{r_{0}}\cap \partial V_{r_{0}}\left( \tau
\right) $ respectively such that we have
\begin{equation*}
\partial \left( \left( \Omega _{1}\left( \tau \right) \right)
_{r_{0}}\smallsetminus V_{r_{0}}\left( \tau \right) \right) =\Gamma
_{1,r_{0}}\left( \tau \right) \cup \Gamma _{2,r_{0}}\left( \tau \right)
\text{ .}
\end{equation*}
We prove (\ref{4.85}) when $i=1$, the case $i=2$ being analogous.

\noindent We have, for every $\tau \in \left[ 0,T\right] $,
\begin{equation}
\int\nolimits_{\Omega _{1}\left( \tau \right) \smallsetminus G\left( \tau
\right) }u_{1}^{2}\left( x,\tau \right) dx  \label{4.140}
\end{equation}
\begin{equation*}
\leq \int\nolimits_{\left( \Omega _{1}\left( \tau \right) \smallsetminus
\left( \Omega _{1}\left( \tau \right) \right) _{r_{0}}\right) \smallsetminus
G\left( \tau \right) }u_{1}^{2}\left( x,\tau \right)
dx+\int\nolimits_{\omega _{r_{0}}\left( \tau \right) }u_{1}^{2}\left(
x,\tau \right) dx\text{ .}
\end{equation*}
By (\ref{4.95}) we have
\begin{equation}
\int\nolimits_{\left( \Omega _{1}\left( \tau \right) \smallsetminus \left(
\Omega _{1}\left( \tau \right) \right) _{r_{0}}\right) \smallsetminus
G\left( \tau \right) }u_{1}^{2}\left( x,\tau \right) dx\leq
CH^{2}R_{0}^{n-1}r_{0}\text{ ,}  \label{4.145}
\end{equation}
where $C$ depends on $E$ and $M$ only.

\noindent By the divergence theorem we have
\begin{equation}
0=\int\nolimits_{\omega _{r_{0}}\left( \tau \right) }\left( div\left(
\kappa _{1}\left( x,t\right) \nabla u_{1}\right) -\partial _{t}u_{1}\right)
u_{1}dxdt  \label{4.147}
\end{equation}%
\begin{eqnarray*}
&=&\int\nolimits_{0}^{\tau }dt\int\nolimits_{\partial \omega
_{r_{0}}\left( t\right) }\left( \kappa _{1}\left( x,t\right) \nabla
u_{1}\cdot \nu \right) u_{1}ds-\int\nolimits_{\omega _{r_{0}}\left( \left(
0,\tau \right) \right) }\kappa _{1}\left( x,t\right) \nabla u_{1}\cdot
\nabla u_{1}dxdt \\
&&-\frac{1}{2}\int\nolimits_{\omega _{r_{0}}\left( \tau \right)
}u_{1}^{2}\left( x,\tau \right) dx-\frac{1}{2}\int\nolimits_{\left(
\partial \omega _{r_{0}}\right) \left( \left( 0,\tau \right) \right)
}u_{1}^{2}\left( N\cdot e_{n+1}\right) dS\text{ ,}
\end{eqnarray*}%
where, as usual, $\left( \partial \omega _{r_{0}}\right) \left( \left(
0,\tau \right) \right) =\bigcup\limits_{t\in \left( 0,\tau \right)
}\partial \omega _{r_{0}}\left( t\right) \times \left\{ t\right\} $, $N$ is
the exterior unit normal to $\omega _{r_{0}}\left( \left( 0,\tau \right)
\right) $, $ds$ and $dS$ are the $\left( n-1\right) $-dimensional and the $n$%
-dimensional surface element respectively.

\noindent Therefore
\begin{equation}
\frac{1}{2}\int\nolimits_{\omega _{r_{0}}\left( \tau \right)
}u_{1}^{2}\left( x,\tau \right) dx  \label{4.150}
\end{equation}
\begin{equation*}
\leq \int\nolimits_{0}^{\tau }dt\int\nolimits_{\partial \omega
_{r_{0}}\left( t\right) }\left( \kappa _{1}\left( x,t\right) \nabla
u_{1}\cdot \nu \right) u_{1}ds+\frac{1}{2}\int\nolimits_{\left( \partial
\omega _{r_{0}}\right) \left( \left( 0,\tau \right) \right) }u_{1}^{2}\left(
N\cdot e_{n+1}\right) dS\text{.}
\end{equation*}
By (\ref{4.150}) and (\ref{4.95}) we have
\begin{equation}
\int\nolimits_{\omega _{r_{0}}\left( \tau \right) }u_{1}^{2}\left( x,\tau
\right) dx\leq \frac{CH}{R_{0}}\int\nolimits_{0}^{\tau
}dt\int\nolimits_{\partial \omega _{r_{0}}\left( t\right) }\left\vert
u_{1}\left( x,t\right) \right\vert ds  \label{4.155}
\end{equation}
\begin{equation*}
=\frac{CH}{R_{0}}\left( \int\nolimits_{0}^{\tau }dt\int\nolimits_{\Gamma
_{1,r_{0}}\left( t\right) }\left\vert u_{1}\left( x,t\right) \right\vert
ds+\int\nolimits_{0}^{\tau }dt\int\nolimits_{\Gamma _{2,r_{0}}\left(
t\right) }\left\vert u_{1}\left( x,t\right) \right\vert ds\right) \text{ ,}
\end{equation*}
where $C$ depends on $\lambda _{0}$ and $E$ only.

Now, let us notice that, if $t\in \left[ 0,\tau \right] $ and $r_{0}\in
\left( 0,\overline{r}_{0}\right] $ then $dist\left( x,\widetilde{\Sigma }%
\right) \geq \frac{R_{1}}{2}$ for every $x\in \Gamma _{1,r_{0}}\left(
t\right) $. Therefore if $t\in \left[ 0,\tau \right] $ and $x\in \Gamma
_{1,r_{0}}\left( t\right) $ then there exists $y\in \partial \Omega
_{1}\left( t\right) \smallsetminus \widetilde{\Sigma }$ such that $%
\left\vert y-x\right\vert =dist\left( x,\partial \Omega _{1}\left( t\right)
\right) =r_{0}$. Since \textit{supp}$f\left( .,t\right) \subset \widetilde{%
\Sigma }$ we have $u_{1}\left( y,t\right) =0$, hence by (\ref{4.95}) we get,
for every $x\in \Gamma _{1,r_{0}}\left( t\right) $, $t\in \left[ 0,\tau %
\right] $,
\begin{equation}
\left\vert u_{1}\left( x,t\right) \right\vert =\left\vert u_{1}\left(
x,t\right) -u_{1}\left( y,t\right) \right\vert \leq \left\Vert \nabla
u_{1}\right\Vert _{L^{\infty }\left( \Omega _{1}\left( t\right) \right)
}\left\vert x-y\right\vert \leq H\frac{r_{0}}{R_{0}}\ \text{.}  \label{4.160}
\end{equation}
Analogously, if $t\in \left[ 0,\tau \right] $ and $x\in \Gamma
_{2,r_{0}}\left( t\right) $ then there exists $y\in \partial \Omega
_{2}\left( t\right) \smallsetminus \widetilde{\Sigma }$ such that $%
\left\vert y-x\right\vert =dist\left( x,\partial \Omega _{2}\left( t\right)
\right) =r_{0}$. Since $u_{2}\left( y,t\right) =0$, by (\ref{4.95}) and (\ref%
{4.135}) we obtain, for every $x\in \Gamma _{2,r_{0}}\left( t\right) $, $%
t\in \left[ 0,\tau \right] $,
\begin{equation}
\left\vert u_{1}\left( x,t\right) \right\vert \leq H\frac{r_{0}}{R_{0}}%
+\left\vert u\left( x,t\right) \right\vert \leq CH\left( \frac{r_{0}}{R_{0}}%
+\left( \frac{\varepsilon }{H}\right) ^{\widetilde{s}_{1}s_{2}^{m}}\right)
\text{ , }  \label{4.165}
\end{equation}
where $C$ depend on $\lambda _{0},\Lambda _{0},E,M,F_{0},\beta $ and $%
R_{0}^{2}T^{-1}$ only and $\widetilde{s}_{1}$, $s_{2}$ are the same of (\ref%
{4.135}).

By (\ref{4.140}), (\ref{4.145}), (\ref{4.155}), (\ref{4.160}) and (\ref%
{4.165}) we have, for every $\tau \in \left[ 0,T\right] $ and $r_{0}\in
\left( 0,\overline{r}_{0}\right] $,
\begin{equation}
\int\nolimits_{\Omega _{1}\left( \tau \right) \smallsetminus G\left( \tau
\right) }u_{1}^{2}\left( x,\tau \right) dx\leq CR_{0}^{n}H^{2}\left( \frac{%
r_{0}}{R_{0}}+\left( \frac{\varepsilon }{H}\right) ^{\widetilde{s}%
_{1}s_{2}^{m}}\right) \text{ ,}  \label{4.170}
\end{equation}
where $C$ depend on $\lambda _{0},\Lambda _{0},E,M,F_{0},\beta $ and $%
R_{0}^{2}T^{-1}$ only.

\noindent Denote
\begin{equation*}
\widetilde{\varepsilon }=\left( \frac{\varepsilon }{H}\right) ^{\widetilde{s}%
_{1}}\text{ ,}
\end{equation*}
\begin{equation*}
\widetilde{\varepsilon }_{1}=\min \left\{ e^{-1},\exp \left( -\exp \left(
2C_{0}\left( \frac{R_{0}}{\overline{r}_{0}}\right) ^{n}\left\vert \log
s_{2}\right\vert \right) \right) \right\} \text{.}
\end{equation*}
If $\widetilde{\varepsilon }<\widetilde{\varepsilon }_{1}$ then we choose
\begin{equation*}
r_{0}=R_{0}\left( \frac{2C_{0}\left\vert \log s_{2}\right\vert }{\log
\left\vert \log \widetilde{\varepsilon }\right\vert }\right) ^{1/n}
\end{equation*}
at the right-hand side of (\ref{4.170}) and we have
\begin{equation}
\int\nolimits_{\Omega _{1}\left( \tau \right) \smallsetminus G\left( \tau
\right) }u_{1}^{2}\left( x,\tau \right) dx\leq CR_{0}^{n}H^{2}\left( \log
\left\vert \log \left( \frac{\varepsilon }{H}\right) ^{\widetilde{s}%
_{1}}\right\vert \right) ^{-\frac{1}{n}}\text{ ,}  \label{4.175}
\end{equation}
where $C$ depends on $\lambda _{0},\Lambda _{0},E,M,F_{0},\beta $ and $%
R_{0}^{2}T^{-1}$ only.

\noindent Otherwise, if $\widetilde{\varepsilon }\geq \widetilde{\varepsilon
}_{1}$ then we have trivially
\begin{equation}
\int\nolimits_{\Omega _{1}\left( \tau \right) \smallsetminus G\left( \tau
\right) }u_{1}^{2}\left( x,\tau \right) dx\leq MR_{0}^{n}H^{2}\leq \frac{%
MR_{0}^{n}H^{2}}{\widetilde{\varepsilon }_{1}}\left( \frac{\varepsilon }{H}%
\right) ^{\widetilde{s}_{1}}\text{ .}  \label{4.180}
\end{equation}
By (\ref{4.175}) and (\ref{4.180}) we have (\ref{4.85}).$\blacksquare $

\bigskip

\textbf{Proof of Proposition \ref{Pr4.3}. }Let us prove (\ref{4.85}) and (%
\ref{4.185}) when $i=1$, the case $i=2$ being analogous. By the divergence
theorem we have, for every $\tau \in \left[ 0,T\right] $
\begin{equation}
\frac{1}{2}\int\nolimits_{\Omega _{1}\left( \tau \right) \smallsetminus
G\left( \tau \right) }u_{1}^{2}\left( x,\tau \right) dx  \label{4.190}
\end{equation}%
\begin{equation*}
\leq \int\nolimits_{0}^{\tau }dt\int\nolimits_{\partial \left( \Omega
_{1}\left( t\right) \smallsetminus G\left( t\right) \right) }\left( \kappa
_{1}\left( x,t\right) \nabla u_{1}\cdot \nu \right) u_{1}ds+\frac{1}{2}%
\int\nolimits_{\left( \partial \left( \Omega _{1}\smallsetminus G\right)
\right) \left( \left( 0,\tau \right) \right) }u_{1}^{2}\left( N\cdot
e_{n+1}\right) dS\text{ .}
\end{equation*}%
Since, for every $t\in \left[ 0,T\right] $,
\begin{equation*}
\partial \left( \Omega _{1}\left( t\right) \smallsetminus G\left( t\right)
\right) \subset \left( \partial \Omega _{1}\left( t\right) \smallsetminus
\widetilde{\Sigma }\right) \cup \left( \left( \partial \Omega _{2}\left(
t\right) \smallsetminus \widetilde{\Sigma }\right) \cap \partial G\left(
t\right) \right)
\end{equation*}%
and since $u_{i}\left( .,t\right) =0$ on $\partial \Omega _{i}\left(
t\right) \smallsetminus \widetilde{\Sigma }$, , for every $t\in \left[ 0,T%
\right] $ and $i=1,2$, by (\ref{4.95}) and (\ref{4.190}) we have, for every $%
\tau \in \left[ 0,T\right] $ ,
\begin{equation}
\int\nolimits_{\Omega _{1}\left( \tau \right) \smallsetminus G\left( \tau
\right) }u_{1}^{2}\left( x,\tau \right) dx\leq CHR_{0}^{n-2}T\left\Vert
u\right\Vert _{L^{\infty }\left( \left( \partial G\right) \left( \left[
0,\tau \right] \right) \right) }\text{ , }  \label{4.195}
\end{equation}%
where $u=u_{1}-u_{2}$ and $C$ depends on $\lambda _{0}$ and $E$ only.

\noindent Now it is simple to check that for any $r\in \left( 0,R_{0}\right]
$, any $t_{0}\in \left[ 0,\tau \right] $ and any $x\in \partial G\left(
t_{0}\right) $ such that $dist\left( x,\partial G\left( t_{0}\right) \right)
\geq r$ we have
\begin{equation*}
dist\left( x,\partial \Omega _{i}\left( t\right) \right) \geq \frac{r}{2}%
\text{ , for every }t\in \left( t_{0}-\frac{r^{2}}{2E+1},t_{0}\right] \text{%
, }i=1,2\text{,}
\end{equation*}%
hence, since we have $\partial G\left( t\right) \subset \partial \Omega
_{1}\left( t\right) \cup \partial \Omega _{2}\left( t\right) $, we obtain
\begin{equation*}
dist\left( x,\partial G\left( t\right) \right) \geq \frac{r}{2}\text{ , for
every }t\in \left( t_{0}-\frac{r^{2}}{2E+1},t_{0}\right] \text{, }i=1,2\text{%
,}
\end{equation*}%
therefore
\begin{equation}
W_{r}\left( x,t_{0}\right) :=B_{\frac{r}{2}}\left( x\right) \times \left(
t_{0}-\frac{r^{2}}{2E+1},t_{0}\right] \subset G\left( \left( -\infty ,t_{0}%
\right] \right) \ \text{. }  \label{4.196}
\end{equation}%
Let $t_{0}\in \left[ 0,\tau \right] $ and $z\in \partial G\left(
t_{0}\right) $ be fixed. Since $\partial G\left( t_{0}\right) $ is of
Lipschitz class with constants $\rho _{0}$, $L$ there exists $\zeta \in
\mathbb{R}^{n}$, $\left\vert \zeta \right\vert =1$, such that $\mathcal{C}%
\left( z,\zeta ,\alpha ,R\right) \subset G\left( t_{0}\right) $, where $%
\alpha =$arctan$\frac{1}{L}$ (recall that $\mathcal{C}\left( z,\zeta ,\alpha
,R\right) $ is defined by (\ref{subs3.2.1})). Let us denote $\overline{\mu }=%
\frac{\rho _{0}}{1+\sin \alpha }$ and, for any $\mu \in \left( 0,\overline{%
\mu }\right) $, $z_{\mu }=z+\mu \zeta $. We have
\begin{equation}
dist\left( z_{\mu },\partial G\left( t_{0}\right) \right) \geq \mu \sin
\alpha \text{ , for any }\mu \in \left( 0,\overline{\mu }\right) \text{ ,}
\label{4.197}
\end{equation}%
Now, denoting with $\overline{\delta }=\frac{\sin \alpha }{\sqrt{2E+1}}$, $%
\overline{\rho }=\overline{\mu }\cos \frac{\alpha }{2}$, by (\ref{4.196})
and (\ref{4.197}) we have
\begin{equation}
\mathcal{S}\left( \left( z,t_{0}\right) ,\zeta ,\frac{\alpha }{2},\overline{%
\delta },\overline{\rho }\right) \subset \bigcup\limits_{\mu \in \left( 0,%
\overline{\mu }\right) }W_{\mu \sin \alpha }\left( z_{\mu },t_{0}\right)
\subset G\left( \left( -\infty ,t_{0}\right] \right)  \label{4.200}
\end{equation}%
(recall that $\mathcal{S}\left( \left( z,t_{0}\right) ,\zeta ,\frac{\alpha }{%
2},\overline{\delta },\overline{\rho }\right) $ is defined by (\ref%
{subs3.2.3})).

Let us denote by $\widetilde{u}$ the trivial extension of $u$ (i.e. $%
\widetilde{u}\left( .,t\right) =0$ if $t<0$), by (\ref{4.95}) we have
\begin{equation}
\left\Vert \widetilde{u}\right\Vert _{L^{\infty }\left( \mathcal{S}\left(
\left( z,t_{0}\right) ,\zeta ,\frac{\alpha }{2},\overline{\delta },\overline{%
\rho }\right) \right) }+\overline{\rho }^{\beta }\left[ \widetilde{u}\left(
.,t_{0}\right) \right] _{\beta ,\mathcal{C}\left( z,\zeta ,\alpha ,\overline{%
\rho }\right) }\leq H\text{ .}  \label{4.201}
\end{equation}
Denote by
\begin{equation*}
\alpha _{1}=\arcsin \left( \min \left\{ \sin \frac{\alpha }{2},\overline{%
\delta }\left( 1-\sin \frac{\alpha }{2}\right) \right\} \right) \text{ ,}
\end{equation*}
\begin{equation*}
\mu _{1}=\frac{\overline{\rho }}{1+\sin \alpha _{1}}\text{ , }w_{1}=z+\mu
_{1}\zeta \text{ , }\rho _{1}=\frac{1}{4}\mu _{1}\eta _{1}\sin \alpha _{1}%
\text{ ,}
\end{equation*}
where $\eta _{1}\in \left( 0,1\right) $ is defined in Theorem \ref{Th3.4}
and depends on $\lambda _{0},\Lambda _{0},E,M,F_{0},\beta $ and $%
R_{0}^{2}T^{-1}$ only.

\noindent We have
\begin{equation}
dist\left( w_{1},\partial G\left( t_{0}\right) \right) \geq \min \left\{
\rho _{0}-\left\vert w_{1}-z\right\vert ,\left\vert w_{1}-z\right\vert \sin
\alpha \right\} =\rho _{0}\widetilde{\eta }\text{ ,}  \label{4.205}
\end{equation}%
where
\begin{equation*}
\widetilde{\eta }=\min \left\{ 1-\frac{\cos \frac{\alpha }{2}}{1+\sin \alpha
}\frac{1}{1+\sin \alpha _{1}},\frac{\sin \alpha }{1+\sin \alpha _{1}}\frac{%
\cos \frac{\alpha }{2}}{1+\sin \alpha }\right\} \text{.}
\end{equation*}%
Now $\left( G\left( t_{0}\right) \right) _{\frac{\rho _{0}\widetilde{\eta }}{%
2}}$ is connected and, by (\ref{4.205}) $w_{1}\in \left( G\left(
t_{0}\right) \right) _{\frac{\rho _{0}\widetilde{\eta }}{2}}$. Therefore
arguing as in the proof of Proposition \ref{Pr4.2}, we get, by iterated
application of the two-sphere one-cylinder inequality and by (\ref{4.115}),
\begin{equation}
\left( \rho _{1}^{-n}\int\nolimits_{B_{\rho _{1}}\left( \overline{x}\right)
}\widetilde{u}^{2}\left( x,t_{0}\right) dx\right) ^{1/2}\leq C\varepsilon
^{s_{3}}H^{1-s_{3}}\text{ ,}  \label{4.210}
\end{equation}%
where $s_{3}$, $s_{3}\in \left( 0,1\right) $, and $C$ depends on $\lambda
_{0},\Lambda _{0},E,M,F_{0},\beta $ and $R_{0}^{2}T^{-1}$ only.

\noindent By (\ref{4.201}), (\ref{4.210}) and (\ref{3.840}) we obtain
\begin{equation}
\left\vert u\left( z,t_{0}\right) \right\vert \leq CH\left\vert \log \left(
\frac{\varepsilon }{H}\right) \right\vert ^{\frac{-1}{C}}\text{ ,}
\label{4.215}
\end{equation}
where $C$, $C>1$, depends on $\lambda _{0},\Lambda _{0},E,M,F_{0},\beta $
and $R_{0}^{2}T^{-1}$ only.

\noindent Finally (\ref{4.195}) and (\ref{4.215}) give (\ref{4.185}).$%
\blacksquare $

\bigskip

\textbf{Proof of Proposition \ref{Pr4.4}. }Let $\tau \in \left( 0,T\right] $
and $P\in \partial \Omega \left( \tau \right) $. There exists a rigid
tranformation of space coordinates under which we have $\left( P,\tau
\right) =\left( 0,\tau \right) $ and

\begin{equation*}
\Omega \left( \left( -\infty ,+\infty \right) \right) \cap \left(
B_{R_{0}}\left( 0\right) \times \left( \tau -R_{0}^{2},\tau
+R_{0}^{2}\right) \right)
\end{equation*}
\begin{equation*}
=\left\{ \left( x,t\right) \in B_{R_{0}}\left( 0\right) \times \left( \tau
-R_{0}^{2},\tau +R_{0}^{2}\right) :x_{n}>\varphi \left( x^{\prime },t\right)
\right\} \text{ ,}
\end{equation*}
where $\varphi \in C^{2,\beta }\left( B_{R_{0}}^{\prime }\left( 0\right)
\times \left( \tau -R_{0}^{2},\tau +R_{0}^{2}\right) \right) $ satisfying $%
\varphi \left( 0,\tau \right) =0$ and $\left\Vert \varphi \right\Vert
_{C^{2,\beta }\left( B_{R_{0}}^{\prime }\left( 0\right) \times \left( \tau
-R_{0}^{2},\tau +R_{0}^{2}\right) \right) }\leq ER_{0}$.

Let us introduce the following notations
\begin{equation*}
\rho _{0}=\frac{R_{0}}{64\sqrt{2+E^{2}}}\text{ , }y_{0}=e_{n}\frac{R_{0}}{4}%
\text{ ,}
\end{equation*}
\begin{equation*}
\Gamma =\left\{ \left( x^{\prime },\varphi \left( x^{\prime },\tau \right)
\right) :x^{\prime }\in B_{R_{0}}^{\prime }\left( 0\right) \right\} \text{ .}
\end{equation*}
We have
\begin{equation}
B_{2\rho _{0}}\left( y_{0}\right) \subset \left( \Omega \left( \tau \right)
\right) _{\rho _{0}}\text{ .}  \label{4.230}
\end{equation}
For $z\in \left( \Omega \left( \tau \right) \right) _{2\overline{\rho }}$
denote
\begin{equation}
\varepsilon _{0}=\left( \overline{\rho }^{-n}\int\nolimits_{B_{\overline{%
\rho }}\left( z\right) }u^{2}\left( x,\tau \right) dx\right) ^{1/2}\text{ }
\label{4.235}
\end{equation}
and let $\rho _{1}=\min \left\{ \overline{\rho },\rho _{0}\right\} $.

Let $y\in B_{2\rho _{0}}\left( y_{0}\right) $ and let $\gamma $ be an arc in
$\left( \Omega \left( \tau \right) \right) _{\rho _{1}}$ joining $z$ to $y$.
Since $dist\left( \gamma ,\partial \Omega \left( \tau \right) \right) \geq
\rho _{1}$, denoting $\sigma _{0}\left( \tau \right) =\min \left\{ \sqrt{%
\tau },\frac{\rho _{1}}{\sqrt{2E+4}}\right\} $, we have
\begin{equation}
B_{\sigma _{0}\left( \tau \right) }\left( x\right) \subset \left( \tau
-\sigma _{0}^{2}\left( \tau \right) ,\tau \right] \subset \Omega \left(
\left( 0,\tau \right] \right) \text{ , for every }x\in \gamma \text{.}
\label{4.240}
\end{equation}%
Let us denote $\sigma _{1}\left( \tau \right) =\eta _{1}\frac{\sigma
_{0}\left( \tau \right) }{6}$, where $\eta _{1}\in \left( 0,1\right) $ is
defined in Theorem \ref{Th3.4} and depends on $\lambda _{0},\Lambda
_{0},E,M,F_{0},\beta $ and $R_{0}^{2}T^{-1}$ only. Arguing as in Proposition %
\ref{Pr4.2} we define $x_{i}\in \gamma $, $i=1,...,m$, such that $x_{1}=z$, $%
\left\vert x_{i+1}-x_{i}\right\vert =2\sigma _{1}\left( \tau \right) $, $%
i=1,...,m$, $\left\vert x_{m}-y\right\vert \leq 2\sigma _{1}\left( \tau
\right) $. Hence $m\leq m\left( \tau \right) :=c_{0}M\left( \frac{R_{0}}{%
2\sigma _{1}\left( \tau \right) }\right) ^{n}$, where $c_{0}$ depends on $n$
only. By an iterated application of inequality (\ref{3.281}) with $r=\sigma
_{1}\left( \tau \right) $, $\rho =3\sigma _{1}\left( \tau \right) $, $R=%
\frac{1}{2}\sigma _{0}\left( \tau \right) $ over the chain of balls $%
B_{\sigma _{1}\left( \tau \right) }\left( x_{i}\right) $ and by (\ref{4.95}%
), (\ref{4.235}) and (\ref{4.240}) we obtain
\begin{equation}
\left( \sigma _{1}^{-n}\left( \tau \right) \int\nolimits_{B_{\sigma
_{1}\left( \tau \right) }\left( y\right) }u^{2}\left( x,\tau \right)
dx\right) ^{1/2}\leq C\varepsilon _{0}^{s_{2}^{m\left( \tau \right)
}}H^{1-s_{2}^{m\left( \tau \right) }}\text{ ,}  \label{4.245}
\end{equation}%
where $s_{2}$, $s_{2}\in \left( 0,1\right) $, and $C$, $C>1$, depend on $%
\lambda _{0},\Lambda _{0},E,M,F_{0},\beta $ and $R_{0}^{2}T^{-1}$ only. By
the interpolation inequality (\ref{4.125}) and by (\ref{4.245}) we get
\begin{equation}
\left\Vert u\left( .,\tau \right) \right\Vert _{L^{\infty }\left( B_{\sigma
_{1}\left( \tau \right) }\left( y\right) \right) }\leq C_{1}\varepsilon
_{0}^{\mu \left( \tau \right) }H^{1-\mu \left( \tau \right) }\text{ ,}
\label{4.250}
\end{equation}%
where $\mu \left( \tau \right) =\frac{2}{n+2}s_{2}^{m\left( \tau \right) }$
and $C_{1}$, $C_{1}>1$, depends on $\lambda _{0},\Lambda
_{0},E,M,F_{0},\beta $ and $R_{0}^{2}T^{-1}$ only. Hence
\begin{equation}
\left\Vert u\left( .,\tau \right) \right\Vert _{L^{\infty }\left( B_{2\rho
_{0}}\left( y_{0}\right) \right) }\leq C_{1}\varepsilon _{0}^{\mu \left(
\tau \right) }H^{1-\mu \left( \tau \right) }\text{ .}  \label{4.255}
\end{equation}

Now, in order to apply Proposition \ref{Pr3.8} let us introduce the
following notation
\begin{equation*}
\alpha =\arcsin \frac{1}{16\sqrt{2+E^{2}}}\text{ , }\delta _{\tau }=\frac{%
\sqrt{\tau }}{R_{0}}\text{ , }R_{1}=\frac{R_{0}}{4}\cos \alpha \text{.}
\end{equation*}
It is simple to check that
\begin{equation}
\mathcal{S}\left( \left( 0,\tau \right) ,e_{n},\alpha ,\delta _{\tau
},R_{1}\right) \subset \Omega \left( \left( 0,\tau \right] \right) \text{ }
\label{4.260}
\end{equation}
and
\begin{equation}
B_{\rho _{2}}\left( w_{1}\right) \subset \mathcal{C}\left( 0,e_{n},\alpha
,R_{1}\right) \text{ ,}  \label{4.265}
\end{equation}
where
\begin{equation*}
w_{1}=\frac{R_{1}}{1+\sin \alpha _{1}\left( \tau \right) }e_{n}\text{ ,}
\end{equation*}
\begin{equation*}
\rho _{2}=\frac{1}{4}\eta _{1}\frac{R_{1}\sin \alpha _{1}\left( \tau \right)
}{1+\sin \alpha _{1}\left( \tau \right) }\text{ ,}
\end{equation*}
\begin{equation*}
\alpha _{1}\left( \tau \right) =\arcsin \left( \min \left\{ \sin \alpha
,\delta _{\tau }\left( 1-\sin \alpha \right) \right\} \right) \text{ .}
\end{equation*}
(recall that $\mathcal{C}\left( 0,e_{n},\alpha ,R_{1}\right) $ and $\mathcal{%
S}\left( \left( 0,\tau \right) ,e_{n},\alpha ,\delta _{\tau },R_{1}\right) $
are defined by (\ref{subs3.2.1}) and (\ref{subs3.2.3}) respectively).

\noindent By (\ref{3.840}), (\ref{4.95}), (\ref{4.255}), (\ref{4.260}) and (%
\ref{4.265}) we have, when $C_{1}\left( \frac{\varepsilon _{0}}{H}\right)
^{\mu \left( \tau \right) }<1$,
\begin{equation}
\left\vert u\left( 0,\tau \right) \right\vert \leq C_{2}\left( \mu \left(
\tau \right) \right) ^{-B\left( \tau \right) }H\left\vert \log \left( \frac{%
\varepsilon _{0}}{eH}\right) \right\vert ^{-B\left( \tau \right) }\text{ ,}
\label{4.270}
\end{equation}
where
\begin{equation}
B\left( \tau \right) =\frac{1}{C_{2}}\left\vert \log \frac{1-\frac{1}{4}\eta
_{1}\sin \alpha _{1}\left( \tau \right) }{1+\frac{1}{4}\eta _{1}\sin \alpha
_{1}\left( \tau \right) }\right\vert \text{ ,}  \label{4.275}
\end{equation}
and $C_{2}$, $C_{2}>1$, depends on $\lambda _{0},\Lambda
_{0},E,M,F_{0},\beta $ and $R_{0}^{2}T^{-1}$ only. Otherwise, if $%
C_{1}\left( \frac{\varepsilon _{0}}{H}\right) ^{\mu \left( \tau \right)
}\geq 1$ then we have trivially
\begin{equation}
\left\vert u\left( 0,\tau \right) \right\vert \leq H\leq C_{1}H\left( \frac{%
\varepsilon _{0}}{H}\right) ^{\mu \left( \tau \right) }\text{ ,}
\label{4.280}
\end{equation}
By (\ref{4.270}) and (\ref{4.280}) we have
\begin{equation}
\left\Vert u\left( .,\tau \right) \right\Vert _{L^{\infty }\left( \partial
\Omega \left( \tau \right) \right) }\leq C_{2}H\left( \frac{\log \left(
eC_{1}\right) }{\mu \left( \tau \right) }\right) ^{B\left( \tau \right)
}\left\vert \log \left( \frac{\varepsilon _{0}}{eH}\right) \right\vert
^{-B\left( \tau \right) }\text{ ,}  \label{4.290}
\end{equation}
By (\ref{4.55}) and (\ref{4.290}) we get
\begin{equation}
F_{1}\left( \tau \right) \leq He^{C\frac{R_{0}^{n-1}}{\min \left\{ \tau
^{\left( n-1\right) /2},\overline{\rho }^{n-1}\right\} }}\left\vert \log
\left( \frac{\varepsilon _{0}}{eH}\right) \right\vert ^{-\frac{\min \left\{
\tau ^{1/2},\overline{\rho }\right\} }{CR_{0}}}\text{ ,}  \label{4.300}
\end{equation}
$C$, $C>1$, depends on $\lambda _{0},\Lambda _{0},E,M,F_{0},\beta $ and $%
R_{0}^{2}T^{-1}$ only. By (\ref{4.300}) we obtain (\ref{4.220}).$%
\blacksquare $

\bigskip

\textbf{Proof of Proposition \ref{Pr4.5}. }The proof of such a Proposition
is the same, with obvious changes, of Corollary 2.3 in \cite{AlBRVe1}.

\bigskip

\textbf{Proof of Proposition \ref{Pr4.6}. }Let $t,\tau \in \left[ 0,T\right]
$ satisfy $\left\vert t-\tau \right\vert \leq R_{0}^{2}\min \left\{ \frac{1}{%
2E},1\right\} $. By (\ref{4.25}) we have that $\Omega _{1}\left( t\right) $
and $\Omega _{1}\left( \tau \right) $, $\Omega _{2}\left( t\right) $ and $%
\Omega _{2}\left( \tau \right) $ are relative graphs and
\begin{equation}
\gamma \left( \Omega _{i}\left( t\right) ,\Omega _{i}\left( \tau \right)
\right) \leq E\frac{\left\vert t-\tau \right\vert }{R_{0}}\text{ , }i=1,2
\label{4.315}
\end{equation}
(recall that $\gamma \left( .,.\right) $ is defined by (\ref{4.304})).

\noindent By Proposition \ref{Pr4.5} and (\ref{4.315}) there exists a $%
C^{1,\beta }$ diffeomorphisms $\Phi _{i}:\mathbb{R}^{n}\rightarrow \mathbb{R}%
^{n}$ such that $\Phi _{i}\left( \Omega _{i}\left( t\right) \right) =\Omega
_{i}\left( \tau \right) $, $i=1,2$, and
\begin{equation}
\left\Vert \Phi _{i}-Id\right\Vert _{L^{\infty }\left( \mathbb{R}^{n}\right)
}\text{ , }\left\Vert \Phi _{i}^{-1}-Id\right\Vert _{L^{\infty }\left(
\mathbb{R}^{n}\right) }\leq CE\frac{\left\vert t-\tau \right\vert }{R_{0}}%
\text{, }i=1,2\text{,}  \label{4.320}
\end{equation}%
where $C$ depends on $E,M$ and $\beta $ and only.

\noindent Now let $\overline{x}\in \overline{\Omega _{1}\left( t\right) }$
be such that
\begin{equation*}
dist\left( \overline{x},\Omega _{2}\left( \tau \right) \right)
=\sup\limits_{x\in \Omega _{1}\left( t\right) }dist\left( x,\Omega
_{2}\left( \tau \right) \right) \text{ .}
\end{equation*}%
By the triangular inequality and by (\ref{4.320}) we have
\begin{equation*}
dist\left( \overline{x},\Omega _{2}\left( \tau \right) \right) \leq
dist\left( \Phi _{1}\left( \overline{x}\right) ,\Omega _{2}\left( \tau
\right) \right) +\left\vert \Phi _{1}\left( \overline{x}\right) -\overline{x}%
\right\vert
\end{equation*}%
\begin{equation*}
\leq \sup\limits_{x\in \Omega _{1}\left( \tau \right) }dist\left( x,\Omega
_{2}\left( \tau \right) \right) +CE\frac{\left\vert t-\tau \right\vert }{%
R_{0}}\text{ ,}
\end{equation*}%
where $C$ depends on $E,M$ and $\beta $ and only. Hence
\begin{equation}
\sup\limits_{x\in \Omega _{1}\left( t\right) }dist\left( x,\Omega _{2}\left(
\tau \right) \right) \leq \sup\limits_{x\in \Omega _{1}\left( \tau \right)
}dist\left( x,\Omega _{2}\left( \tau \right) \right) +CE\frac{\left\vert
t-\tau \right\vert }{R_{0}}\text{ ,}  \label{4.325}
\end{equation}%
where $C$ is the same of above.

\noindent Now we will prove that there exists a constant $C$ depending on $%
E,M$ and $\beta $ and only such that
\begin{equation}
\sup\limits_{x\in \Omega _{1}\left( t\right) }dist\left( x,\Omega _{2}\left(
t\right) \right) \leq \sup\limits_{x\in \Omega _{1}\left( t\right)
}dist\left( x,\Omega _{2}\left( \tau \right) \right) +CE\frac{\left\vert
t-\tau \right\vert }{R_{0}}\text{ .}  \label{4.330}
\end{equation}
Let $x\in \Omega _{1}\left( t\right) $, we distinguish two cases

\noindent a) $x\notin \Omega _{2}\left( \tau \right) $,

\noindent b) $x\in \Omega _{2}\left( \tau \right) $.

\noindent If case a) occurs then there exists $z\in \partial \Omega
_{2}\left( \tau \right) $ such that
\begin{equation}
dist\left( x,\Omega _{2}\left( \tau \right) \right) =dist\left( x,\partial
\Omega _{2}\left( \tau \right) \right) =\left\vert x-z\right\vert \text{ .}
\label{4.335}
\end{equation}%
Up to a rigid transformation of the space coordinates we have $z=\left(
z^{\prime },\varphi \left( z^{\prime },\tau \right) \right) $, where $%
\varphi \in C^{2,\beta }\left( B_{R_{0}}^{\prime }\left( 0\right) \times
\left( \tau -R_{0}^{2},\tau +R_{0}^{2}\right) \right) $ and $\left\Vert
\varphi \right\Vert _{C^{2,\beta }\left( B_{R_{0}}^{\prime }\left( 0\right)
\times \left( \tau -R_{0}^{2},\tau +R_{0}^{2}\right) \right) }\leq ER_{0}$.
Therefore, for any $t\in \left( \tau -R_{0}^{2},\tau +R_{0}^{2}\right) $ we
have
\begin{equation*}
\left\vert x-z\right\vert \geq \left\vert x-\left( z^{\prime },\varphi
\left( z^{\prime },t\right) \right) \right\vert -\left\vert \varphi \left(
z^{\prime },t\right) -\varphi \left( z^{\prime },\tau \right) \right\vert
\geq dist\left( x,\Omega _{2}\left( t\right) \right) -E\frac{\left\vert
t-\tau \right\vert }{R_{0}}\text{ ,}
\end{equation*}%
hence by (\ref{4.335}) we get
\begin{equation*}
dist\left( x,\Omega _{2}\left( t\right) \right) \leq dist\left( x,\Omega
_{2}\left( \tau \right) \right) +E\frac{\left\vert t-\tau \right\vert }{R_{0}%
}\text{ .}
\end{equation*}%
If case b) occurs, by (\ref{4.320}) we have
\begin{equation*}
dist\left( x,\Omega _{2}\left( t\right) \right) \leq \left\vert x-\Phi
_{2}^{-1}\left( x\right) \right\vert \leq CE\frac{\left\vert t-\tau
\right\vert }{R_{0}}\text{ ,}
\end{equation*}%
where $C$ depends on $E,M$ and $\beta $ and only.

\noindent Therefore, for any $x\in \Omega _{1}\left( t\right) $, we have
\begin{equation*}
dist\left( x,\Omega _{2}\left( t\right) \right) \leq dist\left( x,\Omega
_{2}\left( \tau \right) \right) +CE\frac{\left\vert t-\tau \right\vert }{%
R_{0}}\text{ ,}
\end{equation*}
where $C$ depends on $E,M$ and $\beta $ and only and (\ref{4.330}) follows.

\noindent Now by (\ref{4.325}) and (\ref{4.330}) we have
\begin{equation}
\sup\limits_{x\in \Omega _{1}\left( t\right) }dist\left( x,\Omega _{2}\left(
t\right) \right) \leq \sup\limits_{x\in \Omega _{1}\left( \tau \right)
}dist\left( x,\Omega _{2}\left( \tau \right) \right) +CE\frac{\left\vert
t-\tau \right\vert }{R_{0}}\text{ ,}  \label{4.340}
\end{equation}
where $C$ depends on $E,M$ and $\beta $ and only. Analogously we have
\begin{equation}
\sup\limits_{x\in \Omega _{2}\left( t\right) }dist\left( x,\Omega _{1}\left(
t\right) \right) \leq \sup\limits_{x\in \Omega _{2}\left( \tau \right)
}dist\left( x,\Omega _{1}\left( \tau \right) \right) +CE\frac{\left\vert
t-\tau \right\vert }{R_{0}}\text{ ,}  \label{4.345}
\end{equation}
where $C$ is the same of (\ref{4.340}).

\noindent By (\ref{4.340}) and (\ref{4.345}) we obtain
\begin{equation*}
d_{\mathcal{H}}\left( \overline{\Omega _{1}\left( t\right) },\overline{%
\Omega _{2}\left( t\right) }\right) \leq d_{\mathcal{H}}\left( \overline{%
\Omega _{1}\left( \tau \right) },\overline{\Omega _{2}\left( \tau \right) }%
\right) +CE\frac{\left\vert t-\tau \right\vert }{R_{0}}\text{ ,}
\end{equation*}
where $C$ is the same of (\ref{4.340}) and (\ref{4.310}) follows.$%
\blacksquare $

\subsection{Some extensions of Theorem \protect\ref{Th4.1}\label{Subs4.3}}

The extensions of Theorem \ref{Th4.1} that we give in the present section
can be proved with slight changes. We shall give short outline of the proofs.

\bigskip

\textbf{Variant I.} Assume that $\left\{ \Omega _{i}\left( t\right) \right\}
_{t\in \mathbb{R}}$, $\left\{ A_{i}\left( t\right) \right\} _{t\in \left[ 0,T%
\right] }$, $\left\{ I_{i}\left( t\right) \right\} _{t\in \left[ 0,T\right]
} $, $i=1,2$ satisfy the same hypotheses of Theorem \ref{Th4.1}, (in
particular $A_{i}\left( t\right) =A$, $t\in \left[ 0,T\right] $, $i=1,2$)
but instead of $\Omega _{1}\left( 0\right) =\Omega _{2}\left( 0\right) $
assume that $\partial A\cap \partial I_{i}\left( t\right) =\varnothing $ for
$t\in \left[ 0,T\right] $, $i=1,2$ and
\begin{equation}
d_{\mathcal{H}}\left( \overline{\Omega _{1}\left( 0\right) },\overline{%
\Omega _{2}\left( 0\right) }\right) \leq R_{0}\varepsilon \text{ .}
\label{4.2.1}
\end{equation}%
Let $f\in C^{2,\beta }\left( A\times \left[ 0,T\right] \right) $ satisfy (%
\ref{4.50}) and for any $t\in \left[ 0,T\right] $%
\begin{equation}
\left( \frac{1}{\left\vert A\right\vert }\int\nolimits_{A}\left\vert
f\left( x,t\right) -f_{A}\left( t\right) \right\vert ^{2}ds\right)
^{1/2}\geq F_{1}\left( t\right) \text{ ,}  \label{4.2.5}
\end{equation}%
where $f_{A}\left( t\right) =\frac{1}{\left\vert A\right\vert }%
\int\nolimits_{A}f\left( x,t\right) ds$ and $F_{1}$ is a strictly
increasing continuous function on $\left[ 0,T\right] $ such that $%
F_{1}\left( 0\right) =0$. Let $c_{1}$, $c_{2}$ be two numbers satisfying the
conditions
\begin{equation}
\left\vert c_{i}\right\vert \leq F_{0}\text{ , }i=1,2\text{.}  \label{4.2.10}
\end{equation}%
Let $u_{0}\in C^{2,\beta }\left( \mathbb{R}^{n}\right) $ satisfy
\begin{equation}
\left\Vert u_{0}\right\Vert _{C^{2,\beta }\left( \mathbb{R}^{n}\right) }\leq
F_{0}\text{ , }u_{0}\left( .,0\right) =c_{i}\text{, on }I_{i}\left( 0\right)
\text{ , }i=1,2\text{.}  \label{4.2.15}
\end{equation}%
Let $u_{i}\in C^{2,\beta }\left( \Omega _{i}\left( \left[ 0,T\right] \right)
\right) $ be the solution to
\begin{equation}
\left\{
\begin{array}{c}
div\left( \kappa \left( x,t,u_{i}\right) \nabla u_{i}\right) -\partial
_{t}u_{i}=0\text{, \ \ in }\Omega _{i}\left( \left( 0,T\right) \right) \text{%
,} \\
u_{i}=f\text{, \ \ \ \ \ \ \ \ \ \ \ on }A\times \left( 0,T\right] \text{,\ }
\\
u_{i}=c_{i}\text{, \ \ \ \ \ \ \ \ \ \ \ \ \ on }I_{i}\left( \left( 0,T%
\right] \right) \text{, \ \ } \\
u_{i}\left( .,0\right) =u_{0}\text{, \ \ \ \ \ \ \ \ \ \ in }\Omega
_{i}\left( 0\right) \text{ ,}%
\end{array}%
\right.  \label{4.2.20-23}
\end{equation}%
If we have
\begin{equation}
R_{0}\left\Vert \kappa _{1}\nabla u_{1}\cdot \nu -\kappa _{2}\nabla
u_{2}\cdot \nu \right\Vert _{L^{2}\left( \Sigma \times \left( 0,T\right)
\right) }\leq \varepsilon \text{ ,}  \label{4.2.30}
\end{equation}%
where $\kappa _{i}=\kappa \left( x,t,u_{i}\left( x,t\right) \right) $, $%
i=1,2 $, then we have, for every $\tau \in \left( 0,T\right] $ and $%
\varepsilon \in \left( 0,1\right) $,
\begin{equation}
\sup\limits_{t\in \left[ \tau ,T\right] }d_{\mathcal{H}}\left( \overline{%
\Omega _{1}\left( t\right) },\overline{\Omega _{2}\left( t\right) }\right)
\leq R_{0}C_{1}\left( \tau \right) \left\vert \log \varepsilon \right\vert
^{-\frac{1}{C_{2}\left( \tau \right) }}\text{ ,}  \label{4.2.35}
\end{equation}%
where $C_{1}\left( \tau \right) $ and $C_{2}\left( \tau \right) $ are
positive constants depending on $\tau ,\lambda _{0},\Lambda
_{0},E,M,F_{0},\beta ,R_{0}^{2}T^{-1}$ and $F_{1}\left( .\right) $ only.

\bigskip

\textbf{Proof of Variant I. }If $\varepsilon \leq \frac{d_{0}}{R_{0}}$,
where $d_{0}$ is defined in Proposition \ref{Pr4.8}, then we have easily
from (\ref{4.361}), (\ref{4.2.1}) and (\ref{4.2.15})
\begin{equation}
\left\vert c_{1}-c_{2}\right\vert \leq F_{0}\varepsilon \text{.}
\label{4.2.40}
\end{equation}%
By (\ref{4.50}), (\ref{4.2.10}) and (\ref{4.2.15}) we have that (\ref{4.95})
holds true. Now, arguing as in the proof of Proposition \ref{Pr4.2} and
taking into account (\ref{4.2.40}) we have, instead of (\ref{4.85}),
\begin{equation}
\int\nolimits_{\Omega _{i}\left( t\right) \smallsetminus G\left( t\right)
}\left( u_{i}\left( x,t\right) -c_{i}\right) ^{2}dx\leq R_{0}^{n}H^{2}\omega
\left( \frac{\varepsilon }{H}\right) \text{ , }i=1,2\text{, }t\in \left[ 0,T%
\right] \text{ ,}  \label{4.2.45}
\end{equation}%
where $\omega $ is an increasing continuous function on $\left[ 0,+\infty
\right) $ which satisfies (\ref{4.90}). Analogously, arguing as in the proof
of Proposition \ref{Pr4.3} we have that, if $\partial G\left( t\right) $ is
of Lipschitz class with constants $\rho _{0}$, $L$ then (\ref{4.2.45}) holds
true with $\omega $ satisfying (\ref{4.185}). Concerning Proposition \ref%
{Pr4.4}, it is enough to observe that the functions $v_{i}:=u_{i}-c_{i}$ are
solutions to the equation
\begin{equation*}
div\left( \kappa \left( x,t,v_{i}+c_{i}\right) \nabla v_{i}\right) -\partial
_{t}v_{i}=0\text{, \ \ in }\Omega _{i}\left( \left( 0,T\right) \right) \text{%
, }i=1,2\text{,}
\end{equation*}%
so we have that (\ref{4.220}) holds true replacing there $u$ by $u_{i}-c_{i}$
and $\Omega \left( t\right) $ by $\Omega _{i}\left( t\right) $, $i=1,2$, $%
t\in \left[ 0,T\right] $. With the changes indicated above and proceeding as
in the prof of Theorem \ref{Th4.1}, (\ref{4.2.35}) follows.$\blacksquare $

\bigskip

\textbf{Variant II.} Assume that $\left\{ \Omega _{i}\left( t\right)
\right\} _{t\in \mathbb{R}}$, $\left\{ A_{i}\left( t\right) \right\} _{t\in %
\left[ 0,T\right] }$, $\left\{ I_{i}\left( t\right) \right\} _{t\in \left[
0,T\right] }$, $i=1,2$ satisfy the same hypotheses of Variant I except for
assumption (\ref{4.2.1}). Let $f\in C^{2,\beta }\left( A\times \left[ 0,T%
\right] \right) $ satisfy $f\in C^{2,\beta }\left( A\times \left[ 0,T\right]
\right) $ satisfy (\ref{4.45}), (\ref{4.50}), (\ref{4.55}) and $f\left(
.,0\right) =0$ on $\partial A$. Let $u_{i}\in C^{2,\beta }\left( \Omega
_{i}\left( \left[ 0,T\right] \right) \right) $ be the solution to (\ref%
{4.1-4.4}) and assume that $u_{0}\equiv 0$. If (\ref{4.2.30}) is satisfied
then we have for any $\tau \in \left( 0,T\right] $ and $\varepsilon \in
\left( 0,1\right) $%
\begin{equation}
\sup\limits_{t\in \left[ \tau ,T\right] }d_{\mathcal{H}}\left( \overline{%
\Omega _{1}\left( t\right) },\overline{\Omega _{2}\left( t\right) }\right)
\leq R_{0}e^{\left( \frac{CH}{F_{1}\left( \tau \right) }\right)
^{C}}\left\vert \log \varepsilon \right\vert ^{-\frac{1}{C+\left(
F_{1}\left( \tau \right) /H\right) ^{-C}}}\text{ ,}  \label{4.2.50}
\end{equation}%
where $C$ is a positive constant depending on $\tau ,\lambda _{0},\Lambda
_{0},E,M,F_{0},\beta ,R_{0}^{2}T^{-1}$ and $F_{1}\left( .\right) $ only.

\bigskip

\textbf{Proof of Variant II. }Propositions \ref{Pr4.2} and \ref{Pr4.3}
continue to hold in this case. Concerning Proposition \ref{Pr4.4}, it is
simple to improve inequality (\ref{4.220}) using in its proof instead of $u$
the trivial extension $\widetilde{u}$. Thus, setting $\rho _{1}=\min \left\{
\overline{\rho },\rho _{0}\right\} $, $\sigma _{1}=\frac{\eta _{1}\rho _{1}}{%
6\sqrt{2E+4}}$ we have, instead of (\ref{4.255}) the inequality
\begin{equation}
\left\Vert u\left( .,\tau \right) \right\Vert _{L^{\infty }\left( B_{2\rho
_{0}}\left( y_{0}\right) \right) }\leq C_{1}\varepsilon ^{\mu }H^{1-\mu }%
\text{ ,}  \label{4.2.60}
\end{equation}%
where $\mu =\frac{2}{n+2}s_{2}^{m}$, $m\leq c_{0}M\left( \frac{R_{0}}{%
2\sigma _{1}}\right) ^{n}$ and $s_{2}$, $s_{2}\in \left( 0,1\right) $, $%
C_{1} $ depend on $\lambda _{0},\Lambda _{0},E,M,F_{0},\beta $ and $%
R_{0}^{2}T^{-1} $ only.

In order to apply Proposition \ref{Pr3.8} to $\widetilde{u}$ observe that,
with the same notation used in the proof of Proposition \ref{Pr4.4} and
setting $\delta _{1}=\frac{\sin \alpha }{\sqrt{2E+1}}$, we have
\begin{equation*}
\mathcal{S}\left( \left( 0,\tau \right) ,e_{n},\alpha ,\delta
_{1},R_{1}\right) \subset \Omega \left( \left( -\infty ,\tau \right] \right)
\text{ .}
\end{equation*}%
Furthermore, denoting by%
\begin{equation*}
\alpha _{1}=\arcsin \left( \min \left\{ \sin \alpha ,\delta _{1}\left(
1-\sin \alpha \right) \right\} \right) \text{ ,}
\end{equation*}%
\begin{equation*}
B\left( \tau \right) =\frac{1}{C_{2}}\left\vert \log \frac{1-\frac{1}{4}\eta
_{1}\sin \alpha _{1}\left( \tau \right) }{1+\frac{1}{4}\eta _{1}\sin \alpha
_{1}\left( \tau \right) }\right\vert \text{ ,}
\end{equation*}%
where $C_{2}$ depends on $\lambda _{0},\Lambda _{0},E,M,F_{0},\beta $ and $%
R_{0}^{2}T^{-1}$ only, we have, instead of (\ref{4.300}) the inequality
\begin{equation}
F_{1}\left( \tau \right) \leq He^{C\left( \frac{R_{0}}{\overline{\rho }}%
\right) ^{n-1}}\left\vert \log \left( \frac{\varepsilon _{0}}{eH}\right)
\right\vert ^{-\frac{\overline{\rho }}{CR_{0}}}\text{ ,}  \label{4.2.65}
\end{equation}%
where $C$, $C>1$, depends on $\lambda _{0},\Lambda _{0},E,M,F_{0},\beta $
and $R_{0}^{2}T^{-1}$ only. By (\ref{4.2.65}) we have, for every $\overline{%
\rho }>0$, $t\in \left( 0,T\right] $, $z\in \left( \Omega \left( t\right)
\right) _{2\overline{\rho }}$ ,
\begin{eqnarray}
&&\int\nolimits_{B_{\overline{\rho }}\left( z\right) }u^{2}\left(
x,t\right) dx  \label{4.2.70} \\
&\geq &H^{2}\overline{\rho }^{n}\exp \left( -e^{\frac{CR_{0}^{n}}{\overline{%
\rho }^{n}}}\left( \frac{F_{1}\left( t\right) }{H}\right) ^{\frac{CR_{0}}{%
\overline{\rho }}}\right) \text{ ,}  \notag
\end{eqnarray}%
where $C$, $C>1$, depends on $\lambda _{0},\Lambda _{0},E,M,F_{0},\beta $
and $R_{0}^{2}T^{-1}$ only.

In order to prove (\ref{4.2.50}) we argue as in the proof of Theorem \ref%
{Th4.1}, in doing this we use (\ref{4.2.70}) instead of (\ref{4.300}) with $%
\widetilde{u}$ instead of $u$. We get, instead of (\ref{4.375}) and (\ref%
{4.376}) respectively ($\sigma $ being defined by (\ref{4.370}))
\begin{equation}
d_{m}\left( t\right) \leq \frac{CR_{0}}{b_{1}\left( t\right) }\left( \frac{%
\sigma }{MH^{2}}\right) ^{\frac{1}{\log \frac{C}{b_{1}\left( t\right) }}}%
\text{ ,}  \label{4.2.75}
\end{equation}
and
\begin{equation}
d\left( t\right) \leq \frac{CR_{0}}{b_{1}\left( t\right) }\left( \frac{%
\sigma }{MH^{2}}\right) ^{\frac{1}{\log \frac{C}{b_{1}\left( t\right) }}}%
\text{ ,}  \label{4.2.80}
\end{equation}
where
\begin{equation}
b_{1}\left( t\right) =\exp \left( -\left( \frac{e^{-1}F_{1}\left( t\right) }{%
H}\right) ^{-C_{1}}\right) \text{ ,}  \label{4.2.85}
\end{equation}
and $C$, $C>1$, $C_{1}$ depend on $\lambda _{0},\Lambda _{0},E,M,F_{0},\beta
$ and $R_{0}^{2}T^{-1}$ only.

Now, let us denote
\begin{equation*}
\widetilde{\psi }\left( \mu \right) =\left( e^{-1}F_{1}\left( T_{0}\mu
\right) \right) ^{-C_{1}}+\log \left( Ce\right) \text{ , }\mu \in \left( 0,1%
\right]
\end{equation*}
and let $\widetilde{\Psi }_{0}=\inf\limits_{\left( 0,1\right] }\widetilde{%
\psi }$, notice $\widetilde{\Psi }_{0}\geq 1$.

\noindent By (\ref{4.2.80}) and Proposition \ref{Pr4.6} we have
\begin{equation}
\sup\limits_{t\in \left[ 0,T\right] }d\left( t\right) \leq R_{0}\left(
Ce^{\left( \widetilde{\psi }\left( \mu \right) \right) }\left( \frac{\sigma
}{MH^{2}}\right) ^{\frac{1}{\widetilde{\psi }\left( \mu \right) }}+\mu
\right) \ \text{,}\ \text{for every }\mu \in \left( 0,1\right] \text{,}
\label{4.2.90}
\end{equation}
where $C$, $C>1$, depends on $\lambda _{0},\Lambda _{0},E,M,F_{0},\beta $
and $R_{0}^{2}T^{-1}$ only.

\noindent By (\ref{4.2.90}) and Proposition \ref{Pr4.2} we have
\begin{equation*}
\sup\limits_{t\in \left[ 0,T\right] }d_{\mathcal{H}}\left( \overline{\Omega
_{1}\left( t\right) },\overline{\Omega _{2}\left( t\right) }\right) \leq
R_{0}\widetilde{\phi }\left( \frac{1}{M}\left( \log \left\vert \log \frac{%
\varepsilon }{H}\right\vert \right) ^{-1/n}\right) \text{ ,}
\end{equation*}
for every $\varepsilon \in \left( 0,e^{-1}H\right) $, where
\begin{equation*}
\widetilde{\phi }\left( s\right) =\left\{
\begin{array}{c}
\left( \widetilde{\psi }^{-1}\left( \left\vert \log \left( \frac{s}{MH^{2}}%
\right) \right\vert \right) \right) ^{1/4}+\exp \left( -\widetilde{m}%
_{0}\left\vert \log \left( \frac{s}{MH^{2}}\right) \right\vert ^{3/4}\right)
\text{, if }s\in \left( 0,\widetilde{a}\right) \text{,} \\
C_{2}\widetilde{a}^{-1}s\text{ , if }s\in \left[ \widetilde{a},+\infty
\right) ,%
\end{array}
\right.
\end{equation*}
$\widetilde{a}=e^{-2\widetilde{\Psi }_{0}}MH^{2}$ , $\widetilde{m}%
_{0}=1-\left( 2\widetilde{\Psi }_{0}\right) ^{-1/2}$ and $C_{2}$ depends on $%
M$ and $E$ only. Now let us define

\begin{equation*}
\widetilde{\varepsilon }_{1}=\sup \left\{ \tau \in \left( 0,e^{-1}\right) :%
\widetilde{\phi }\left( \frac{1}{M}\left( \log \left\vert \log \tau
\right\vert \right) ^{-1/n}\right) \leq \frac{d_{0}}{R_{0}}\right\} \text{ ,}
\end{equation*}
where $d_{0}$ is defined in Proposition \ref{Pr4.8}. From such a Proposition %
\ref{Pr4.8} we have that if $\varepsilon \in \left( 0,\widetilde{\varepsilon
}_{1}H\right] $ then there exist $\rho _{0}\in \left( 0,R_{0}\right] $ and $%
L>0$, $\frac{\rho _{0}}{R_{0}}$ and $L$ depend on $E$ only, such that $%
\partial G\left( t\right) $ is of Lipschitz class of constants $\rho _{0}$, $%
L$. Therefore by Proposition \ref{Pr4.3} and (\ref{4.376}) we have
\begin{equation*}
d\left( t\right) \leq \frac{CR_{0}}{b_{1}\left( t\right) }\left( C\left\vert
\log \frac{\varepsilon }{H}\right\vert ^{-1/C}\right) ^{\frac{1}{\log \frac{C%
}{b_{1}\left( t\right) }}}\text{ ,}
\end{equation*}
where $C$, $C>1$, depends on $\lambda _{0},\Lambda _{0},E,M,F_{0},\beta $
and $R_{0}^{2}T^{-1}$ only. Otherwise, if $\varepsilon >\widetilde{%
\varepsilon }_{1}H$ we have trivially $d\left( t\right) \leq C^{\prime }%
\frac{\varepsilon }{H\widetilde{\varepsilon }_{1}}$, where $C^{\prime }$
depends on $E$ and $M$ only and (\ref{4.2.50}) follows.$\blacksquare $

\bigskip

\textbf{Variant III.} The stability results stated in Theorem \ref{Th4.1},
Variant I and Variant II continue to hold if the first equation in (\ref%
{4.1-4.4}) is replaced by the equation
\begin{equation}
div\left( \kappa \left( x,t,u\right) \nabla u\right) +B\left( x,t,u,\nabla
u\right) -\partial _{t}u=0\text{,}\ \text{in }\Omega \left( \left(
0,T\right) \right) \text{,}  \label{4.2.150}
\end{equation}%
where the matrix $\kappa \left( x,t,u\right) $ satisfies the same hypotheses
of Theorem \ref{Th4.1} and $B:\mathbb{R}^{2n+2}\rightarrow \mathbb{R}$
satisfies the following conditions
\begin{equation}
B\left( x,t,0,0\right) =0\text{, for every }\left( x,t\right) \in \mathbb{R}%
^{n+1}\text{ ,}  \label{4.2.155}
\end{equation}%
\begin{equation}
R_{0}^{2}\left\Vert \nabla _{x}B\right\Vert _{L^{\infty }\left( \mathbb{R}%
^{2n+2}\right) }+R_{0}^{2}\left\Vert \partial _{u}B\right\Vert _{L^{\infty
}\left( \mathbb{R}^{2n+2}\right) }+R_{0}\left\Vert \nabla _{p}B\right\Vert
_{L^{\infty }\left( \mathbb{R}^{2n+2}\right) }\leq \Lambda _{0}\text{ .}
\label{4.2.170}
\end{equation}

\bigskip

\textbf{Proof of Variant III. }To obtain the proof of Variant III we need to
change the proof of Theorem \ref{Th4.1} only for what concerns the proofs of
inequalities (\ref{4.155}) and (\ref{4.195}). Actually such inequalities can
be achieved multiplying the new equations by $u_{1}e^{-\gamma t}$, with $%
\gamma $ large enough, and then proceeding in the same way of Propositions %
\ref{Pr4.2} and \ref{Pr4.3} .$\blacksquare $

\bigskip

\textbf{Variant IV.} If the first equation in (\ref{4.1-4.4}) is replaced by
a linear equation then we can get Theorem \ref{Th4.1} and the stability
results of Variant I, II under relaxed a priori assumptions. In what follows
we will discuss in more detail the case in which the first equation in (\ref%
{4.1-4.4}) is replaced by the equation
\begin{equation}
div\left( \kappa \left( x,t\right) \nabla u\right) +b\left( x,t\right) \cdot
\nabla u+c\left( x,t\right) u-\partial _{t}u=0\text{,}\ \text{in }\Omega
\left( \left( 0,T\right) \right) \text{ ,}  \label{4.2.190}
\end{equation}%
where $\kappa \left( x,t\right) =\left\{ \kappa ^{ij}\left( x,t\right)
\right\} _{i,j=1}^{n}$ is a symmetric $n\times n$ matrix, $b\left(
x,t\right) $ and $c\left( x,t\right) $ are a vector valued measurable
function and a measurable function respectively. We assume that $\kappa
\left( x,t\right) $ satisfies the following conditions for every $\xi \in
\mathbb{R}^{n}$ and every $\left( x,t\right) ,\left( y,s\right) \in \mathbb{R%
}^{n+1}$,
\begin{equation}
\lambda _{0}\left\vert \xi \right\vert ^{2}\leq
\sum\limits_{i,j=1}^{n}\kappa ^{ij}\left( x,t\right) \xi _{i}\xi _{j}\leq
\lambda _{0}^{-1}\left\vert \xi \right\vert ^{2}\text{ }  \label{4.2.195}
\end{equation}%
and
\begin{equation}
\left( \sum\limits_{i,j=1}^{n}\left( \kappa ^{ij}\left( x,t\right) -\kappa
^{ij}\left( y,\tau \right) \right) ^{2}\right) ^{1/2}\leq \frac{\Lambda _{0}%
}{R_{0}}\left( \left\vert x-y\right\vert ^{2}+\left\vert t-s\right\vert
\right) ^{1/2}\text{ .}  \label{4.2.200}
\end{equation}%
Concerning $b\left( x,t\right) $ and $c\left( x,t\right) $ we assume that
\begin{equation}
R_{0}\left\Vert b\right\Vert _{L^{\infty }\left( \mathbb{R}^{n+1}\right)
}+R_{0}^{2}\left\Vert c\right\Vert _{L^{\infty }\left( \mathbb{R}%
^{n+1}\right) }\leq \Lambda _{0}\text{.}  \label{4.2.205}
\end{equation}%
Assume that $\left\{ \Omega _{i}\left( t\right) \right\} _{t\in \mathbb{R}}$%
, $i=1,2$, are two families of domains satisfying (\ref{4.10}) which
(instead of (\ref{4.25})), satisfy
\begin{equation}
\partial \Omega _{i}\left( \left[ 0,T\right] \right) \text{ is a portion of }%
\partial \Omega _{i}\left( \left( -\infty ,+\infty \right) \right)
\label{4.2.210}
\end{equation}%
\begin{equation*}
\text{of class }C^{1,1}\text{ with constants }R_{0}\text{, }E\text{.}
\end{equation*}%
For $i=1,2$ and for any $t\in \left[ 0,T\right] $ let $A_{i}\left( t\right) $%
, $I_{i}\left( t\right) $ satisfy (\ref{4.20}). Assume that, for any $t\in %
\left[ 0,T\right] $, $A_{1}\left( t\right) =A_{2}\left( t\right) =A$ and
also\ (instead of \ref{4.26})
\begin{equation}
A\times \left[ 0,T\right] \text{, }I\left( \left[ 0,T\right] \right) \text{
are portions of }\partial \Omega \left( \left( -\infty ,+\infty \right)
\right)  \label{4.2.220}
\end{equation}%
\begin{equation*}
\text{of class }C^{1,1}\text{ with constants }R_{0}\text{, }E\text{.}
\end{equation*}%
Let $u_{0}\in L^{\infty }\left( \mathbb{R}^{n}\right) $, $f\in H^{\frac{1}{2}%
,\frac{1}{4}}\left( A\times \left( 0,T\right) \right) $ satisfy (instead of (%
\ref{4.38})-(\ref{4.50}))
\begin{equation}
\left\Vert u_{0}\right\Vert _{L^{\infty }\left( \mathbb{R}^{n}\right)
}+\left\Vert f\right\Vert _{H^{\frac{1}{2},\frac{1}{4}}\left( A\times \left(
0,T\right) \right) }\leq F_{0}\text{ ,}  \label{4.2.225}
\end{equation}%
\begin{equation}
\text{\textit{supp}}f\subset \widetilde{\Sigma }\times \left[ 0,T\right]
\text{ ,}  \label{4.2.226}
\end{equation}%
\begin{equation}
\widetilde{F}_{1}\left( \tau \right) \leq \frac{\left\Vert f\right\Vert
_{L^{2}\left( A\times \left( \tau ,T\right) \right) }}{\left\Vert
u_{0}\right\Vert _{L^{\infty }\left( \mathbb{R}^{n}\right) }+\left\Vert
f\right\Vert _{H^{\frac{1}{2},\frac{1}{4}}\left( A\times \left( 0,T\right)
\right) }}\text{, for }\tau \in \left[ 0,T\right] \text{,}  \label{4.2.230}
\end{equation}%
where $\widetilde{F}_{1}\left( \tau \right) $ is a given strictly increasing
function on $\left[ 0,T\right] $ such that $\widetilde{F}_{1}\left( 0\right)
=0$.

\noindent Let $u_{i}$, $i=1,2$, be the weak solution to
\begin{equation*}
\left\{
\begin{array}{c}
div\left( \kappa \nabla u_{i}\right) +b\cdot \nabla u_{i}+cu_{i}-\partial
_{t}u_{i}=0\text{, \ \ in }\Omega _{i}\left( \left( 0,T\right) \right) \text{%
,} \\
u_{i}=f\text{, \ \ \ \ \ \ \ \ \ \ \ on }A\times \left( 0,T\right] \text{,\ }
\\
u_{i}=c_{i}\text{, \ \ \ \ \ \ \ \ \ \ \ \ \ on }I_{i}\left( \left( 0,T%
\right] \right) \text{, \ \ } \\
u_{i}\left( .,0\right) =u_{0}\text{, \ \ \ \ \ \ \ \ \ \ in }\Omega
_{i}\left( 0\right) \text{ ,}%
\end{array}
\right.
\end{equation*}
Assume that
\begin{equation*}
R_{0}\left\Vert \kappa \nabla u_{1}\cdot \nu -\kappa \nabla u_{2}\cdot \nu
\right\Vert _{L^{2}\left( \Sigma \times \left( 0,T\right) \right) }\leq
\varepsilon \text{ ,}
\end{equation*}
and assume that one of the following conditions is satisfied
\begin{equation*}
d_{\mathcal{H}}\left( \overline{\Omega _{1}\left( 0\right) },\overline{%
\Omega _{2}\left( 0\right) }\right) \leq R_{0}\varepsilon \text{ }
\end{equation*}
or
\begin{equation}
u_{0}\equiv 0\text{,}  \label{4.2.242}
\end{equation}
then we have, for any $\tau \in \left( 0,T\right] $ and $\varepsilon \in
\left( 0,1\right) $,
\begin{equation*}
\sup\limits_{t\in \left[ \tau ,T\right] }d_{\mathcal{H}}\left( \overline{%
\Omega _{1}\left( t\right) },\overline{\Omega _{2}\left( t\right) }\right)
\leq R_{0}C_{1}\left( \tau \right) \left\vert \log \varepsilon \right\vert
^{-\frac{1}{C_{2}\left( \tau \right) }}\text{ ,}
\end{equation*}
where $C_{1}\left( \tau \right) $ and $C_{2}\left( \tau \right) $ are
positive constants depending on $\tau ,\lambda _{0},\Lambda
_{0},E,M,F_{0},\beta ,R_{0}^{2}T^{-1}$ and $\widetilde{F}_{1}\left( .\right)
$ only. In addition, if (\ref{4.2.242}) is satisfied then
\begin{equation*}
C_{1}\left( \tau \right) =\exp \frac{C}{\left( \widetilde{F}_{1}\left( \tau
\right) \right) ^{C}}\text{ , }C_{2}\left( \tau \right) =C+\left( \widetilde{%
F}_{1}\left( \tau \right) \right) ^{-C}\text{ ,}
\end{equation*}
where $C$, $C>1$, depends on $\lambda _{0},\Lambda _{0},E,M,F_{0},\beta $
and $R_{0}^{2}T^{-1}$ only.

\section{Exponential instability\label{INST}}

In the present section we prove that the Dirichlet inverse problem studied
above is exponentially unstable, that is, the stability estimate for such
problem cannot be better than logarithmic. In order to prove this
instability property we consider two Dirichlet inverse problems: in the
first problem the space dimension $n$ is equal to $1$ and the unknown
boundary is a curve $x=s\left( t\right) $, in the second problem the unknown
boundary is not time varying. Since the second problem has been studied in
\cite{DcRVe} and it requires more technical proofs than the first one, here
we give the detailed proofs for the first problem only.

The common tool used to prove the instability character of the above inverse
problems is an abstract theorem of which we give the statement below and
refer for the proof to \cite{DcR}.

Let $\left( X,d\right) $ be a metric space and let $H$ be a separable
Hilbert space with scalar product $\left( .,.\right) $. We denote by $%
H^{\prime }$ the dual of $H$ and, for any $v\in H$ and any $v\prime \in
H^{\prime }$, we denote by $\left\langle v^{\prime },v\right\rangle $ the
duality pairing between $v^{\prime }$ and $v$. We denote by $\mathcal{L}%
\left( H,H^{\prime }\right) $ the space of bounded linear operators between $%
H$ and $H^{\prime }$ with the usual operator norm. We shall also fix $\gamma
:H\smallsetminus \left\{ 0\right\} \rightarrow \left[ 0,+\infty \right] $
such that
\begin{equation*}
\gamma \left( \lambda v\right) =\gamma \left( v\right) \text{, for every }%
v\in H\smallsetminus \left\{ 0\right\} \text{ and }\lambda \in \mathbb{R}%
\smallsetminus \left\{ 0\right\} \text{ .}
\end{equation*}%
Let us remark that the function $\gamma $ may attain both the values $0$ and
$+\infty $ and can be thought of as a suitable Rayleigh quotient.

Let $F$ be a function from $X$ to $\mathcal{L}\left( H,H^{\prime }\right) $,
that is for any $x\in X$, $F\left( x\right) $ is a linear and bounded
operator between $H$ and $H^{\prime }$. We also fix a reference operator $%
F_{0}\in \mathcal{L}\left( H,H^{\prime }\right) $ and a reference point $%
x_{0}\in X$ (notice that it is not required any connection between $F_{0}$
and $F\left( x_{0}\right) $). For any $\varepsilon >0$ denote $%
X_{\varepsilon }=\left\{ x\in X:d\left( x,x_{0}\leq \varepsilon \right)
\right\} $.

\begin{definition}
\label{DefInst}Let $\left( Y,d_{Y}\right) $ be a metric space. Let $%
\varepsilon $ be a positive number. We shall say that a subset $W$ of $Y$ is
a $\varepsilon $-discrete set if for any two distinct points $w_{1},w_{2}\in
W$ we have $d_{Y}\left( w_{1},w_{2}\right) \geq \varepsilon $.
\end{definition}

\begin{theorem}
\label{ThabsInst}Let us assume that the following conditions are satisfied.

\noindent i) There exist positive constants $\varepsilon _{0}$, $C_{1}$ and $%
\alpha _{1}$ such that for every $\varepsilon \in \left( 0,\varepsilon
_{0}\right) $ there exists an $\varepsilon $-discrete set $Z$ contained in $%
X_{\varepsilon }$ with at least $\exp \left( C_{1}\varepsilon ^{-\alpha
_{1}}\right) $ elements.

\noindent ii) There exist three positive constants $p$, $C_{2}$ and $\alpha
_{2}$ and an orthonormal basis of $H$, $\left\{ \psi _{k}\right\} _{k\in
\mathbb{N}}$ such that the following condition hold.

For any $k\in \mathbb{N}$, we have that $\gamma \left( \psi _{k}\right)
<+\infty $, and for any $q\in \mathbb{N}$,
\begin{equation}
\#\left\{ k\in \mathbb{N:}\gamma \left( \psi _{k}\right) \leq q\right\} \leq
C_{2}\left( 1+q\right) ^{p}\text{,}  \label{inst3}
\end{equation}

where $\#$ denotes the number of elements.

For any $x\in X$ and $\left( k,l\right) \in \mathbb{N\times N}$ we have
\begin{equation}
\left\vert \left\langle \left( F\left( x\right) -F_{0}\right) \psi _{k},\psi
_{l}\right\rangle \right\vert \leq C_{2}\exp \left( -\alpha _{2}\max \left\{
\gamma \left( \psi _{k}\right) ,\gamma \left( \psi _{l}\right) \right\}
\right) \text{.}  \label{inst4}
\end{equation}
Then there exists a positive constant $\varepsilon _{1}$, depending on $%
\varepsilon _{0}$, $C_{1}$, $C_{2}$, $\alpha _{1}$, $\alpha _{2}$ and $p$
only, such that for every $\varepsilon \in \left( 0,\varepsilon _{1}\right) $
we can find $x_{1}$ and $x_{2}$ satisfying
\begin{equation}
x_{1},x_{2}\in X_{\varepsilon }\text{ , }d\left( x_{1},x_{2}\right) \geq
\varepsilon \text{,}  \label{inst5a}
\end{equation}
\begin{equation}
\left\Vert F\left( x_{1}\right) -F\left( x_{2}\right) \right\Vert _{\mathcal{%
L}\left( H,H^{\prime }\right) }\leq 2\exp \left( -\varepsilon ^{-\alpha
_{1}/2\left( p+1\right) }\right) \text{.}  \label{inst5b}
\end{equation}
\end{theorem}

\subsection{Exponential instability of the Dirichlet inverse problem with
time-varying unknown boundary\label{INSTT}}

We shall use the following Hilbert space. The space $H=H^{\frac{3}{4}}\left(
\left( \frac{\pi }{2},\frac{3\pi }{2}\right) \right) $, its dual $H^{\prime
}=H_{1}=H^{-\frac{3}{4}}\left( \left( \frac{\pi }{2},\frac{3\pi }{2}\right)
\right) $ and $H_{0}=H^{\frac{1}{4}}\left( \left( \frac{\pi }{2},\frac{3\pi
}{2}\right) \right) $. We consider the interpolation spaces between $H_{0}$
and $H_{1}$, \cite{LiMa}. For any $\theta $, $0\leq \theta \leq 1$, we
define $H_{\theta }$ as $\left[ H_{0},H_{1}\right] _{\theta }$, where the
latter denotes the interpolation at level $\theta $ between $H_{0}$ and $%
H_{1}$. The norm in $H_{0}$ and $H_{1}$ will be denoted by $\left\Vert
.\right\Vert _{\theta }$. Recall that for any $\theta \in \left[ 0,1\right] $
there exists a constant $C_{\theta }$, depending on $\theta $ only, such
that the following interpolation inequality holds true
\begin{equation}
\left\Vert \psi \right\Vert _{\theta }\leq C_{\theta }\left\Vert \psi
\right\Vert _{0}^{1-\theta }\left\Vert \psi \right\Vert _{1}^{\theta }\text{%
, for every }\psi \in H_{0}\text{.}  \label{inst13}
\end{equation}%
By using the interpolation properties of fractional order Sobolev spaces,
see \cite{LiMa}, we can characterize $H_{\theta }$ as follows
\begin{equation}
H_{\theta }=\left\{
\begin{array}{cc}
H^{\frac{1}{4}-\theta }\left( \left( \frac{\pi }{2},\frac{3\pi }{2}\right)
\right) \text{, } & \text{if }0\leq \theta \leq \frac{1}{4}\text{,} \\
H^{-\left( \frac{1}{4}-\theta \right) }\left( \left( \frac{\pi }{2},\frac{%
3\pi }{2}\right) \right) \text{,} & \text{if }\frac{1}{4}\leq \theta \leq 1%
\text{, }\theta \neq \frac{3}{4}\text{.}%
\end{array}%
\right.  \label{inst14}
\end{equation}%
Let us notice that the interesting case $\theta =\frac{1}{4}$, where we have
$H_{\theta }=L^{2}\left( \left( \frac{\pi }{2},\frac{3\pi }{2}\right)
\right) $ and that $H_{\frac{3}{4}}$ does not coincide with $H^{-\frac{1}{2}%
}\left( \left( \frac{\pi }{2},\frac{3\pi }{2}\right) \right) $.

Let $m$ be a given positive integer number and let $b$, $\delta $ be
positive numbers as above, we define
\begin{equation*}
X_{mb\delta }=
\begin{array}{c}
\left\{ s\in C^{m}\left( \left[ \frac{\pi }{2},\frac{3\pi }{2}\right]
\right) :\left\Vert s\right\Vert _{C^{m}\left( \left[ \frac{\pi }{2},\frac{%
3\pi }{2}\right] \right) }\leq b\text{, }\right. \\
\left. 1\leq s\leq 1+\delta \text{, }s\left( \frac{\pi }{2}\right) =s\left(
\frac{3\pi }{2}\right) =1\right\} \text{.}%
\end{array}%
\end{equation*}
Let us consider the metric space $\left( X,d\right) :=\left( X_{mb\delta
},d_{0}\right) $ where $d_{0}\left( s_{1},s_{2}\right) =\left\Vert
s_{1}-s_{2}\right\Vert _{L^{\infty }\left( \left[ \frac{\pi }{2},\frac{3\pi
}{2}\right] \right) }$.

For any $s\in X$ we denote
\begin{equation*}
Q_{s}=\left\{ \left( x,t\right) \in \mathbb{R}^{2}:t\in \left( \frac{\pi }{2}%
,\frac{3\pi }{2}\right) \text{, }0<x<s\left( t\right) \right\} \text{,}
\end{equation*}
if $s\left( t\right) =s_{0}\left( t\right) =1$ for every $t\in \left[ \frac{%
\pi }{2},\frac{3\pi }{2}\right] $ we set $Q_{0}=Q_{s_{0}}$. Notice that
there exists a constant $C$, $C\geq 1$, such that for any $s_{1},s_{2}\in X$
we have $C^{-1}d_{\mathcal{H}}\left( \overline{Q_{s_{1}}},\overline{Q_{s_{2}}%
}\right) \leq d_{0}\left( s_{1},s_{2}\right) \leq C^{-1}d_{\mathcal{H}%
}\left( \overline{Q_{s_{1}}},\overline{Q_{s_{2}}}\right) $.

\noindent For any $s\in X$ we consider the linear operator $\mathcal{D}%
_{s}:H\rightarrow H_{0}$ which is defined as follows. For any $\psi \in H$,
let $u\in H^{2,1}\left( Q_{s}\right) $ be the solution to
\begin{equation}
\left\{
\begin{array}{c}
\partial _{t}u-\partial _{x}^{2}u=0\text{, in }Q_{s}\text{,} \\
u\left( x,\frac{\pi }{2}\right) =0\text{, \ \ \ \ for }x\in \left[ 0,1\right]
\text{,} \\
u\left( s\left( t\right) ,t\right) =0\text{, for }t\in \left[ \frac{\pi }{2},%
\frac{3\pi }{2}\right] \text{,} \\
u\left( 0,t\right) =\psi \left( t\right) \text{, for }t\in \left[ \frac{\pi
}{2},\frac{3\pi }{2}\right] \text{,}%
\end{array}%
\right.  \label{inst15}
\end{equation}%
then, for any $\psi \in H$, we set
\begin{equation}
\mathcal{D}_{s}\psi =-\partial _{x}u\left( 0,t\right) \text{, }t\in \left(
\frac{\pi }{2},\frac{3\pi }{2}\right) \text{.}  \label{inst16}
\end{equation}%
By the trace theorem we have that the linear operator $\mathcal{D}_{s}$ is
bounded. We can also consider $\mathcal{D}_{s}$ as an operator belonging to $%
\mathcal{L}\left( H,H^{\prime }\right) $, by setting, for any $\psi ,\phi
\in H$,
\begin{equation}
\left\langle \mathcal{D}_{s}\psi ,\phi \right\rangle _{H^{\prime
},H}=\left\langle -\partial _{x}u\left( 0,.\right) ,\phi \right\rangle
_{H^{\prime },H}=\int\nolimits_{\frac{\pi }{2}}^{\frac{3\pi }{2}}-\partial
_{x}u\left( 0,t\right) \phi \left( t\right) dt\text{,}  \label{inst17}
\end{equation}%
where $u$ solves (\ref{inst15}) and $\left\langle .,.\right\rangle
_{H^{\prime },H}$ is the duality pairing between $H^{\prime }$ and $H$. We
refer to the operator $\mathcal{D}_{s}$ to as the \textit{%
Dirichlet-to-Neumann }map.

\begin{theorem}
\label{Thinst1}Let $m\in \mathbb{N}$ and $b$, $b>0$, be fixed. Then there
exists a positive constant $\delta _{1}$ depending on $m$ and $b$ only, such
that for every $\delta \in \left( 0,\delta _{1}\right) $ we can find $%
s_{1},s_{2}\in X$ satisfying
\begin{equation}
d_{0}\left( s_{i},s_{0}\right) \leq \delta \text{, }i=1,2\text{; }%
d_{0}\left( s_{1},s_{2}\right) \geq \delta  \label{inst18}
\end{equation}
and for any $\theta \in \left[ 0,1\right] $%
\begin{equation}
\left\| \mathcal{D}_{s_{1}}-\mathcal{D}_{s_{2}}\right\| _{\mathcal{L}\left(
H,H_{\theta }\right) }\leq C\exp \left( -\theta \delta ^{-\frac{1}{6m}%
}\right) \text{,}  \label{inst19}
\end{equation}
where $C$ is a constant depending on $m$, $b$ and $\theta $ only.
\end{theorem}

In order to prove Theorem \ref{Thinst1}, we need to introduce some notation.
For any $s\in X$ let $\widetilde{s}$ be its extension to $1$ to $\left[
0,2\pi \right] $ and denote $\widetilde{Q}_{s}=\left\{ \left( x,t\right) \in
\mathbb{R}^{2}:t\in \left( 0,2\pi \right) \text{, }0<x<\widetilde{s}\left(
t\right) \right\} $. Let us consider the Hilbert space $L^{2}\left( \left(
0,2\pi \right) \right) $ endowed with the scalar product
\begin{equation*}
\left( \psi ,\phi \right) _{0}=\int\nolimits_{0}^{2\pi }\psi \left(
t\right) \phi \left( t\right) dt\text{, for any }\psi ,\phi \in L^{2}\left(
\left( 0,2\pi \right) \right) \text{.}
\end{equation*}%
We have that the set
\begin{equation}
\left\{ \psi _{n}=\frac{1}{\sqrt{\pi }}\sin \left( \frac{n}{2}t\right) :n\in
\mathbb{N}\right\} \text{,}  \label{inst22}
\end{equation}%
is an orthonormal basis of $L^{2}\left( \left( 0,2\pi \right) \right) $ with
respect to the scalar product $\left( .,.\right) _{0}$.

For any $\sigma \in \left( 0,1\right] $ consider the Sobolev space $%
H_{,0}^{\sigma }\left( \left( 0,2\pi \right) \right) $. We have the
following properties. First, $H_{,0}^{\sigma }\left( \left( 0,2\pi \right)
\right) =H^{\sigma }\left( \left( 0,2\pi \right) \right) $ if and only if $%
\sigma \leq \dfrac{1}{2}$. Then, for any $\sigma \neq \dfrac{1}{2}$ we can
endow $H^{\sigma }\left( \left( 0,2\pi \right) \right) $ with the following
scalar product
\begin{equation*}
\left( \psi ,\phi \right) _{\sigma }=\sum\limits_{n=1}^{\infty }\left(
1+\left( \frac{n}{2}\right) ^{2\sigma }\right) \left( \psi ,\psi _{n}\right)
_{0}\left( \phi ,\psi _{n}\right) _{0}\text{ , for every }\psi ,\phi \in
H^{\sigma }\left( \left( 0,2\pi \right) \right) \text{,}
\end{equation*}%
with respect to which the set
\begin{equation}
\left\{ \widetilde{\psi }_{n}=\frac{\psi _{n}}{\sqrt{1+\left( \frac{n}{2}%
\right) ^{2\sigma }}}:n\in \mathbb{N}\right\} \text{ ,}  \label{inst24}
\end{equation}%
is an othonormal basis of $H_{,0}^{\sigma }\left( \left( 0,2\pi \right)
\right) $.We define $\gamma :\widetilde{H}\smallsetminus \left\{ 0\right\}
\rightarrow \left[ 0,+\infty \right] $, such that
\begin{equation*}
\gamma \left( \widetilde{\psi }\right) =\frac{\left\Vert \widetilde{\psi }%
\right\Vert _{\widetilde{H}}}{\left\Vert \widetilde{\psi }\right\Vert _{%
\widetilde{H}_{0}}}\text{, for every }\widetilde{\psi }\in \widetilde{H}%
\smallsetminus \left\{ 0\right\} \text{,}
\end{equation*}%
Notice that
\begin{equation*}
\gamma \left( \widetilde{\psi }_{n}\right) =\sqrt{\frac{n}{2}}\text{ , for
every }n\in \mathbb{N}
\end{equation*}%
and let us call for every positive integer $q$, $N_{1}\left( q\right)
=\#\left\{ n\in \mathbb{N:}\gamma \left( \widetilde{\psi }_{n}\right) \leq
q\right\} $. We have
\begin{equation}
N_{1}\left( q\right) \leq 2\left( 1+q\right) ^{2}\text{, for every }q\in
\mathbb{N}\text{.}  \label{inst25}
\end{equation}

Let $\widetilde{H}=H_{,0}^{\frac{3}{4}}\left( \left( 0,2\pi \right) \right) $%
, $\widetilde{H}^{\prime }=\widetilde{H}_{1}=H^{-\frac{3}{4}}\left( \left(
0,2\pi \right) \right) $ and $\widetilde{H}_{0}=H^{\frac{1}{4}}\left( \left(
0,2\pi \right) \right) $. For any $s\in X$, let us define in the same way as
before the Dirichlet-to-Neumann map associated to $\widetilde{Q}_{s}$, that
is the linear and bounded operator $\widetilde{\mathcal{D}}_{s}:\widetilde{H}%
\rightarrow \widetilde{H}_{0}$ such that $\widetilde{\psi }\in \widetilde{H}$%
, we set
\begin{equation*}
\widetilde{\mathcal{D}}_{s}\widetilde{\psi }=-\partial _{x}\widetilde{u}%
\left( 0,t\right) \text{, }t\in \left( 0,2\pi \right) \text{,}
\end{equation*}%
where $\widetilde{u}$ solves
\begin{equation}
\left\{
\begin{array}{c}
\partial _{t}\widetilde{u}-\partial _{x}^{2}\widetilde{u}=0\text{, in }%
\widetilde{Q}_{s}\text{,} \\
\widetilde{u}\left( x,0\right) =0\text{, \ \ \ \ for }x\in \left[ 0,1\right]
\text{,} \\
\widetilde{u}\left( \widetilde{s}\left( t\right) ,t\right) =0\text{, for }%
t\in \left[ 0,2\pi \right] \text{,} \\
\widetilde{u}\left( 0,t\right) =\widetilde{\psi }\left( t\right) \text{, for
}t\in \left[ 0,2\pi \right] \text{.}%
\end{array}%
\right.  \label{inst26}
\end{equation}%
Let us also define the following two linear and bounded operators $%
\widetilde{G},\widetilde{G}^{\ast }:\widetilde{H}\rightarrow \widetilde{H}$
such that for any $\widetilde{\psi }\in \widetilde{H}$ we have
\begin{eqnarray}
\left( \widetilde{G}\widetilde{\psi }\right) \left( t\right) &=&\widetilde{%
\psi }\left( t\right) e^{-1/t}t^{-3/2}\text{,}  \label{inst27} \\
\left( \widetilde{G}^{\ast }\widetilde{\psi }\right) \left( t\right) &=&%
\widetilde{\psi }\left( t\right) e^{-1/\left( 2\pi -t\right) }\left( 2\pi
-t\right) ^{-3/2}\text{.}  \notag
\end{eqnarray}%
We have that there exists a constant $C_{1}$ such that
\begin{equation}
\left\Vert \widetilde{G}\right\Vert _{\mathcal{L}\left( \widetilde{H},%
\widetilde{H}\right) }\text{,}\left\Vert \widetilde{G}^{\ast }\right\Vert _{%
\mathcal{L}\left( \widetilde{H},\widetilde{H}\right) }\leq C_{1}\text{.}
\label{inst28}
\end{equation}%
For any $s\in X$, we define the following linear and bounded operator $%
\widetilde{\mathcal{D}}_{1,s}:\widetilde{H}\rightarrow \widetilde{H}^{\prime
}$ as follows
\begin{equation}
\left\langle \widetilde{\mathcal{D}}_{1,s}\widetilde{\psi },\widetilde{\phi }%
\right\rangle _{\widetilde{H}^{\prime },\widetilde{H}}=\left\langle
\widetilde{\mathcal{D}}_{s}\widetilde{G}\widetilde{\psi },\widetilde{G}%
^{\ast }\widetilde{\phi }\right\rangle _{\widetilde{H}^{\prime },\widetilde{H%
}}\text{, for any }\widetilde{\psi },\widetilde{\phi }\in \widetilde{H}\text{%
.}  \label{inst29}
\end{equation}%
The proof of Theorem \ref{Thinst1} is an immediate consequence of the
following two propositions whose proofs we postpone for a while.

\begin{proposition}
\label{Prinst1}There exists a constant $C_{2}$, depending on $m$ and $b$
only, such that for any $s\in X$ and any $\widetilde{\psi }\in \widetilde{H}$
we have
\begin{equation}
\left\Vert \left( \widetilde{\mathcal{D}}_{s}-\widetilde{\mathcal{D}}%
_{s_{0}}\right) \widetilde{\psi }\right\Vert _{\widetilde{H}_{0}}\leq
C_{2}\left\Vert \widetilde{u}_{0}\right\Vert _{H^{2,1}\left( \left(
0,1/2\right) \times \left( 0,2\pi \right) \right) }\text{,}  \label{inst30}
\end{equation}
where $\widetilde{u}_{0}$ solves (\ref{inst26}) with $s=s_{0}$.
\end{proposition}

\begin{proposition}
\label{Prinst2}There exists a positive constant $\delta _{1}$, depending on $%
m$ and $b$ only, such that for any $\delta \in \left( 0,\delta _{1}\right) $
we can find $s_{1},s_{2}\in X$ satisfying (\ref{inst18}) such that
\begin{equation}
\left\Vert \widetilde{\mathcal{D}}_{1,s_{1}}-\widetilde{\mathcal{D}}%
_{1,s_{2}}\right\Vert _{\mathcal{L}\left( \widetilde{H},\widetilde{H}%
^{\prime }\right) }\leq 2\exp \left( -\delta ^{-\frac{1}{6m}}\right) \text{ .%
}  \label{inst31}
\end{equation}
\end{proposition}

\textbf{Proof of Theorem \ref{Thinst1}.} By the trace inequality and
standard estimates for parabolic equations we have, for any $\widetilde{\psi
}\in \widetilde{H}$,
\begin{equation}
\left\Vert \widetilde{\mathcal{D}}_{s_{0}}\widetilde{\psi }\right\Vert _{%
\widetilde{H}_{0}}\leq C_{3}\left\Vert \widetilde{\psi }\right\Vert _{%
\widetilde{H}}\text{ ,}  \label{inst31b}
\end{equation}%
where $C_{3}$ depends on $m$ and $b$ only.

\noindent By Proposition \ref{Prinst1} and by (\ref{inst31b}) we have
\begin{equation}
\left\Vert \widetilde{\mathcal{D}}_{s}\right\Vert _{\mathcal{L}\left(
\widetilde{H},\widetilde{H}_{0}\right) }\leq C_{4}\text{, for every }s\in X%
\text{,}  \label{inst32}
\end{equation}
where $C_{4}$ depends on $m$ and $b$ only.

Now, for any $\psi \in H$, let $\widetilde{\psi }\leq \widetilde{H}$ be its
extension to $0$ to $\left( 0,2\pi \right) \smallsetminus \left( \pi /2,3\pi
/2\right) $. We have that $J:H\rightarrow \widetilde{H}$, where $J\left(
\widetilde{\psi }\right) =\psi $ for any $\psi \in H$, is a linear isometry.
Furthermore, we have that the linear operators $G,G^{\ast }:H\rightarrow H$
such that
\begin{equation*}
G\left( \psi \right) =\widetilde{G}\left( \widetilde{\psi }\right) _{\mid
\left( \pi /2,3\pi /2\right) }\text{ , }G^{\ast }\left( \psi \right) =%
\widetilde{G}^{\ast }\left( \widetilde{\psi }\right) _{\mid \left( \pi
/2,3\pi /2\right) }\text{ ,}
\end{equation*}%
are invertible and
\begin{equation}
\left\Vert G\right\Vert _{\mathcal{L}\left( H,H\right) }\text{ , }\left\Vert
G^{\ast }\right\Vert _{\mathcal{L}\left( H,H\right) }\text{ , }\left\Vert
G^{-1}\right\Vert _{\mathcal{L}\left( H,H\right) }\text{ , }\left\Vert
\left( G^{\ast }\right) ^{-1}\right\Vert _{\mathcal{L}\left( H,H\right)
}\leq C_{5}\text{,}  \label{inst33}
\end{equation}%
where $C_{5}$ is a constant.

\noindent For any $s\in X$ and for any $\psi ,\phi \in H$, we have
\begin{equation}
\mathcal{D}_{s}\psi =\widetilde{\mathcal{D}}_{s}\widetilde{\psi }_{\mid
\left( \pi /2,3\pi /2\right) }  \label{inst34}
\end{equation}
and
\begin{equation}
\left\langle \mathcal{D}_{s}\psi ,\phi \right\rangle _{H^{\prime
},H}=\left\langle \widetilde{\mathcal{D}}_{s}J\widetilde{\psi },J\widetilde{%
\phi }\right\rangle _{\widetilde{H}^{\prime },\widetilde{H}}=\left\langle
\widetilde{\mathcal{D}}_{1,s}JG^{-1}\widetilde{\psi },J\left( G^{\ast
}\right) ^{-1}\widetilde{\phi }\right\rangle _{\widetilde{H}^{\prime },%
\widetilde{H}}\text{.}  \label{inst35}
\end{equation}
By (\ref{inst32}) and (\ref{inst34}) we have
\begin{equation}
\left\Vert \mathcal{D}_{s_{1}}-\mathcal{D}_{s_{2}}\right\Vert _{\mathcal{L}%
\left( H,H_{0}\right) }\leq 2C_{4}\text{, for every }s_{1},s_{2}\in X\text{.}
\label{inst36}
\end{equation}
By (\ref{inst33}) and (\ref{inst35}) we derive that, for every $%
s_{1},s_{2}\in X$, we have
\begin{equation}
\left\Vert \mathcal{D}_{s_{1}}-\mathcal{D}_{s_{2}}\right\Vert _{\mathcal{L}%
\left( H,H^{\prime }\right) }\leq C_{5}^{2}\left\Vert \widetilde{\mathcal{D}}%
_{1,s_{1}}-\widetilde{\mathcal{D}}_{1,s_{2}}\right\Vert _{\mathcal{L}\left(
\widetilde{H},\widetilde{H}^{\prime }\right) }\text{.}  \label{inst37}
\end{equation}
By Proposition \ref{Prinst2}, (\ref{inst36}), (\ref{inst37}) and
interpolation inequality (\ref{inst13}) Theorem \ref{Thinst1} follows.$%
\blacksquare $

\bigskip

\textbf{Proof of Proposition \ref{Prinst1}.} Let $\chi \in C^{\infty }\left( %
\left[ 0,+\infty \right) \right) $ be a function such that $0\leq \chi \leq
1 $, $\chi \left( x\right) =0$, for $x\in \left[ 0,\frac{1}{4}\right] $ and $%
\chi \left( x\right) =1$, for $x\in \left[ \frac{1}{2},+\infty \right) $.
Let $v=\widetilde{u}-\left( 1-\chi \right) \widetilde{u}_{0}$, where $%
\widetilde{u}$ is the solution to (\ref{inst26}). We have
\begin{equation}
\left\{
\begin{array}{c}
\partial _{t}v-\partial _{x}^{2}v=f\text{, in }\widetilde{Q}_{s}\text{,} \\
v\left( x,0\right) =0\text{, \ \ \ \ for }x\in \left[ 0,1\right] \text{,} \\
v\left( \widetilde{s}\left( t\right) ,t\right) =0\text{, for }t\in \left[
0,2\pi \right] \text{,} \\
v\left( 0,t\right) =0\text{, for }t\in \left[ 0,2\pi \right] \text{,}%
\end{array}%
\right.  \label{inst38}
\end{equation}%
where $f=-\chi ^{\prime \prime }\widetilde{u}_{0}-2\chi ^{\prime }\widetilde{%
u}_{0}^{\prime }$. Since \textit{supp\thinspace }$f\subset \left[ \frac{1}{4}%
,\frac{1}{2}\right] $ we have, for a constant $C_{6}$,
\begin{equation}
\left\Vert f\right\Vert _{L^{2}\left( \widetilde{Q}_{s}\right) }\leq
C_{6}\left\Vert \widetilde{u}_{0}\right\Vert _{H^{2,1}\left( \left(
0,1/2\right) \times \left( 0,2\pi \right) \right) }\text{.}  \label{inst39}
\end{equation}%
Now denote by $w$, $f_{1}$ the functions defined as follows
\begin{equation*}
w\left( \xi ,t\right) =v\left( \xi \widetilde{s}\left( t\right) ,t\right)
\text{ , }f_{1}\left( \xi ,t\right) =f\left( \xi \widetilde{s}\left(
t\right) ,t\right) \text{ , for every }\left( \xi ,t\right) \in \left[ 0,1%
\right] \times \left[ 0,2\pi \right] \text{.}
\end{equation*}%
We have that $w$ solves the following problem
\begin{equation*}
\left\{
\begin{array}{c}
\partial _{t}w-\left( \frac{1}{\widetilde{s}^{2}\left( t\right) }\partial
_{\xi }^{2}w+\frac{\xi \widetilde{s}^{\prime }\left( t\right) }{\widetilde{s}%
\left( t\right) }\partial _{\xi }w\right) =f_{1}\text{, in }\left(
0,1\right) \times \left( 0,2\pi \right) \text{,} \\
w\left( \xi ,0\right) =0\text{, \ \ \ \ for }\xi \in \left[ 0,1\right] \text{%
,} \\
w\left( 1,t\right) =0\text{, for }t\in \left[ 0,2\pi \right] \text{,} \\
w\left( 0,t\right) =0\text{, for }t\in \left[ 0,2\pi \right] \text{,}%
\end{array}%
\right.
\end{equation*}%
Therefore, by standard regularity estimates for parabolic equations, (\ref%
{inst39}) and a simple change of variable we have
\begin{equation}
\left\Vert w\right\Vert _{H^{2,1}\left( \left( 0,1\right) \times \left(
0,2\pi \right) \right) }\leq C_{8}\left\Vert f\right\Vert _{L^{2}\left(
\widetilde{Q}_{s}\right) }\text{,}  \label{inst40b}
\end{equation}%
where $C_{8}$ depends on $b$ and $m$ only.

\noindent On the other side, for every $t\in \left( 0,2\pi \right) $ we have
\begin{equation*}
\partial _{\xi }w\left( 0,t\right) =-\widetilde{s}\left( t\right) \partial
_{x}v\left( 0,t\right) =-\widetilde{s}\left( t\right) \left( \partial _{x}%
\widetilde{u}-\partial _{x}u_{0}\right) \left( 0,t\right) =\widetilde{s}%
\left( t\right) \left( \widetilde{\mathcal{D}}_{s}\widetilde{\psi }-%
\widetilde{\mathcal{D}}_{s_{0}}\widetilde{\psi }\right) \left( 0,t\right)
\text{,}
\end{equation*}%
hence
\begin{equation}
\left( \widetilde{\mathcal{D}}_{s}\widetilde{\psi }-\widetilde{\mathcal{D}}%
_{s_{0}}\widetilde{\psi }\right) \left( 0,t\right) =-\frac{1}{\widetilde{s}%
\left( t\right) }\partial _{\xi }w\left( 0,t\right)  \label{inst40c}
\end{equation}%
and by the trace theorem there exists a constant $C_{9}$ such that
\begin{equation*}
\left\Vert \partial _{\xi }w\left( 0,.\right) \right\Vert _{\widetilde{H}%
_{0}}\leq C_{9}\left\Vert w\right\Vert _{H^{2,1}\left( \left( 0,1\right)
\times \left( 0,2\pi \right) \right) }\text{.}
\end{equation*}%
By the just obtained inequality (\ref{inst39}), (\ref{inst40b}) and (\ref%
{inst40c}) we get (\ref{inst30}).$\blacksquare $

\bigskip

In order to prove Proposition \ref{Prinst2} we need further auxiliary
propositions and new notation. The following proposition has been proved in
\cite{Man}. Let us fix $m\in \mathbb{N}$ and a positive numer $b$. Let $%
\varepsilon $ be a fixed positive number, we denote by $X_{\varepsilon }$
the set $\left\{ s\in X:d_{0}\left( s,s_{0}\right) \leq \varepsilon \right\}
$. Notice that $X_{\varepsilon }=X_{mb\varepsilon }$.

\begin{proposition}
\label{Mandache}There exists a positive constant $\varepsilon _{0}$,
depending on $m$ and $b$ only, such that for every $\varepsilon \in \left(
0,\varepsilon _{0}\right) $ there exists $W\subset X_{\varepsilon }$
satifying the following properties. $W$ is $\varepsilon $-discrete with
respect to the distance $d_{0}$ and $W$ has at least $\exp \left(
2^{-1}\varepsilon _{0}^{\frac{1}{m}}\varepsilon ^{-\frac{1}{m}}\right) $
elements.
\end{proposition}

\bigskip

For any $n\in \mathbb{N}$ denote
\begin{equation*}
y_{n}\left( t\right) =\left\{
\begin{array}{c}
\frac{1}{\sqrt{\pi }}e^{-1/t}t^{-3/2}\sin \left( \frac{n}{2}\right) t\text{,
if }t>0\text{,} \\
0\text{, if }t\leq 0\text{.}%
\end{array}
\right.
\end{equation*}
Notice that $y_{n}\in C^{\infty }\left( \mathbb{R}\right) $ and $y_{n}\in
H^{m}\left( \mathbb{R}\right) $ for any $n,m\in \mathbb{N}$. Let $U_{n}$ be
the solution to the following boundary value problem
\begin{equation}
\left\{
\begin{array}{c}
\partial _{t}U_{n}-\partial _{x}^{2}U_{n}=0\text{, in }\left( 0,1\right)
\times \left( 0,+\infty \right) \text{,} \\
U_{n}\left( x,t\right) =0\text{, in }\left[ 0,1\right] \times \left( -\infty
,0\right] \text{,} \\
U_{n}\left( 1,t\right) =y_{n}\left( t\right) \text{, for }t\in \mathbb{R}%
\text{,} \\
U_{n}\left( 0,t\right) =0\text{, for }t\in \mathbb{R}\text{.}%
\end{array}
\right.  \label{inst41}
\end{equation}
The following decay estimate of exponential type will be crucial.

\bigskip

\begin{proposition}
\label{Prinst3}For any fixed $n\in \mathbb{N}$, let $U_{n}$ solve (\ref%
{inst41}). Let us fix $\rho _{0}\in \left( 0,1\right) $. Then there exist
positive constants $K_{1}$ and $k_{1}$ depending on $\rho _{0}$ only, such
that
\begin{equation}
\left\Vert U_{n}\right\Vert _{L^{2}\left( B_{\rho _{0}}\times \mathbb{R}%
\right) }\leq K_{1}\exp \left( -k_{1}\gamma \left( \widetilde{\psi }%
_{n}\right) \right) \text{. }  \label{inst42}
\end{equation}
\end{proposition}

\textbf{Proof. }Let us denote
\begin{equation*}
\widehat{U}_{n}\left( x,\xi \right) =\int\nolimits_{-\infty }^{+\infty
}e^{-i\xi t}U_{n}\left( x,t\right) dt\text{ ,}
\end{equation*}
\begin{equation*}
\widehat{y}_{n}\left( \xi \right) =\int\nolimits_{-\infty }^{+\infty
}e^{-i\xi t}y_{n}\left( t\right) dt\text{.}
\end{equation*}
By (\ref{inst41}) we get
\begin{equation}
\left\{
\begin{array}{c}
\partial _{x}^{2}\widehat{U}_{n}-i\xi \widehat{U}_{n}=0\text{, in }\left(
0,1\right) \times \mathbb{R}\text{,} \\
\widehat{U}_{n}\left( 1,\xi \right) =\widehat{y}_{n}\left( \xi \right) \text{%
, for }\xi \in \mathbb{R}\text{,} \\
\widehat{U}_{n}\left( 0,\xi \right) =0\text{, for }\xi \in \mathbb{R}\text{.}%
\end{array}
\right.  \label{inst46}
\end{equation}
Therefore
\begin{equation}
\widehat{U}_{n}\left( x,\xi \right) =\widehat{y}_{n}\left( \xi \right) \frac{%
\sinh \sigma \left( 1-x\right) }{\sinh \sigma }\text{,}  \label{inst48}
\end{equation}
where $\sigma =\sqrt{\frac{\left\vert \xi \right\vert }{2}}\left( i-\text{%
sign}\xi \right) $. Notice that $\sigma ^{2}=i\xi $. Hence there exists a
constant $C$ such that
\begin{equation}
\left\vert \widehat{U}_{n}\left( x,\xi \right) \right\vert =C\left\vert
\widehat{y}_{n}\left( \xi \right) \right\vert e^{-x\sqrt{\frac{\left\vert
\xi \right\vert }{2}}}  \label{inst48a}
\end{equation}
By using classical tables of integral transforms, see for instance \cite{Erd}%
, we have that, for any $\xi \in \mathbb{R}$,
\begin{eqnarray*}
\widehat{y}_{n}\left( \xi \right) &=&\frac{e^{-\sqrt{2\left\vert \xi
-n/2\right\vert }}}{2i}\left( \cos \left( \sqrt{2\left\vert \xi
-n/2\right\vert }\right) -\text{sign}\left( \xi -n/2\right) i\sin \left(
\sqrt{2\left\vert \xi -n/2\right\vert }\right) \right) \\
&&-\frac{e^{-\sqrt{2\left\vert \xi +n/2\right\vert }}}{2i}\left( \cos \left(
\sqrt{2\left\vert \xi +n/2\right\vert }\right) -\text{sign}\left( \xi
+n/2\right) i\sin \left( \sqrt{2\left\vert \xi +n/2\right\vert }\right)
\right) \text{, }
\end{eqnarray*}
hence the following estimate holds true
\begin{equation}
\left\vert \widehat{y}_{n}\left( \xi \right) \right\vert \leq e^{-\sqrt{%
2\left\vert \xi -n/2\right\vert }}+e^{-\sqrt{2\left\vert \xi +n/2\right\vert
}}\text{, for every }\xi \in \mathbb{R}\text{.}  \label{inst47}
\end{equation}
By Plancherel Theorem, (\ref{inst48a}) and (\ref{inst47}) we have, for any $%
\rho _{0}\in \left( 0,1\right) $ and for any $n\in \mathbb{N}$,
\begin{equation*}
\int\nolimits_{\left( 0,\rho _{0}\right) \times \left( 0,+\infty \right)
}\left\vert U\left( x,t\right) \right\vert ^{2}dxdt=\frac{1}{2\pi }%
\int\nolimits_{0}^{\rho _{0}}dx\int\nolimits_{-\infty }^{+\infty }\left\vert
\widehat{U}_{n}\left( x,\xi \right) \right\vert ^{2}d\xi
\end{equation*}
\begin{equation*}
\leq C\int\nolimits_{0}^{\rho _{0}}dx\int\nolimits_{-\infty }^{+\infty
}\left\vert \widehat{y}_{n}\left( \xi \right) \right\vert ^{2}e^{\sqrt{%
2\left\vert \xi \right\vert }x}d\xi \leq C\rho _{0}\int\nolimits_{-\infty
}^{+\infty }e^{-2\sqrt{2\left\vert \xi -n/2\right\vert }}e^{\sqrt{%
2\left\vert \xi \right\vert }\rho _{0}}d\xi
\end{equation*}
\begin{equation*}
C\rho _{0}\left( \int\nolimits_{\left\vert \xi \right\vert \leq n/4}e^{-2%
\sqrt{2\left\vert \xi -n/2\right\vert }}e^{\sqrt{2\left\vert \xi \right\vert
}\rho _{0}}d\xi +\int\nolimits_{\left\vert \xi \right\vert \geq n/4}e^{-2%
\sqrt{2\left\vert \xi -n/2\right\vert }}e^{\sqrt{2\left\vert \xi \right\vert
}\rho _{0}}d\xi \right) \leq K_{1}e^{-k_{1}\sqrt{\frac{n}{2}}}\text{,}
\end{equation*}
where $K_{1}$, $k_{1}$ depend on $\rho _{0}$ only. Recalling that $\gamma
\left( \widetilde{\psi }_{n}\right) =\sqrt{\frac{n}{2}}$ , for every $n\in
\mathbb{N}$, by the inequality proved above we get (\ref{inst42}).$%
\blacksquare $

\bigskip

\begin{lemma}
\label{Leinst}There exist positive constants $K_{2}$ and $k_{2}$ depending
on $m$ and $b$ only, such that for any $n\in \mathbb{N}$ and any $s\in X$,
we have
\begin{equation}
\left\Vert \left( \widetilde{\mathcal{D}}_{1,s}-\widetilde{\mathcal{D}}%
_{1,s_{0}}\right) \widetilde{\psi }_{n}\right\Vert _{\widetilde{H}^{\prime
}}\leq K_{2}\exp \left( -k_{2}\gamma \left( \widetilde{\psi }_{n}\right)
\right) \text{.}  \label{inst43}
\end{equation}
\end{lemma}

\textbf{Proof. }By (\ref{inst28}) and (\ref{inst29}) we have
\begin{equation}
\left\Vert \left( \widetilde{\mathcal{D}}_{1,s}-\widetilde{\mathcal{D}}%
_{1,s_{0}}\right) \widetilde{\psi }_{n}\right\Vert _{\widetilde{H}^{\prime
}}\leq C_{1}\left\Vert \left( \widetilde{\mathcal{D}}_{s}-\widetilde{%
\mathcal{D}}_{s_{0}}\right) \widetilde{G}\widetilde{\psi }_{n}\right\Vert _{%
\widetilde{H}^{\prime }}\text{.}  \label{inst43a}
\end{equation}%
Let $U_{n}$ solve (\ref{inst41}). Then the restriction of $U_{n}$ to the
time interval $\left( 0,2\pi \right) $ solves (\ref{inst26}) with $s=s_{0}$
and boundary data at $x=0$ and $x=1$ equal respectively to $\widetilde{G}%
\widetilde{\psi }_{n}$ and $0$. Therefore, by Proposition \ref{Prinst1} and (%
\ref{inst43a}), we have
\begin{equation*}
\left\Vert \left( \widetilde{\mathcal{D}}_{1,s}-\widetilde{\mathcal{D}}%
_{1,s_{0}}\right) \widetilde{\psi }_{n}\right\Vert _{\widetilde{H}^{\prime
}}\leq K_{3}\left\Vert U_{n}\right\Vert _{H^{2,1}\left( \left( 0,\frac{1}{2}%
\right) \times \left( 0,2\pi \right) \right) }\text{.}
\end{equation*}%
Then the conclusion follows by Proposition \ref{Prinst3} and a Caccioppoli
type inequality of the form
\begin{equation*}
\left\Vert U_{n}\right\Vert _{H^{2,1}\left( \left( 0,\rho _{1}\right) \times
\left( 0,2\pi \right) \right) }\leq K_{4}\left\Vert U_{n}\right\Vert
_{L^{2}\left( \left( 0,\rho _{0}\right) \times \left( 0,2\pi \right) \right)
}\text{,}
\end{equation*}%
where $0<\rho _{1}<\rho _{0}<1$ and $K_{4}$ depends on $\rho _{0}$ and $\rho
_{1}$ only. For similar Caccioppoli type inequalities we refer to Section 6,
Chapter III, of \cite{LSU}.$\blacksquare $

\bigskip

\textbf{Proof of Proposition \ref{Prinst2}. }First we introduce the operator
which is the adjoint to $\widetilde{\mathcal{D}}_{1,s}$. For any $s\in X$,
let us define the linear and bounded operator $\widetilde{\mathcal{D}}%
_{s}^{\ast }:\widetilde{H}\rightarrow \widetilde{H}_{0}$ such that for any $%
\widetilde{\phi }\in \widetilde{H}$, we have
\begin{equation*}
\widetilde{\mathcal{D}}_{s}^{\ast }\widetilde{\phi }=-\partial _{x}%
\widetilde{v}\left( 0,t\right) \text{ ,}
\end{equation*}
where $\widetilde{v}$ solves
\begin{equation}
\left\{
\begin{array}{c}
\partial _{t}\widetilde{v}+\partial _{x}^{2}\widetilde{v}=0\text{, in }%
\widetilde{Q}_{s}\text{,} \\
\widetilde{v}\left( x,2\pi \right) =0\text{, \ \ \ \ for }x\in \left[ 0,1%
\right] \text{,} \\
\widetilde{v}\left( \widetilde{s}\left( t\right) ,t\right) =0\text{, for }%
t\in \left[ 0,2\pi \right] \text{,} \\
\widetilde{v}\left( 0,t\right) =\widetilde{\phi }\text{, for }t\in \left[
0,2\pi \right] \text{.}%
\end{array}
\right.  \label{inst44}
\end{equation}
Then we define $\widetilde{\mathcal{D}}_{1,s}^{\ast }:\widetilde{H}%
\rightarrow \widetilde{H}^{\prime }$ so that
\begin{equation*}
\left\langle \widetilde{\mathcal{D}}_{1,s}^{\ast }\widetilde{\phi },%
\widetilde{\psi }\right\rangle _{\widetilde{H}^{\prime },\widetilde{H}%
}=\left\langle \widetilde{\mathcal{D}}_{s}^{\ast }\widetilde{G}^{\ast }%
\widetilde{\phi },\widetilde{G}\widetilde{\psi }\right\rangle _{\widetilde{H}%
^{\prime },\widetilde{H}}\text{ , for every }\widetilde{\phi },\widetilde{%
\psi }\in \widetilde{H}\text{.}
\end{equation*}
By using the weak formulation of (\ref{inst26}) and (\ref{inst44}) it is
easy to show that the following adjointness property holds true for any $%
s\in X$%
\begin{equation}
\left\langle \widetilde{\mathcal{D}}_{1,s}\widetilde{\psi },\widetilde{\phi }%
\right\rangle _{\widetilde{H}^{\prime },\widetilde{H}}=\left\langle
\widetilde{\mathcal{D}}_{1,s}^{\ast }\widetilde{\phi },\widetilde{\psi }%
\right\rangle _{\widetilde{H}^{\prime },\widetilde{H}}\text{ , for every }%
\widetilde{\psi },\widetilde{\phi }\in \widetilde{H}\text{.}  \label{inst44a}
\end{equation}
It is not difficult to show that, by the change of variable in time, (\ref%
{inst43}) holds true if we replace there $\widetilde{\mathcal{D}}_{1,s}$
with $\widetilde{\mathcal{D}}_{1,s}^{\ast }$. Therefore,using (\ref{inst43a}%
), we have
\begin{equation}
\left\vert \left\langle \left( \widetilde{\mathcal{D}}_{1,s}-\widetilde{%
\mathcal{D}}_{1,s_{0}}\right) \widetilde{\psi }_{n},\widetilde{\psi }%
_{n^{\prime }}\right\rangle _{\widetilde{H}^{\prime },\widetilde{H}%
}\right\vert \leq K_{5}\exp \left( -k_{2}\max \left\{ \gamma \left(
\widetilde{\psi }_{n}\right) ,\gamma \left( \widetilde{\psi }_{n^{\prime
}}\right) \right\} \right) \text{,}  \label{inst45}
\end{equation}
where $K_{5}$ depends on $m$ and $b$ only.

Finally, recalling Proposition \ref{Mandache}, the inequality (\ref{inst25})
and the decay estimate (\ref{inst45}) we notice that we are exactly in the
position of applying Theorem \ref{ThabsInst} and Proposition \ref{Prinst2}
follows.$\blacksquare $

\subsection{Stability properties of the Dirichlet inverse problem with
unknown boundary independent on time\label{INSTNT}}

In this section, for the reader's convenience, first we present the
formulation of the Dirichlet inverse problem when the unknown boundary is
independent on time and the equation is linear, then we state, without
proofs, the stability properties (i.e. stability estimate and exponential
instability) for such a problem. Such a proof has the same structure of the
proof of Theorem \ref{Thinst1} but it presents additional difficulties with
respect to this Theorem. For the more details we refer to \cite{DcRVe}.

Referring to the problem formulated at beginning of Section \ref{Sec4}, in
this section $\Omega \left( t\right) =\Omega $ and $I\left( t\right) =I$ for
every $t\in \mathbb{R}$, where $\Omega $ is a bounded domain in $\mathbb{R}%
^{n}$ with a sufficiently smooth boundary $\partial \Omega $ and $I$ is an
inaccessible part of $\partial \Omega $. Given a nontrivial function $f$ on $%
A\times \left( 0,T\right) $, let us consider the following Cauchy-Dirichlet
(direct) problem
\begin{equation}
\left\{
\begin{array}{c}
div\left( \kappa \left( x,t\right) \nabla u\right) -\partial _{t}u=0\text{,
\ \ in }\Omega \times \left( 0,T\right) \text{,} \\
u=f\text{, \ \ \ \ \ \ \ \ \ \ \ on }A\times \left( 0,T\right] \text{,\ } \\
u=0\text{, \ \ \ \ \ \ \ \ \ \ \ \ \ on }I\times \left( 0,T\right] \text{, \
\ } \\
u\left( .,0\right) =0\text{, \ \ \ \ \ \ \ \ \ \ in }\Omega \text{ ,}%
\end{array}
\right.  \label{DcRV1}
\end{equation}
where $\kappa \left( x,t\right) =\left\{ \kappa ^{ij}\left( x,t\right)
\right\} _{i,j=1}^{n}$ denotes a known symmetric matrix which satisfies a
hypothesis of uniform ellipticity and some smothness conditions that we
shall specify below.

Given an open portion $\Sigma $ of $\partial \Omega $ such that $\Sigma
\subset A$, we consider the inverse problem of determinig $I$, from the
knowledge of $\kappa \left( x,t\right) \nabla u\cdot \nu $ on $\Sigma \times %
\left[ 0,T\right] $.

Let $\lambda _{0}$, $\Lambda _{0}$, $E$, $M$, $R_{0}$, $\widetilde{F}$ and $%
\beta $ be given positive numbers with $\lambda _{0}$, $\beta \in \left( 0,1%
\right] $.

\textit{i) A priori information on the domain }$\Omega $.

We assume
\begin{equation}
\left\vert \Omega \right\vert \leq MR_{0}^{n}\text{.}  \label{DcRV2}
\end{equation}%
Concerning the regularity of $\partial \Omega $ we assume that
\begin{equation}
\partial \Omega \text{ is of }C^{0,1}\text{ class a with constants }R_{0}%
\text{, }E\text{.}  \label{DcRV3}
\end{equation}%
In addition we assume that there exist open portions $\Sigma $, $\widetilde{%
\Sigma }$ within $A$ such that $\Sigma \subset \widetilde{\Sigma }$ and
there exists a point $P_{0}\in \Sigma $ such that
\begin{equation}
\partial \Omega \cap B_{R_{0}}\left( P_{0}\right) \subset \Sigma \text{ , }
\label{DcRV5a}
\end{equation}%
and
\begin{equation}
\widetilde{\Sigma }\cap \left( I\right) ^{R_{0}}=\varnothing \text{ .}
\label{DcRV5b}
\end{equation}%
and
\begin{equation}
I\cup A=\partial \Omega \text{ , Int}I\cap \text{Int}\left( A\right)
=\varnothing \text{ , }I\cap A=\partial A=\partial I\text{ .}  \label{DcRV4}
\end{equation}%
Recall that interior and boundaries are intented in the relative topology of
$\partial \Omega $.

\textit{ii) A priori information on the unknown part of the boundary.} We
suppose that we have
\begin{equation}
I\text{ is of }C^{1,\beta }\text{ class with constant }R_{0}\text{, }E\text{.%
}  \label{DcRV6}
\end{equation}

\textit{iii) A priori information about prescribed boundary datum.}

We shall assume that

\begin{equation}
f\in H^{\frac{1}{2},\frac{1}{4}}\left( A\times \left( 0,T\right) \right)
\text{,}  \label{DcRV7a}
\end{equation}
\begin{equation}
\text{\textit{supp}}f\subset \widetilde{\Sigma }\text{ }  \label{DcRV7b}
\end{equation}
and
\begin{equation}
\frac{\left\| f\right\| _{H^{1/2,1/4}\left( A\times \left( 0,T\right)
\right) }}{\left\| f\right\| _{L^{2}\left( A\times \left( 0,T\right) \right)
}}\leq \widetilde{F}\text{.}  \label{DcRV7c}
\end{equation}

\bigskip

\textit{iii) Assumptions about the thermal conductivity }$\kappa $.

The thermal conductivity $\kappa $ is assumed to be a given function from $%
\mathbb{R}^{n}\times \mathbb{R}$ with values in $n\times n$ matrices
satisfying the following conditions. When $\xi \in \mathbb{R}^{n}$ and $%
\left( x,s\right) ,\left( y,\tau \right) \in \mathbb{R}^{n+1}$ we assume
that
\begin{equation}
\lambda _{0}\left\vert \xi \right\vert ^{2}\leq \kappa \left( x,t\right) \xi
\cdot \xi \leq \lambda _{0}^{-1}\left\vert \xi \right\vert ^{2}\text{, }
\label{DcRV8}
\end{equation}%
\begin{equation}
\left\vert \kappa \left( x,t\right) -\kappa \left( y,\tau \right)
\right\vert \leq \Lambda _{0}\left( \frac{\left\vert x-y\right\vert }{R_{0}}+%
\frac{\left\vert t-s\right\vert }{R_{0}}\right) \text{.}  \label{DcRV9a}
\end{equation}%
In the next theorem we shall use the following notation
\begin{equation*}
W\left( \Omega \times \left( 0,T\right) \right) =\left\{ v:v\in L^{2}\left(
\left( 0,T\right) ;H^{1}\left( \Omega \right) \right) \text{, }\partial
_{t}v\in L^{2}\left( \left( 0,T\right) ;H^{-1}\left( \Omega \right) \right)
\right\} \text{.}
\end{equation*}

\begin{theorem}
\label{ThDcRVe}Let $\Omega _{1}$, $\Omega _{2}$ be domains satisfying (\ref%
{DcRV2}) (\ref{DcRV3}). For any $i=1,2$ let $A_{i}$, $I_{i}$ satisfying (\ref%
{DcRV4}), be the accessible and inaccessible part of $\partial \Omega _{i}$
respectively. Assume that $A_{1}=A_{2}=A$ and that $\Omega _{1}$ and $\Omega
_{2}$ lie on the same side of $A$. Let us take $\Sigma $, $\widetilde{\Sigma
}$ within $A$ such that $\Sigma \subset \widetilde{\Sigma }$ satisfying (\ref%
{DcRV5a}) (\ref{DcRV5b}). Finally, we suppose that, for any $i=1,2$, $I_{i}$
satisfies the a priori information (\ref{DcRV6}).

Let us assume that (\ref{DcRV7a}), (\ref{DcRV7b}), (\ref{DcRV7c}), (\ref%
{DcRV8}) and (\ref{DcRV9a}) are also satisfied and let $u_{i}\in W\left(
\Omega _{i}\times \left( 0,T\right) \right) $ be the weak solution to (\ref%
{DcRV1}) when $\Omega =\Omega _{i}$, $i=1,2$. If, given $\varepsilon >0$, we
have
\begin{equation}
R_{0}\left\| k\nabla u_{1}\cdot \nu -k\nabla u_{2}\cdot \nu \right\|
_{L^{2}\left( \Sigma \times \left( 0,T\right) \right) }\leq \varepsilon
\text{ ,}  \label{DcRV10}
\end{equation}
then we have, for every $\varepsilon \in \left( 0,1\right) $,
\begin{equation}
d_{\mathcal{H}}\left( \overline{\Omega _{1}},\overline{\Omega _{2}}\right)
\leq R_{0}C_{1}\left| \log \varepsilon \right| ^{-\frac{1}{C_{2}}}\text{ ,}
\label{DcRV11}
\end{equation}
where $C_{1}$, $C_{2}$ are positive constants depending on $\lambda _{0}$, $%
\Lambda _{0}$, $E$, $M$, $R_{0}$, $\widetilde{F}$, $\beta $ and $%
R_{0}^{2}T^{-1}$ only.
\end{theorem}

In order to analyse the instability of the inverse problem we study what
happens in the most favourable situation, that is we suppose to have strong
a priori information on $I$, $A$ and $\kappa $ as simple as possible. In
addition we take into account all possible measurements.

More precisely in the instability analysis we deal within the following
framework. Let us assume $\Omega =B_{1}\smallsetminus D$ where $D$ is a
compact subset of $B_{1}$ ($D$ represents an unknown cavity). Let $%
A=\partial B_{1}$, $I=\partial D$. We also set $\Sigma =\widetilde{\Sigma }%
=A $. We set $Q=B_{1}\times \left( \frac{\pi }{2},\frac{3\pi }{2}\right) $
and $\Gamma =A\times \left( \frac{\pi }{2},\frac{3\pi }{2}\right) $. In
addition we assume that $\kappa $ is equal to the identity matrix. Let us
fix a positive integer $m$. To each $C^{m}$ regular cavity $D$ we associate
its Dirichlet-to-Neumann map $\mathcal{D}_{D}$, that is the operator that
maps each prescribed temperature $f$ on $A\times \left( \frac{\pi }{2},\frac{%
3\pi }{2}\right) $ into the corresponding heat flux $\frac{\partial u}{%
\partial \nu }_{\mid A\times \left( \frac{\pi }{2},\frac{3\pi }{2}\right) }$%
. We shall give below a formal definition of $\mathcal{D}_{D}$.

We shall use the following Hibert space. The space $H=H_{,0}^{\frac{3}{2},%
\frac{3}{4}}\left( \Gamma \right) $, its dual $H^{\prime }=H_{1}=H^{-\frac{3%
}{2},-\frac{3}{4}}\left( \Gamma \right) $ and $H_{0}=H^{\frac{1}{2},\frac{1}{%
4}}\left( \Gamma \right) $. We consider the interpolation spaces between $%
H_{0}$ and $H_{1}$, \cite{LiMa}. For any $\theta $, $0\leq \theta \leq 1$,
we define $H_{\theta }$ as $\left[ H_{0},H_{1}\right] _{\theta }$, where the
latter denotes the interpolation at level $\theta $ between $H_{0}$ and $%
H_{1}$. By using the interpolation properties of fractional order Sobolev
spaces, see \cite{LiMa}, we can characterize $H_{\theta }$ as follows
\begin{equation*}
H_{\theta }=\left\{
\begin{array}{cc}
H^{2\left( \frac{1}{4}-\theta \right) ,\frac{1}{4}-\theta }\left( \left(
\frac{\pi }{2},\frac{3\pi }{2}\right) \right) \text{, } & \text{if }0\leq
\theta \leq \frac{1}{4}\text{,} \\
H^{-2\left( \frac{1}{4}-\theta \right) ,-\left( \frac{1}{4}-\theta \right)
}\left( \left( \frac{\pi }{2},\frac{3\pi }{2}\right) \right) \text{,} &
\text{if }\frac{1}{4}\leq \theta \leq 1\text{, }\theta \neq \frac{3}{4}\text{%
.}%
\end{array}
\right.
\end{equation*}
Let us notice that the interesting case $\theta =\frac{1}{4}$, where we have
$H_{\theta }=L^{2}\left( \Gamma \right) $ and that $H_{\frac{3}{4}}$ do not
coincide with $H^{-1,-\frac{1}{2}}\left( \Gamma \right) $.

Let us fix an integer $m\geq 2$ and positive constants $\delta $, $b$ and $r$%
.

For a strictly positive function $g$ defined on $\partial B_{1}$ we recall
the definition of its radial subgraph,
\begin{equation*}
\text{subgraph}_{rad}\left( g\right) =\left\{ y\in \mathbb{R}^{n}:y=\rho
\omega \text{, }0\leq \rho \leq g\left( r\omega \right) \text{, }\omega \in
\partial B_{1}\right\} \text{.}
\end{equation*}
Then we define
\begin{eqnarray*}
X_{mb\delta }\left( \overline{B}_{r}\right) &=&\left\{ \text{subgraph}%
_{rad}\left( g\right) :g\in C^{\infty }\left( \partial B_{r}\right) \text{,}%
\right. \\
&&\left. \left\| g\right\| _{C^{m}\left( \partial B_{r}\right) }\leq b\text{
and }r\leq g\leq r+\delta \right\} \text{.}
\end{eqnarray*}
Let us consider the metric space $\left( X,d\right) =\left( X_{mb\frac{1}{4}%
}\left( \overline{B}_{1/2}\right) ,d_{\mathcal{H}}\right) $, where $d_{%
\mathcal{H}}$ denotes the Hausdorff distance. Let us notice that every $D\in
X$ is closed, is starshaped with respect to the origin and satisfies $%
\overline{B}_{1/2}\subset D\subset \overline{B}_{3/4}$.

For any $D\in X$, we set $Q\left( D\right) =\left( B_{1}\smallsetminus
D\right) \times \left( \frac{\pi }{2},\frac{3\pi }{2}\right) $, $\Gamma
\left( D\right) =\partial D\times \left( \frac{\pi }{2},\frac{3\pi }{2}%
\right) $. If $D=\emptyset $ then we set $Q\left( D\right) =Q$ and $\Gamma
\left( D\right) =\emptyset $.

For any $D\in X\cup \left\{ \emptyset \right\} $, we consider the operator $%
\mathcal{D}_{D}:H\rightarrow H_{0}$ which is defined as follows. For any $%
\psi \in H$, let $u\in H^{2,1}\left( Q\left( D\right) \right) $ be the
solution to
\begin{equation}
\left\{
\begin{array}{c}
\Delta u-\partial _{t}u=0\text{, \ \ in }Q\left( D\right) \text{,} \\
u=\psi \text{, \ \ \ \ \ \ \ \ \ \ \ on }\Gamma \text{,\ } \\
u=0\text{, \ \ \ \ \ \ \ \ \ \ \ \ \ on }\Gamma \left( D\right) \text{, \ \ }
\\
u\left( .,\frac{\pi }{2}\right) =0\text{, \ \ \ \ \ \ \ \ \ \ in }\Omega
\smallsetminus D\text{ ,}%
\end{array}
\right.  \label{DcRV15}
\end{equation}
Then, for any $\psi \in H$, we set
\begin{equation*}
\mathcal{D}_{D}\psi =\frac{\partial u}{\partial \nu }_{\mid \Gamma }\text{ ,
}u\text{ solution to (\ref{DcRV15}).}
\end{equation*}

We have that Theorems 4.3 and 6.2 in Chapter 4 of \cite{LiMa} imply,
respectively, existence and uniqueness of a solution $u\in H^{2,1}\left(
Q\left( D\right) \right) $ to (\ref{DcRV15}) and its continuous dependence
from the boundary datum $\psi \in H$. In addition, by trace theorem (Chapter
4 Theorem 2.1 of \cite{LiMa}) we can conclude that, for any $D\in X\cup
\left\{ \emptyset \right\} $ the operator $\mathcal{D}_{D}:H\rightarrow
H_{0} $ is linear and bounded. We can also consider $\mathcal{D}_{D}$ as an
operator belonging to $\mathcal{L}\left( H,H^{\prime }\right) $, by setting,
for any $\psi ,\phi \in H$,
\begin{equation*}
\left\langle \mathcal{D}_{D}\psi ,\phi \right\rangle _{H^{\prime
},H}=\left\langle \frac{\partial u}{\partial \nu }_{\mid \Gamma },\phi
\right\rangle _{H^{\prime },H}=\int\nolimits_{\Gamma }\frac{\partial u}{%
\partial \nu }\phi dS\text{,}
\end{equation*}
where $u$ solves (\ref{DcRV15}) and $\left\langle .,.\right\rangle
_{H^{\prime },H}$ is the duality pairing between $H^{\prime }$ and $H$.

\begin{theorem}
\label{ThinstDcRV}Let us fix an integer $m\geq 2$ and a positive constant $b$%
. Let $\left( X,d\right) =\left( X_{mb\frac{1}{4}}\left( \overline{B}%
_{1/2}\right) ,d_{\mathcal{H}}\right) $. Then there exists a positive
constant $\overline{\delta }$, depending on $m$ and $b$ only, such that for
every $\delta \in \left( 0,\overline{\delta }\right) $ we can find $%
D_{1},D_{2}\in X$ satisfying
\begin{equation*}
d\left( D_{i},\overline{B}_{1/2}\right) \leq \delta \text{, for any }i=1,2%
\text{; }d\left( D_{1},D_{2}\right) \geq \delta
\end{equation*}%
and for any $\theta \in \left[ 0,1\right] $%
\begin{equation*}
\left\Vert \mathcal{D}_{D_{1}}-\mathcal{D}_{D_{2}}\right\Vert _{\mathcal{L}%
\left( H,H_{\theta }\right) }\leq C\exp \left( -\theta \delta ^{-\frac{1}{6m}%
}\right) \text{,}
\end{equation*}%
where $C$ is a constant depending on $m$, $b$ and $\theta $ only.
\end{theorem}

\textbf{Acknowledgments.} I would like to express my gratitude to Professor
Elena Beretta and Professor Edi Rosset for useful and stimulating discussion
about the problem.

\bigskip

Sergio Vessella,

DiMaD, Universit\`{a} di Firenze,

via Lombroso 6/17

50134 Firenze (Italy).

\textit{E-mail}: sergio.vessella@dmd.unifi.it

\end{document}